\documentclass[12pt,twoside,leqno,showkeys]{article}
\usepackage{amsmath}
\usepackage{amssymb}
\usepackage{amsxtra}
\usepackage{amscd}
\usepackage{amsthm}
\usepackage[mathscr]{eucal}

\usepackage{graphicx}
\usepackage{color}

\theoremstyle{plain}

\newtheorem{prop}[subsection]{Proposition}

\newtheorem{sbthm}[subsubsection]{Theorem}
\newtheorem{sbprop}[subsubsection]{Proposition}
\newtheorem{sbcor}[subsubsection]{Corollary}
\newtheorem{sblem}[subsubsection]{Lemma}

\theoremstyle{definition}

\newtheorem{sbrem}[subsubsection]{Remark}

\newtheorem{sbpara}[subsubsection]{}

\newenvironment{pf}{\proof[\proofname]}{\endproof}

\usepackage{fancyhdr}
\begin{document}

\title{Classifying spaces of degenerating mixed Hodge structures, V:
Extended period domains and algebraic groups}

\author
{Kazuya Kato, Chikara Nakayama, Sampei Usui}

\maketitle
\renewcommand{\mathbb}{\bold}

\newcommand\Cal{\mathcal}
\newcommand\define{\newcommand}
\define\gp{\mathrm{gp}}%
\define\fs{\mathrm{fs}}%
\define\an{\mathrm{an}}%
\define\mult{\mathrm{mult}}%
\define\add{\mathrm{add}}%
\define\Ker{\mathrm{Ker}\,}%
\define\Coker{\mathrm{Coker}\,}%
\define\Hom{\mathrm{Hom}\,}%
\define\Ext{\mathrm{Ext}\,}%
\define\rank{\mathrm{rank}\,}%
\define\gr{\mathrm{gr}}%
\define\cHom{\Cal{Hom}}
\define\cExt{\Cal Ext\,}%

\define\cC{\Cal C}
\define\cD{\Cal D}
\define\cO{\Cal O}
\define\cS{\Cal S}
\define\cM{\Cal M}
\define\cG{\Cal G}
\define\cH{\Cal H}
\define\cE{\Cal E}
\define\cF{\Cal F}
\define\cN{\Cal N}
\define\cQ{\Cal Q}
\define\fF{\frak F}
\define\fg{\frak g}
\define\fh{\frak h}
\define\Dc{\check{D}}
\define\Ec{\check{E}}

\newcommand{\N}{{\mathbb{N}}}
\newcommand{\Q}{{\mathbb{Q}}}
\newcommand{\Z}{{\mathbb{Z}}}
\newcommand{\R}{{\mathbb{R}}}
\newcommand{\C}{{\mathbb{C}}}
\newcommand{\bN}{{\mathbb{N}}}
\newcommand{\bQ}{{\mathbb{Q}}}
\newcommand{\bF}{{\mathbb{F}}}
\newcommand{\bZ}{{\mathbb{Z}}}
\newcommand{\bP}{{\mathbb{P}}}
\newcommand{\bR}{{\mathbb{R}}}
\newcommand{\bC}{{\mathbb{C}}}
\newcommand{\bS}{{\bold{S}}}
\newcommand{\bbQ}{{\bar \mathbb{Q}}}
\newcommand{\ol}[1]{\overline{#1}}
\newcommand{\too}{\longrightarrow}
\newcommand{\respect}{\rightsquigarrow}
\newcommand{\compatible}{\leftrightsquigarrow}
\newcommand{\upc}[1]{\overset {\lower 0.3ex \hbox{${\;}_{\circ}$}}{#1}}
\newcommand{\Gmlog}{\bG_{m, \log}}
\newcommand{\Gm}{\bG_m}
\newcommand{\ep}{\varepsilon}
\newcommand{\Spec}{\operatorname{Spec}}
\newcommand{\val}{{\mathrm{val}}} 
\newcommand{\n}{\operatorname{naive}}
\newcommand{\bs}{\operatorname{\backslash}}
\newcommand{\Gal}{\operatorname{{Gal}}}
\newcommand{\gal}{{\rm {Gal}}({\bar \Q}/{\Q})}
\newcommand{\galp}{{\rm {Gal}}({\bar \Q}_p/{\Q}_p)}
\newcommand{\gall}{{\rm{Gal}}({\bar \Q}_\ell/\Q_\ell)}
\newcommand{\wep}{W({\bar \Q}_p/\Q_p)}
\newcommand{\wel}{W({\bar \Q}_\ell/\Q_\ell)}
\newcommand{\Ad}{{\rm{Ad}}}
\newcommand{\BS}{{\rm {BS}}}
\newcommand{\even}{\operatorname{even}}
\newcommand{\End}{{\rm {End}}}
\newcommand{\odd}{\operatorname{odd}}
\newcommand{\GL}{\operatorname{GL}}
\newcommand{\np}{\text{non-$p$}}
\newcommand{\g}{{\gamma}}
\newcommand{\G}{{\Gamma}}
\newcommand{\Lam}{{\Lambda}}
\newcommand{\La}{{\Lambda}}
\newcommand{\lam}{{\lambda}}
\newcommand{\la}{{\lambda}}
\newcommand{\uL}{{{\hat {L}}^{\rm {ur}}}}
\newcommand{\uQp}{{{\hat \Q}_p}^{\text{ur}}}
\newcommand{\sel}{\operatorname{Sel}}
\newcommand{\dt}{{\rm{Det}}}
\newcommand{\Sig}{\Sigma}
\newcommand{\fil}{{\rm{fil}}}
\newcommand{\SL}{{\rm{SL}}}
\newcommand{\spl}{{\rm{spl}}}
\newcommand{\st}{{\rm{st}}}
\newcommand{\Isom}{{\rm {Isom}}}
\newcommand{\Mor}{{\rm {Mor}}}
\newcommand{\bg}{\bar{g}}
\newcommand{\id}{{\rm {id}}}
\newcommand{\cone}{{\rm {cone}}}
\newcommand{\al}{a}
\newcommand{\ChL}{{\cal{C}}(\La)}
\newcommand{\Image}{{\rm {Image}}}
\newcommand{\toric}{{\operatorname{toric}}}
\newcommand{\torus}{{\operatorname{torus}}}
\newcommand{\Aut}{{\rm {Aut}}}
\newcommand{\Qp}{{\mathbb{Q}}_p}
\newcommand{\barQp}{{\mathbb{Q}}_p}
\newcommand{\Qpur}{{\mathbb{Q}}_p^{\rm {ur}}}
\newcommand{\Zp}{{\mathbb{Z}}_p}
\newcommand{\Zl}{{\mathbb{Z}}_l}
\newcommand{\Ql}{{\mathbb{Q}}_l}
\newcommand{\Qlur}{{\mathbb{Q}}_l^{\rm {ur}}}
\newcommand{\F}{{\mathbb{F}}}
\newcommand{\eps}{{\epsilon}}
\newcommand{\epsLa}{{\epsilon}_{\La}}
\newcommand{\epsLaVxi}{{\epsilon}_{\La}(V, \xi)}
\newcommand{\epsOLaVxi}{{\epsilon}_{0,\La}(V, \xi)}
\newcommand{\Qplin}{{\mathbb{Q}}_p(\mu_{l^{\infty}})}
\newcommand{\otimesQplin}{\otimes_{\Qp}{\mathbb{Q}}_p(\mu_{l^{\infty}})}
\newcommand{\galFl}{{\rm{Gal}}({\bar {\Bbb F}}_\ell/{\Bbb F}_\ell)}
\newcommand{\gallur}{{\rm{Gal}}({\bar \Q}_\ell/\Q_\ell^{\rm {ur}})}
\newcommand{\galFF}{{\rm {Gal}}(F_{\infty}/F)}
\newcommand{\galFv}{{\rm {Gal}}(\bar{F}_v/F_v)}
\newcommand{\galF}{{\rm {Gal}}(\bar{F}/F)}
\newcommand{\epsVxi}{{\epsilon}(V, \xi)}
\newcommand{\epsOVxi}{{\epsilon}_0(V, \xi)}
\newcommand{\plim}{\lim_
{\scriptstyle 
\longleftarrow \atop \scriptstyle n}}
\newcommand{\sig}{{\sigma}}
\newcommand{\ga}{{\gamma}}
\newcommand{\del}{{\delta}}
\newcommand{\Vss}{V^{\rm {ss}}}
\newcommand{\Bst}{B_{\rm {st}}}
\newcommand{\Dpst}{D_{\rm {pst}}}
\newcommand{\Dcrys}{D_{\rm {crys}}}
\newcommand{\DdR}{D_{\rm {dR}}}
\newcommand{\Fin}{F_{\infty}}
\newcommand{\Kla}{K_{\lambda}}
\newcommand{\Ola}{O_{\lambda}}
\newcommand{\Mla}{M_{\lambda}}
\newcommand{\Det}{{\rm{Det}}}
\newcommand{\Sym}{{\rm{Sym}}}
\newcommand{\LaSa}{{\La_{S^*}}}
\newcommand{\cX}{{\cal {X}}}
\newcommand{\MHG}{{\frak {M}}_H(G)}
\newcommand{\tauMla}{\tau(M_{\lambda})}
\newcommand{\Fvur}{{F_v^{\rm {ur}}}}
\newcommand{\Lie}{{\rm {Lie}}\,}
\newcommand{\cB}{{\cal {B}}}
\newcommand{\cL}{{\cal {L}}}
\newcommand{\cW}{{\cal {W}}}
\newcommand{\fq}{{\frak {q}}}
\newcommand{\cont}{{\rm {cont}}}
\newcommand{\SC}{{SC}}
\newcommand{\Om}{{\Omega}}
\newcommand{\dR}{{\rm {dR}}}
\newcommand{\crys}{{\rm {crys}}}
\newcommand{\hatSig}{{\hat{\Sigma}}}
\newcommand{\rdet}{{{\rm {det}}}}
\newcommand{\ord}{{{\rm {ord}}}}
\newcommand{\Alb}{{{\rm {Alb}}}}
\newcommand{\BdR}{{B_{\rm {dR}}}}
\newcommand{\BdRO}{{B^0_{\rm {dR}}}}
\newcommand{\Bcrys}{{B_{\rm {crys}}}}
\newcommand{\Qw}{{\mathbb{Q}}_w}
\newcommand{\barkappa}{{\bar{\kappa}}}
\newcommand{\cP}{{\Cal {P}}}
\newcommand{\cZ}{{\Cal {Z}}}
\newcommand{\oppLa}{{\Lambda^{\circ}}}
\newcommand{\bG}{{\mathbb{G}}}
\newcommand{\br}{{{\bold r}}}
\newcommand{\triv}{{\rm{triv}}}
\newcommand{\sub}{{\subset}}
\newcommand{\LD}{{D^{\star,\mild}_{\SL(2)}}}
\newcommand{\LbD}{{D^{\star}_{\SL(2)}}}
\newcommand{\dbDv}{{D^{\star}_{\SL(2),\val}}}
\newcommand{\nspl}{{{\rm nspl}}}
\newcommand{\lval}{{[\val]}}
\newcommand{\mild}{{{\rm{mild}}}}
\newcommand{\lan}{\langle}
\newcommand{\ran}{\rangle}
\newcommand{\rar}{{{\rm {rar}}}}
\newcommand{\rac}{{{\rm {rac}}}}
\newcommand{\Rep}{{\operatorname{Rep}}}
\newcommand{\LMH}{{\rm {LMH}}}
\newcommand{\MHS}{{\rm {MHS}}}
\newcommand{\red}{{{\rm red}}}
\newcommand{\cV}{{\Cal {V}}}
\newcommand{\HS}{{\rm {HS}}}
\newcommand{\all}{{\rm {all}}}
\newcommand{\bi}{{\rm {\mathbb i}}}
\newcommand{\bc}{{\rm {\mathbb c}}}
\newcommand{\zero}{{\rm {\mathbb 0}}}
\newcommand{\pa}{{\rm {W}}}

\def\ts{\textstyle\sum}
\def\tp{\textstyle\prod}
\def\t{\tilde}
\def\tsqrt{\textstyle\sqrt}
\def\vf{\varphi}
\def\Int{\operatorname{Int}}
\def\bm{\bold m}
\let\x=\times
\def\a{\alpha}
\def\b{\beta}

\begin{abstract} 
For a linear algebraic group $G$ over $\Q$, we consider the period domains $D$ classifying $G$-mixed Hodge structures, and construct the extended 
period domains $D_{\BS}$, $D_{\SL(2)}$, and $\Gamma \bs D_{\Sig}$. 
  In particular, we give toroidal partial compactifications of 
mixed Mumford--Tate domains. 
\end{abstract}

\renewcommand{\thefootnote}{\fnsymbol{footnote}}



\section*{Contents}

\noindent \S\ref{s:intro}. Introduction

\noindent \S\ref{s:D}. The period domain $D$

\ref{ss:1} Definition of the period domain $D$

\ref{ss:1.2} Relation to \cite{De0} Section 5 and \cite{De} 1.5--1.8 of Deligne

\ref{realanstrD} The real analytic structure of $D$

\ref{ss:cxanastrD} The complex analytic structure of $D$

\ref{ss:pol} Polarizability

\ref{ss:Dold} Relations with usual period domains and Mumford--Tate domains

\noindent \S\ref{s:DBS}. The space of Borel--Serre orbits

\ref{ss:cL} Real analytic manifolds with corners

\ref{ss:BSli} Borel--Serre liftings

\ref{ss:BS1} Review of Borel--Serre theory

\ref{ss:DBS} The set $D_{\BS}$

\ref{ss:BSan} The real analytic structure of $D_{\BS}$

\ref{ss:BSproperty} Global properties of $D_{\BS}$

\noindent \S\ref{s:DSL}. The space of SL(2)-orbits

\ref{ss:red} The set $D_{\SL(2)}$ when $G$ is reductive

\ref{ss:DSL} The sets $D_{\SL(2)}$ and $D^{\star}_{\SL(2)}$

\ref{SL2an0} Weight filtrations associated to $\SL(2)$-orbits

\ref{SL2an} The real analytic structures of $D_{\SL(2)}$ and $D^{\star}_{\SL(2)}$

\ref{pfan} Proofs for Section \ref{SL2an}

\ref{Fan} The fan of parabolic subgroups 

\ref{s:star_to_BS} Relation of $D^{\star}_{\SL(2)}$ and $D_{\BS}$

\ref{s:Shim} Case of Shimura varieties

\ref{ss:star_to_II} Relation of $D^{\star}_{\SL(2)}$ and $D^{II}_{\SL(2)}$

\ref{ss:val1} Valuative spaces, I

\ref{ss:SL2gl} Global properties of $D^I_{\SL(2)}$, $D^{II}_{\SL(2)}$, $D^{\star}_{\SL(2)}$ etc.

\noindent \S\ref{s:DSig}. The space of nilpotent orbits

\ref{ss:setDSig} The sets $D_{\Sig}$ and $D_{\Sig}^{\sharp}$

\ref{ss:Esig} $E_{\sig}$ and the spaces of nilpotent orbits

\ref{ss:ratio} The space of ratios

\ref{ss:val2} Valuative spaces, II

\ref{ss:CKS} Nilpotent orbits and $\SL(2)$-orbits

\ref{ss:property} Properties of the spaces of nilpotent orbits

\ref{ss:val3} Valuative spaces, III

\ref{ss:mild} Mild degeneration

\ref{ss:fund} The fundamental diagram in examples

\ref{ss:functoriality} Functoriality in $G$

\ref{ss:GLMH} $G$-log mixed Hodge structures

\ref{ss:mGLMH} Moduli of $G$-log mixed Hodge structures and period maps

\ref{ss:inf} Infinitesimal study

\ref{ss:gen1} Generalizations, I

\ref{ss:gen2} Generalizations, II

A Appendix

Correction to \cite{KNU2} Part IV

Correction to \cite{KNU3}

List of notation

References

\renewcommand{\thefootnote}{\arabic{footnote}}

\setcounter{section}{-1}
\section{Introduction}\label{s:intro}

In this paper, we generalize the theories of Mumford--Tate domains (see Green--Griffiths--Kerr's book \cite{GGK}) 
 and their toroidal partial compactifications by Kerr--Pearlstein (\cite{KP}) to mixed Hodge theory. We also construct the corresponding Borel--Serre spaces and the spaces of $\SL(2)$-orbits.
In \cite{KNU3}, we described this generalization briefly with an application to construct the toroidal partial compactifications of  higher Albanese manifolds. In this paper, we give the details.

  In Section \ref{s:D}, for a linear algebraic group $G$, we define the period domain $D$ as a space of $G$-mixed Hodge structures. Here a $G$-mixed Hodge structure means an exact $\otimes$-functor from the category of linear representations of $G$ endowed with weight filtrations 
  to the category of mixed Hodge structures. 
  In the case where $G$ is reductive, our period domain $D$ is essentially the Mumford--Tate domain studied in \cite{GGK}. 
  
  In Sections \ref{s:DBS}, \ref{s:DSL}, and \ref{s:DSig},
we construct the extended 
period domains $D_{\BS}$, $D_{\SL(2)}$, and $\Gamma \bs D_{\Sig}$, the space of Borel--Serre orbits, the space of $\SL(2)$-orbits, and the space of nilpotent orbits, respectively. 
   In the case where $G$ is reductive, 
   $\Gamma \bs D_{\Sigma}$ is essentially the space of Kerr--Pearlstein (\cite{KP}). 

   Our method to construct $D_{\BS}$, $D_{\SL(2)}$, and $\Gamma \bs D_{\Sig}$ is similar to that for the usual period domain developed in \cite{KU2} and in the 
preceding parts of this series of papers \cite{KNU2}. 
  Actually, the present paper is regarded as a generalization of our pervious works in Parts  I--IV of \cite{KNU2}.
  We prove similar results (the Hausdorffness etc.) for these spaces to what we have proved so far for the corresponding spaces in the case of usual period domain. We obtain a fundamental diagram 
\medskip
\newline
\setbox0=\vbox{\parindent=0pt
$
{
\begin{matrix}
&&&&&&&&&& D^{\star}_{\SL(2),\val}&&\overset {\eta^{\star}}
\longrightarrow &&D_{\BS,\val}\\
&&&&&&&&&\swarrow&\downarrow &\searrow&&& \downarrow\\
&&&&D^{\sharp}_{\Sigma,\lval}&\overset{\psi}\to & D^I_{\SL(2),\val}&\to& D^{II}_{\SL(2),\val}&&D^{\star,+}_{\SL(2)}&&D^{\star,\pa}_{\SL(2)}&\to&D_{\BS}\\
&&&&\downarrow &&\downarrow&&\downarrow&\swarrow&&\searrow&\downarrow&&\\
 \Gamma \bs D_{\Sig,\val}& \leftarrow & D_{\Sig,\val}^{\sharp}&\to &D^{\sharp}_{\Sig,[:]} & \overset {\psi} \to &D^I_{\SL(2)}&\to& D^{II}_{\SL(2)}&&&&D^{\star}_{\SL(2)}&&\\
\downarrow &&&&\downarrow&&&&&&&&&&\\
 \Gamma \bs D_{\Sig}&&\longleftarrow &&D_{\Sig}^{\sharp}&&&&&&&&&&
\end{matrix}
}
$}
\noindent\scalebox{.85}[.85]{\box0}
\medskip
\newline 
of the same style as the fundamental diagram in Part IV of \cite{KNU2}, which is commutative and in which the maps respect the structures of the spaces. 
The good properties of these spaces are proved starting from the spaces on the right-hand-side and then moving to the left-hand-side, similarly as in \cite{KU2} and \cite{KNU2}. 
  \medskip

  The authors thank Teruhisa Koshikawa for helpful discussions. 
  The authors thank Katsutoshi Shinohara for a help to complete Example 2 in Remark \ref{r:actonD}. 
  K.\ Kato was 
partially supported by NFS grants DMS 1001729, DMS 1303421, DMS 1601861, and DMS 2001182.
C.\ Nakayama was 
partially supported by JSPS Grants-in-Aid for Scientific Research (C) 22540011, (B) 23340008, (C) 16K05093, and (C) 21K03199.
S.\ Usui was 
partially supported by JSPS Grants-in-Aid for Scientific Research (B) 23340008 and (C) 17K05200.

\medskip

\section{The period domain $D$}\label{s:D}

  Let $G$ be a linear algebraic group over $\Q$.
  Let $G_u$ be the unipotent radical of $G$, and let $G_{\red}=G/G_u$ be the reductive quotient  of $G$. 
  Let $\mathrm{Rep}(G)$ be the category of 
finite-dimensional linear representations of $G$ over $\Q$.

\subsection{Definition of the period domain $D$}\label{ss:1}

\begin{sbpara} 

In this Section \ref{ss:1}, we give the definition of our period domain $D$ considered in this paper.  As we will see in Section 1.2, this is the period domain $D$ in our previous work \cite{KNU3}. In this Section 1.1, we introduce  $D$ from the following point of view: We want to have a period domain endowed with an action of $G$. More precisely, we want to have a complex analytic period domain $D$ for $\Q$-mixed Hodge structures which is endowed with a real analytic action of the subgroup $G(\R)G_u(\C)$ of $G(\C)$. Then we naturally come to the definition of our $D$.

\end{sbpara}

\begin{sbpara}\label{GMHS} 
For a subfield $E$ of $\R$, let  $(E\MHS)$  be the category of $E$-mixed Hodge structures. 

By a {\it $G$-mixed Hodge structure} ($G$-MHS for short),  we mean an exact $\otimes$-functor   (\cite{De3} 2.7) $H:\Rep(G) \to (\Q\MHS);V\mapsto(H(V)_\Q,W,F)$.  (This definition of $G$-MHS is slightly different from the definition in our previous work \cite{KNU3}.)

\end{sbpara}

\begin{sbpara}
\label{movact} Let $D_{\all}(G)$ be the set of all isomorphism classes of $G$-MHS  $H:\Rep(G)\to(\bQ\MHS)$ which preserves the underlying $\Q$-vector spaces (that is, for any $V\in \Rep(G)$, $H(V)_\Q=V$ as a $\Q$-vector space) and 
which satisfy the following conditions (i) and (ii). 

(i) For every $V\in \Rep(G)$, the weight filtration $W$ of $H(V)$ is stabilized by the action of $G$ on $V$.

(ii) For every $V\in \Rep(G)$, $G_u$ acts trivially on $\gr^WV$.

Then the subgroup $G(\R)G_u(\C)$ of $G(\C)$ acts on $D_{\all}(G)$: 
  For $H\in D_{\all}(G)$ and $g\in G(\R)G_u(\C)$, $gH\in D_{\all}(G)$ is defined as follows. For $V\in \Rep(G)$, $(gH(V))_\Q=V$, the weight filtration $W$ of $gH(V)$ is that of $H(V)$, and the Hodge filtration $F$ of $gH(V)$ is given by $F^p(gH(V)):= gF^pH(V)\subset V_\C:=\C\otimes_\Q V$, where $F^pH(V)$ is the Hodge filtration of $H(V)$. 
\end{sbpara}

\begin{sbpara}\label{perd} By a {\it period domain of $G$-MHS}, we mean a $G(\R)G_u(\C)$-orbit in $D_{\all}(G)$. 

This is the definition of the period domain of this paper, which we denote by $D$. 

\end{sbpara}

\begin{sbpara}
As will be shown in \ref{Dalltop}, for a reasonable topology on $D_{\all}(G)$, each period domain is open and closed in $D_{\all}(G)$ and hence as a topological space, $D_{\all}(G)$ is a disjoint union of all period domains in $D_{\all}(G)$. 

\end{sbpara}

\begin{sbrem} We think this definition of $D$ is a   natural way to have a period domain for $\Q$-MHS with an action of $G(\R)G_u(\C)$. 

To consider a $G$-MHS is a natural way to have an action of $G$ on the period domain because each $V\in \Rep(G)$ has an action of $G$. 

The group $G(\R)$ acts on $D_{\all}(G)$ by the condition (i) in \ref{movact}. We remark how this condition (i) is important. Note that if $(V, W, F)$ is a $\Q$-MHS with $V\in \Rep(G)$, $W$ the weight filtration, and $F$  the Hodge filtration, then for $g\in G(\R)$,  we have an $\R$-MHS $(V_\R, gW_\R, gF)$. But  if $gW_\R$ is rational for all $g\in G(\R)$ (we need this to have the action of $G(\R)$ on the set of $\Q$-MHS), $W$ must be stabilized by the connected component of $G(\R)$ containing $1$. (That is, rational weight filtrations can not move continuously.) Hence it is natural to put the condition that $G$ stabilizes $W$.

We have the action of $G_u(\C)$ on $D_{\all}(G)$ by the condition (ii) in \ref{movact}. 

\end{sbrem}

\begin{sbrem} The condition $G$ stabilizes $W$ may be natural also for the following reason. In the definition of a period domain of mixed Hodge structures, it is natural to fix the weight filtration and move the Hodge filtration. This is because
in a variation of mixed Hodge structure, the weight filtration is put on the local system and the Hodge filtration is put on the vector bundle,  and hence the weight filtration is  constant locally and Hodge filtration varies continuously. 

\end{sbrem}

\begin{sbpara}\label{GHS}
For a subfield $E$ of $\R$, 
by an $E$-HS, we will mean an $E$-MHS which is a direct sum of pure $E$-Hodge structures. (They are also called {\it split} $E$-mixed Hodge structure.) Let $(E\HS)$ be the category of $E$-HS. 
By $G$-HS, we mean a $G$-MHS which has values in  ($\Q$HS)$\subset$ ($\Q\MHS$).

\end{sbpara}

\begin{sbprop}\label{movred}
  Assume that $G$ is reductive and let $H\in D_{\all}(G)$. Then $H$ is a $G$-HS. For $V\in \Rep(G)$ and for $w\in \Z$, if  $V_w$ denotes the underlying $\Q$-vector space of the part of $H(V)$ of weight $w$, $V_w$ is $G$-stable in $V$.
\end{sbprop}

\begin{pf}
Since  any representation of a reductive group is  semisimple, the weight filtration on $V$ has a splitting which is compatible with the action of $G$. Hence there is a decomposition $V=\bigoplus_w V_w$ such that $V_w$ is $G$-stable and $W_wV_w=V_w$, $W_{w-1}V_w=0$. We have $H(V)=\bigoplus_w H(V_w)$ and $H(V_w)$ is a pure Hodge structure of weight $w$.
\end{pf}

\subsection{Relation to \cite{De0} Section 5 and \cite{De} 1.5--1.8 of Deligne}\label{ss:1.2}

\begin{sbpara} \label{SCR0}
  Let $S_{\C/\R}$ be the Weil restriction of ${\mathbb G}_m$ from $\C$ to $\R$. 
    It represents the functor $A\mapsto (\C\otimes_\R A)^\times$ for commutative rings $A$ over $\R$.
  We have  $S_{\C/\R}(\R)=\C^\times$, which is understood as $\C^\times$ regarded as an algebraic group over $\R$.

 Let $w: {\mathbb G}_{m,\R}\to S_{\C/\R}$ be the homomorphism induced from the natural maps $A^\times\to (\C\otimes_\R A)^\times$ 
 for commutative rings $A$ over $\R$.
 
 In \cite{De0} and \cite{De}, Deligne related $S_{\C/\R}$ to the theory of Hodge structures as follows. 
\end{sbpara}

\begin{sbpara}\label{SCR1} The following (1) and (2) are identified. 

(1) An $\R$-HS.

(2)  A (finite-dimensional) linear representation of $S_{\C/\R}$ over $\R$. 

In fact, a  linear representation of $S_{\C/\R}$ over $\R$ is equivalent to a finite-dimensional $\R$-vector space $V$ endowed with a decomposition
$$V_\C:=\C\otimes_\R V=\bigoplus_{p,q\in \Z} \; V_\C^{p,q}$$
such that for any $p,q$, $V_\C^{q,p}$ coincides with the complex conjugate of $V_{\C}^{p,q}$ (that is, the image of $V_\C^{p,q}$ under $\C\otimes_\R V\to \C\otimes_\R V\;;\;z\otimes v\mapsto {\bar z}\otimes v$). For a linear representation $V$ of $S_{\C/\R}$, the corresponding decomposition is defined by 
$$V_\C^{p,q}=\{v\in V_\C\;|\; [z] v= z^p{\bar z}^qv\;\text{for}\;z\in \C^\times\}.$$
Here $[z]$ denotes $z$ regarded as an element of $S_{\C/\R}(\R)$. By taking $V_\C^{p,q}$ as the $(p,q)$-Hodge component of $V_\C$, this is understood as an $\R$-HS with underlying $\R$-vector space $V$. 

The weight $m$-part of this $\R$-HS is the part on which $z\in \R^\times={\mathbb G}_{m,\R}(\R)$ acts as multiplication by $z^m$ via $w: {\mathbb G}_{m,\R} \to S_{\C/\R}$. 
\end{sbpara}

\begin{sbpara}\label{SCR2}

(1) The following (1.1) and (1.2) are identified.

(1.1) An exact $\otimes$-functor $H:\Rep(G) \to (\R\HS)$ 
 with the underlying $\R$-vector space $H(V)_\R=\R \otimes_\Q V$ ($V\in \Rep(G)$). 

(1.2) A homomorphism $h: S_{\C/\R} \to G_\R$.

\medskip

(2) The following (2.1) and (2.2) are identified.

(2.1)
A $G$-HS $H:\Rep(G) \to (\Q\HS)$ (\ref{GHS}) preserving the underlying $\Q$-vector spaces.

(2.2) A homomorphism $h:S_{\C/\R}\to G_\R$
 such that $h\circ w:{\mathbb G}_{m,\R}\to G_\R$  is defined over $\Q$.

We explain (1).
If we are given a homomorphism $S_{\C/\R}\to G_\R$ in (1.2), for any $V\in \Rep(G)$, $H(V)_\R=\R \otimes_{\Q} V$ has an action of $S_{\C/\R}$ and is regarded as an $\R$-HS. 
Conversely, if we have a functor $H:\Rep(G) \to (\R\HS)$ having the property in (1.1), then  by \ref{SCR1} and by the theory of Tannakian categories, we have a homomorphism $S_{\C/\R}\to G_\R$.  This shows (1). 

Taking account of $\Q$-structures, we get (2) from (1).

\end{sbpara}

\begin{sbpara}\label{hatDa}

Let ${\hat D}_{\all}(G)$ be the set of all isomorphism classes of exact $\otimes$-functors $H:\Rep(G) \to (\R \MHS)$ with the underlying $\R$-vector space $H(V)_\R=\R \otimes_\Q V$ ($V\in \Rep(G)$) such that for every $V \in \Rep(G)$, the action of $G$ on $V$ stabilizes the weight filtration $W$ on $V_\R=\R\otimes_\Q V$ and the action of $G_u$ on $\gr^WV_\R$ is trivial.
(Note that $D_{\all}(G)$ in \ref{movact} is identified with the set of those functors whose weight filtrations are $\Q$-rational.)

Let ${\hat D}'_{\all}(G)$ be  the set of all pairs $(h, \delta)$, where $h$ is a homomorphism $S_{\C/\R}\to G_\R$ and $\delta$ is an element of 
$\Lie(G)_\R$ satisfying the following condition (i). Consider the adjoint action of $G$ on $\Lie(G)$ and the induced action of $S_{\C/\R}$ on $\Lie(G)_{\bR}$ via $h$.  By this, $\Lie(G)_{\bR}$ becomes an $\R$-HS. For $m,n\in \Z$, let $H^{m,n}$ be the Hodge $(m,n)$-component of  the weight $(m+n)$-part of $\Lie(G)_{\bC}$. 

(i) $\delta \in \bigoplus_{m<0,n<0} H^{m,n}$.

 \end{sbpara}

\begin{sbprop}\label{SCR4}  We have a bijection
$${\hat D}'_{\all}(G)\to {\hat D}_{\all}(G)\;;\;(h,\delta)\mapsto e^{i\delta}H,$$ where  $H$ is the functor $\Rep(G) \to (\R\HS)$  corresponding to $h$ in $\ref{SCR2}$ $(1)$.
\end{sbprop}

\begin{pf}
  This follows from  \cite{CKS} Proposition (2.20). 
\end{pf}

\begin{sblem}\label{Lie}

Let $H\in {\hat D}_{\all}(G)$ ($\ref{hatDa}$). Consider the adjoint action of $G$ on $\Lie(G)$ and consider the weight filtration of $\Lie(G)_{\bR}$ and the Hodge filtration of $\Lie(G)_{\bC}$ defined by $H$. 

$(1)$ Let $\ell\in \Lie(G)_{\bR}$ and $w\in \Z$. Then  $\ell\in W_w\Lie(G)_{\bR}$ if and only if, for any $V\in \Rep(G)$ and any $k\in \Z$, the Lie action $\Lie(G)_{\bR} \times V_\R \to V_\R$ satisfies $\ell W_kV_\R\subset W_{k+w}V_\R$. 

$(2)$ Let  $\ell\in \Lie(G)_{\bC}$ and $p\in \Z$. Then $\ell\in F^p\Lie(G)_{\bC}$ if and only if, for any $V\in \Rep(G)$ and any $r\in \Z$, the Lie action $\Lie(G)_{\bC} \times V_\C \to V_\C$ satisfies $\ell F^rV_\C\subset F^{r+p}V_\C$.
\end{sblem}

\begin{pf} 
  The only if parts are clear. 

  Take a $V\in \Rep(G)$ such that the map $\Lie(G) \to \End(V)$ induced by the Lie action is injective. 
  Then since this is a homomorphism of $\R$-MHS, we have $W_w\Lie(G)_{\bR} = \Lie(G)_{\bR} \cap W_w\End(V_\R) = \{h\in\Lie(G)_{\bR}\;|\; hW_kV_\R\subset W_{k+w}V_\R$ for any $k\in \Z\}$, and a similar relation of Hodge filtrations. This proves the if parts. 
\end{pf}

\begin{sbpara}\label{kandh} (1) Let $\Psi_W(G)$ be the set of all homomorphisms $k: \bG_m \to G_{\red}$ satisfying the following conditions (i) and (ii).

(i) The image of $k$ is contained in the center of $G_{\red}$.

(ii) For some  homomorphism (and hence for any homomorphism) $\tilde k : \bG_m\to G$ which lifts $k$, the adjoint action of $\bG_m$ on $\Lie(G_u)$ via $\tilde k$ is of weight $\le-1$.

\medskip

Note that in the case where $G$ is reductive, the condition (ii) is automatically satisfied. 

\medskip

(2) Let $\Psi_H(G)$ be the set of all homomorphisms $h: S_{\C/\R}\to G_{\red,\R}$ satisfying the following condition: There exists $k\in \Psi_W(G)$ such that $h\circ w: \bG_{m,\R}\to G_{\red,\R}$ comes from $k$ (see \ref{SCR0} for $w$).

We have a canonical map $\Psi_H(G) \to \Psi_W(G)\;;\;h\mapsto k$. 
\end{sbpara}

\begin{sbpara}
\label{wt} 
Let $k\in \Psi_W(G)$.

  Then, for any $V \in \Rep(G)$, the action of ${\mathbb G}_m$  on $V$ via a lifting $\tilde k$ of $k$  defines 
a  rational increasing filtration $W$ on $V$ called the {\it weight filtration}, which is independent of the lifting. For any $V\in \Rep(G)$, $W_wV$ is $G$-stable for any $w\in \Z$ and  the action of $G_u$ in $\gr^WV$ is trivial. 

\end{sbpara}

\begin{sbprop}
\label{DPsi1}

$(1)$ Let $H\in {\hat D}_{\all}(G)$ ($\ref{hatDa}$). 
Then the  following conditions {\rm (i)} and {\rm (ii)} are equivalent.

{\rm (i)} $H\in D_{\all}(G)$ ($\ref{movact}$).

{\rm (ii)} The restriction of $H$ to $\Rep(G_{\red})$ is a $G$-HS ($\ref{GHS}$), and the corresponding homomorphism $h: S_{\C/\R}\to G_{\red, \R}$ ($\ref{SCR2}$ $(2)$) belongs to $\Psi_H(G)$ ($\ref{kandh}$ $(2)$). 

$(2)$ Let $H\in {\hat D}_{\all}(G)$ ($\ref{hatDa}$). 
If the equivalent conditions in $(1)$ are satisfied,  the weight filtrations of $H$ are given by the image $k$ of $h$ in $\Psi_W(G)$ as in $\ref{wt}$.

$(3)$ If $G$ is reductive, the map $D_{\all}(G) \to \Psi_H(G)\;;\;H\mapsto h$ is bijective. 
\end{sbprop}

\begin{pf}
 In the case where $G$ is reductive, (1) follows from Proposition \ref{movred}. 
 
 (3) follows from the case of (1) where $G$ is reductive. 

 We prove (1) and (2).

  Assume (i). 
  Then, by Proposition \ref{movred}, the restriction is a $G$-HS. 
  Further, the induced homomorphism $h\circ w: \bG_{m,\R}\to G_{\red,\R}$ is defined over $\Q$ and its image is in the center of $G_{\red,\R}$. We prove that $\Lie(G_u)$ is of weight $\leq -1$, which implies (ii). 
  Since $G_u$ acts on $\gr^WV$ trivially  for any $V\in \Rep(G)$, the Lie action of $\Lie(G_u)$ on $V$ 
    induces the zero action on $\gr^WV$. By Lemma \ref{Lie}, this proves that $\Lie(G_u)$ has weights $\leq -1$. 
    
  Assume (ii). Since $H$ comes from ${\hat D}'_{\all}(G)$ (Proposition \ref{SCR4}),  the weight filtration $W$ of $H$ is  given by the image $k$ of $h$ in $\Psi_W(G)$ as in \ref{wt}.  Hence the weight filtration of $H$ is rational, that is, $H\in D_{\all}(G)$.
  Thus we have proved (1) and (2).
\end{pf}

By Proposition \ref{DPsi1} (1), we have a map $$D_{\all}(G)\to \Psi_H(G)\;;\; H \mapsto h.$$

\begin{sbcor}\label{ucom}
If $D_{\all}(G)$ is not empty, $G_u$ is contained in the commutator subgroup of $G$. 
\end{sbcor}

\begin{pf} If $D_{\all}(G)$ is not empty, via the maps $D_{\all}(G)\to \Psi_H(G) \to \Psi_W(G)$, we see that $\Psi_W(G)$ is not empty. For a lifting $\tilde k :\bG_m \to G$ of an element $k:\bG_m \to G_{\red}$ of $\Psi_W(G)$, since $\Lie(G_u)$ is of weight $\le -1$ for  $\tilde k$, we have $G_u=[\tilde k(\bG_m), G_u]$. 
\end{pf}

\begin{sbprop}
\label{DPsi}
  Consider the maps $D_{\all}(G)\to \Psi_H(G)\to \Psi_W(G)$. 

$(1)$ Let $H\in D_{\all}(G)$ and let $h$ be the image of $H$ in $\Psi_H(G)$. Then 
for $g_1\in G(\R)$ and $g_2\in G_u(\C)$, the image of $g_1g_2H\in D_{\all}(G)$ in $\Psi_H(G)$ is $g_{1,\red}hg_{1,\red}^{-1}: S_{\C/\R}\to G_{\red,\R}$.
Here $g_{1,\red}$ is the image of $g_1$ in $G_{\red}(\R)$.

$(2)$ Let $h\in \Psi_H(G)$ and let $k$ be the image of $h$ in $\Psi_W(G)$. Then for $g\in G_{\red}(\R)$, the image of $ghg^{-1}$ in $\Psi_W(G)$ is $k$.

$(3)$ The map $D_{\all}(G) \to \Psi_H(G)$ is surjective. The action of $G_u(\C)$ on each fiber of this  map is transitive.
\end{sbprop}

\begin{pf}
(1) and (2) are straightforwards. 

 We prove  
(3).
Let $h\in \Psi_H(G)$.  
  By Proposition \ref{SCR4}, the fiber of $h$ in $D_{\all}(G)$ consists of the images $e^{i \delta}H$ of $(\tilde h, \delta)\in {\hat D}'_{\all}(G)$ in ${\hat D}_{\all}(G)$, where $\tilde h: S_{\C/\R}\to G_\R$ are liftings of $h$ and $H$ denotes the functor $\Rep(G) \to (\R\HS)$ corresponding to $\tilde h$.  The fiber is not empty because a  lifting $\tilde h$ of $h$ exists. 
Let $(\tilde h, \delta), (\tilde h', \delta')\in {\hat D}'_{\all}(G)$, where $\tilde h, \tilde h'$ are liftings of $h$ and let $H$ and $H'$ be the functors corresponding to $\tilde h$ and $\tilde h'$, respectively. For $V\in \Rep(G)$, the weight filtration $W$ on $V$ given by $e^{i\delta}H$ and that given by $e^{i\delta'}H'$ coincide because they are given by  the image of $h$ in $\Psi_W(G)$ as in \ref{wt}. Let $s_V: \gr^WV_\R \to  V_\R$ and $s'_V: \gr^WV_\R \to V_\R$ be the splittings of $W$ over $\R$ defined by $\tilde h\circ w$ and $\tilde h' \circ w$, respectively. Then $(s'_V \circ s_V^{-1}: V \to V)_V$ comes from an element $u$ of $G_u(\R)$ and $\tilde h'= u\tilde h u^{-1}$. We have  $H'=uH$ and  hence $e^{i\delta'}H'=ge^{i\delta}H$, where $g= e^{i\delta'}ue^{-i\delta}\in G_u(\C)$. 
\end{pf} 

We will prove a more precise result on the fibers in Theorem \ref{fiber} below. 

\begin{sbpara}\label{kDall}
Fix $k\in \Psi_W(G)$. Then the inverse image of $k$ in $D_{\all}(G)$ under $D_{\all}(G)\to \Psi_W(G)$ is identified with the set of all $G$-MHS which preserve the underlying $\Q$-vector spaces and whose weight filtrations are given by $k$ as in \ref{wt}.

\end{sbpara}

\begin{sbpara}\label{D}  Fix a homomorphism  $$h_0: S_{\C/\R}\to G_{\red,\R}$$
which belongs to $\Psi_H(G)$ (\ref{kandh}). 
Then, by Proposition \ref{DPsi}, there is a unique $G(\R)G_u(\bC)$-orbit $D(G, h_0)$ in $D_{\all}(G)$ whose image in $\Psi_H(G)$ is the set of all $G(\R)$-conjugates of $h_0$. We call this period domain of $G$-MHS  (\ref{perd}) the 
{\it period domain associated to $h_0$}. By \ref{kDall},   $D(G, h_0)$ coincides with the period domain associated to $h_0$ defined in our previous paper \cite{KNU3}.  

In the case where $G$ is reductive, $D(G, h_0)$ is identified with the set of all $G(\R)$-conjugates of $h_0$ (Proposition \ref{DPsi1} (3)). 
In this case, the definition of the period domain $D(G, h_0)$ as the set of $G(\R)$-conjugates of $h_0$ appears in Section 5 of \cite{De0} and 
 1.5 of \cite{De} of Deligne. We borrow the notation $h_0$ from \cite{De}.

\end{sbpara}

 \begin{sbpara}\label{D2} If $G$ is reductive and $h_0$ in \ref{D} satisfies the Shimura data in 1.5 of \cite{De}, then as in 1.8 of ibid., 

\medskip

(1)  $G(\Q) \bs (D(G, h_0) \times G({\bf A}_\Q^{\infty}))/K$, where ${\bf A}_\Q^{\infty}$  is the adele ring of $\Q$ without $\infty$-component and $K$ is an open compact subgroup of $G({\bf A}^{\infty}_\Q)$,

\medskip
\noindent
is a {\it Shimura variety over $\C$ associated to $h_0$}.

We expect that the set (1) in general (it is a complex analytic space as in Section \ref{ss:cxanastrD} below) is also important in number theory. \cite{Ka} is a trial of the study in this direction.

 \end{sbpara}

 We show that each period domain is isolated in $D_{\all}(G)$. 
 
 \begin{sbprop}\label{Dalltop}
 $(1)$ If $G$ is reductive,   each $G(\R)$-orbit  is open in $D_{\all}(G)$, where we endow $D_{\all}(G)$ with the topology induced  by the compact-open topology of   $\Hom_{\mathrm{cont}}(\C^\times, G(\R))$ via the injection $D_{\all}(G)\simeq \Psi_H(G)\to \Hom_{\mathrm{cont}}(\C^\times, G(\R))$.
 
$(2)$ For a general $G$, for any topology of $D_{\all}(G)$ such that the map $D_{\all}(G) \to D_{\all}(G_{\red})$ is continuous for the above topology of $D_{\all}(G_{\red})$, each $G(\R)G_u(\C)$-orbit is open in $D_{\all}(G)$. That is, as a topological space, $D_{\all}(G)$ is the disjoint union of period domains. 
 
 \end{sbprop}
 
 \begin{pf} (1) The set of homomorphisms from $\bG_{m,\R}$ to the center of $G_\R$ is discrete for the topology induced by the compact-open topology of $\Hom_{\text{cont}}(\R^\times, G(\R))$. Hence we are reduced to proving that for $S^{(1)}_{\C/\R}:= \text{Ker}(\text{norm}: S_{\C/\R}\to \bG_{m, \R})$, every $G(\R)$-conjugacy class in $\Hom(S^{(1)}_{\C/\R}, G_\R)$ is open for the topology induced by the compact-open topology of $\Hom_{\text{cont}}(S^{(1)}_{\C/\R}(\R), G(\R))$. But this follows from the case $K=S_{\C/\R}^{(1)}(\R)=\{z\in \C^\times\;|\;|z|=1\}$ and $L= G(\R)$ of the result of Lee and Wu \cite{LW} that for a compact group $K$ and for a locally compact group $L$, each $L$-conjugacy class in $\Hom_{\text{cont}}(K, L)$ is open for the compact-open topology.

   (2) follows from (1)  because $D$ is the inverse image of its image $D_{\red}$ in $D_{\all}(G_{\red})$ and $D_{\red}$ is open in $D_{\all}(G_{\red})$ by (1). 
 \end{pf}
 
\begin{sbpara} We compare the above period domain with the Griffiths period domain \cite{Gr} and its generalization \cite{U} to MHS. (A more precise comparison is given in Section \ref{ss:Dold}.)

The period domain in \cite{U} classifies MHS with a fixed weight filtration and fixed Hodge numbers of each $\gr^W_w$, and there the Hodge filtrations move. 

In the definition of the period domain of $G$-MHS in the present paper, fixing $W$ in ibid.\ corresponds to fixing $k_0\in \Psi_W(G)$ as in \ref{wt}, and fixing Hodge numbers of each $\gr^W_w$ in ibid.\ corresponds to the fact that we fix the $G_{\red}(\C)$-conjugacy class of the composition $\bG_{m,\C}\to \bG_{m,\C}\times \bG_{m,\C}=S_{\C/\R, \C} \overset{h_0}\to G_{\red,\C}$, where the first arrow is $z\mapsto (z,1)$. Moving the Hodge filtration in ibid.\ corresponds to moving $H\in D$ by $G(\R)G_u(\C)$.

\end{sbpara}

\subsection{The real analytic structure of $D$} 
\label{realanstrD}

Let $h_0: S_{\C/\R}\to G_{\red,\R}$ be as in \ref{D} and let $D=D(G, h_0)$ be the associated period domain (\ref{D}). We consider the real analytic structure of $D$.

\begin{sbpara}\label{ranD}
$D$ is regarded as $(G(\R)G_u(\C))/I_x$ for the isotropy subgroup $I_x$ of $G(\R)G_u(\C)$ at $x\in D$. 
  Furthermore, $I_x$ is a real algebraic subgroup of $G(\R)G_u(\C)$. This gives a real analytic structure on $D$, and it is independent of the choice of $x\in D$.

\end{sbpara}

\begin{sbpara}\label{splW0}
The image of the composite map $D\to \Psi_H(G) \to \Psi_W(G)$ is a one point $k_0\in \Psi_W(G)$.  Hence the weight filtration $W$ on $V\in \Rep(G)$ given by $x\in D$ is independent of $x$ and it is defined by $k_0$ as in \ref{wt}.

 Let $\spl(W)$ be the set of all isomorphisms of $\otimes$-functors from $\Rep(G)$ to the category of $\R$-vector spaces
$$(V\mapsto \gr^WV_\R) \overset{\sim}\to (V \mapsto V_\R)\quad\text{preserving the weight filtrations}.$$
(The notation  $\spl(W_\R)$ may be better because we are considering splittings over $\R$, but we use the notation $\spl(W)$ for simplicity and for the compatibility with our notation in \cite{KNU2} Parts I--IV.)  Then  $\spl(W)$ is a $G_u(\R)$-torsor.  Hence it is regarded as a real analytic manifold. 

\end{sbpara}

\begin{sbpara}\label{Dred} 
Let $D_{\red}=D(G_{\red}, h_0)$. This set is identified with the set of all $G(\R)$-conjugates of $h_0$ in $\Hom(S_{\C/\R}, G_{\red,\R})$.
We have a canonical surjective projection $D\to D_{\red}$. 
For $x\in D$, let $x_{\red}$ be the image of $x$ in $D_{\red}$. 

\end{sbpara}

\begin{sbpara}\label{CKSdelta}
Let $\cL=\bigoplus_{w\leq -2} \gr^W_w\Lie(G)_\R$. 

Consider the adjoint action of $G$ on $\Lie(G)$. 
Since the weight filtration is stable under the action of $G$ and since $G_u$ acts trivially on 
 $\gr^W\Lie(G)$ (\ref{wt}), $G_{\red}$ acts on $\gr^W\Lie(G)$.  
For $p\in D_\red$,  define $\cL(p)= \gr^W\Lie(G)_{\bR} \cap (\bigoplus_{m<0, n<0} H^{m,n})\subset \cL$, where $H^{m,n}$ denotes the Hodge $(m,n)$-component of  $\gr^W_{m+n} \Lie(G)_{\bC}$ with respect to $p$. 
 
 By Proposition \ref{SCR4}, any element $x$ of $D$ is written uniquely as  $$x=s(e^{i\delta}p): V\mapsto s(e^{i\delta}p(\gr^W(V))),$$ where $p\in D_{\red}$, $s\in \spl(W)$, $\delta\in \cL(p)$. In fact, this is the understanding of $x$ as the image of $(\tilde h, s\delta s^{-1})\in {\hat D}'_{\all}(G)$ in ${\hat D}_{\all}(G)$, where $\tilde h: S_{\C/\R}\to G_\R$  is the lifting of the homomorphism $h:S_{\C/\R}\to G_{\red,\R}$, corresponding to $p$, defined as $\tilde h(z)= sh(z)s^{-1}$ ($z\in S_{\C/\R}$). 
 
  We denote this $\delta\in \cL(p)$ by $\delta(x)$. 

\end{sbpara}

 \begin{sbpara}\label{splW} 
We have a canonical real analytic map 
$$\spl_W:D\to\spl(W)\;;\;x\mapsto\spl_W(x)$$ which is a modification $\spl_W(x)=s\circ\exp(\zeta)$ 
of the real analytic map $D\to \spl(W)\;;\;x\mapsto s$ in \ref{CKSdelta}
by an element $\zeta\in \cL(p)$ explained below.
This splitting $\spl_W(x)$ is called the {\it canonical splitting of $W$ at $x$}.

We explain $\zeta\in\cL(p)$.
For each $V\in \Rep(G)$, $\zeta_V\in \cL_{p(V)}$ is defined as a universal Lie polynomial of the Hodge $(j,k)$-components $\delta_{V, j, k}$ of $\delta_V$ (\cite{CKS} (3.28), (6.60); see also \cite{KNU1} 1.4, Appendix, \cite{KNU2} Part II 1.2). 
We can show $\zeta_{V \otimes V'}= \zeta_V \otimes \text{id}_{V'}+\text{id}_V \otimes \zeta_{V'}$, and we have $\zeta\in \cL(p)$.

\end{sbpara}

\begin{sbprop}\label{Dandgr}
We have a canonical  isomorphism of real analytic manifolds 
$$D\overset{\sim}\to  \{(p, s, \delta)\in D_{\red} \times \spl(W)\times \cL\;|\; \delta\in \cL(p)\}$$ 
given by $x\mapsto (x_{\red}, \spl_W(x), \delta(x))$.

\end{sbprop}

We have an isomorphism of the same form  even if we replace the map $\spl_W: D\to \spl(W)$ by the map $D\to \spl(W)\;;\;x\mapsto s$ of \ref{CKSdelta}, but the isomorphism in Proposition \ref{Dandgr} behaves better in degeneration (see Remark \ref{splrem1} below).

\begin{sbpara}\label{Lp=Lp'}
If $p, p'\in D_{\red}$, $p'=gp$ for some $g\in G_{\red}(\R)$, and we have $\cL(p')= \Ad(g)\cL(p)$. Hence all $\cL(p)$ ($p\in D_{\red}$) are isomorphic as  graded $\R$-vector spaces. 

Let $L=\cL(p)$ for some $p\in D_{\red}$. 

By Proposition \ref{Dandgr}, we have

\end{sbpara}

  \begin{sbcor}
  $D$ is an $L$-bundle over $D_{\red}\times \spl(W)$.
  
  \end{sbcor}

\begin{sbpara}\label{Dspl}
Let $D_{\spl}$ be the part of $D$ consisting of exact $\otimes$-functors $\Rep(G)\to(\Q\MHS)$ such that the image of the composition $\Rep(G) \to(\Q\MHS)\to(\R\MHS)$ is contained in $(\R\HS)$ (that is, such that the images are $\R$-split mixed Hodge structures). Then $$D_{\spl}= \{s(p)\;|\; s\in \spl(W), p\in D_{\red}\}=\{x\in D\;|\; \delta(x)=0\}$$
and $D_{\spl}$ is a closed real analytic submanifold of $D$. 
  Here $s(p): \Rep(G) \to (\Q\MHS)$ sends a $V \in \Rep(G)$ to the $\Q$-MHS on the underlying $\Q$-vector space of $V$ which is induced by $p(\gr^WV)$ and $s(V)\colon \gr^WV_{\bR} \simeq V_{\bR}$. 
  Let $D_{\nspl}= D\smallsetminus D_{\spl}$. 
\end{sbpara}

In the rest of this Section \ref{realanstrD}, we explain how $\spl_W: D\to \spl(W)$ is important and why we prefer this map to the map $D\to \spl(W)\;;\;x \mapsto s$ in \ref{CKSdelta}.

\begin{sbrem}\label{splrem1}

When we consider degeneration along a nilpotent orbit, the canonical splitting $\spl_W:D\to\spl(W)$ behaves better than the splitting $D\to\spl(W)$ in \ref{CKSdelta}.
We explain this.
Assume that $(N_1,\dots,N_n,F)$ generates a  nilpotent orbit as in \ref{nilp2} below.
For $y=(y_j)_{1\le j\le n}\in\R^n$ with $y_j$ are sufficiently large for all $j$, let $\spl_W(y)$ and $s(y)$ be the above splittings associated to the mixed Hodge structure $(W,\exp(\sum_{j=1}^{n}iy_jN_j)F)$.
Then, $\spl_W(y)$ converges in $\spl(W)$ when $y_j/y_{j+1}\to\infty$ ($1\le j\le n$, $y_{n+1}$ denotes $1$) (\cite{KNU1} Theorem 0.5 (1), \cite{KNU2} Part II 2.4.2 (i)), whereas $s(y)$ can diverge.
For this convergence, the term $\zeta$ in the canonical splitting $\spl_W(y)$ plays a crucial role (\cite{KNU1} Example 13.3).
\end{sbrem}

The canonical splitting of $W$ of MHS has a special importance and a characterization related to the theory of $\SL(2)$-orbits
as in Remark \ref{splrem2} and Remark \ref{splrem3} below  (\cite{KNU1} 8.7).

\begin{sbrem}\label{splrem2}  Assume that we are given a nilpotent orbit $(H_\R, \langle \cdot,\cdot\rangle, N, F)$ of weight $w$ as in \cite{KU2} 5.4.1.
   Let $W'$ be the $-w$ twist of the monodromy filtration $M$ of $N$. (That is, $W'$ is the twist of $M$ such that the central graded quotient 
of $W'$  is of weight $w$ whereas the central graded quotient of $M$ is of weight $0$.)  Then $(W', F)$ is a MHS. The canonical splitting of $W'$ of this MHS is explained as follows.

For $y\gg 0$, $(H_\R, \langle \cdot,\cdot\rangle, \exp(iyN)F)$ is a polarized Hodge structure. Let $s_{\BS}(y)$ be the unique splitting of $W'$  such that $s_{\BS}(y)(\gr^{W'}_w)$ and $s_{\BS}(y)(\gr^{W'}_{w'})$ with $w\neq w'$ are orthogonal with respect to the Hodge metric of 
 $(H_\R, \langle \cdot,\cdot\rangle, \exp(iyN)F)$. (This splitting of $W'$ is treated  in the theory of Borel--Serre lifting. See \ref{orth} below.)  As is proved in  \cite{KNU1}, when $y\to \infty$, $s_{\BS}(y)$ converges to the canonical splitting $\spl_{W'}(F)$ of $W'$ associated to the MHS $(W', F)$. 
 
 In the theory of $\SL(2)$-orbits, we have a homomorphism of algebraic groups $\rho:\SL(2)_\R\to  \Aut_\R(H_\R)$ over $\R$ associated to $(N, F)$. It  is the unique homomorphism such that $\rho\begin{pmatrix} 1/t & 0\\ 0&t\end{pmatrix}$ acts on the part of $H_\R$ of weight $j$  for  $\spl_{W'}(F)$  as $t^{j-w}$ and such that the Lie algebra homomorphism $\Lie(\rho) : {\frak {sl}}(2)_\R\to \End_\R(H_\R)$ sends the matrix $\begin{pmatrix} 0& 1 \\ 0&0 \end{pmatrix}$ to $N$. That is, the action of the diagonal part of $\SL(2)$ in the theory of $\SL(2)$-orbits gives the canonical splitting $\spl_{W'}(F)$. 
  This $\spl_{W'}(F)$ is denoted as $\spl_{W'}^{\BS}(F)$ in \ref{BSspl}.
\end{sbrem}

\begin{sbrem}\label{splrem3} The canonical splitting of the weight filtration is functorial for MHS. Any MHS $F'$ is embedded in a MHS $F$ which appears in Remark \ref{splrem2} (\cite{KNU1}). Hence Remark \ref{splrem2} characterizes the canonical splitting.  

\end{sbrem}

\subsection{The complex analytic structure of $D$}
\label{ss:cxanastrD}
We consider the complex analytic structure of $D$.

Fix a homomorphism $h_0: S_{\C/\R}\to G_{\red,\R}$ as in \ref{D}.

\begin{sbpara}\label{Y}
Let $Y$ be  the set of all isomorphism classes of exact $\otimes$-functors from $\Rep(G)$ to the following category $\cC$ preserving the underlying vector spaces and the weight 
filtrations. 

$\cC$ is the category of triples $(V, W, F)$, where $V$ is a finite-dimensional $\Q$-vector space, $W$ is an increasing filtration on $V$ 
(called the weight filtration), and $F$ is a decreasing filtration on $V_\C$ (called the Hodge filtration).  

Then $G(\C)$ acts on $Y$ by changing the Hodge filtration $F$. We have $D\subset Y$ and $D$ is stable in $Y$ under the action of $G(\R)G_u(\C)$. 

Let $$\Dc:=G(\C)D \subset Y.$$ Since the action of $G(\C)$ on $\Dc$ is transitive and the isotropy group of each point of $\Dc$ is an algebraic subgroup of $G(\C)$, $\Dc$ has a natural structure of a complex analytic manifold as a $G(\bC)$-homogeneous space.
\end{sbpara}

\begin{sbprop}\label{tan} For $x\in \Dc$, the tangent space of $\Dc$ at $x$ is canonically isomorphic to $\Lie(G)_{\bC}/F(x)^0\Lie(G)_{\bC}$, where $F(x)$ denotes the Hodge filtration of $x$ on $\Lie(G)_{\bC}$ defined by the adjoint action of $G$ on $\Lie(G)$. 

The tangent bundle of $\Dc$ is canonically isomorphic to $\Lie(G)_{\cO}/F^0\Lie(G)_{\cO}$, where $\Lie(G)_{\cO}:=\cO\otimes\Lie(G)$ with $\cO=\cO_{\Dc}$ the sheaf of holomorphic functions on $\Dc$.
\end{sbprop}

\begin{pf}
By definition of $F(x)$, $F(x)^0\Lie(G)_{\bC}$ is the Lie algebra of the isotropy subgroup of $G(\C)$ at $x$ under the action of $G(\C)$ on $\check D$ in \ref{Y}.
The assertions of this proposition follow.
\end{pf}

\begin{sbprop}\label{open}
$D$ is open in $\Dc$. 
\end{sbprop}

  \begin{pf}  Let $x\in D$.
Since the Hodge filtration $F(x)^{\bullet}\Lie(G_{\red})_{\C}$ is pure of weight $0$, the 
map $\Lie(G_{\red})_\R \to \Lie(G_{\red})_{\C}/F(x)^0\Lie(G_{\red})_\C$ is surjective.  Hence the map $\Lie(G)_{\bR}+\Lie(G_u)_\C\to \Lie(G)_{\bC}/F(x)^0\Lie(G)_{\bC}$ is surjective. 
  Since 
$\Lie(G)_{\C}/F(x)^0\Lie(G)_{\bC}$  is the tangent space of $\Dc$ at $x$ (Proposition \ref{tan}), the last surjectivity shows that $G(\R)G_u(\C) x $ is a neighborhood of $x$ in $\Dc$.
\end{pf}

\begin{sbcor}
  $D$ is a complex analytic manifold. 
\end{sbcor}

\begin{sbrem}This Proposition \ref{open} is  Proposition 3.2.7 of \cite{KNU3}. The proof of it given there is wrong.

\end{sbrem}

The real analytic structure of $D$ given in \ref{ranD} coincides with the one induced by this complex analytic structure. 

\begin{sbthm}\label{fiber} The map $D\to D_{\red}$ is smooth and surjective. For $x\in D$, the fiber over $x_{\red}\in D_{\red}$ ($\ref{Dred}$) in $D$  is isomorphic to $G_u(\C)/F(x)^0G_u(\C)$. Here $F(x)^0G_u(\C)$ is the $\C$-valued points of the algebraic subgroup of $G_{u,\C}$ whose Lie algebra is $F(x)^0\Lie(G_u)_{\C})$. 
\end{sbthm}

\begin{pf} The smoothness follows from the smoothness of $D$ and $D_{\red}$ and the surjectivity of the map $T_xD\to T_{x_{\red}}(D_{\red})$ of tangent spaces which follows, by Proposition \ref{tan}, from the surjectivity of $\Lie(G)\to \Lie(G_{\red})$. Consider the action of $G_u(\C)$ on the fiber. The action is transitive by Proposition \ref{DPsi} (3). We prove that the  isotropy subgroup of $G_u(\C)$ at $x$ is $F(x)^0G_u(\C)$. Let $g\in G_u(\C)$. Then $g\in F(x)^0G_u(\C)$ if and only if $\log(g) \in F(x)^0\Lie(G_u)_{\C})$, 
that is, if and only if $\log(g) \in F(x)^0\Lie(G)_{\bC}$.
By Lemma \ref{Lie} (2), the last condition is equivalent to the condition that $\log(g)F(x)^pV_\C\subset F(x)^pV_\C$ for any $V\in \Rep(G)$ 
for the Lie action. This condition is equivalent to the condition that $g F(x)^pV_\C = F(x)^pV_\C$ for any $V\in \Rep(G)$,
that is, $g$ fixes $x$. 
\end{pf}

\begin{sbpara}\label{faith} 
In this paper, we will often use the following fact (see \cite{Mi2} Theorem 4.14): If $\cG$ is a linear algebraic group over a field $E$ and if $V_1$ is a finite-dimensional  faithful representation of $\cG$, $V_1$ generates the $\otimes$-category $\Rep_k(\cG)$ of all finite-dimensional representations of $\cG$ over $E$. That is, 
all $V\in \Rep_k(\cG)$ can be constructed from $V_1$ by taking $\otimes$, $\oplus$, the dual, and subquotients.
\end{sbpara}

\begin{sblem}\label{DV} 
For $V\in\Rep(G)$, define $D(V)$ (resp.\ $\Dc(V)$) as the set $\{FH(V)\;|\; H\in D\}$ (resp.\ $\{FH(V)\;|\; H\in \Dc\}$) of decreasing filtrations on $V_\C$.

Assume that $V\in \Rep(G)$ is faithful. Then the map $D\to D(V)$ (resp.\ $\Dc\to \Dc(V)$)$; H\mapsto FH(V)$ is a bijection. 
If $H\in D$ (resp.\ $\Dc$), $D(V)$ (resp.\ $\Dc(V)$) coincides with the $G(\R)G_u(\C)$-orbit (resp.\ $G(\C)$-orbit) in the set of decreasing filtrations on $V_\C$ containing $FH(V)$.
\end{sblem}

\begin{pf}  Since $V$ is faithful, the map $D\to D(V)$ (resp.\ $\Dc\to \Dc(V)$)$; H\mapsto FH(V)$ is injective by \ref{faith}, and hence bijective by definition of $D(V)$ (resp.\ $\Dc(V)$).
The action of $G(\R)G_u(\C)$ (resp.\ $G(\C)$) on $D$ (resp.\ $\Dc$) is transitive. 
The second assertion follows.
\end{pf}

\subsection{Polarizability} \label{ss:pol}
  For a linear algebraic group $G$, let $G':=[G,G]$ be the commutator algebraic subgroup.
  Note that $G_u\subset G'$ if $D_{\all}(G)$ is non-empty (\ref{ucom}).

\begin{sbpara}\label{defpol}  

There are two formulations of polarization of a Hodge structure: 
the \lq\lq classical formulation'' (\cite{Gr} I 2) 
and the formulation by Deligne (\cite{De}). 
  We adopted the former in Part I--Part IV of this series of papers \cite{KNU2}. 

  Let $H$ be a $\Q$-Hodge structure of weight $w$. Then a polarization in the \lq\lq classical sense'' is a $\Q$-bilinear form
$\langle \cdot,\cdot\rangle: H_\Q\times H_\Q\to \Q$ which is symmetric if $w$ is even and anti-symmetric if $w$ is odd, satisfying 
$\langle F^p, F^{w+1-p}\rangle =0$ for any $p$, 
where $F$ is the Hodge filtration, and the condition that the Hermitian form $(\cdot, \cdot): H_\C\times H_\C\to \C$ defined by $(x, y)=\langle x, i^{p-q}\bar y\rangle$ for $x\in H_\C$ and $y\in H^{p,q}_\C$ is positive definite. (This positive definite Hermitian form  $(\cdot,\cdot)$ is called the {\it Hodge metric} of the polarization. Note that the restriction of the Hodge metric $(\cdot,\cdot)$ of the polarization  to $H_\R\times H_\R$ is a positive definite symmetric bilinear form $H_\R\times H_\R\to \R$ and it is written as $(x, y)\mapsto \langle x, h(i)y\rangle$ for the action $h$ of $S_{\C/\R}$ on $H_\R$ (\ref{SCR2}).)  

  On the other hand, a polarization in the sense of Deligne is a homomorphism $p: H_\Q\otimes H_\Q\to \Q(-w)=\Q \cdot (2\pi i)^{-w}$ of Hodge structures of weight $2w$ such that the induced $\Q$-bilinear form $H_\Q\times H_\Q \to \Q\;;\; 
(x, y)\mapsto (2\pi i)^wp(x\otimes y)$ is a polarization in the above classical sense. 

To keep the consistency with Parts I--IV of this series of papers, we adopt in this Part V the formulation of polarization in the classical sense. When we use the formulation of Deligne, we will say that it is a polarization in the sense of Deligne. 

We formulate polarizations of $\R$-Hodge structures in the same way. 
\end{sbpara}

\begin{sbpara}\label{pol} 
  Let $h_0: S_{\C/\R}\to G_{\red,\R}$ be as in \ref{D}. 
  Let $C:=h_0(i)$  (Weil operator) be the image of $i\in \C^\times = S_{\C/\R}(\R)$ by $h_0$ in $G_{\red}(\R)$. 
 
  We say that $h_0$ is {\it $\R$-polarizable} if 
$\{a\in (G_{\red})'(\R)\;|\; Ca=aC\}$ is a maximal compact subgroup of $(G_{\red})'(\R)$. 

That is,  $h_0$ is $\R$-polarizable if and only if $\Ad(C)$ on  $\Lie(G'_{\red})_{\bR})$ is a Cartan involution.
\end{sbpara}

In the following lemma, which  is a variant of \cite{De2} Section 2, 
we compare several polarizabilities (the above $\R$-polarizability is put as the condition (4.0)).

\begin{sblem}\label{pol2}  Let $h_0: S_{\C/\R}\to G_{\red,\R}$ be as in $\ref{D}$. Then, for $a=1, 2, 3$, the following conditions $(a.1)$ and $(a.2)$ are equivalent. Furthermore, the conditions $(4.0)$, $(4.1)$, and $(4.2)$ are equivalent. For the conditions $(a):=(a.1)\Leftrightarrow (a.2)$ 
($a=1,2,3$), and the condition $(4):= (4.0)\Leftrightarrow (4.1) \Leftrightarrow (4.2)$, we have the implications $(1) \Rightarrow (2) \Rightarrow (4)$ and $(1) \Rightarrow (3) \Rightarrow (4)$.

$(1.1)$ (resp.\ $(1.2)$). There is a homomorphism $t: G_{\red}\to \bG_m$ such that $t(h_0(w(x)))=x^{-2}$ ($x\in \bG_m$). Furthermore, if we consider the action of $G$ on $\Q\cdot (2\pi i)^r$ ($r\in \Z$) via $t^r: G_{\red}\to \bG_m$ and identify $H(\Q\cdot (2\pi i)^r)$ for $H\in D=D(G,h_0)$ with the Hodge structure $\Q(r)$, then, for every $H\in D$ and every $V\in \Rep(G)$ (resp.\ for some $H\in D$ and some faithful representation $V\in \Rep(G)$) and for each $w\in \Z$, there exists a homomorphism  $\gr^W_w(V)\otimes \gr^W_w(V) \to \Q\cdot (2\pi i)^{-w}$ in $\Rep(G)$  which polarizes the Hodge structure $\gr^W_wH(V)$ of weight $w$ in the sense of Deligne.

$(2.1)$ (resp.\ $(2.2)$).  There is a homomorphism $t: G_{\red,\R}\to \bG_{m,\R}$ such that $t(h_0(w(x)))=x^{-2}$ ($x\in \bG_{m,\R}$) and such that, for every $H\in D$ and every $V\in \Rep(G)$ (resp.\ for some $H\in D$ and some faithful representation $V\in \Rep(G)$) and for each $w\in \bZ$, there exists an $\R$-bilinear form $\langle \cdot,\cdot\rangle: \gr^W_w(V)_\R\times \gr^W_w(V)_\R \to \R$ satisfying $\langle gx, gy\rangle= t(g)^{-w}\langle x, y\rangle$  ($g\in G_{\red,\R}$) which polarizes the $\R$-Hodge structure $\gr^W_wH(V)_\R$ of weight $w$.

$(3.1)$ (resp.\ $(3.2)$). For every $H\in D$ and every $V\in \Rep(G)$ (resp.\ For some $H\in D$ and some faithful representation $V\in \Rep(G)$) and for each $w\in \bZ$, there exists a $\Q$-bilinear form $\langle \cdot,\cdot\rangle: \gr^W_w(V)\times \gr^W_w(V) \to \Q$ satisfying $\langle gx, gy\rangle= \langle x, y\rangle$  for $g\in G'_{\red}$ which polarizes the Hodge structure $\gr^W_wH(V)$ of weight $w$.

$(4.0)$  The  homomorphism $h_0:S_{\C/\R}\to G_{\red,\R}$ is $\R$-polarizable in the sense of $\ref{pol}$.

$(4.1)$ (resp.\ $(4.2)$). For every $H\in D$ and every $V\in \Rep(G)$ (resp.\ For some $H\in D$ and some faithful representation $V\in \Rep(G)$) and for each $w\in \bZ$, there exists an $\R$-bilinear form $\langle \cdot,\cdot\rangle: \gr^W_w(V)_\R\times \gr^W_w(V)_\R \to \R$ satisfying $\langle gx, gy\rangle= \langle x, y\rangle$  for $g\in G'_{\red,\R}$ which polarizes the $\R$-Hodge structure $\gr^W_wH(V)_\R$ of weight $w$.

\end{sblem}

\begin{pf} For $a=1,2,3,4$, the implication ($a.1$) $\Rightarrow$ ($a.2$) is clear and the implication ($a.2$) $\Rightarrow$ ($a.1$) follows from \ref{faith}. 
For $b=1, 2$, the implications ($1.b$) $\Rightarrow$ ($2.b$) $\Rightarrow$ ($4.b$) and the implications ($1.b$) $\Rightarrow$ ($3.b$) $\Rightarrow$ ($4.b$) are clear.

We prove the equivalence of  (4.0) and (4.2).

  To see the equivalence, by taking $\gr^{W}$, we may assume that $G$ is reductive. 
  Then the equivalence is an analogue of \cite{De2} lemme 2.8, and proved as follows.
  We assume that $G$ is reductive.  

  Assume (4.0). We prove that (4.2) is satisfied.

 Let $K:=\{a \in G'(\bR)\,|\,Ca=aC\}$. 
 Let $S_{\C/\R}^{(1)}$ be the kernel of the norm map $S_{\C/\R}\to \bG_{m,\R}$ and let $K_1\subset G(\R)$ be the image of $S_{\C/\R}^{(1)}(\R)=\{z\in \C^\times\;|\;|z|=1\}$ under $h_0$. 
 We first show that there is a compact subgroup $K_2$ of $G(\R)$ which contains  
 both $K$ and $K_1$. Let $J= \{a\in G_\R\;|\; Ca=aC\}$ and let $K_2$ be a maximal compact subgroup of $J(\R)$ which contains $K_1$. Then, since $K_2$ contains some conjugate of $K$ in $J(\R)$ and since $K$ is normal in $J(\R)$, we have $K\subset K_2$.
 
  By \cite{Mo}, there is a finite-dimensional faithful representation $V$ of $G_{\bR}$ and a positive definite symmetric $\R$-bilinear form $(\cdot, \cdot): V\times V\to \R$ which is fixed by $K_2$ such that $G_\bR$ is stable in $\text{Aut}_\R(V)$ under the transpose $g\mapsto {}^tg$ with respect to $(\cdot,\cdot)$. 
  Note that the last condition implies that $G'_{\bR}$ is also stable under the transpose. 

\smallskip

\noindent 
{\bf Claim 1.}  The Cartan involution $\theta_K: G'_\R \to G'_\R$ associated to $K$ is $g \mapsto {}^t g^{-1}$.

\smallskip

  Note that this claim is also used in \cite{BS}. 

\noindent
{\it Proof of Claim 1.} This is an algebraic homomorphism and its set of fixed points is compact and contains $K$.  Since $K$  is a maximal compact subgroup of $G'(\bR)$ by the assumption, $K$ coincides with the set of fixed points of $\theta_K$. 
This proves Claim 1.

\smallskip

On the other hand, we know that $g \mapsto C^{-1}gC$ is the Cartan involution of $G'_\R$ associated to $K$. 
Hence, by Claim 1, we have $C^{-1}gC= {}^t g^{-1}$. 

  Put $\langle x, y \rangle := (x, C^{-1}y)$. 
  We show that it is $G'_\R$-invariant. 
  Let $g$ be in $G'_\R$. 
  Then we have $\langle gx,gy \rangle=(gx, C^{-1}gy)=(gx,{}^tg^{-1}C^{-1}y)
= (g^{-1}gx, C^{-1}y) = (x, C^{-1}y)=\langle x, y \rangle$. 

Let  $V_w$ for $w\in \Z$ be the part of $V$ of weight $w$  with respect to $h_0$.
Let $c\in \R^\times\subset \C^\times=S_{\C/\R}(\R)$. We prove that ${}^th_0(c)v=c^wv$ for $v\in V_w$. Since ${}^t h_0(c)$ belongs to $G(\R)$ in $\text{Aut}_\R(V)$, we have ${}^th_0(c)V_w=V_w$.  For every $v'\in V_w$, we have $({}^th_0(c)v, v')= (v, h_0(c)v')= (v, c^wv')=(c^wv, v')$. Since $(\cdot, \cdot): V_w\times V_w\to \R$ is non-degenerate, we have ${}^th_0(c)v=c^wv$. 
We prove $\langle V_w, V_{w'}\rangle=0$ unless $w=w'$. 
For $v\in V_w$ and $v'\in V_{w'}$, $c^w(v, v')= (h_0(c)v, v')= (v, {}^th_0(c)v')= (v, c^{w'}v')=c^{w'}(v, v')$. Hence if $w\neq w'$, then $(V_w, V_{w'})=0$ and hence $\langle V_w, V_{w'}\rangle=0$.

Let $\langle\cdot,\cdot\rangle_w: V_{w,\R}\times V_{w,\R}\to \R$ be the pairing induced by $\langle \cdot,\cdot\rangle$. 
We prove that $\langle\cdot, \cdot\rangle_w$ is a polarization on $V_w$. Let $(\cdot,\cdot)_w: V_{w, \C}\times V_{w, \C} \to \C$ be the  positive definite Hermitian form induced by $(\cdot,\cdot)$.  Let $H_w^{p,q}$ ($p+q=w$) be the $(p,q)$-Hodge component of $V_{w, \C}$ with respect to $h_0$. It is sufficient to show that $(H_w^{p,q}, H_w^{p',q'})=0$ unless $p=p'$. Let $u\in K_1=S_{\C/\R}^{(1)}(\R)$. Then, since $(\cdot, \cdot)_w$ is $K_1$-invariant, we have, for $v\in H_w^{p,q}$ and $v'\in H_w^{p',q'}$, $(v, v')= ([u]v, [u]v')= u^{p-q}u^{q'-p'}(v,v')= u^{2(p-p')}(v,v')$. Hence $p=p'$. This proves that the condition (4.2) is satisfied.

  Assume (4.2).
  We prove that (4.0) is satisfied. Take a faithful representation $V\in \Rep(G)$. By our assumption, 
  there is a $G'_\R$-invariant bilinear form $\langle\cdot,\cdot\rangle$ on $V_\R$ such that $(x,y):=
\langle x, Cy\rangle$ is positive definite and symmetric. 

\smallskip

\noindent 
{\bf Claim 2.}  For $g\in G'_\R$, we have $C^{-1}gC ={}^t g^{-1}$, where the transpose is with respect to $(\cdot,\cdot)$. 

\smallskip

\noindent {\it Proof.} This follows from the fact that $\langle\cdot, \cdot\rangle$ is $G'_\R$-stable. 

The following Claim 3 is well-known. 

\smallskip
\noindent
{\bf Claim 3.} Let $U$ be a finite-dimensional $\R$-vector space endowed with a positive definite symmetric bilinear form $(\cdot, \cdot)$. Let $\cG$ be an algebraic subgroup of $\GL_U$ which is stable under $g\mapsto {}^tg$, where the transpose is with respect to $(\cdot, \cdot)$. Then $\{g\in \cG(\R)\;|\; {}^tg=g^{-1}\}$ is a maximal compact subgroup of $\cG(\R)$. 

  By Claim 2 and Claim 3, $\{g \in G'(\bR)\,|\, C^{-1}gC=g\}$ is a maximal compact subgroup of $G'(\bR)$.
  Thus (4.2) implies (4.0). 
\end{pf}

\begin{sbpara}\label{example} The conditions (1)--(4) in Lemma \ref{pol2} 
are different from each other as the following examples show. 

{\it Example} 1. Let $E$ be a cubic extension field of $\Q$ having one real place and one complex place and let  $G=E^\times$ regarded as a torus over $\Q$ of dimension $3$. Let $h_0 :S_{\C/\R}\to G_\R$ be the homomorphism such that the induced map $S_{\C/\R}(\R)=\C^\times \to G(\R)= (E\otimes_\Q \R)^\times= \R^\times \times \C^\times$ sends  $z\in \C^\times$ to $(z\bar z, z^2)$. This example satisfies the conditions (2) and (4), but does not satisfy the conditions (1), (3). In fact, (2) is satisfied because we have $t: G_\R\to \bG_{m,\R}$ which sends $(r,z)\in G(\R)= \R^\times\times \C^\times$ to $r(z\bar z)^{-1}$, and for $V=E\in \Rep(G)$, we have an $\R$-bilinear form on $E_\R= \R\times \C$ with values in $\R\cdot (2\pi i)^{-2}=\R$ 
given by $((x_1, y_1), (x_2, y_2))\mapsto x_1x_2 - y_1\bar y_2-\bar y_1 y_2$ ($x_j\in \R$, $y_j\in \C$) which polarizes the Hodge structure  of $E$ of weight $2$ associated to $h_0$. 
But the condition (3) is not satisfied because there is no bilinear form on the $\Q$-vector space $E$ which polarizes the Hodge structure associated to $h_0$.

{\it Example} 2. Let $E$ be an imaginary quadratic field over $\Q$ and let $G=E^\times$ regarded as a torus over $\Q$ of dimension $2$. Let  $h_0: S_{\C/\R}\to G_\R$ be the homomorphism such that the induced map $S_{\C/\R}(\R)=\C^\times \to G(\R)= (E\otimes_\Q \R)^\times= \C^\times$ sends  $z\in \C^\times$ to $z/\bar z$. This example  satisfies the conditions (3) and (4), but does not satisfy the conditions (1), (2). In fact,   there is no homomorphism $t:G_\R\to \bG_{m,\R}$ such that $t(h_0(w(x)))=x^{-2}$. The condition (3) is satisfied because $E\times E\to \Q\;;\; (x_1,x_2)\mapsto -x_1\bar x_2-\bar x_1 x_2$ ($x_j\in E$) is the polarization of the Hodge structure of $E$ of weight $0$ associated to $h_0$.

\end{sbpara}

  We thank Teruhisa Koshikawa for his advice on $\R$-polarizability and for showing the above Example 1 to us.

  \begin{sbpara} \label{X}
  Let $\cG$ be a semisimple algebraic group over $\R$. 
  
Let $\frak X$ be the set of all maximal compact subgroups of $\cG(\R)$. Then $\frak X$ is not empty. The group $\cG(\R)$ acts on $\frak X$ by conjugation, and this action is transitive. See \cite{Mo} Theorem 3.1,  \cite{BH} Proposition 1.12. For $K\in \frak X$, we have a bijection $\cG(\R)/K\to \frak X\;;\;g\mapsto gKg^{-1}$. Via this bijection, we regard $\frak X$ as a real analytic manifold. This real analytic structure is independent of the choice of $K$.
This $\frak X$ is called the {\it symmetric space associated to $\cG$}.

\end{sbpara}

 \begin{sbpara}\label{DtoX}

 Consider the commutator subgroup  $G'_{\red}=[G_{\red},G_{\red}]$ 
 of $G_{\red}=G/G_u$. This is a semisimple algebraic group. 
Let $\frak X$ be the symmetric space associated to $G'_{\red,\R}$. 

Let $h_0:S_{\C/\R}\to G_{\red,\R}$ be $\R$-polarizable as in \ref{pol}.
From $D_{\red}=D(G_{\red},h_0)$ (\ref{Dred}), we have a canonical map 
$$D_{\red}\to \frak X$$ 
which sends $h\in D_{\red}\subset \Hom(S_{\C/\R}, G_{\red,\R})$ to a maximal compact subgroup $K=\{g\in G'_{\red}(\R)\;|\; gh(i)=h(i)g\}$ of $G'_{\red}(\R)$ associated to the Weil operator $h(i)$.

 \end{sbpara}

\begin{sbpara}\label{Gamma} Let $\Gamma$ be a subgroup of $G(\Q)$. 
  We call $\Gamma$ an {\it arithmetic subgroup} (resp.\ {\it semi-arithmetic subgroup}) of $G(\Q)$ if the following condition $(A)$ (resp.\ $(SA)$) is satisfied.

$(A)$  There are $n\geq 1$ and  an injective homomorphism $\rho: G\to\GL(n)$ such that $\Gamma$ is a subgroup of $\{g \in G(\Q)\;|\; \rho(g)\in \GL(n,\Z)\}$ of finite index.

$(SA)$   There are $n\geq 1$ and  an injective homomorphism $\rho: G\to\GL(n)$ such that $\rho(\Gamma)\subset \GL(n,\Z)$.
That is, there is a faithful representation $V \in \Rep(G)$ and a $\bZ$-lattice $L$ in $V$ such that $L$ is stable under the action of $\Gamma$.

The terminology arithmetic subgroup is used by many people. We hope the terminology semi-arithmetic group is acceptable. 
\end{sbpara}

 The next two lemmas are straightforward. 

\begin{sblem}\label{Gamma2}  

$(1)$  The condition $(A)$ is equivalent to the following condition $('A)$.

$('A)$  For every $n\geq 1$ and every homomorphism $\rho: G\to \GL(n)$, $\rho(\Gamma) \cap \GL(n,\Z)$ is of finite index  in $\rho(\Gamma)$ and in $\rho(G(\Q)) \cap \GL(n,\Z)$.  

$(2)$ The condition $(SA)$ is equivalent to the following condition $('SA)$.

$('SA)$ For every $n\geq 1$ and every homomorphism $\rho: G\to \GL(n)$, $\rho(\Gamma) \cap \GL(n,\Z)$ is of finite index  in $\rho(\Gamma)$.  

\end{sblem}

\begin{sblem}\label{redss2general} 
  Let $f: G_1\to G_2$ be a homomorphism of linear algebraic groups over $\Q$, let $\Gamma_1$ be a subgroup of $G_1(\Q)$ and let $\Gamma_2$ be
  the image of $\Gamma_1$ in $G_2(\Q)$.
  
  $(1)$  If $\Gamma_1$ is a semi-arithmetic subgroup of $G_1(\Q)$,   $\Gamma_2$ is a semi-arithmetic subgroup of $G_2(\Q)$.
  
  $(2)$ If  $\Gamma_1$ is an arithmetic subgroup of $G_1(\Q)$ and if $f$ is surjective, $\Gamma_2$ is an arithmetic subgroup of $G_2(\Q)$.
    \end{sblem}

\begin{sbrem} 

If $\Gamma$ is a semi-arithmetic subgroup of $G(\Q)$, $\Gamma$ is discrete in $G(\R)$. (The converse is not valid: A subgroup of $G(\Q)$ which is discrete in $G(\R)$ need not be semi-arithmetic.  For example, let $G=\SL(2)$, and let $\Gamma$ be the subgroup of $G(\Q)$ consisting of diagonal matrices with diagonal entries $(2^n, 2^{-n})$ ($n\in \Z$). Then $\Gamma$ is discrete in $G(\R)=\SL(2,\R)$ but $\Gamma$ is not semi-arithmetic.  Another example is $G=\Gm$ and $\Gamma=\{2^n\,|\,n \in \bZ\}$.)
\end{sbrem}

\begin{sbprop}\label{actonD}
  Let $h_0: S_{\C/\R}\to G_{\red,\R}$ be as in $\ref{D}$.
  Assume that $h_0$ is $\R$-polarizable ($\ref{pol}$).
  Let $\Gamma$ be a semi-arithmetic subgroup of $G'(\Q)$  ($\ref{Gamma}$).
Then the following holds. 

$(1)$ The action of $\Gamma$ on $D$ is proper and the quotient space $\Gamma \bs D$ is Hausdorff. 

$(2)$ If $\Gamma$ is torsion-free, the action of $\Gamma$ on $D$ is free (that is, if $\gamma\in \Gamma$ and if   $\gamma p=p$ for some $p\in D$, then $\gamma=1$), and
the projection $D\to \Gamma \bs D$ is a local homeomorphism.
\end{sbprop}

This Proposition \ref{actonD}  follows from its stronger version Theorem \ref{BSgl}. (We will not use this \ref{actonD} before we prove \ref{BSgl}.)

\begin{sbrem} 
\label{r:actonD}
 Proposition \ref{actonD} for a semi-arithmetic subgroup of $G(\Q)$ need  not be true as is shown in the following examples.

{\it Example} 1. Let $E$ be a number field (a finite extension of $\Q$), and let $G$ be the algebraic group $E^\times$ over $\Q$ (that is, $G(R)=(E\otimes_\Q R)^\times$ for any commutative ring $R$ over $\Q$). Then the unit group $O_E^\times$ of $O_E$ is an arithmetic  subgroup of $G(\Q)$. Take any homomorphism $h_0 : S_{\C/\R} \to G_\R$ (for example, the trivial homomorphism). Then $h_0$ is $\R$-polarizable and $D$ consists of one point. The action of $\Gamma:=O_E^\times$  on the one-point set $D$ is proper if and only if $\Gamma$ is finite. But $O_E^\times$ need not be finite, e.g., for a real quadratic extension field $E$ over $\Q$.

{\it Example} 2. Let $E$ be a real quadratic field, and 
let $G$ be the algebraic group 
$$\left\{\left.
\begin{pmatrix}
a & b \\
0 & 1 
\end{pmatrix}
\right|a \in E^\times,\, b \in E\right\}
$$ over $\Q$ (that is, $G(R)=
\left\{\left.
\left(
\begin{matrix}
a & b \\
0 & 1 
\end{matrix}
\right)
\right|a \in (E\otimes_\Q R)^{\times},\, b \in E\otimes_\Q R\right\}$ 
for any commutative ring $R$ over $\Q$). 
  Then the group 
$\Gamma:=\left\{\left.
\left(
\begin{matrix}
a & b \\
0 & 1 
\end{matrix}
\right)
\right|a \in O_E^\times,\, b \in O_E\right\}$ is an arithmetic  subgroup of $G(\Q)$.
  Let $h_0 : S_{\C/\R} \to G_{\red,\R}$ be the composite 
$S_{\C/\R} \to \bG_{m,\R} \to E^{\times}_{\bR} = G_{\red,\R}$, 
where the first homomorphism is the norm inverse and the second is the natural inclusion, that is, the one such that the induced map $S_{\C/\R}(\R)=\C^\times \to (E\otimes \bR)^{\times} = G_{\red}(\R)$
sends $z\in \C^\times$ to 
$|z|^{-2}$.  
Then $h_0$ is $\R$-polarizable (simply because $(G_{\red})'$ is trivial). Let $H\in D$ be the element corresponding to the composition $S_{\C/\R}\overset{h_0}\to E_\R^\times \subset G_\R$. By Theorem \ref{fiber}, or by \ref{DV} applied to $V=E^2$, 
we have $D=G_u(\bC)=E \otimes_{\Q}\C$, where $g\in G_u(\C)$ is identified with $gH\in D$.  
 Via this identification, 
  the action of $\Gamma$ on $D$ comes from the adjoint action of $G$ on $G_u$, and described as 
$$\left(\begin{matrix}
a & b \\
0 & 1 
\end{matrix}\right)\cdot
(x \otimes w)=(ax+b) \otimes w\qquad (a \in O_E^{\times},\, b \in O_E,\, x \in E,\, w \in \bC).$$
  The subspace $\Gamma \bs (E \otimes_{\Q} \bR)$ of the quotient $\Gamma \bs D$ is not Hausdorff because it is homeomorphic to 
the quotient of the real torus $O_E \bs (E \otimes_{\bQ} \bR)$ by the action of $O_E^{\times}$, and the last action has a dense orbit (\cite{R}). 
  Hence, $\Gamma \bs D$ is also not Hausdorff.
\end{sbrem}

\begin{sbprop}\label{A=A*}
Assume that the condition $(1)$ of Lemma $\ref{pol2}$ is satisfied. Then  for a semi-arithmetic  subgroup $\Gamma$ of $G(\Q)$, $\Gamma \cap G'(\Q)$ is of finite index in $\Gamma$.

\end{sbprop}

\begin{pf} 
By $G_u\subset G'$ (\ref{ucom}), we have $G/G'\overset{\sim}\to G_{\red}/G'_{\red}$. Hence by replacing $\Gamma$ by the image of $\Gamma$ in $G_{\red}(\Q)$,  we are reduced to the case $G$ is reductive.

Assume that $G$ is reductive. Let $\Gamma_0=\Gamma \cap Z(G)(\Q)$ and $\Gamma_1=\Gamma \cap G'(\Q)$. Since $Z(G) \times G'\to G$ is an isogeny, the image of $\Gamma_0\times \Gamma_1 \to \Gamma$ is of finite index. Hence it is sufficient to prove that $\Gamma_0$ is finite. We prove this. 

The image of $\Gamma$ under $t: G(\Q)\to {\bf G}_m(\Q)=\Q^\times$ is contained in $\{\pm 1\}$. 
 Hence in the faithful representation $V$ of $G$ in (1.2) in Lemma \ref{pol2}, for $H\in D(G, h_0)$, the action of some subgroup of $\Gamma_0$ of finite index preserves the Hodge filtration and the polarization of $H(\gr^W_wV)$ and hence preserves the Hodge metric for every $w$. The elements of $G(\R)$ which preserve these Hodge metrics for all $w$ form a compact subgroup and $\Gamma_0$ is discrete, and hence $\Gamma_0$ is  finite.  
\end{pf}

\subsection{Relations with usual period domains and Mumford--Tate domains}\label{ss:Dold}

In this section, in \ref{Dusual}--\ref{clEx}, we  explain that the classical Griffiths domains \cite{Gr} and their mixed Hodge generalization in \cite{U} are 
essentially regarded as 
special cases of the period domains in this paper. 
  In this case, our partial compactifications essentially coincide with those in \cite{KNU2} Part III. Next in \ref{relMT}, we explain that the Mumford--Tate domains studied in \cite{GGK} are regarded as important  cases of our period domains for reductive $G$. For Mumford--Tate domains, our toroidal  partial compactifications essentially coincide with those in \cite{KP}. 

\begin{sbpara}\label{Dusual}

  Let $\Lambda=(H_{0, \Q}, W, (\langle \cdot,\cdot\rangle_w)_w, (h^{p,q})_{p,q})$ be as in \cite{KNU2} Part III 2.1.1. That is, $H_{0,\Q}$ is a finite-dimensional $\Q$-vector space,  $W$ is a finite  increasing filtration on $H_{0, \Q}$,  $\langle\cdot,\cdot\rangle_w$ for each $w\in \Z$ is a non-degenerate $\Q$-bilinear form $\gr^W_w\times \gr^W_w\to \Q$ which is symmetric if $w$ is even and anti-symmetric if $w$ is odd, $h^{p,q}$ are non-negative integers given for each $(p,q)\in \Z^2$ such that $\sum_{p,q} h^{p,q}=\dim_\Q H_{0,\Q}$, $\sum_{p+q=w} h^{p,q}= \dim_\Q \gr^W_w$ for every $w\in \Z$, and $h^{p,q}=h^{q,p}$ for all $(p,q)$.

  Let $G$ be the subgroup of $\Aut(H_{0,\Q}, W)\times {\bf G}_m$ consisting of all elements $(g,t)$ such that 
  $\langle gx,gy\rangle_w=t^w\langle x,y\rangle_w$ for all $w$ and for all $x,y\in \gr^W_wH_{0,\Q}$.

  Let $D(\Lambda)$ be the period domain of \cite{U}.
 As a set, it is the set of all decreasing filtrations $F$ on $H_{0,\C}=\C\otimes_\Q H_{0,\Q}$ such that $\dim_\C(\gr^p_F(\gr^W_wH_{0,\C}))=h^{p,w-p}$ for all $w,p\in \Z$ and such that $(\gr^W_w, F(\gr^W_w), \langle\cdot,\cdot\rangle_w)$ is a polarized Hodge structure for any $w\in \Z$. 

  Let $D^{\pm}(\La)$ be the set of all decreasing filtrations on $H_{0,\C}$ such that $\dim_\C(\gr^p_F(\gr^W_w))=h^{p,w-p}$ for all $w,p\in \Z$ and such that either $(\gr^W_w, F(\gr^W_w), \langle\cdot,\cdot\rangle_w)$ is a polarized Hodge structure for any $w\in \Z$ or $(\gr^W_w, F(\gr^W_w), (-1)^w\langle\cdot,\cdot\rangle_w)$ is a polarized Hodge structure for any $w\in \Z$. Then $D^{\pm}(\La)=D(\La)$ if and only if $\gr^W_wH_0=0$ for all odd $w$. If $\gr^W_w\neq 0$ for some odd $w$, there is a $(g,t)\in G(\Q)$ such that $t<0$ and that $D^{\pm}(\La)= D(\La) \coprod gD(\La)$. 
 \end{sbpara}
 
 \begin{sbpara}
 
Assume that $D(\Lam)$ is not empty and fix an $\br\in D(\Lambda)$.
  Then the Hodge decomposition of  $\br(\gr^W)$ induces a homomorphism 
$h_0: S_{\C/\R} \to G_{\red,\R}$. We have
$\langle h_0(z)x, h_0(z)y\rangle_w= |z|^{2w}\langle x, y\rangle_w$ for $z\in \C^\times=S_{\C/\R}(\R)$ and $x,y\in\gr^W_w$. The condition (1) of Lemma \ref{pol2} is satisfied. 
 \end{sbpara} 
 
\begin{sbprop} We have an isomorphism $$D\overset{\sim}\to D^{\pm}(\La)\;;\;H\mapsto H(H_{0,\Q}),$$ where $D=D(G,h_0)$ is the period domain of the present paper. 
\end{sbprop}  

\begin{pf}
  This is seen by Lemma \ref{DV}. 
\end{pf}

\begin{sbpara} 
  In the above situation, the extended period domains $D_{\BS}$, $D_{\SL(2)}$, and $\Gamma \bs D_{\Sig}$ in this paper generalize those in \cite{KNU2} Part I--Part IV.
\end{sbpara}

\begin{sbpara}\label{clEx} 
(Classical) {\it Example}. Take $\La$ in \ref{Dusual} as follows. $H_{0,\Q}=\Q^2=\Q e_1+\Q e_2$, $W_1=H_{0,\Q}$, $W_0=0$, $\langle \;,\;\rangle_1$ is the anti-symmetric form characterized by $\langle e_2, e_1\rangle_1=1$, $h^{1,0}=h^{0,1}=1$, and other $h^{p,q}$ are $0$. 
Identify $\GL(2)$ with the subgroup $G$ in \ref{Dusual} by $g\mapsto(g,\det g)$.
Let $h_0: S_{\C/\R}\to G_\R$ be the homomorphism whose homomorphism of $\R$-valued points is $\C^\times \to \GL(2,\R)\;;\;z\mapsto \langle z\rangle$, where 
$$\langle z\rangle= \begin{pmatrix}  a & -b\\ b & a\end{pmatrix}  \quad \text{for}  \; z=a+bi  \; (a, b\in \R, (a,b) \neq (0,0)).$$
Then we have unique isomorphisms $$%
\frak H^{\pm} \simeq D:=D(G,h_0), \quad \bP^1(\C) \simeq \Dc,$$
which send $i\in \frak H\subset \frak H^{\pm}$ to $h_0$ such that the former is %
$\GL(2,\R)$-equivariant and the latter is %
$\GL(2,\C)$-equivariant. 
  Here $%
\frak H^{\pm}$ denotes %
the disjoint union of the upper half plain and the lower half plain, which are interchanged by 
$\begin{pmatrix}0&1\\1&0\end{pmatrix}$.

In fact, consider the natural action of $G=\GL(2)$ on $V=\Q^2$.
The eigenvalues of %
$h_0(z)$ are $z$, $\bar z$, and the eigenspace decomposition is $V_\C=\C\begin{pmatrix}i\\1\end{pmatrix}\oplus\C\begin{pmatrix}-i\\1\end{pmatrix}$, which yields the Hodge decomposition corresponding to $i\in\frak H$.
\end{sbpara}

\begin{sbpara}\label{relMT} Let $H$ be a polarized $\Q$-Hodge structure. 

Let $\cC$ be the Tannakian category of $\Q$-Hodge structures generated by $H$ and $\Q(1)$.  Let $M$ be the Tannakian group of $\cC$. This means that $\cC$ is identified with the Tannakian category  $\Rep(M)$. This $M$ is reductive. Let $h_0: S_{\C/\R}\to M_\R$ be the associated homomorphism. Then our period domain $D(M, h_0)$ coincides with the Mumford--Tate domain in \cite{GGK} associated to $H$. 

Take  $\La=(H_{0,\Q}, W, (\langle\cdot,\cdot\rangle_w)_w, (h^{p,q})_{p,q})$ in \ref{Dusual} as follows. $H_{0,\Q}$ is the $\Q$-structure $H_\Q$ of $H$. 
 $W$ is the weight filtration on $H_\Q$ (that is, if $w_0$ denotes the weight of $H_\Q$, $W_w=H_\Q$ if $w\geq w_0$ and $W_w=0$ if $w<w_0$). $\langle \cdot,\cdot \rangle_{w_0}$ is the polarization of $H$ (times $(2\pi i)^{w_0}$). $h^{p,q}$ is the dimension of the $(p,q)$-Hodge component of $H$. Let $L$ be the algebraic group $G$ in \ref{Dusual}. Then $M$  is identified with the smallest algebraic closed subgroup $M$  of $L$ defined over $\Q$ such that $M_\R$  contains the image of the homomorphism  $S_{\C/\R}\to L_\R$ associated to $H$.  The Mumford--Tate domain $D(M, h_0)$ is identified with the $M(\R)$-orbit in $D(\La)$ containing the class of $H$. 
 
 This $h_0$ satisfies the condition (1) in Lemma \ref{pol2}. 
 
 The period domain $D(M/Z, \bar h_0)$, where $Z$ is the center of $M$ and $\bar h_0$ denotes the composition $S_{\C/\R}\overset{h_0}\to M_\R\to (M/Z)_\R$, for the semisimple algebraic group $M/Z$ is also considered as a Mumford--Tate domain, for example, as in \cite{KP}. The period domain $D(M, h_0)$ is identified with an open and closed subspace of the period domain $D(M/Z, \bar h_0)$ (Proposition \ref{redss3}).

\end{sbpara}

\section{The space of Borel--Serre orbits}\label{s:DBS}

We define and study the space $D_{\BS}\supset D$ of Borel--Serre orbits. This is the $G$-MHS version of the space $D(\La)_{\BS}\supset D(\La)$ for the classical period domain $D(\La)$ (\ref{ss:Dold}) defined and studied in \cite{KU1} (in the pure case) \cite{KNU2} Part I (in the mixed case).

For an algebraic group $G$, $G^{\circ}$ denotes the connected component of $G$ as an algebraic group which contains $1\in G$. 
  A parabolic subgroup $P$ of $G$ is a closed algebraic subgroup $P\subset G^{\circ}$  such that $G/P$ is a projective variety (see e.g.\ \cite{B2} IV (11.2)).

  The organization of this Section 2 is as follows. 
  After preliminaries in Sections \ref{ss:cL} and \ref{ss:BSli} and a review of Borel--Serre space $\frak X_{\BS}$ in Section \ref{ss:BS1}, we define our space $D_{\BS}$ as a set in Section \ref{ss:DBS}, which we endow with a real analytic structure (precisely speaking, a structure of a real analytic manifold with corners) in Section \ref{ss:BSan}. In Section \ref{ss:BSproperty}, we prove the nice properties 
  of $D_{\BS}$ (e.g.\ Hausdorffness of the quotient by a  semi-arithmetic subgroup of $G'(\Q)$).

\subsection{Real analytic manifolds with corners}\label{ss:cL}

As will be explained in Section \ref{ss:BSan}, our space $D_{\BS}$ is a real analytic manifold with corners.

In this Section \ref{ss:cL}, we review this notion real analytic manifold with corners (\cite{BS} Appendix by A.\ Douady and L.\ Herault) and consider spherical compactifications as examples 
of real analytic manifolds with corners.

\begin{sbpara}
Let $m,n\geq 0$ and consider the topological space $S=\R^m \times \R^n_{\geq 0}$ which is endowed with the inverse image $\cO_S$ of the sheaf of ($\R$-valued) real analytic functions on $\R^{m+n}$. That is, $\cO_S$ is the sheaf of functions which are locally extendable to real analytic functions on an open subset of $\R^m \times \R^n$. 

A {\it real analytic manifold with corners} is a locally ringed space over $\R$ which has an open covering whose each member is isomorphic to an open set of the above  $S=\R^m \times \R^n_{\geq 0}$  with the restriction of $\cO_S$ for some $m,n\geq 0$. 
\end{sbpara}

\begin{sbpara}\label{cgv} 
Consider 
 a finite-dimensional graded real vector space $V=\bigoplus_{w\le-1}V_w$. In the rest of this Section \ref{ss:cL},
 we review the compactification $\overline V$ of $V$ defined as a real analytic manifold with corners in \cite{KNU2} Part I Section 7. We call it the {\it spherical compactification} of $V$ because as a topological space, it coincides with the spherical compactification of $V$ in \cite{BeS} Definition 2.1. We will use this $\overline V$ in Section \ref{ss:DBS}. 

As in \cite{KU2} and \cite{KNU2}, the property compact includes Hausdorff in our terminology. 
\end{sbpara}

 \begin{sbpara}\label{cgv1}  
 Consider the action of the multiplicative group $\R_{>0}$ on $V$ given by $t\circ (\sum_w v_w):=\sum_w t^wv_w$ ($t\in \R_{>0}$, $v_w\in V_w$).

Let $$\overline V:=V\times^{\R_{>0}} \R_{\geq 0} \smallsetminus \{(0,0)\},$$ which we endow with the natural topology. 

Recall that for a group $H$, for a set $X$ on which $H$ acts from the right, and for a set $Y$ on which $H$ acts from the left, $X\times^H Y$ denotes the quotient of $X\times Y$ by the equivalence relation $(xh,y)\sim (x, hy)$ ($x\in X$, $y\in Y$, $h\in H$). When we use the notation $X\times^H Y$ in this paper, $H$ is a commutative group (as above) and hence right or left in the action doesn't matter.  

Embed $V$ into $\overline V$ by $v\mapsto \text{class}(v,1)$ ($v\in V$). 
 We have $\overline{V}\smallsetminus V=\{\text{class}(v,0)\;|\;v\in V \smallsetminus \{0\}\}$, and for $v\in V\smallsetminus \{0\}$, when $t\in \R_{>0}$ converges to $0$, $t\circ v=\text{class}(tv,1)=\text{class}(v,t)$ converges to $\text{class}(v,0)$.  This space $\overline{V}$ is covered by the two open sets,  $V$ and the complement $\overline{V}\smallsetminus \{0\}$ of $0\in V\subset\overline{V}$.

 \end{sbpara}

 \begin{sbpara}\label{cgv2}  
We regard $\overline V$  as a real analytic manifold with corners as follows. 

There is a real analytic map $f: V\smallsetminus \{0\} \to \R_{>0}$ such that 
\smallskip

(1) $f(t\circ v)= tf(v)$ for any $t\in \R_{>0}$ and $v\in V\smallsetminus \{0\}$. 
\smallskip

For example, 
 taking a base $(e_{w,j})_j$ of $V_w$ for each $w$ and taking an integer $m<0$ satisfying $m/w\in 2\Z$ for any $w$ such that  $V_w\neq 0$, let $f(\sum_{w,j} x_{w,j}e_{w,j}):= (\sum_{w,j}  x_{w,j}^{m/w})^{1/m}$. 
 Then, this $f$ satisfies the condition (1).

  Let $f: V\smallsetminus \{0\}\to \R_{>0}$ be a real analytic map satisfying (1). Let $V^{(1)}=f^{-1}(1)$. Then $V^{(1)}\times \R_{>0}\overset{\sim}\to V\smallsetminus \{0\}\;;\;(v, t)\mapsto t\circ v$. 
The inverse map is $v\mapsto(f(v)^{-1}\circ v,f(v))$.  
  We have a canonical homeomorphism  $V^{(1)}\times \R_{\geq 0} \to  \overline V \smallsetminus \{0\}$. We endow $\overline V \smallsetminus \{0\}$ with the structure of a real analytic manifold with corners via this homeomorphism. This structure is independent of the choice of $f$. This is because if $f': V\smallsetminus \{0\}\to \R_{>0}$ is also a real analytic map satisfying (1) and if $V^{(1)'}:=(f')^{-1}(1)$, the isomorphism $V^{(1)}\times \R_{\geq 0} \overset{\sim}\to V^{(1)'}\times \R_{\geq 0}\;;\;(v, t) \mapsto (f'(v)^{-1}\circ v, tf'(v))$ of real analytic manifolds with corners is compatible with the homeomorphisms to $\bar V \smallsetminus \{0\}$.

  Furthermore, the restriction of this structure to $V\smallsetminus \{0\}$ coincides with the natural structure of it as a real analytic manifold. Hence there is a unique structure on $\overline V$ of a real analytic manifold with corners whose restriction to $\overline V \smallsetminus \{0\}$ is the structure which we just defined and whose restriction to $V$ is the natural structure of $V$ as a real analytic manifold.

 \end{sbpara}

 \begin{sbpara}\label{cgv4}  
The map $V^{(1)}\times \R_{\geq 0} \to \overline V\smallsetminus \{0\}$ in \ref{cgv2} extends to a surjective continuous map $V^{(1)}\times[0,\infty]\to\overline V$ which sends all points $(v,\infty)$ ($v\in V^{(1)}$) to $0\in\overline V$. Via this map, $\overline V$ is homeomorphic to the quotient of the compact space $V^{(1)}\times  [0, \infty]$  obtained by identifying all $(v, \infty)$ ($v \in V^{(1)}$). Hence $\overline V$ is compact. 

 \end{sbpara}

\subsection{Borel--Serre liftings}\label{ss:BSli}

  \begin{sbpara} {\it Borel--Serre lifting.}\label{BSli}  
Let $\cG$ be a semisimple algebraic group over $\R$. Let $\cP$ be a parabolic subgroup of $\cG$, and let 
$\cS_{\cP}$ be the largest $\R$-split torus in $\cP_{\red}=\cP/\cP_u$, where $\cP_u$ is the unipotent radical of $\cP$. Let 
$K$ be a maximal compact subgroup of $\cG(\R)$.

  Then we have a unique homomorphism $\cS_{\cP}\to \cP\;;\; a\mapsto a_K$ which lifts the inclusion  map $\cS_{\cP}\hookrightarrow \cP_{\red}$ and which satisfies $\theta_K(a_K)=a_K^{-1}$ for any $a\in \cS_{\cP}(\R)$, where $\theta_K$ is the Cartan involution of $\cG(\R)$ associated to $K$ (\cite{BS} Proposition 1.6).

  We call this $a_{K}$ the {\it Borel--Serre lifting of $a\in \cS_{\cP}$ at $K$}.

    \end{sbpara}

 \begin{sbrem}\label{orth} We remark that the Borel--Serre lifting is understood as the splitting of a filtration by taking the orthogonal complements.
 
 By \cite{Mo} (see also \cite{BH} Section 1), there is a finite-dimensional faithful representation $V$ of $\cG$ and a positive definite  symmetric $\R$-bilinear form $(\cdot, \cdot): V\times V\to \R$ which is fixed by $K$ such that $\cG$ is stable under the transpose $g\mapsto {}^tg$ for $(\cdot,\cdot)$. Furthermore, there are $\R$-subspaces of $V$ such that  $0=V_0\subset V_1 \subset \dots \subset V_m=V$ which are stable under $\cP$ satisfying $(g-1)V_j\subset V_{j-1}$ for all $g\in \cP_u$ and $1\leq j\leq m$. For example, we can take $V_j= \{v\in V\;|\; (g_1-1)\cdots (g_j-1)
 v=0\;\text{for all}\;g_1,\dots, g_j\in \cP_u\; \}$ (then $V_j$ are $\cP$-stable and $V_m=V$ for some $m$). 
  
  We have a commutative diagram
 $$\begin{matrix}    \cS_{\cP}(\R) & \subset & \cP_{\red}& \to & \prod_{j=1}^m \text{Aut}_\R(V_j/V_{j-1})\\
 \downarrow &&&& \downarrow\\
 \cP(\R) &&\overset{\subset}\longrightarrow && \text{Aut}_\R(V),\end{matrix}$$
 where the left vertical arrow is the Borel--Serre lifting at $K$ and the right vertical arrow is given by the orthogonal decomposition of the filtration $(V_j)_j$. That is, we have a unique decomposition $V=\bigoplus_{j=1}^m  V_{[j]}$ such that $V_j=\bigoplus_{k\leq j} V_{[k]}$ for all $j$ and such that $V_{[j]}$ and $V_{[k]}$ are orthogonal for $(\cdot,\cdot)$ if $j\neq k$. The right vertical arrow is given by $V_j/V_{j-1}\simeq V_{[j]}$.

  We prove that the above diagram is commutative, that is, $a_K$ for $a\in \cS_{\cP}(\R)$ preserves $V_{[j]}$. Note that  the Cartan involution of $\cG(\R)$ associated to $K$ is given by $\theta_K(g)=({}^tg)^{-1}$ for $g\in \cG(\R)$. Hence ${}^t(a_K)=\theta_K(a_K)^{-1}=a_K$. Since $V_{[j]}= (\bigoplus_{k \geq  j} V_{[k]})\cap  V_j$ and $a_K$ preserves $V_j$, it is sufficient to prove that $a_K$ preserves $\bigoplus_{k\geq j} V_{[k]}$. Note that $\bigoplus_{k\geq j} V_{[k]}$ is the 
 annihilator of $V_{j-1}$ for  $(\cdot,\cdot)$.  For $x\in \bigoplus_{k\geq j} V_{[k]}$ and $y\in V_{j-1}$, we have $(a_Kx, y)= (x, {}^t(a_K)y)= (x, a_Ky) =0$.
Hence $a_Kx\in \bigoplus_{k\geq j} V_{[k]}$.

\end{sbrem}

\subsection{Review of Borel--Serre theory}\label{ss:BS1}

  Let $G$ be a semisimple algebraic group over $\Q$. 
  Let $\frak X$ be the associated symmetric space as in \ref{X}.

  In this Section \ref{ss:BS1}, we review how the Borel--Serre space $\frak X_{\BS}\supset \frak X$ is constructed in the paper of Borel--Serre \cite{BS}.

  \begin{sbpara} \label{Xrt}
  Let $P$ be a parabolic subgroup of $G$. Let $S_P$ be the largest $\Q$-split torus in the center of $P_{\red}:=P/P_u$. 
  Let $X(S_P)$ be the character group of $S_P$. Let $X(S_P)^+$ be the submonoid of $X(S_P)$ generated by roots. Here a {\it root} means an element $\chi\in X(S_P)$ such that for some (equivalently, for any)  lifting $s: S_P\to P$ of the embedding $S_P\hookrightarrow P_{\red}$, there is a non-zero element $v$ of $\Lie(P_u)$ such that $\Ad(s(t))v= \chi(t)^{-1}v$ for $t\in S_P$. Then $X(S_P)^+$ is a free 
  monoid, that is, $X(S_P)^+\simeq \N^n$ for some $n\geq 0$. The basis $\Delta(P)$ of the monoid $X(S_P)^+$ is called the set of {\it fundamental roots} (or, of {\it simple roots}).

  Let $$A_P:=\Hom(X(S_P)^+, \R^{\mult}_{>0})=\R_{>0}^{\Delta(P)}\subset \bar A_P:=\Hom(X(S_P)^+, \R^{\mult}_{\geq 0})= \R_{\geq 0}^{\Delta(P)}.$$
  Since $X(S_P)^+$ generates a subgroup of $X(S_P)$ of finite index, the identification $S_P(\R)=\Hom(X(S_P), \R^\times)$ induces an isomorphism $S_P(\R)^{\circ}\overset{\sim}\to A_P$, where $S_P(\R)^{\circ}$ denotes the connected component of the topological group $S_P(\R)$ containing $1$.

  \end{sbpara}

  \begin{sbpara} {\it Borel--Serre action} ({\it geodesic action}, in the terminology of \cite{BS} 3.3).  Let $P$ be as in \ref{Xrt}. We have the action $\circ$ of the group $A_P$ on $\frak X$, which we call the Borel--Serre action,  defined as follows. 
  For $K\in \frak X$ and $a\in A_P$, let $a_K\in P(\R)$ be the Borel--Serre lifting of $a$ at $K$ obtained by applying \ref{BSli} to
  $\cG=G_\R$, $\cP=P_\R$ and $\cS_{\cP}\supset S_{P,\R}$. 
   Define $$a\circ K= a_K K a_K^{-1}.$$

  \end{sbpara}

  \begin{sbpara}\label{XBSs} 
   As a set, the Borel--Serre space $\frak X_{\BS}$ is defined by 
  $\frak X_{\BS}:= \{(P, Z)\}$, 
  where $P$ runs over parabolic subgroups of $G$ and $Z$ is an $A_P$-orbit for the Borel--Serre action. 
  \end{sbpara}

  \begin{sbpara}\label{XBSp} 
For a parabolic subgroup $P$ of $G$, let
$$
\frak X_{\BS}(P):=\{(Q,Z)\in\frak X_{\BS}\,|\,Q\supset P \}.
$$
Then we have a bijection
$$
\frak X_{\BS}(P)\simeq\frak X\times^{A_P} \bar A_P,\quad (Q,Z)\leftrightarrow(K,a),
$$
defined as follows. 

Let $Q$ be a parabolic subgroup of $G$ such that $Q\supset P$.
Then, $\Delta(Q)$ is regarded as a subset of $\Delta(P)$, $A_Q$ is regarded as a subgroup of $A_P$, and $Q\mapsto \Delta(Q)$ is a bijection from the set of all parabolic subgroups of $G$ such that $Q\supset P$ to the set of all subsets of $\Delta(P)$. 

For $(Q,Z)\in\frak X_{\BS}(P)$, $K$ is any element of $Z$ and $a\in\bar A_P$ is defined by $a(\chi)=0$ if $\chi\in\Delta(Q)$ and $a(\chi)=1$ if $\chi\notin\Delta(Q)$.
For $(K,a)\in\frak X \times \bar A_P$, $Q$ is the parabolic subgroup of $G$ containing $P$ such that $\Delta(Q)=\{\chi\in\Delta(P)\,|\,a(\chi)=0\}$ and $Z:=\{a'\circ K\,|\,a'\in A_P, \chi(a')=a(\chi)\;\text{for any}\;\chi\in\Delta(P)-\Delta(Q)\}$.
Note that, when $Q\supset P$, at a common $K$, the Borel--Serre action of $A_Q$ for $Q$ coincides with its action for $P$ regarding $A_Q$ as a subgroup of $A_P$.

  \end{sbpara}

  \begin{sbpara}\label{XBSr}
  
The set $\frak X_{\BS}$ has a structure of a real analytic manifold with corners defined as follows. 
  \medskip

 For a parabolic subgroup $P$ of $G$, there exists a real analytic map $f: \frak X\to A_P$  satisfying
\medskip

(1)  $f(a\circ K)=af(K)$ for all $a\in A_P$ and all $K\in\frak X$.
\medskip

We sketch the proof of this assertion.
  Let ${}^\circ P:=\bigcap_\chi\Ker(\chi^2:P\to{\mathbb G}_m)$, where $\chi$ runs over all homomorphisms $P\to\R^\times$ of algebraic groups over $\Q$.
Then, $P_u\sub\,{}^\circ P$ and the composition $A_P\to P(\R)/P_u(\R)\to P(\R)/{}^\circ P(\R)$ is an isomorphism by \cite{BS} 1.2, which we use by taking $P_{\red}$ as $G$ there.
Let $\pi:P(\R)\to P(\R)/{}^\circ P(\R)\simeq A_P$.
Fix $K\in\frak X$.
Since $G(\R)=P(\R)K$ (see \cite{B1} \S11) and since $\pi$ kills the compact subgroup $P(\R)\cap K$, there exists a unique map $G(\R)\to A_P$ sending $pk$ to $\pi(p)$ $(p\in P(\R), k\in K)$, which factors through $f:\frak X\simeq G(\R)/K\to A_P$.
Since the action of $a\in A_P$ on $G(\R)$ is $a\circ pk=pa_Kk$, $f$ satisfies the property $f(a\circ K)=af(K)$. 
\medskip

 The set  $\frak X_{\BS}(P)$ is regarded as a real analytic manifold with corners as follows. 
  
Let $\frak X^{(1)}:=f^{-1}(1)$.
Then we have

\medskip

(2)  $\frak X^{(1)}\times A_P \overset{\sim}\to \frak X\;;\; (K, a) \mapsto a \circ K$, an isomorphism of real analytic manifolds,
\medskip

(3) $\frak X^{(1)}\times \bar A_P \to \frak X\times^{A_P} \bar A_P = \frak X_{\BS}(P)$, a bijection of sets. 
\medskip

We regard  $\frak X_{\BS}(P)$ as a real analytic manifold with corners via the bijection (3). This structure of $\frak X_{\BS}(P)$ does not depend on the choice of $f$ because if $f': \frak X\to A_P$ satisfies the same condition as $f$ and if $Y:= \frak X^{(1)}=f^{-1}(1)$ and $Y':= (f')^{-1}(1)$, the map (3) and the map $Y'\times \bar A_P\to \frak X_{\BS}(P)$ are compatible with the isomorphism $Y\times \bar A_P\overset{\sim}\to Y'\times \bar A_P\;;\; (K, a)\mapsto (f'(K)^{-1}\circ K, f'(K)a)$ whose inverse is $(K, a)\mapsto (f(K)^{-1}K, f(K)a)$. 

  $\frak X_{\BS}$ has a unique structure of a real analytic manifold with corners such that for any parabolic subgroup $P$ of $G$,  $\frak X_{\BS}(P)$ is open in $\frak X_{\BS}$ and the structure of $\frak X_{\BS}(P)$ as a real analytic manifold with corners  defined above coincides with the restriction of that of $\frak X_{\BS}$. 
This follows from the following (4) and  (5). 
\medskip

(4)  Let $Q$ be a parabolic subgroup of $G$ such that $Q\supset P$.  Then $\frak X_{\BS}(Q)$ is an open subset of $\frak X_{\BS}(P)$ and the structure of $\frak X_{\BS}(Q)$ as a real analytic manifold with corners coincide with the restriction of that of $\frak X_{\BS}(P)$.
\medskip

(4) is shown as follows. Let $f: \frak X\to A_P$ be a real analytic map satisfying (1) and let $Y=f^{-1}(1)$. Let $p$ be the projection $A_P\to A_Q$. Then $f':= h \circ f: \frak X\to A_Q$ satisfies the condition on $f$ with $Q$ replacing $P$. Let $Y'=(f')^{-1}(1)$.  Let $T$ be the kernel of $p$. Then the map  (3) and the map $Y' \times \bar A_Q\to \frak X_{\BS}(Q)$ are compatible with the open immersion $Y' \times \bar A_Q \simeq Y \times T \times A_Q \overset{\subset}\to Y\times \bar A_P$, where the first isomorphism is induced by $Y\times T\overset{\sim}\to Y'\;;\; (K, a) \mapsto a \circ K$ and the second map is induced by $T\times \bar A_Q\to \bar A_P\;;\;(t, a) \mapsto ta$. 

(5)  For parabolic subgroups $P$ and $Q$ of $G$,  $\frak X_{\BS}(P)\cap \frak X_{\BS}(Q)=\frak X_{\BS}(P*Q)$. 
  Here $P*Q$ is the algebraic subgroup of $G$ generated by $P$ and $Q$, which is a parabolic subgroup of $G$. 
\end{sbpara}

\begin{sbpara}\label{proper} We will often use the following basic things about proper actions. Let $H$ be a Hausdorff topological group.

(1) If $H$ acts properly on a topological space $X$, the quotient space $H\bs X$ is Hausdorff.

(2) If $H$ is discrete and if $H$ acts properly on a Hausdorff space $X$, and if this action is free (that is, $hx=x$ with $h\in H$ and $x\in X$ implies $h=1$), the projection $X \to  H\bs X$ is a local homeomorphism.

(3) Assume that $H$ acts on topological spaces $X$ and $Y$, and let $f:X\to Y$ be an equivariant continuous map. 

(3.1) If $H$ acts properly on $Y$ and if $X$ is Hausdorff, then $H$ acts properly on $X$.

(3.2) If $H$ acts on $X$ properly and if the map $f$ is proper and surjective, then $H$ acts properly on $Y$.

Here in (3.2) and throughout this paper, as in \cite{KU2} and \cite{KNU2}, for a continuous map $f: X\to Y$ of topological spaces, $f$ is proper means that it satisfies the following (a) and (b). (a) For any topological space $Z$, the induced map $X\times_Y Z \to Z$ is a closed map. (b) The map $X \to X\times_Y X$ is a closed map. 

(4) Let $H_1$ be a closed normal subgroup of $H$ and let $H_2:= H/H_1$. Assume that $H$ (resp.\ $H_2$) acts on a topological space $X_1$ (resp.\ $X_2$) continuously. Assume that for $j=1,2$, the action of $H_j$ on $X_j$ is proper and free. Assume that there are a neighborhood $U$ of $1$ in $H_2$ and a continuous map $U\to H$ which lifts the inclusion map $U\to H_2$. Then the diagonal action of $H$ on $X_1\times X_2$ is proper and free.

(5) Let $H_1$ be a closed subgroup of $H$ of finite index. Assume that $H$ acts on a topological space $X$ continuously. If the action of $H_1$ on $X$ is proper, then the action of $H$ on $X$ is proper. 

\medskip

For proofs of (1), (2), (3), see \cite{Bou} Ch.3 \S4 no.2 Proposition 3, ibid.\ Ch.3 \S4 no.4 Corollary, ibid.\ Ch.3 \S4 no.2 Proposition 5, respectively. The proof of (4) is given in \cite{KNU2} Part III Definition 4.2.4. (5) is proved in an elementary manner. 
\end{sbpara}
  
\begin{sbpara}
\label{onXBS}

 Let $\Gamma$ be a semi-arithmetic subgroup of $G(\Q)$(\ref{Gamma}). Then:

$(1)$ The action of $\Gamma$ on $\frak X_{\BS}$ is proper, and the quotient space $\Gamma\bs \frak X_{\BS}$ is Hausdorff.

$(2)$ If $\Gamma$ is torsion-free,  the action of $\Gamma$ on $\frak X_{\BS}$ is free, and the map $\frak X_{\BS}\to \Gamma\bs \frak X_{\BS}$ is a local homeomorphism. 

$(3)$ If $\Gamma$ is an arithmetic subgroup of $G(\Q)$,  $\Gamma\bs \frak X_{\BS}$ is compact. 

\medskip

These are given in \cite{BS} 9.3 Theorem and 9.5  for arithmetic subgroups $\Gamma$ of $G(\Q)$. Since a semi-arithmetic subgroup of $G(\Q)$ is a subgroup of an arithmetic subgroup of $G(\Q)$,  the properness in (1) is reduced to the case of an arithmetic subgroup. The Hausdorffness in (1) follows from the properness in (1) (\ref{proper} (1)). We prove (2). We may assume that $\Gamma$ has a subgroup $\Gamma_1$ of finite index which is a subgroup of a torsion-free arithmetic subgroup. Assume that $\gamma \in \Gamma$ fixes some point of $\frak X_{\BS}$. Take $n\geq 1$ such that $\gamma^n\in \Gamma_1$. By the result for an arithmetic subgroup, we have $\gamma^n=1$. Since $\Gamma$ is torsion-free, we have $\gamma=1$. The rest of (2) follows from the properness in (1) and from this free property (\ref{proper} (2)).

\end{sbpara}
    
\subsection{The set $D_{\BS}$}\label{ss:DBS}
  
  Let $G$ be a linear algebraic group over $\Q$ and let $h_0: S_{\C/\R}\to G_{\red,\R}$ be a homomorphism as in \ref{D} which is $\R$-polarizable. 
   
  \begin{sbpara} \label{para}
   Let $G'_{\red}=[G_{\red},G_{\red}]$ be the commutator subgroup of $G_{\red}=G/G_u$. This is a semisimple algebraic group. 

  We have bijections between the sets
  $$(\text{parabolic subgroups of} \; G)\leftrightarrow(\text{parabolic subgroups of}\; G_{\red})\leftrightarrow(\text{parabolic subgroups of}\;G'_{\red})$$
  given as follows. The bijection from the second set to the first set is to take the inverse image in $G$. The bijection from the second set to the third set is to take the intersection with $G'_{\red}$. 
  
  \end{sbpara}
  
  \begin{sbpara}\label{Drt}
  Let $P$ be a parabolic subgroup of $G_{\red}$.

 Let $P':= P\cap G'_{\red}$. Then $P'$ is a parabolic subgroup of $G'_{\red}$. We will denote $A_{P'}$ by $A_P$. This $A_P$ is also described as follows.

  Let $S_P$ be the largest $\Q$-split torus in the center of $P_{\red}$.  Let  $X(S_P)^+$ be the submonoid of the character group $X(S_P)$ generated by the inverses of characters which appear in the adjoint action on $\Lie(P_u)$ of a lifting of $S_P$ in $P$ (see \ref{Xrt}). Then the canonical map $X(S_P)\to X(S_{P'})$ induces an isomorphism $X(S_P)^+\overset{\sim}\to X(S_{P'})^+$ and 
  a bijection $\Delta(P)\to \Delta(P')$ between the bases. Hence we can write 
  $$A_P = \Hom(X(S_P)^+,\R^{\mult}_{>0})= \R_{>0}^{\Delta(P)} \subset \bar A_P=\Hom(X(S_P)^+, \R^{\mult}_{\geq 0})= \R_{\geq 0}^{\Delta(P)}.$$

     \end{sbpara}

  \begin{sbpara} \label{Aact} 
  We define the {\it Borel--Serre action of $A_P$ on $D$} as follows. Let $a\in A_P$ and $x\in D$. Let $(p, s, \delta)\in D_{\red}\times \spl(W) \times \cL$, with $\delta \in \cL(p)$, be the element corresponding to $x$ (Proposition \ref{Dandgr}). 
  Let $K(p)\in \frak X$ be the image of $p$ under $D_{\red}\to \frak X$ (\ref{DtoX}). We have the Borel--Serre lifting $A_P\to P'(\R)\;;\; a\mapsto a_{K(p)}$ associated to $K(p)$. Then $a \circ x$ is defined to be the element of $D$ whose image in $D_{\red}\times \spl(W)\times \cL$ is $(a_{K(p)}p, s, \Ad(a_{K(p)})\delta)$. 
  
  \end{sbpara}
  
  \begin{sbpara}\label{Bact}
  Let $B_P= \R_{>0} \times A_P$. Then we define the action of $B_P$ on $D_{\nspl}$ (see \ref{Dspl}) as follows. Let $b=(t,a)\in \R_{>0}\times A_P$. Then, for $x\in D_{\nspl}$, $b\circ x:= a \circ (t \circ x)$, where  $t\circ x\in D$ corresponds to the element $(p, s, t \circ \delta)$ of $D_{\red}\times \spl(W) \times \cL$ (see \ref{cgv1}). 
  Here $(p,s,t)$ is the element corresponding to $x$ by Proposition \ref{Dandgr}.
\end{sbpara}
  
  \begin{sbpara}\label{DBSs}
  Let $D_{\BS}$ be the set of all pairs $(P, Z)$, where $P$ is a parabolic subgroup of $G_{\red}$ and $Z$ is either an $A_P$-orbit in $D$ or a $B_P$-orbit in $D_{\nspl}$ for the Borel--Serre action. 
  
  We denote by $D^{\mild}_{\BS}$ the part of $D_{\BS}$ consisting of $A_P$-orbits and by $D_{\nspl,\BS}$ the subset of $D_{\BS}$ consisting of all elements of the form $(P,Z)$ such that $Z\subset D_{\nspl}$.
  ($D^{\mild}_{\BS}$ and $D_{\nspl,\BS}$ correspond to $D^{(A)}_{\BS}$ and $D^{(B)}_{\BS}$ in the notation of \cite{KNU2} Part I 8.1, respectively.)
  \end{sbpara}

\subsection{The real analytic structure of $D_{\BS}$}\label{ss:BSan}

We are in the same setting as in \ref{ss:DBS}.
We endow $D_{\BS}$ with a structure of a real  analytic manifold with corners. 

 \begin{sblem}\label{BSf1}  Let $P$ be a parabolic subgroup of $G$. 
 Then there is a real analytic map $f: D\to A_P$ such that $f(a \circ x)= af(x)$ for any $a\in A_P$ and $x\in D$.

 \end{sblem}  
 
 \begin{pf}  Take a real analytic map $f_{\frak X}: \frak X\to A_P$ such that $f_{\frak X}(a\circ K)=af_{\frak X}(K)$ for any $a\in A_P$ and $K\in \frak X$ (\ref{XBSr}), and define $f$ to be the composite map $D\to D_{\red}\to \frak X\overset{f_{\frak X}}\to A_P$, where the second arrow is as in \ref{DtoX}.  
  \end{pf}
  \begin{sblem}\label{BSf2} Let $P$ be a parabolic subgroup of $G$. Then there is a real analytic map $f: D_{\nspl} \to B_P$ such that $f(b \circ x)=bf(x)$ for any $b\in B_P$ and $x\in D_{\nspl}$. 
 
 \end{sblem}
 
 \begin{pf}
  By Lemma \ref{BSf1}, there is a real analytic map $f_D: D_{\red}\to A_P$ such that $f_D(a\circ p)=af_D(p)$ for any $a\in A_P$ and any $p\in D_{\red}$.
 
  Let $\cL=W_{-2}\gr^W\Lie(G)_{\bR}$ be as in \ref{CKSdelta}. It is a graded $\R$-vector space with weights $\le-2$. By \ref{cgv2}, there is a real analytic map
  $f_{\cL}:\cL\smallsetminus\{0\}\to\R_{>0}$  such that $f_{\cL}(t \circ \delta)= tf_{\cL}(\delta)$ for any $t\in\R_{>0}$ and $\delta\in \cL\smallsetminus \{0\}$. 
  Define $f:D_{\nspl}\subset  D_{\red}\times \spl(W) \times (\cL\smallsetminus \{0\})\to B_P=A_P\times \R_{>0}$ by
 $f(p, s,\delta) = (f_D(p), f_{\cL}(\mathrm{Ad}(f_D(p)_p)^{-1}\delta))$ $(p\in D_{\red}, s \in \spl(W), \delta\in \cL(p))$, 
where $f_D(p)_p$ is the Borel--Serre lifting of  $ f_D(p) \in A_P$ at $p$.  
  This $ f $ satisfies $f(b \circ x)= bf(x)$. 
 \end{pf}
 
 \begin{sbpara}\label{DBSr} 
   
   Let $P$ be a parabolic subgroup of $G$.
Let $D_{\BS}(P):=\{(Q,Z)\in D_{\BS}\,|\,Q\supset P\}$.

Let $D_{\BS}^{\mild}(P):=D_{\BS}(P)\cap D_{\BS}^{\mild}$ and $D_{\nspl,\BS}(P):=D_{\BS}(P)\cap D_{\nspl,\BS}$.
Then, in the same way as in the case of $\frak X_{\BS}(P)$, we have
$$
D_{\BS}^{\mild}(P)=D\times^{A_P} \bar A_P,\quad D_{\nspl, \BS}(P)= D\times^{B_P} \bar B_P.
$$

   We endow $D_{\BS}(P)$ with the structure of a real analytic manifold in the following way. 
 
 First,   the set  $U_1:=D^{\mild}_{\BS}(P)$ is regarded as a real analytic manifold with corners as follows. 
  By Lemma \ref{BSf1}, there is a real analytic map $f_A: D\to A_P$  satisfying

\medskip

(1) $f_A(a\circ x)=af_A(x)$ for all $a\in A_P$ and all $x\in D$.
\medskip

\noindent
Let $D_A^{(1)}=f_A^{-1}(1)$. Then 
\medskip

(2) $D_A^{(1)} \times A_P \overset{\sim}\to D\;;\; (x, a) \mapsto a \circ x$.
\medskip

\noindent
This  map  (2) induces a bijection
\medskip

(3) $D_A^{(1)}\times \bar A_P \to D\times^{A_P} \bar A_P = D^{\mild}_{\BS}(P)$. 
\medskip

\noindent
We regard  $U_1=D^{\mild}_{\BS}(P)$ as a real analytic manifold with corners via the bijection  (3). This structure of $D_{\BS}^{\mild}(P)$ does not depend on the choice of $f_A$. 

Next, the set $U_2:= D_{\nspl,\BS}(P)$ is regarded as a real analytic manifold with corners as follows. By Lemma \ref{BSf2}, there is a real analytic map $f_B: D_{\nspl}\to B_P$  satisfying

\medskip

(4) $f_B(b\circ x)=bf_B(x)$ for all $b\in B_P$ and all $x\in D_{\nspl}$.
\medskip

\noindent
Let $D_B^{(1)}=f_B^{-1}(1)$. Then 
\medskip

(5) $D_B^{(1)}\times B_P \overset{\sim}\to D_{\nspl}\;;\; (x, b) \mapsto b \circ x$.
\medskip

\noindent
This  map  (5) induces a bijection
\medskip

(6) $D_B^{(1)}\times \bar B_P \to D\times^{B_P} \bar B_P = D_{\nspl, \BS}(P)$. 
\medskip

\noindent
We regard $U_2:=D_{\nspl, \BS}(P)$ as a real analytic manifold with corners via the bijection  (6). This structure of $D_{\nspl,\BS}(P)$ does not depend on the choice of $f_B$.

It is easy to prove that for $j=1,2$ and for the structure of a real analytic manifold with corners on $U_j$,  the intersection $U_1\cap U_2$ is an open set of $U_j$ and that  the restriction of the structure of the real analytic manifold with corners of $U_1$ to $U_1\cap U_2$ coincides with that of $U_2$. Hence there is a unique structure of a real analytic manifold with corners on $U_1\cup U_2=D_{\BS}(P)$ for which $U_1$ and $U_2$ are open sets and whose restriction to $U_j$ coincides with that of $U_j$ for $j=1,2$.

\end{sbpara}

\begin{sbprop}\label{DBSr2}  The set $D_{\BS}$ has a unique structure of a real analytic manifold with corners for which $D_{\BS}(P)$ is an open set and whose restriction to $D_{\BS}(P)$ coincides with that of $D_{\BS}(P)$ defined in $\ref{DBSr}$ for any parabolic subgroup $P$ of $G$. 

\end{sbprop}

\begin{pf} It is easy to see that for parabolic subgroups $P$ and $Q$ of $G$ such that $Q\supset P$, $D_{\BS}(Q)$ is an open set of $D_{\BS}(P)$ and the restriction of the structure of $D_{\BS}(P)$ as a real analytic manifold with corners to $D_{\BS}(Q)$ coincides with that of $D_{\BS}(Q)$. Furthermore, $D_{\BS}(P)\cap D_{\BS}(Q) = D_{\BS}(P*Q)$. This proves Proposition \ref{DBSr2}. 
  \end{pf}

  \begin{sbrem}\label{ZBJ0}
  When $G$ is semisimple, the Borel--Serre partial compactification $\overline{G(\R)}_{\BS}$ in \cite{BJ2} Proposition 6.3 coincides with our $D(G,h_0)_{\BS}$ in Proposition \ref{DBSr2} as real analytic manifolds with corners.
\end{sbrem}

\begin{sbpara} We have a morphism   $$D_{\BS}\to \frak X_{\BS}$$ of real analytic manifolds with corners induced by $D_{\red,\BS}\to \frak X_{\BS}\;;\; (P, Z) \mapsto (P, Z')$, where $Z'$ is the image of $Z\subset D_{\red}$ under $D_{\red}\to \frak X$ (\ref{DtoX}). 

\end{sbpara}      
 \begin{sbprop}\label{DXBS1}  If $G$ is reductive, the map $D_{\BS}\to \frak X_{\BS}$ is proper. 
  \end{sbprop}
 
 \begin{pf} 
Since $G$ is reductive, we have $D_\BS(P)=D_\BS^\mild(P)$.
Hence it is sufficient to prove that $D_{\BS}(P)\to \frak X_{\BS}(P)$ is proper. 
Since $D_{\BS}(P)=D_A^{(1)}\times \bar A_P$ (\ref{DBSr} (3)) and $\frak X_{\BS}(P)= \frak X^{(1)}\times \bar A_P$ (\ref{XBSr} (3)), the properness of $D_{\BS}(P)\to \frak X_{\BS}(P)$ is reduced to that of $D_A^{(1)}\to \frak X^{(1)}$.
The latter is reduced to the properness of $D\to \frak X$ by \ref{DBSr} (2) and \ref{XBSr} (2).
Take a point $p\in D$, let $K'$ be the isotropic subgroup of $G(\R)$ at $p$, and let $K$ be the isotropic subgroup of $G(\R)$ at the image of $p$ in $\cX$.
Then $K$ is a maximal compact subgroup of $G(\R)$ and $K'\sub K$ is a compact subgroup.
Hence $D=G(\R)/K'\to\frak X=G(\R)/K$ is proper.
\end{pf}

\begin{sbpara}\label{BSpsd}

Let $P$ be a parabolic subgroup of $G$.  Fix $f: D_{\red}\to A_P$, which is induced from Lemma \ref{BSf1}, let $Y=f^{-1}(1)\subset D_{\red}$ and let $b_f$ be the composition $D_{\red,\BS}(P)=Y \times \bar A_P\to Y\subset D_\red$ (see \ref{para}). Consider the isomorphism of real analytic manifolds  $$D \overset{\sim}\to \{(p, s, \delta)\in D_{\red}\times \spl(W) \times \cL\;|\; \delta\in \cL(b_f(p))\}\;;\;x\mapsto (p, \spl_W(x), \Ad(f(p)_p)^{-1} \delta(x)),$$ where $p=x_{\red}$ and $f(p)_p$ is the Borel--Serre lifting of $f(p)$ to $P$ at $p$. 

\begin{sbprop} This isomorphism  extends uniquely to an isomorphism of real analytic manifolds with corners
$$D_{\BS}(P)\overset{\sim}\to \{(p, s, \delta)\in  D_{\red,\BS}(P)\times \spl(W) \times \overline{\cL}\;|\;\delta\in \overline{\cL}(b_f(p))\}.$$
\end{sbprop}

  This map sends $x=(Q,Z)$ to $(p, s, \delta)$ defined as follows.
$p=(Q, Z_{\red})$. $s=\spl_W(y)$ for any $y\in Z$, which is independent of the choice of $y$. $\delta= \Ad(f(z)_z)^{-1}(\delta(y))$, where  $y\in Z$ and $z=y_{\red}$, which is independent of the choice of $y$. 

  By this map, the image of $x\in D_{\BS}(P)$ in $\overline{\cL}$ belongs to $\cL$ if and only if $x$ is an $A_P$-orbit.
\end{sbpara}

 \begin{sbprop}\label{BSbarL}
 The map $D_{\BS}\to D_{\red, \BS}\times \spl(W)$ is proper. It is an $\bar L$-bundle.  
Note that $\bar L$ here is $\overline{\cL}(p)$ for any $p$.
 \end{sbprop}

We hope that our notation $D_{\red, \BS}$ is not confusing with the reductive Borel--Serre space in \cite{Z}.

\subsection{Global properties of $D_{\BS}$}
\label{ss:BSproperty}

We are still in the setting in Section \ref{ss:DBS}.

\begin{sbthm}\label{BSgl} Let $\Gamma$ be a semi-arithmetic subgroup of $G'(\Q)$ ($\ref{Gamma}$). 

$(1)$ The action of $\Gamma$ on $D_{\BS}$ is proper and the quotient space $\Gamma\bs D_{\BS}$ is Hausdorff.

$(2)$ If $\Gamma$ is torsion-free,  the action of $\Gamma$ on $D_{\BS}$ is free and  the map $D_{\BS}\to \Gamma\bs D_{\BS}$ is a local homeomorphism. 

$(3)$ If $\Gamma$ is an arithmetic subgroup of $G'(\Q)$,  $\Gamma\bs D_{\BS}$ is compact. 
\end{sbthm} 

\begin{sbpara}\label{p:DBSHaus} We first prove 
that  $D_{\BS}$ is Hausdorff, that is, the case $\Gamma=\{1\}$ of \ref{BSgl} (1).  

 The map $D_{\BS}\to \frak X \times \spl(W)$ is proper by Propositions \ref{DXBS1} and \ref{BSbarL}. Since $\frak X$ is Hausdorff (\ref{onXBS}) as well as $\spl(W)$ which is isomorphic to $G_u(\R)$, we have that $D_{\BS}$ is Hausdorff.

\end{sbpara}
  
  \begin{sbpara}\label{neat}
We say that a subgroup $\Gamma$ of $G(\Q)$ is {\it neat} if it satisfies the following condition (1).

(1) There is a faithful representation $V\in \Rep(G)$ such that for every element $\gamma\in \Gamma$, the subgroup of $\C^\times$ generated by all the eigenvalues of  the action of $\gamma$ on $V_\C$ is torsion-free.

By \ref{faith}, the condition (1) is equivalent to the following condition (1$^{\prime}$). 

(1$^{\prime}$) For every  $\gamma\in \Gamma$, every $n \geq 1$, and every homomorphism $\rho: G\to \GL(n)$, the subgroup of $\C^\times$ generated by all the eigenvalues of  $\rho(\gamma)$  is torsion-free.

A neat subgroup of $G(\Q)$ is torsion-free.

Every  semi-arithmetic subgroup $\Gamma$  (\ref{Gamma}) has a neat subgroup of finite index (\cite{B1} 17.4).
\end{sbpara}
  
\begin{sbpara}\label{pfBSgl}
  We prove Theorem \ref{BSgl}.
   Let $\Gamma$ be a  semi-arithmetic subgroup of $G'(\Q)$ (\ref{Gamma}).  

  The proof is similar to that of \cite{KNU2} Part I Theorem 9.1.  

By \ref{proper} (5), the properness in (1) is reduced to the case where $\Gamma$ is neat. We prove the properness in (1) and the free property in (2) assuming $\Gamma$ is neat. We apply \ref{proper} (4) to $H=\Gamma$, $H_1=\Gamma_u:= \Gamma\cap G_u(\Q)$, $X_1= \spl(W)$, $X_2=\frak X$. The action of $\Gamma/\Gamma_u$ on $\frak X$ is proper and free by (1) and (2) of \ref{onXBS}, and the action of $\Gamma_u$ on $\spl(W)$ is proper and free because $\spl(W)\simeq G_u(\R)$ on which $\Gamma_u$ acts through the inclusion $\Gamma_u\to G_u(\R)$. Hence the action of $\Gamma$ on $\spl(W) \times \frak X$ is proper and free. By using the canonical continuous equivariant map $D_{\BS}\to \spl(W) \times \frak X$, we have that the action of $\Gamma$ on $D_{\BS}$ is proper and free by \ref{p:DBSHaus} and \ref{proper} (3.1).

The Hausdorffness in (1) follows from the properness in (1)  (\ref{proper} (1)) and the 
 local homeomorphism property in (2) follows from it and from the properness in (1)  (\ref{proper} (2)).

   (3) follows from the compactness of $(\Gamma/\Gamma_u)\bs \frak X_{\BS}$, where $\Gamma_u=\Gamma\cap G_u(\Q)$,  the compactness of $\Gamma_u\bs G_u(\R)$, and Proposition \ref{BSbarL}. 
   \end{sbpara}

   In Theorem \ref{BSgl}, we can use  an arithmetic subgroup and a semi-arithmetic subgroup of $G(\Q)$ (not of $G'(\Q)$) 
  in the following situations in Remark \ref{BSSA1} and Proposition \ref{BSSA2}.

\begin{sbrem}\label{BSSA1} 

 If either $G$ is semisimple or  the condition (1) in Lemma \ref{pol2} is satisfied,   \ref{BSgl} holds for a semi-arithmetic subgroup $\Gamma$ of $G(\Q)$. In fact, $\Gamma \cap G'(\Q)$ is of finite index (cf. Proposition \ref{A=A*} for the latter case). Hence by \ref{proper} (5), we can replace $\Gamma$ by the semi-arithmetic subgroup $\Gamma\cap G'(\Q)$ of $G'(\Q)$. 

\end{sbrem}

\begin{sbprop}\label{BSSA2}
  Assume that $G$ is reductive. Let $\Gamma$ be a subgroup of $G(\Q)$. Then$:$
 
$(1)$ If $\Gamma$ is a semi-arithmetic subgroup of $G(\Q)$, $\Gamma \bs D_{\BS}$ is Hausdorff.

$(2)$ Let $Z$ be the center of $G$. If $\Gamma$ is a semi-arithmetic subgroup of $G(\Q)$ and the image of $\Gamma$ in $(G/Z)(\Q)$ is torsion-free, then the map $D_{\BS}\to \Gamma\bs D_{\BS}$ is a local homeomorphism.

$(3)$ If $\Gamma$ is an arithmetic subgroup of $G(\Q)$,  $\Gamma \bs D_{\BS}$ is compact.
\end{sbprop}

See \ref{redpf} for the proof. 

\section{The space of SL(2)-orbits}\label{s:DSL}
  Let $G$ be a linear algebraic group defined over $\Q$.
  We define and study the extended period domains $D_{\SL(2)}$ and  $D^{\star}_{\SL(2)}$ of $\SL(2)$-orbits which contain $D$. This $D_{\SL(2)}$ (resp.\ $D^{\star}_{\SL(2)}$) is the  $G$-MHS version of the space $D(\La)_{\SL(2)}$ 
(resp.\  $D(\La)^{\star}_{\SL(2)}$) 
for the classical period domain $D(\La)$  (\ref{ss:Dold}) defined and studied in \cite{KU1} (in the pure case) and  \cite{KNU2} Part II (in the mixed case) (resp.\ in \cite{KNU2} Part IV). 

  Assume that we are given an element $h_0:S_{\C/\R}\to G_{\red,\R}$ of $\Psi_H(G)$ (\ref{kandh}). We assume the  $\R$-polarizability (\ref{pol}). 
  We denote by $k_0 \in \Psi_W(G)$ the homomorphism $\bG_m\to G_{\red}$ (defined over $\Q$) which induces $h_0\circ w: \bG_{m,\R}\to G_{\red,\R}$.

\subsection{The set $D_{\SL(2)}$ when $G$ is reductive}
\label{ss:red}

In this Section \ref{ss:red}, we assume that $G$ is reductive.

\begin{sbpara}\label{red1} Let $n\geq 0$. By an {\it $\SL(2)$-orbit in $n$ variables}, we mean a pair  $(\rho, \varphi)$  $$\rho: \SL(2)_\R^n\to G_\R, \quad \varphi: \bP^1(\C)^n\to \check D,$$
where $\rho$ is a homomorphism of algebraic groups over $\R$ and $\varphi$ is a holomorphic map  satisfying the following conditions (i)--(iii).

(i)   $\varphi(gz)= \rho(g)\varphi(z)$  \; ($g\in \SL_2(\C)^n,  \; z\in \bP^1(\C)^n$).

(ii)   $\varphi(\frak H^n)\sub D$.

(iii) For $z\in \frak H^n$, the homomorphism $\Lie(\rho): \frak{sl}(2,\R)^n \to \Lie(G_\R)=\Lie(G)_\R$ of Lie algebras induced by $\rho$ is a homomorphism of $\R$-Hodge structures if we endow $\frak{sl}(2,\Q)^n$ with the $\Q$-Hodge structure of weight $0$ associated to $z$ (see below) and endow $\Lie(G)$ with the $\Q$-Hodge structure of weight $0$ associated to $\varphi(z)\in D$ and the adjoint action of $G$ on $\Lie(G)$. 

   Here the $\Q$-Hodge structure of $\frak{sl}(2,\Q)^n$ associated to $z=(z_j)_{1\leq j\leq n}$ is the direct sum of the $\Q$-Hodge structure of the $j$-th $\frak{sl}(2,\Q)$ 
induced by the $\Q$-Hodge structure of $V_\C=\C^2$ corresponding   to $z_j$ 
 (e.g.\ \cite{KU2} 1.2.3) and the adjoint action of $\GL(2)_{\Q}$ on $\frak{sl}(2,\Q)$.
\end{sbpara} 
   
\begin{sbpara}\label{clEx2}  (Classical) {\it Example}.
Assume $G=\GL(2)$ and let $h_0$ be the homomorphism $z\mapsto \langle z\rangle$ in \ref{clEx}. Let $\rho: \SL(2)_\R\to G_\R$ be the inclusion map and let $\varphi: \bP^1(\C)\to \Dc$ be the identity map in \ref{clEx}. Then $(\rho, \varphi)$ is an $\SL(2)$-orbit in one variable. 
\end{sbpara}

\begin{sblem}\label{redp1} Let $n\geq 0$. Then the following sets {\rm (i)}, {\rm (ii)}, {\rm (iii)}, {\rm (iii)}${}'$, and {\rm (iv)} can be identified. 

{\rm (i)} The set of all $\SL(2)$-orbits $(\rho, \varphi)$ in $n$ variables. 

{\rm (ii)} The set of all pairs $(\rho, \varphi)$, where $\rho$ is a homomorphism $\SL(2)_\R^n \to G_\R$ of algebraic groups over $\R$ and $\varphi$ is a holomorphic map  $\frak H^n \to D$ satisfying the condition {\rm (i)} in $\ref{red1}$ for $g\in \SL_2(\R)^n$ and $z \in \frak H^n$ and satisfying
the condition {\rm (iii)} in $\ref{red1}$.

{\rm (iii)} The set of  all pairs $(\rho, \br)$,  where $\rho$ is a homomorphism $\SL(2)_\R^n \to G_\R$ of algebraic groups over $\R$ and $\br\in D$ satisfying the following condition.

The homomorphism $\Lie (\rho) : \frak{sl}(2, \R)^n\to  \Lie (G_\R)$  is a homomorphism of $\R$-Hodge structures with respect to the Hodge structure of $\frak{sl}(2, \R)^n$ associated to $\bi$ and the Hodge structure of $\Lie(G_\R)$ associated to $\br$. 

{\rm (iii)}${}'$ The set of all pairs $(\rho, \br)$, where $\rho$ is a homomorphism $\SL(2)_\R^n \to G_\R$ of algebraic groups over $\R$ and $\br\in D$ satisfying the following condition.

Let $S^{(1)}_{\C/\R}$ be the kernel of norm $S_{\C/\R} \to \bG_{m,\R}; z\mapsto z\bar z$. Let $\xi_1: S_{\C/\R}^{(1)}\to G_\R$ be the homomorphism defined by
$$\xi_1(u):= \br(u)\rho(\langle u\rangle, \dots, \langle u\rangle)^{-1},$$
where we regard $\br$ as a homomorphism $S_{\C/\R}\to G_\R$ and $\langle u\rangle\in \mathrm{SO}(2)_{\R}$ is as in $\ref{clEx}$. 
Then $\xi_1(u)\rho(g)=\rho(g)\xi_1(u)$ for any $u\in S^{(1)}_{\C/\R}$ and $g\in \SL(2)_\R^n$. 

{\rm (iv)} The set of all homomorphisms $\SL(2)_\R^n \times S^{(1)}_{\C/\R}\to G_\R$ satisfying the following conditions {\rm (iv.1)} and {\rm (iv.2)}. Let $\rho$ (resp.\ $\xi_1$) be the restriction of this homomorphism to 
$\SL(2)_\R^n$ (resp.\ $S^{(1)}_{\C/\R}$). 

{\rm (iv.1)} $k_0(-1)= \xi_1(-1)\rho(-1,\dots, -1)$. Here $-1$ inside $\rho(-)$ denotes the scaler matrix $-1$ in $\SL(2)_\R$. 

{\rm (iv.2)} Define the homomorphism $\br: S_{\C/\R}\to G_\R$ by
$\br(t u):= k_0(t) \xi_1(u)\rho(\langle u\rangle, \dots, \langle u\rangle)$,
where $t\in \bG_{m,\R}\subset S_{\C/\R}$ and $u\in S^{(1)}_{\C/\R}$. Then $\br \in D$. 
\end{sblem}

\begin{pf}
  We will prove that the first three sets are identified. 
  In fact, the bijectivities of the evident maps (i) $\to$ (ii) $\to $ (iii) are well-known 
(for their 1-variable cases, see \cite{Sc}, \cite{U2} Remark (2.2), and \cite{Satake}), but we give a proof here. 
  To prove that (i)$\to$ (ii) $\to $ (iii) are bijective, 
since they are injective, 
it is sufficient to prove that (i) $\to$ (iii) is surjective.  Let $(\rho, \br)$ be an element of the set (iii). We prove that if $g\in H:=\{g\in\SL(2,\C)^n\;|\; g\bi = \bi\}$,  then $\rho(g)\br=\br$. This will imply that we have a map $\varphi: \bP^1(\C)^n \to \check{D}\;;\; g\bi \mapsto \rho(g)\br$ having the desired properties. 
Since $H$ is connected, it is sufficient to prove that for $X\in \Lie(H)=F^0_{\bi}{\frak {sl}}(2, \C)^n$, $\Lie(\rho)(X)$ respects the Hodge filtration of $\br$. But this follows from $\Lie(\rho)(X)\in F^0_{\br} \Lie(G)_{\bC}$. 

We prove (iii) $=$ (iii)${}'$.  Since the Hodge structures on the Lie algebras in the condition in the definition of (iii) are of weight $0$,
the condition in (iii) is equivalent to the condition that for $u\in S^{(1)}_{\C/\R}(\R)$,  the action  $\Ad(\langle u\rangle,\dots, \langle u\rangle)$ on $\frak{sl}(2,\R)^n$ and the adjoint action of $\br(u)$ on $\Lie(G_\R)$ are compatible via $\Lie(\rho)$. This is equivalent to the condition that the action $\Ad(\xi_1(u))$ on $\Lie(\rho)(\frak{sl}(2,\R)^n)$ is trivial for every $u$ and hence to the condition that 
$\xi_1(u)$ and $\rho(g)$ commute for every $u$ and $g\in \SL(2)^n_\R$. 

The evident map (iii)${}'$ $\to$ (iv) is bijective because the inverse map is given as in (iv.2). 
Thus we can identify the sets (i), (ii), (iii), (iii)${}'$, and (iv). 
\end{pf}

\begin{sbprop}\label{redp2} Let $(\rho, \varphi)$ be an $\SL(2)$-orbit in $n$ variables and let $\br=\varphi(\bi): S_{\C/\R}\to G_\R$ and $\xi_1:S^{(1)}_{\C/\R}\to G_\R$ be as in Lemma $\ref{redp1}$. 
  Define a homomorphism $\xi: S_{\C/\R}\to G_\R$ of algebraic groups over $\R$ 
by $$\xi(t):=k_0(t) \rho(\mathrm{diag}(1/t, t), \dots, \mathrm{diag}(1/t, t))\quad\text{for}\; t\in \bG_{m,\R}\subset S_{\C/\R},$$
$$\xi(u):=\xi_1(u)\quad \text{for}\;  u\in S_{\C/\R}^{(1)},$$
 where $\xi_1$ is as in Lemma $\ref{redp1}$ {\rm (iv)}.  Then the Hodge filtration of $\xi$ is $\varphi(\zero)$, where $\zero=(0, \dots, 0)\in \bP^1(\C)^n$.

\end{sbprop}

\begin{pf}

 To relate $\xi$ and $\br$,  we use the Cayley element 
 $$
\bc:=\frac{1}{\sqrt{2}}\begin{pmatrix}1&i\\i&1\end{pmatrix}\in\SL(2,\C).
$$
We use the following properties (i) and (ii) of $\bc$.

(i) By the action of $\SL_2(\C)$ on $\bP^1(\C)$, $\bc$ sends $0\in \bP^1(\C)$ to $i\in \bP^1(\C)$.

(ii) $\langle u\rangle = \bc\, \text{diag}(u^{-1},u)\,\bc^{-1}$ for $u\in S^{(1)}_{\C/\R}$. 

Let $V\in \Rep(G)$, let $V_{\C,\xi}$ (resp.\ $V_{\C,\br}$) be $V_\C$ endowed with the action of $S_{\C/\R}$ via $\xi$ (resp.\ $\br$), and let $V_{\C,\xi}^{p,q}$ (resp.\ $V_{\C,\br}^{p,q}$) be the Hodge $(p,q)$-component of $V_{\C, \xi}$ (resp.\ %
$V_{\C,\br}$). 

\smallskip

{\bf Claim.}  $\rho(\bc, \dots, \bc)(\bigoplus_q V_{\C, \xi}^{p,q})= \bigoplus_q V_{\C, \br}^{p,q}$ in $V_\C$.

\smallskip

We prove Claim. For $(a,b,c)\in \Z^3$, let $V_\C^{(a,b, c)}$ be the part of $V_\C$ on which $k_0(t)$ for $t\in \C^\times$ acts as %
$t^a$, $\xi_1(u)$ for $u\in S^{(1)}_{\C/\R}$ acts as 
 $u^b$, and $\rho(\text{diag}(1/t, t), \dots, \text{diag}(1/t, t))$ for $t\in \C^\times$ acts as $t^c$. 
 Then $V_\C$ is the direct sum of these $V_\C^{(a,b,c)}$. On $V_\C^{(a,b,c)}$, $\xi(tu)$ ($t\in \bG_{m,\R}\subset S_{\C/\R}$, $u\in S^{(1)}_{\C/\R}$) acts as $t^{a+c}u^b$ and hence $V_\C^{(a,b,c)}$ has the Hodge type $((a+b+c)/2, (a-b+c)/2)$. On the other hand, by the above property (ii) of $\bc$, $\br(tu)$ acts on $\rho(\bc, \dots, \bc)V_\C^{(a,b,c)}$ as $t^au^{b+c}$ and hence $\rho(\bc, \dots, \bc)V_\C^{(a,b,c)}$ has the Hodge type $((a+b+c)/2, (a-b-c)/2)$. This proves Claim.
 
 Since $\varphi(\zero)= \rho(\bc, \dots, \bc)^{-1}\varphi(\bi)$ by the property %
 (i) of $\bc$, Claim shows that $\xi$ gives $\varphi(\zero)$. 
  \end{pf}

The fact that the Cayley element  relates $\varphi(\zero)$ and $\varphi(\bi)$ was used also in \cite{U2}.

  \begin{sbpara}\label{redas}

  Let $(\rho, \varphi)$ be an $\SL(2)$-orbit in $n$ variables. Then the {\it rank} of $(\rho, \varphi)$ is defined to be the number of $j$ ($1\leq j\leq n$) such that the $j$-th component  $\rho_j :\SL(2)_\R\to G_\R$ of $\rho$ is a non-trivial homomorphism.

  Let $r$ be the rank of $(\rho, \varphi)$ and let $(\rho', \varphi')$ be the $\SL(2)$-orbit in $r$ variables defined as follows. Let $\{s(1), \dots, s(r)\}$ ($s(1)< \dots <s(r)$) be the set of all $j$ ($1\leq j\leq n$) such that $\rho_j$ is non-trivial. Define $\rho': \SL(2)_\R^r \to G_\R$ to be the unique homomorphism satisfying $\rho(g_1, \dots, g_n)= \rho'(g_{s(1)}, \dots, g_{s(r)})$ and let $\varphi'$ be the unique map $\bP^1(\C)^r \to  \Dc$ such that $\varphi(z_1, \dots, z_n)= \varphi'(z_{s(1)}, \dots, z_{s(r)})$.

  We call $(\rho', \varphi')$ the $\SL(2)$-orbit in $r$ variables of rank $r$ {\it associated to $(\rho, \varphi)$}.

  \end{sbpara}

  \begin{sbpara}\label{redW}

Let $(\rho, \varphi)$ be an $\SL(2)$-orbit in $n$ variables of rank $n$. For $1\leq j\leq n$, define homomorphisms $\tau_j^{\star}, \tau_j: \bG_{m,\R}\to G_\R$ as follows. 
  $$\tau^{\star}_j(t)=\rho(g_1, \dots, g_n) \quad (t \in \bR^{\times}),  \quad \text{where}$$ 
  $$g_k= \begin{pmatrix}   1/t & 0\\ 0 & t\end{pmatrix}  \;\text{if} \; 1\leq k\leq j, \quad g_k= 1 \;\text{if} \; j< k\leq n.$$
  $$\tau_j(t):= \tau_j^{\star}(t)k_0(t)\quad (t \in \bR^{\times}),$$
where $k_0\in\Psi_W(G)$ is as in the beginning of Section \ref{s:DSL}.

  \end{sbpara}

  \begin{sbpara}\label{redeq} We define the equivalence relation between $\SL(2)$-orbits. 
  
   An $\SL(2)$-orbit in $n$ variables and an $\SL(2)$-orbits in $n'$-variables are equivalent if and 
only if they have the same rank  $r$ and their associated $\SL(2)$-orbits of rank $r$ are equivalent in the
following sense.
  
  For $\SL(2)$-orbits in $n$ variables of rank $n$, the equivalence relation is given as follows: $(\rho, \varphi) \sim (\rho', \varphi')$ if and only if there is a $t=(t_j)_{1\leq j \leq n}\in \R_{>0}^n$ such that
  $$\rho'(g)=\tau(t)\rho(g)\tau(t)^{-1}, \quad \varphi'(z)=\tau(t)\varphi(z)$$
  for any $g\in \SL(2)^n_\R$ and $z\in \frak H^n$, where $\tau(t)$ denotes $ \prod_{j=1}^n \tau_j(t_j)$ associated to $(\rho, \varphi)$. 
  \medskip
  
  In terms of Lemma \ref{redp1} (iv), $(\rho, \xi_1)\sim (\rho', \xi'_1)$ if and only if  there is a $t=(t_j)_{1\leq j \leq n}\in \R_{>0}^n$ such that
  $\rho'(g)=\tau(t)\rho(g)\tau(t)^{-1}$
  for every $g\in \SL(2)^n_\R$ and $\xi_1=\xi'_1$. 
  \end{sbpara}
  
\begin{sbrem} 
We have the same equivalence relation $\sim$ even if we replace $\tau_j$ by $\tau^{\star}_j$. This is because $\tau(t)\rho(g)\tau(t)^{-1}= \tau^{\star}(t)\rho(g)\tau^{\star}(t)^{-1}$ and the actions of $\tau(t)$ and $\tau^{\star}(t)$ on $D$ are the same. Here $\tau^{\star}(t)= \prod_j \tau^{\star}_j(t_j)$ which is associated to $(\rho, \varphi)$. 
\end{sbrem}

  \begin{sbpara} Let $(\rho, \varphi)$ be an $\SL(2)$-orbit in $n$ variables. 
  
  Let $1\leq j\leq n$. 
  For  $V\in \Rep(G)$, let  $W^{(j)}_{\bullet}V_\R$ be the increasing filtration on $V_\R$ defined by the action of $\bG_{m,\R}$ via $\tau_j$.  That is, for $w\in \Z$, $W^{(j)}_wV_\R$ is the part of $V_\R$ on which the action of $\bG_{m,\R}$ via $\tau_j$ is of weights $\leq w$.

 For $1\leq j\leq n$, define  $N_j\in \Lie(G_\R)$ by
     $$N_j= \Lie(\rho) (g_1, \dots, g_n) \quad \text{with}\; g_j= \begin{pmatrix} 0 & 1 \\ 0& 0\end{pmatrix}, \; g_k=0 \;\text{for}\;  k\neq j,$$
     where $\Lie(\rho)$ is the homomorphism ${\frak{sl}}(2,\R)^n \to \Lie(G_\R)$ induced by $\rho$. Let $W^{(0)}=W$. 
      Then for $0\leq j \leq k \leq n$ and for any $a_l\in \R_{>0}$ ($j< l \leq k$) and any $V\in \Rep(G)$, the filtration $W^{(k)}_{\bullet}V_\R$ is the relative weight filtration of the nilpotent operator $\sum_{j<l\leq k} a_lN_l:V_\R\to V_\R$ with respect to $W^{(j)}_{\bullet}V_\R$.
               
For $1\leq j\leq n$, the following conditions (i)--(iv) are equivalent. (i) $\rho_j$ is trivial. (ii) $\tau_j$ is trivial. (iii) $W^{(j)}=W^{(j-1)}$. (iv) $N_j=0$.

  We say that the  weight filtrations of $(\rho, \varphi)$ are {\it rational} if the weight filtrations $W^{(j)}$ ($1\leq j\leq n$)  of the associated $\SL(2)$-orbit in $n$ variables of rank $n$ has the following property: 
 For any $V\in \Rep(G)$, 
the filtrations $W^{(j)}_{\bullet}V_\R$ ($1\leq j\leq n$) on $V_\R$ 
are $\Q$-rational. This rationality depends only on the equivalence class of $(\rho, \varphi)$.

\end{sbpara}
  
  \begin{sbpara}\label{redSL}
   We define the set $D_{\SL(2)}$ as the set of all equivalence classes of $\SL(2)$-orbits whose weight filtrations are rational.

     \end{sbpara}
    
\begin{sbrem} The rationality condition on the weight filtrations will become important to have the good properties (Hausdorff property etc.) of the quotient space $\Gamma \bs D_{\SL(2)}$ for a semi-arithmetic subgroup $\Gamma$ of $G'(\Q)$  (cf.\ Section \ref{ss:gen1}).

\end{sbrem}

  \begin{sbpara}\label{red2} 
Let $p\in D_{\SL(2)}$.
  Let $(\rho, \varphi)$ be the $\SL(2)$-orbit in $n$ variables of rank $n$ whose class is $p$.

  Let $$\tau_p, \tau_p^{\star}: \bG_{m,\R}^n \to G_\R$$ be the homomorphisms $\tau$ and $\tau^{\star}$ associated to $(\rho, \varphi)$ as in \ref{redeq}, respectively.
 Let $\tau_{p,j}$ and $\tau^{\star}_{p, j}$ ($1\leq j\leq n$) be the $j$-th component $\bG_{m,\R} \to G_\R$ of $\tau_p$ (resp.\ $\tau^{\star}_p$).

 We will call $$\{\varphi(iy_1, \dots, iy_n)\in D\;|\; y_j\in \R_{>0} \;(1\leq j\leq n)\}= \tau_p(\R_{>0}^n) \varphi(\bi)=\tau^{\star}_p(\R_{>0}^n)\varphi(\bi)$$ the {\it torus orbit} associated to $p$.  
 
        \end{sbpara}

\subsection{The sets $D_{\SL(2)}$ and $D^{\star}_{\SL(2)}$}\label{ss:DSL}

When we do not assume that $G$ is reductive, unlike the case where $G$ is reductive, there are several sets of $\SL(2)$-orbits. Here we define the most important two sets $D_{\SL(2)}$ and $D^{\star}_{\SL(2)}$. Other sets of $\SL(2)$-orbits $D^{\star,+}_{\SL(2)}$, $D^{\star,\pa}_{\SL(2)}$ will be defined later. 

 \begin{sbpara}\label{DSL}
Let $D_{\red}=D(G_{\red},h_0)$ be as in \ref{Dred}.
Let $D_{\red, \SL(2)}=D(G_{\red}, h_0)_{\SL(2)}$ be $D_{\SL(2)}$ for $(G_{\red}, h_0)$ in \ref{redSL}.
This  $D_{\red, \SL(2)}$ is a $G$-MHS version of the space $D_{\SL(2)}(\gr^W)^{\sim}$ in \cite{KNU2} Part II 3.5.1.  

For $p\in D_{\red, \SL(2)}$ and let 
  $\tau_p, \tau^{\star}_p:\bG_{m,\R}^n \to G_{\red,\R}$ be the homomorphism in \ref{red2} associated to $p$, where 
$n$ is the rank of $p$. 

 Let $D_{\SL(2)}$ (resp.\ $D^{\star}_{\SL(2)}$) be the set of pairs $(p, Z)$, where $p$ is an element of $D_{\red, \SL(2)}$ and $Z$ is a subset of $D$ whose image in $D_{\red}$ coincides with the torus orbit (\ref{red2}) of $p$, satisfying one of the following conditions (A) and (B) for every (equivalently,  for all) $z\in Z$. Let $s= \spl_W(z)$.

 (A)  Let $h: {\bf G}_{m,\R}^n \to G_\R$ be the homomorphism defined by
 $$h(t):= s\tau_p(t)s^{-1}  \quad (\text{resp.}\; s\tau^{\star}_p(t)s^{-1}).$$
 Then $Z$ is a $h(\R_{>0}^n)$-orbit.
 
 (B)  Let $h: {\bf G}_m\times {\bf G}_{m,\R}^n\to G$ be the homomorphism defined by
 $$h(t_0, (t_j)_{1\leq j\leq n}):= s( k_0(t_0) \prod_{j=1}^n \tau_p(t_j))s^{-1}\quad (\text{resp.}\; s( k_0(t_0) \prod_{j=1}^n \tau^{\star}_p(t_j))s^{-1}).$$
 Then $Z$ is an $h(\R_{>0}\times \R_{>0}^n)$-orbit contained in $D_{\nspl}$. 
 
 \end{sbpara}

 \begin{sbpara}\label{SL2AB} Let $x=(p, Z)\in D_{\SL(2)}$ 
(resp.\ $D^{\star}_{\SL(2)}$). 
  We call $x$ an $A$-orbit in the case (A), and a $B$-orbit in the case (B).

 We define the {\it mild part of $D^{\star}_{\SL(2)}$} as the part consisting of all $A$-orbits and denote it by $D^{\star,\mild}_{\SL(2)}$.

We denote by $D_{\SL(2),\spl}$ the subset of $D_{\SL(2)}$ consisting of $(p,Z)$ whose $Z$ is contained in $D_{\spl}$, and by $D_{\SL(2),\nspl}$ its complement. 
 
       \end{sbpara}
    
      \begin{sbpara}\label{taux0} For $x=(p, Z)\in D_{\SL(2)}$ (resp.\  $D^{\star}_{\SL(2)}$), the homomorphism $h$ in \ref{DSL}  is independent of the choice of $z\in Z$. In the case $x\in D_{\SL(2)}$ (resp.\  $D^{\star}_{\SL(2)}$), we denote this homomorphism $h$  by $\tau_x$ (resp.\  $\tau^{\star}_x$).

      \end{sbpara}

\subsection{Weight filtrations associated to $\SL(2)$-orbits}\label{SL2an0}

In this Section \ref{SL2an0}, we give preliminary definitions which are used in Section \ref{SL2an}.

\begin{sbpara}\label{Wlem} Let $E$ be a field  of characteristic $0$ and let $\cG$ be a linear algebraic group over $E$. Let $\Rep_E(\cG)$ be the category of finite-dimensional representations of $\cG$ over $E$. Let $W'$ be an increasing filtration on the functor $V\mapsto V$ from  $\Rep_E(\cG)$ to the category of finite-dimensional $E$-vector spaces. Then, by \cite{SR} Chap.\ IV 2.4, the following two conditions (i) and (ii) are equivalent. 

(i) There is a homomorphism $a: \bG_{m,E} \to \cG$ which defines $W'$. That is, for $V\in \Rep_E(\cG)$, $W'_wV$ is the part of $V$ on which the action of $\bG_{m,E}$ via $a$ is of weights $\leq w$. 

(ii) $W'$ is an exact $\otimes$-filtration
 in the sense of \cite{SR} Chapter IV 2.1. That is, %
$\gr^{W'}_w$ is an exact functor for every $w\in \Z$, and %
$W'_w(V_1\otimes V_2)= \sum_{j+k=w} \;W'_jV_1 \otimes_E W'_kV_2$ for every $V, V_2\in \Rep_E(\cG)$ and $w\in \Z$.

We denote by  $\frak W(\cG)$  the set of all $W'$ satisfying the above equivalent conditions (i) and (ii). 

If $E\subset E'$, via the injection $\frak W(\cG)\to \frak W(\cG_{E'})$ 
induced by the homomorphism $\Hom(\bG_{m,E},\cG) \to \Hom(\bG_{m,E'},\cG_{E'})$, we often regard $\frak W(\cG)$ as a subset of $\frak W(\cG_{E'})$. 

Actually only the cases $E=\Q$ and $E=\R$ are important to us. 
\end{sbpara}

    \begin{sbpara}\label{taux} Let $x=(p, Z)\in D_{\SL(2)}$ and let $n$ be the rank of $p=x_{\red}$. In the case where $x$ is an $A$ (resp.\ $B$)-orbit,  we defined in \ref{taux0} a homomorphism $$\tau_x: \bG_{m,\R}^n=\bG_{m,\R}^{\{1,\dots,n\}} \to G_\R  \quad \; (\text{resp.}\;  \bG_{m,\R}\times \bG^n_{m,\R}=\bG_{m,\R}^{\{0,\dots, n\}}\to G_\R).$$ 
    For $1\leq j\leq n$ (resp.\  $0\leq j \leq n$), let
     $W_x^{(j)}$ be the increasing filtration 
     on the functor $V\mapsto V_\R=\R\otimes_\Q V$ 
     given by $\tau_{x,j}:{\bf G}_{m,\R}\to G$. For $1\leq j\leq n$, we have $W_x^{(j)}= s(W_p^{(j)})$. That is, for $V\in \Rep(G)$, $W_{x,w}^{(j)}V_\R= \bigoplus_{k\in \Z} \; s(W^{(j)}_{p,w}\gr^W_k V_\R)$. In the case where $x$ is a $B$-orbit, $W^{(0)}_x=W$. These $W_x^{(j)}$ ($1\leq j\leq n$) (resp.\ ($0\leq j\leq n$)) are called the {\it weight filtrations associated to $x$}.

For $x\in D_{\SL(2)}$, let $$\Phi(x)\subset \frak W(G_\R)$$ be the set of weight filtrations associated to $x$. 

We have $W\in \Phi(x)$ if and only if $x$ is a $B$-orbit.  If $G$ is reductive, $W\notin \Phi(x)$ and we have  $\Phi(x)\subset \frak W(G)$.

    \end{sbpara}

  \begin{sbprop}\label{PhiZ} $(p,Z)\in D_{\SL(2)}$ (resp.\ $D^{\star}_{\SL(2)}$) is determined by $(\Phi(p), Z)$.
      \end{sbprop}

\begin{pf}
When $G$ is reductive, the proof of \cite{KU1} Lemma 3.10 for $\Q$-polarization works for $\R$-polarization. 
  The general non-reductive case follows easily from this.
\end{pf}

  This will be used in the proof of the injectivity in Propositions \ref{emb1} (2), \ref{emb2} (2), \ref{emb3} (2) in Section \ref{SL2an}.

  \begin{sbprop}\label{variance} 
  
   Let $x\in D_{\SL(2)}$ and let $n$ be the rank of $x_{\red}\in D_{\red,\SL(2)}$.  Let $V\in \Rep(G)$, $V\neq 0$. Let $0\leq j\leq n$. Let $W_x^{(0)}:=W$ also in the case where $x$ is an $A$-orbit. Define the mean value $\mu_j\in \Q$ 
  and the variance $\sig_j^2\in \Q$ of $W_x^{(j)}$ by 
  $$\mu_j := \sum_{w\in \Z} \dim_\R(\gr^{W_x^{(j)}}_wV_\R)w/\dim_\R(V_\R), \quad  \sig^2_j:=  \sum_{w\in \Z} \dim_\R(\gr^{W_x^{(j)}}_wV_\R)(w-\mu_j)^2/\dim_\R(V_\R).$$
  Then $\mu_j$ is independent of $j$ and $\sig_j^2\leq \sig_k^2$ if $j\leq k$. Furthermore, if $V$ is a faithful representation of $G$, then $\sig^2_j<\sig^2_k$ for $0\leq j<k\leq n$.

  \end{sbprop}

\begin{pf}
This follows from  \cite{KNU2} Part II Proposition 2.1.12.
\end{pf}

\begin{sbpara}\label{cW(G)}

We denote by $\cW(G)$ the set of all  subsets $\Phi$ of $\frak W(G_\R)$ such that $\Phi=\Phi(x)$ for some $x\in D_{\SL(2)}$.

For $\Phi\in \cW(G)$, let  $\Phi_{\red}\in \cW(G_{\red})$ be $\{\gr^W(W')\;|\; W'\in \Phi, \;W'\neq W\}$. If $x=(p, Z)\in D_{\SL(2)}$, then $\Phi(x)_{\red}=\Phi(p)$.

Let  $\Phi\in \cW(G)$ and let $n$ be the order of $\Phi_{\red}$. In the case $W\notin \Phi$ (resp.\ $W\in \Phi$), we often identify $\Phi$ with 
the set $\{1,\dots, n\}$ (resp.\ 
$\{0, 1, \dots, n\}$) 
in the unique way that if $W', W''\in \Phi$ correspond to $i, j\in \{1, \dots, n\}$ (resp.\ $\{0, 1, \dots, n\}$), respectively, then 
$i\leq j$ if and only if the variance of $W'$ $\leq $ the variance of $W''$ (\ref{variance}).

\end{sbpara}

\begin{sbpara} 
For $\Phi\in \cW(G_{\red})$, let $D_{\SL(2)}(\Phi)$ (resp.\ $D^{\star}_{\SL(2)}(\Phi)$)  be the set of all points $(p, Z)$ of $D_{\SL(2)}$ (resp.\ $D^{\star}_{\SL(2)}$) such that $\Phi(p)\subset \Phi$. 

For $\Phi\in \cW(G)$, let $D_{\SL(2)}(\Phi)=\{x\in D_{\SL(2)}\;|\; \Phi(x)\subset \Phi\}$.

Then, in the next Section \ref{SL2an}, 
$$D^{\star}_{\SL(2)}= \bigcup_{\Phi\in \cW(G_{\red})} \; D^{\star}_{\SL(2)}(\Phi)$$
will be an important open covering of %
$D^{\star}_{\SL(2)}$, 
$$D_{\SL(2)}= \bigcup_{\Phi\in \cW(G_{\red})}\;  D_{\SL(2)}(\Phi)$$
will be an important open covering of $D^{II}_{\SL(2)}$, and 
$$D_{\SL(2)}= \bigcup_{\Phi\in \cW(G)} \; D_{\SL(2)}(\Phi)=\Bigl(\bigcup_{\Phi\in \cW(G),  W\notin \Phi} \; D_{\SL(2)}(\Phi)\Bigr)\cup \Bigl(\bigcup_{\Phi\in \cW(G), W\in \Phi} \; D_{\SL(2),\nspl}(\Phi)\Bigr)$$
(here $D_{\SL(2),\nspl}(\Phi):= D_{\SL(2),\nspl}\cap D_{\SL(2)}(\Phi)$) will be an important open covering of $D^I_{\SL(2)}$. 

\end{sbpara}

\begin{sbpara}
\label{P_is_para} Let the notation be as in \ref{Wlem}. 
Assume that $\cG$ is reductive and let $W'\in \frak W(\cG)$. Then, by \cite{SR} Chap.\ IV 2.2.5, 
   the stabilizer $P=\cG^{\circ}_{W'}\subset \cG^{\circ}$ of $W'$ in $\cG^{\circ}$ is a parabolic subgroup of $\cG$. 
   
   We thus have a map 
$$\Hom(\bG_{m,E},\cG)\to \{\text{parabolic subgroup of $\cG$}\}.$$
This map is surjective (cf.\ \cite{Mi2} Theorem 25.1, \cite{Sp} 15.1.2 (ii)). 

Let $W'\in \frak W(\cG)$ and let $P$ be the associated parabolic subgroup of $\cG$. 

If a homomorphism  $a:\bG_{m,E} \to \cG$  defines $W'$, then the image of $a$ in $\cG$ is contained in $P$. The composition $\bG_{m,E} \overset{a}\to P \to P_{\red}$ is independent of the choice of $a$ and has the image in the center of $P_{\red}$. The adjoint action of $\bG_{m,E}$ on $\Lie(P_u)$ via $a$ is of weights $\leq -1$.

Let $\spl(W')$ be the set of all isomorphisms of $\otimes$-functors from $\Rep_E(\cG)$ to the category of $E$-vector spaces
$$(V\mapsto \gr^{W'}V) \overset{\sim}\to (V \mapsto V)\quad\text{preserving the increasing filtrations defined by $W'$}.$$
Then $\spl(W')$ is a $P_u(E)$-torsor. 
   We have a natural bijection from the set of all homomorphisms $\bG_{m,E}\to \cG$  which define $W'$ to $\spl(W')$.

   \end{sbpara}

   \begin{sbpara}\label{BSspl}
Assume that $G$ is reductive.

Let ${\frak W}_0(G_\R)$ be
the set of all $W'\in {\frak W}(G_\R)$ satisfying the following condition.

($\ast$) For some (equivalently, for every) homomorphism $a: \bG_{m,\R}\to G_\R$ which defines $W'$, the image of the homomorphism $a^{\star}: \bG_{m,\R}\to G_\R\;;\;t\mapsto a(t)k_0(t)^{-1}$ belongs to the commutator subgroup $G'$ of $G$, where $k_0$ is as in the beginning of this Section \ref{s:DSL}.

  For a point $p\in D_{\SL(2)}$,   we have $\Phi(p) \subset \frak W(G)\cap  \frak W_0(G_\R)$.  
  In fact, the condition ($\ast$) is satisfied by $W'\in \Phi(p)$ because $\SL(2)^n$ coincides with its commutator subgroup and hence for any linear algebraic group $H$, the image of any homomorphism $\SL(2)^n \to H$ is contained in the commutator subgroup $H'$.

 Let $W'\in \frak W_0(G_\R)$. 
We have a canonical real analytic map $$\spl^{\BS}_{W'}: D\to \spl(W'):=\spl_\R(W')$$   
defined as follows. Let $\iota: \bG_m \to G^{\circ}_{W',\red}$ be the canonical homomorphism in \ref{P_is_para}. 
Let $x\in D$. We have a Borel--Serre lifting  $a^{\star}: \bG_{m,\R} \to (G')^{\circ}_{W',\R}$ of $\bG_m \to (G')^{\circ}_{W',\red}\;;\;t\mapsto \iota(t)k_0(t)^{-1}$ associated to the image of $x$ under $D\to \frak X$. We define $\spl^{\BS}_{W'}(x)$ to be the $\R$-splitting of $W'$ associated to $\bG_{m,\R} \to G^{\circ}_{W'}\;;\;t\mapsto a^{\star}(t)k_0(t)$.

  Note that this splitting $\spl^{\BS}_{W'}(x)$ is also defined as follows.  
  Let $V \in \Rep(G)$. Then for $H=x$, $V=H(V)=\bigoplus_w V_w$ where $V_w$ is the part of $V$ of weight $w$, and we have a polarization $\langle\;,\;\rangle_w$ on $V_{w,\R}=H(V_w)_\R$ in  the condition (4.1) in Lemma \ref{pol2}. 
   Then 
   the splitting $\spl^{\BS}_{W'}(x)$ of $W'_{\bullet}V_\R$  is the direct sum for $w\in \Z$ of the orthogonal decomposition of $W'_{\bullet}V_{w,\R}$ with respect to the Hodge metric associated to $\langle\;,\;\rangle_w$. 
  See \ref{orth}.

  An important fact is that, 
for $p \in D_{\SL(2)}$ of rank $n$,  if $\br$ is an element of the torus orbit of $p$ and  $1\leq j\leq n$, the splitting of $W^{(j)}_p$ given by $\tau_{p,j}:\bG_{m,\R}\to G_\R$ coincides with $\spl_{W_p^{(j)}}^{\BS}(\br)$. See \cite{KU1} Lemma 3.9. This fact is in the basis of the relation 
 of $\SL(2)$-orbits and Borel--Serre orbits considered in this paper and in our former papers \cite{KU2}, \cite{KNU2}. 
  
\end{sbpara}

\begin{sbpara}\label{sPhi}  Let $\Phi\in \cW_\R(G)$. Let $G_{\R, \Phi}=\bigcap_{W'\in \Phi} G_{\R, W'}$, where  
$G_{\R,W'}$ denotes the stabilizer subgroup of $W'$ in $G_\R$. 
By a {\it splitting} of $\Phi$, we mean a homomorphism $$\alpha: \bG_{m, \R}^{\Phi}\to G_{\Phi, \R}$$ such that for each $W'\in \Phi$, the $W'$-component $\bG_{m,\R}\to G_{\R}$ of $\alpha$ is a splitting of $W'$. 

A splitting of $\Phi$ exists. In fact, if $\Phi=\Phi(x)$ for $x\in D_{\SL(2)}$, $\tau_x$ (\ref{taux}) is a splitting of $\Phi$. 

In the case $W\notin \Phi$ (resp.\ $W\in \Phi$), by a {\it distance to $\Phi$-boundary}, we mean a real analytic map $$\beta: D\to \R_{>0}^{\Phi} \quad (\text{resp.} \; D_{\nspl}\to \R^{\Phi}_{>0}) $$ such that $\beta(\alpha(t)x)= t\beta(x)$ for all splittings $\alpha$ of $\Phi$, all $t\in \R^{\Phi}_{>0}$,  and  all $x\in D$ (resp.\ $x\in D_{\nspl}$).

\end{sbpara}

\begin{sbprop}\label{emb30}

 A distance to $\Phi$-boundary exists.
\end{sbprop}

This is a $G$-MHS version of \cite{KNU2} Part II 3.2.5. The proof given below is similar to that of loc.\ cit.

\begin{pf} We first prove the case where $G$ is reductive of \ref{emb30}. 
Take a submonoid $V$ of $X(\bG_m^{\Phi})$ such that $V\cup V^{-1}= X(\bG_m^{\Phi})$, $V \supset X(\bG_m^{\Phi})_+$ and such that $V^\times \cap X(\bG_m^{\Phi})_+=\{1\}$. Let $P=P_V$ be the parabolic subgroup of $G'$ associated to $V$. Then $(G')^{\circ}_{\Phi}\subset P$ %
and, for every splitting $\alpha$ of $\Phi$, the map $\bG_{m,\R}^{\Phi}\to P_{\red,\R}$ is injective and independent of the choice of $\alpha$, and  the image of this homomorphism is contained in $S_P$. Let $\R_{>0}^{\Phi}\to A_P$ be the induced injective map. Take a homomorphism of Lie groups $h: A_P\to \R_{>0}^{\Phi}$ which splits the last injection. Let $\pi: P(\R) \to A_P$ be the canonical map defined in \ref{XBSr}, and let $f:\frak X\to A_P$ be the map there induced by $\pi$ by fixing $K\in \frak X$. Then the composition $D\to \frak X \overset{f}\to A_P\overset{h}\to \R_{>0}^{\Phi}$ is a distance to $\Phi$-boundary. 

Now we consider the general case. Let $n$ be the order of $\Phi_{\red}$. 
Let $\beta_{\red}: D_{\red}\to \R^n_{>0}$ be a distance to the $\Phi_{\red}$-boundary.

In the case $W\notin \Phi$, the composition $D\to D_{\red} \overset{\beta_{\red}}\longrightarrow \R^n_{>0}$ is a distance to $\Phi$-boundary.

Next we prove the case $W\in \Phi$ of \ref{emb30}.  
Fix $p_0\in D_{\red}$. Let $Z$ be the center of $G_{\red}(\R)$. By the $\R$-polarizability, $K=\{g\in G_{\red}(\R)/Z\; |\;gp_0=p_0\}$ is compact. For each integer $w\leq -2$, take a $K$-invariant positive symmetric $\R$-bilinear form $(\cdot, \cdot)_w$ on the part $L_w$ of $L:=\cL(p_0)$. 
Define a map $f :L\smallsetminus \{0\} \to \R_{>0}$ by 
$f(\delta):=(\sum_{w\leq -2}\; (\delta_w,\delta_w)^{-1/w})^{-1/2}$, 
 where $\delta_w$ denotes the
component of $\delta$ of weight $w$. For $p \in D_{\red}$, if $g$ is an element of $G_{\red}(\R)$ 
such that $p = gp_0$, then we have an isomorphism $\Ad(g) : L=\cL(p_0) \overset{\sim}\to \cL(p)$. The map
$f_p :\cL(p)\smallsetminus \{0\}\to \R_{>0}, \delta \mapsto f(\Ad(g)^{-1}\delta)$ is independent of the choice of $g$. This
is because $(g')^{-1}g\in   K$ if $g,g'\in G_{\red}(\R)$ and $gp_0=g'p_0$. Define $\gamma' :D_{\nspl} \to
  \R_{>0}$ by sending an element corresponding to $(p,s,\delta)$ ($p\in D_{\red}$, $s\in \spl(W)$, $\delta\in \cL(p)\smallsetminus \{0\}$) to $ f_p (\delta)$.  For $x\in D_{\nspl}$, define $\gamma(x)=\gamma'(x)\dot \prod_{j=1}^n  \beta_{\red,j}(x_{\red})^{-1}$. 
 Then $(\gamma,\beta_{\red}): D_{\nspl}\to \R_{>0}\times \R^n_{>0}$ is a distance to $\Phi$-boundary. 
\end{pf}

\subsection{The real analytic structures of $D_{\SL(2)}$ and $D^{\star}_{\SL(2)}$}
\label{SL2an}

In this Section \ref{SL2an}, we endow $D_{\SL(2)}$ and $D^{\star}_{\SL(2)}$ with the sheaves of real analytic functions and with the log structures with sign. In fact, like in Part II of \cite{KNU2}, $D_{\SL(2)}$ is endowed with two kinds of such structures, $D^I_{\SL(2)}$ and $D^{II}_{\SL(2)}$. 

A basic property of these structures is that in the case where $G$ is reductive, if $(\rho, \varphi)$ is an $\SL(2)$-orbit in $n$ variables with class $p$ in $D_{\SL(2)}$, $\varphi(iy_1, \dots, iy_n)$ converges to $p$ when $y_j/y_{j+1}\to \infty$ for $1\leq j\leq n$ ($y_{n+1}$ denotes $1$).

      We first review the category $\cB'_\R(\log)$ in \ref{BR1}--\ref{BR3}. The above spaces with the above structures become objects of $\cB'_\R(\log)$

\begin{sbpara}\label{BR1}  The categories $\cB'_\R$ and $\cB'_\R(\log)$ (see \cite{KNU2} Part IV 1.3).

Let $\cB'_\R$ be the category of locally ringed spaces $S$ over $\R$ satisfying the following
condition (i) locally on $S$. Endow $\R^n$ with the sheaf $\cO_{\R^n}$ of real analytic functions. 

(i) There are $n \geq  0$ and a morphism $\iota : S \to \R^n$ of locally ringed spaces over $\R$ such that $\iota$ is injective, the topology of $S$ coincides with the topology induced from that of $\R^n$, and the map $\iota^{-1}(\cO_{\R^n})\to \cO_S$  is surjective.

For an object $S$ of $\cB'_\R$ , we often call the structural sheaf $\cO_S$ the {\it sheaf of real analytic functions} on $S$ (though $S$ need not be a real analytic space).

For an object $S$ of $\cB'_\R$, a {\it log structure with sign} on $S$ means a log structure $M$ on $S$ endowed with 
a submonoid sheaf $M_{>0}$  of $M$ satisfying the following (i) and (ii).

(i) $M_{>0}\supset  \cO^\times_{S,>0}$ . Here $\cO^\times_{S,>0}$  denotes the subgroup sheaf of $\cO^\times_S$  consisting of all local 
sections whose values are $> 0$.

(ii) The map $M_{>0} \times \{\pm 1\}\to  M \;;\;  (f,\epsilon)  \mapsto \epsilon f$ is an isomorphism of sheaves. Here we
regard $\{\pm 1\} \subset \cO_S^\times  \subset  M$.

Let $\cB'_\R(\log)$ be the category of objects of $\cB'_\R$ endowed with an fs log structure with sign.

 Real analytic  manifolds  with corners are regarded as objects of $\cB'_\R(\log)$ (\cite{KNU2} Part IV 1.3.8 (2)).

The category $\cB'_\R(\log)$ has fiber products (ibid.\ Proposition 1.3.11 (1)).

\end{sbpara}

\begin{sbpara}\label{BR2} For an object $X$ of $\cB'_\R(\log)$, for an object $Y$ of $\cB'_\R$  and for a morphism $Y\to X$ in $\cB'_\R$, we have an fs  log structure with sign  $(M_Y, M_{Y,>0})$ on $Y$, called the {\it inverse image} of the log structure with sign $(M_X, M_{X,>0})$ of $X$,  defined as follows. The log structure $M_Y$ is the inverse image of the log structure  $M_X$, that is, $M_Y$ is the push-out of $\cO_Y^\times \leftarrow f^{-1}(\cO_X^\times) \to 
f^{-1}(M_X)$ in the category of sheaves of commutative monoids on $Y$, and $M_{Y,>0}$ is the push-out of $\cO^\times_{Y,>0}\leftarrow f^{-1}(\cO^\times_{X,>0})\to f^{-1}(M_{X,>0})$ in the category of sheaves of commutative monoids on $Y$. 
\end{sbpara}

\begin{sbpara} \label{BR3}  (\cite{KNU2} Part IV 1.3.16.) Let $X$ be an object of $\cB'_\R(\log)$ satisfying the following condition (C).
\smallskip

(C) The canonical map from $\cO_X$ to the sheaf of $\R$-valued functions on $X$ is injective. 

\smallskip

Let $Y$ be a subset of $X$. 
Then $Y$ has a structure of an object of $\cB'_\R(\log)$ satisfying (C) such that for any object $S$ of $\cB'_\R(\log)$ satisfying (C), $\text{Mor}_{\cB'_\R(\log)}(S,Y)$ is identified with $\{f\in \text{Mor}_{\cB'_\R(\log)}(S, X)\;|\;f(S)\subset Y\}$. This structure on $Y$ is  defined as follows. The topology is the one as a subspace of $X$. $\cO_Y$ is the image of the canonical map from the pullback of $\cO_X$ on $Y$ to the sheaf of $\R$-valued functions on $Y$. The  log structure with sign of $Y$ is the inverse image of that of $X$ (\ref{BR2}).

We will say this is the {\it structure of $Y$ as an object of $\cB'_\R(\log)$ induced by the injection $Y\to X$} (or by the embedding of $Y$ in  $X$). 
\end{sbpara}

In the rest of this Section \ref{SL2an}, we  state definitions and properties of the structures $D_{\SL(2)}^I$ and $D^{II}_{\SL(2)}$ of $D_{\SL(2)}$ and  the structure of $D^{\star}_{\SL(2)}$ as objects of $\cB'_\R(\log)$. The almost all proofs will be given in Section \ref{pfan}.  

We first consider the case where $G$ is reductive. In this case, $D^I_{\SL(2)}=D^{II}_{\SL(2)}= D^{\star}_{\SL(2)}$.

\begin{sbprop}\label{emb1}  Assume that $G$ is reductive. Let $\Phi\in \cW(G)$ ($\ref{cW(G)}$). Take a splitting $\alpha$ of $\Phi$  and a distance $\beta$ to $\Phi$-boundary ($\ref{sPhi}$).

$(1)$  There is a unique map 
$$\nu_{\alpha, \beta}: D_{\SL(2)}(\Phi) \to \R^{\Phi}_{\geq 0} \times D \times \prod_{W'\in \Phi} \spl(W')$$
satisfying the following conditions {\rm (i)} and {\rm (ii)}. 

{\rm (i)} For $p\in D$, $\nu_{\alpha,\beta}(p)= (\beta(p), \alpha\beta(p)^{-1}p, (\spl^{\BS}_{W'}(p))_{W'\in \Phi})\in \R^{\Phi}_{> 0}\times D \times  \prod_{W'\in \Phi} \spl(W')$.

{\rm (ii)} For $p\in D_{\SL(2)}(\Phi)$, $\nu_{\alpha,\beta}(p)$   is the limit of $\nu_{\alpha,\beta}(\tau_p(t)r)$  for a fixed point $r$ of the torus orbit ($\ref{red2}$) of $p$ and for  $t\in \R_{>0}^{\Phi(p)}$ which tends to $(0,\dots,0)\in\R_{\geq 0}^{\Phi(p)}$.

$(2)$  This map $ \nu_{\alpha, \beta}$ is injective.

$(3)$ The structure of $D_{\SL(2)}(\Phi)$ as an object of $\cB'_\R(\log)$ induced via the injection $\nu_{\alpha, \beta}$  (\ref{BR3}) is independent of the choice of $(\alpha,\beta)$.

\end{sbprop}

 We will denote this map $\nu_{\alpha,\beta}$ as $p\mapsto (\beta(p), b_{\alpha,\beta}(p), (\spl_{W'}(p))_{W'})$.

  Next we consider $D^{II}_{\SL(2)}$ and $D^{\star}_{\SL(2)}$ for a general linear algebraic group $G$ over $\Q$.

\begin{sbprop}\label{emb2} 
Let $\Phi\in \cW(G_{\red})$. Let $D_{\SL(2)}(\Phi)$ (resp.\ $D^{\star}_{\SL(2)}(\Phi)$) be the inverse image of $D_{\red, \SL(2)}(\Phi)$ in $D_{\SL(2)}$ (resp.\ $D^{\star}_{\SL(2)}$). 
 
  Let $\alpha$ be a splitting of $\Phi$  and let $\beta$ be a distance to $\Phi$-boundary ({\rm \ref{sPhi}}). 
 
 $(1)$ There is a  unique map
 $$\nu_{\alpha, \beta}: D_{\SL(2)}(\Phi)\to  D_{\red,\SL(2)}(\Phi)\times \spl(W) \times \overline{\cL}$$ 
$$\text{(resp.}\;\;  \nu^{\star}_{\alpha, \beta}: D^{\star}_{\SL(2)}(\Phi)\to  D_{\red,\SL(2)}(\Phi)\times \spl(W)\times \overline{\cL}\text{)}$$
satisfying the following conditions {\rm (i)} and {\rm (ii)}. 

{\rm (i)} For $x\in D$, set $p=x_{\red}$.
Then, $$\nu_{\alpha,\beta}(x)=(p, \spl_W(x), \Ad(\alpha\beta(p))^{-1}\delta(x))$$ 
$${\text (resp.\ }
 \nu^{\star}_{\alpha, \beta}(x)=(p, \spl_W(x), \Ad(\alpha^\star\beta(p))^{-1}\delta(x)){\text ).}$$

{\rm (ii)} Let $x=(p,Z)\in D_{\SL(2)}(\Phi)$ (resp.\ $D^\star_{\SL(2)}(\Phi)$) ($\ref{DSL}$). Then  I$\nu_{\alpha,\beta}(x)$ (resp.\ $\nu^{\star}_{\alpha, \beta}(x)$)  is the limit of $\nu_{\alpha,\beta}(\tau_x(t)r)$ (resp.\ $\nu^{\star}_{\alpha, \beta}(\tau_x^{\star}(t)r)$) for a fixed point $r$ of $Z$, where  $t\in \R_{>0}^{\Phi(x)}$ and $t\to (0, \dots,0)\in \R^{\Phi(x)}_{\geq 0}$.

$(2)$ This map  $\nu_{\alpha,\beta}$ (resp.\ $\nu^{\star}_{\alpha, \beta}$) induces the bijection 
 $$D_{\SL(2)}(\Phi)\overset{\sim}\to  \{(p, s, \delta)\in  D_{\red,\SL(2)}(\Phi)\times \spl(W) \times \overline{\cL}\;|\;\delta\in \overline{\cL}(b_{\alpha,\beta}(p))\}$$ 
$$(\text{resp.}\;\;  D^{\star}_{\SL(2)}(\Phi)\overset{\sim}\to  \{(p, s, \delta)\in  D_{\red,\SL(2)}(\Phi)\times \spl(W) \times \overline{\cL}\;|\;\delta\in \overline{\cL}(b_{\alpha,\beta}(p))\}).$$
Here $b_{\alpha,\beta}(p)$ is as after Proposition $\ref{emb1}$.
The ones for $\nu_{\alpha,\beta}$ and for $\nu^\star_{\alpha,\beta}$ coincide. %

$(3)$ The structure of $D_{\SL(2)}(\Phi)$ (resp.\ $D^{\star}_{\SL(2)}(\Phi)$) as an object of $\cB'_\R(\log)$ induced via the bijection $\nu_{\alpha, \beta}$ (resp.\ $\nu^{\star}_{\alpha, \beta}$) is independent of the choice of $(\alpha, \beta)$. 
\end{sbprop}

 In the bijection in the above (2), the image of $x\in D_{\SL(2)}(\Phi)$ (resp.\ 
$x\in D^{\star}_{\SL(2)}(\Phi)$) belongs to $\cL$ if and only if $x$ is an $A$-orbit.

For $D_{\SL(2)}$ (resp.\ $D^{\star}_{\SL(2)}$), this Proposition \ref{emb2} is a $G$-MHS version of  \cite{KNU2} Part II Proposition 3.2.6 (ii), Proposition 3.2.7 (ii), and the situation (c) of \cite{KNU2} Part IV Proposition 2.3.9 (resp.\ of the situation (a) of \cite{KNU2} Part IV Proposition 2.3.9) which treated the extended period domain $D_{\SL(2)}\supset D(\La)$ (resp.\ $D^{\star}_{\SL(2)}\supset D(\La)$) for a classical period domain $D(\La)$ (\ref{ss:Dold}).

\begin{sbrem} 
  This is a remark on the description of a related part of \cite{KNU2}.
  There are fourteen \lq\lq resp.,''s in \cite{KNU2} Part IV 2.3.8, which are divided into two groups, and correspond only to the ones belonging to the same group.   
  Precisely, the fourth, the fifth, the seventh, and the last 
\lq\lq resp.,''s correspond, and the other ten \lq\lq resp.,''s correspond. 
  But four and ten do not correspond. 
\end{sbrem}

\begin{sbprop}\label{str2} $D_{\SL(2)}$ (resp.\ $D^{\star}_{\SL(2)}$) has  a unique structure of an object of $\cB'_\R(\log)$ such that for every $\Phi\in \cW(G_{\red})$, $D_{\SL(2)}(\Phi)$ (resp.\ $D^{\star}_{\SL(2)}(\Phi)$) is open in $D_{\SL(2)}$ and the induced structure on it coincides with the structure in Proposition $\ref{emb2}$ $(3)$. 
\end{sbprop}

The $D_{\SL(2)}$ as an object $\cB'_\R(\log)$ with the above structure is denoted by $D_{\SL(2)}^{II}$.

\medskip

For  $D_{\SL(2)}$ (resp.\ $D^{\star}_{\SL(2)}$), this is the $G$-MHS version of  the situation (c) (resp.\ (a)) of \cite{KNU2} Part IV Proposition 2.3.10 for $D_{\SL(2)}$ (resp.\ $D^{\star}_{\SL(2)}$) extending the classical period domain. For $D_{\SL(2)}$ of ibid.,  \cite{KNU2} Part II Theorem 3.2.10 also gives the same structure of $D^{II}_{\SL(2)}$.

  \begin{sbprop} \label{SLbdl}  
$(1)$ The map  $D^{II}_{\SL(2)}\to D_{\red, \SL(2)}\times \spl(W)$ is proper. 
  It is an $\bar L$-bundle.
  
$(2)$ The map  $D^{\star}_{\SL(2)}\to D_{\red,\SL(2)}\times  \spl(W)$ is proper. It is an $\bar L$-bundle. 
  
  \end{sbprop}

  For $D^{II}_{\SL(2)}$ (resp.\ $D^{\star}_{\SL(2)}$), this is a $G$-MHS version of \cite{KNU2} Part II Theorem 3.5.15 (resp.\ the situation (a) of 
 \cite{KNU2} Part IV Proposition 2.3.16). 
     
     We give a proof of Proposition \ref{SLbdl} assuming Proposition \ref{str2}. Let $p\in D_{\red, \SL(2)}(\Phi)$, $r:=b_{\alpha, \beta}(p)\in D_{\red}$. Take an open neighborhood $V$ of $r$ in $D_{\red}$ and a real analytic map $g: V\to G_{\red}(\R)$ such that $g(r)=1$ and such that $v=g(v)r$ for all $v\in V$. Let  $U\subset D_{\red, \SL(2)}(\Phi)$ be the inverse of $V$ under $b_{\alpha, \beta}: D_{\red,\SL(2)}(\Phi)\to D_{\red}$ and let $\tilde U$ be the inverse image of $U$ under the projection $D_{\SL(2)}(\Phi) \to D_{\red, \SL(2)}(\Phi)$ (resp.\  $D^{\star}_{\SL(2)}(\Phi)\to D_{\red,\SL(2)}(\Phi)$). Then we have an isomorphism $$\tilde U \overset{\sim}\to U \times \spl(W) \times \overline{\cL}(r)\;;\; x \mapsto (x_{\red}, \spl_W(x), \Ad(g(b_{\alpha, \beta}(x_{\red})))^{-1}\delta(b_{\alpha, \beta}(x_{\red}))).$$

     We next consider the structure $D_{\SL(2)}^I$ of $D_{\SL(2)}$. 

\begin{sbprop}\label{emb3}  
Let $\Phi\in \cW(G)$. 

$(1)$ Let $\alpha$ be a splitting of $\Phi$ and let $\beta$ be a distance to $\Phi$-boundary. 

 If $W\notin \Phi$ (resp.\ $W\in \Phi$),there is a  unique map 
$$\mu_{\alpha, \beta}: D_{\SL(2)}(\Phi)\to D   \quad (\text{resp.} \; D_{\SL(2),\nspl}(\Phi) \to D)$$ 
satisfying the following conditions {\rm (i)} and {\rm (ii)}. 

{\rm (i)} If $x\in D$ (resp.\ $x\in D_{\nspl}$), 
$$\mu_{\alpha,\beta}(x)=\alpha(\beta(x))^{-1}x.$$  

{\rm (ii)} Assume $W\notin \Phi$ (resp.\ $W\in \Phi$) and let $x=(p, Z)$ be an element of $D_{\SL(2)}(\Phi)$ (resp.\ $D_{\SL(2),\nspl}(\Phi)$).
Let $z\in Z$. Then 
$\mu_{\alpha,\beta}(x)$ is the limit of $\mu_{\alpha,\beta}(\tau_x(t)z)$ where $t\in \R_{>0}^{\Phi(x)}$ and $t\to (0,\dots,0)\in \R_{\geq 0}^{\Phi}$.

$(2)$  If $W\notin \Phi$ (resp.\ $W\in \Phi$),  the structure of $D_{\SL(2)}(\Phi)$ (resp.\ $D_{\SL(2), \nspl}(\Phi)$) as an object of $\cB'_\R(\log)$ induced by the injection $$D_{\SL(2)}(\Phi) \;(\text{resp.} \; D_{\SL(2),\nspl}(\Phi))\to D^{II}_{\SL(2)}\times D\;;\; x\mapsto (x, \mu_{\alpha,\beta}(x))$$
($\ref{BR3}$) is independent of the choice of $(\alpha,  \beta)$. 
 
\end{sbprop}

This is a $G$-MHS version of 
 \cite{KNU2} Part II Proposition 3.2.6 (i), Proposition 3.2.7 (i).

\begin{sbprop}\label{str3}  $D_{\SL(2)}$ has  a unique structure %
of an object of $\cB'_\R(\log)$ such that for every $\Phi\in \cW(G)$ such that $W\notin \Phi$ (resp.\ $W\in \Phi$),  $D_{\SL(2)}(\Phi)$ (resp.\  $D_{\SL(2),\nspl}(\Phi)$) is open in $D_{\SL(2)}$ and the induced structure on this coincides with that in Proposition $\ref{emb3}$ $(2)$.

\end{sbprop}
  This generalizes \cite{KNU2} Part II Theorem 3.2.10 (i).

  The structure of $D_{\SL(2)}$ in Proposition \ref{str3} as an object of $\cB'_\R(\log)$ is denoted by $D_{\SL(2)}^I$.

 The identity map of $D_{\SL(2)}$ is a morphism $D^I_{\SL(2)}\to D^{II}_{\SL(2)}$ of $\cB'_\R(\log)$. The log structure with sign of 
 $D^I_{\SL(2)}$ is the inverse image of that of $D^{II}_{\SL(2)}$.

To prove the existence of the global structures  Propositions \ref{str2} and \ref{str3}, we need the  following Proposition \ref{U(p)} and the local descriptions in Theorems \ref{SL2loc}, \ref{SL2loc2}, and \ref{SL2loc3} which are proved by using Proposition \ref{U(p)}.

\begin{sbpara}\label{U(p)def}
  Assume that $G$ is reductive. 
  Let $p \in D_{\SL(2)}$.
  We define the subset $U(p)$ of $D_{\SL(2)}$ in the same way as in \cite{KU2} Section 10.2, as follows.  
  Let $(\rho,\varphi)$ be an $\SL(2)$-orbit belonging to $p$.
  For each $J \subset \Phi(p)=\{1,\dots, n\}$, we define $p_J \in D_{\SL(2)}$ as the class of the $\SL(2)$-orbit $(\rho_J, \varphi_J)$ in $m$-variables, where $m$ is the cardinality of $J$, defined as follows.  
$\rho_J(g_1,...,g_m) := \rho(h_1,...,h_n)$, $\varphi_J(z_1,...,z_m) := \varphi(w_1,...,w_n)$,
where $h_j$ and $w_j$ ($1\leq j\leq n$) are as follows. Write $J= \{s_1, \dots, s_m\}\subset \{1, \dots,n\}$, $s_1<\dots <s_m$. If $j \leq s_k$  for some $k$, define $h_j :=g_k$ and $w_j := z_k$ for the smallest integer $k$ with $j \leq  s_k$. Otherwise, $h_j := 1$ and $w_j := i$. 
We have $\Phi(p_J)=J$.

  Let 
$$U(p):=\bigcup_{J \subset \Phi(p)} G'_J(\bR)\cdot p_J \subset 
D_{\SL(2)}(\Phi(p)),$$
where $G'_J$ is the stabilizer of the weight filtrations $W^{(j)}$ $(j\in J)$
in $G'$.
\end{sbpara}
The following proposition  is a generalization of \cite{KU2} Theorem 10.2.2.

\begin{sbprop}\label{U(p)}

  The set $U(p)$ is open in $D_{\SL(2)}$.  
  
  \end{sbprop}
  
  In the case $p\in D$, Proposition \ref{U(p)} follows from the case where $G$ is reductive of (1) of the following Lemma \ref{lemopen1}.

\begin{sblem}\label{lemopen1}
$(1)$ $D$ is a finite disjoint union of $G'(\R)G_u(\C)$-orbits which are open and closed.

$(2)$ $\Dc$ is a  finite disjoint union of $G'(\C)$-orbits which are open and closed.
\end{sblem}

\begin{pf} These are reduced to the case where $G$ is reductive. 
  Assuming that $G$ is reductive, let $Z$ be the center of $G$. 
  Then these are reduced to the facts that 
$G'(\R)Z(\R)$ is an open (and hence closed) subgroup of $G(\R)$ of finite index and $G'(\C)Z(\C)$ is an open (and hence closed) subgroup of $G(\C)$ of finite index.
\end{pf}

  \begin{sbpara}\label{preloc} 
  Assume that $G$ is reductive. Let $p\in D_{\SL(2)}$, let $n$ be the rank of $p$, and take a representative $(\rho, \varphi)$ of $p$. Let $\br=\varphi(\bi)\in D$. 
  We have $\Lie(G'_\R)= \bigoplus_{m\in \Z^n} \; \Lie(G'_\R)_m$, where $$\Lie(G'_\R)_m=\{x\in \Lie(G'_\R)\;|\; \Ad(\tau_p(t))x=\prod_{j=1}^n t_j^{m(j)}x \;\text{for all}\; t\in \R_{>0}^{\Phi}\}.$$ 
For $x\in \Lie(G'_\R)$, let $x_m\in \Lie(G'_\R)_m$ be the $m$-component of $x$. 

Let $K_{\br}$ be the maximal compact subgroup of $G'(\R)$ associated to $\br$. Let $'K_{\br}$ be the stabilizer of $\br$ in $K_{\br}$. 
Take $\R$-subspaces $R \subset \Lie(G'_{\bR})$ and $S \subset \Lie(K_{\br})$ 
such that
$$\Lie(G'_\R)=\Lie(\R_{>0}^n) \oplus \Lie(K_{\br})\oplus R, \quad \Lie(K_{\br})=\Lie('K_{\br})\oplus S,$$
where $\R^n_{>0}$ is embedded in $G'(\R)$ via $\tau_p$, and such that
$$R=\sum_{m\in \Z^{\Phi}} \; R\cap (\Lie(G'_\R)_m \oplus \Lie(G'_\R)_{-m}). $$
  Such $R$ and $S$ exist. 

Let $Y$ be the subset of $\R^n_{\geq 0}\times \Lie(G'_\R) \times \Lie(G'_\R) \times \Lie(G'_\R) \times S$ 
consisting of all elements $(t, f, g, h, k)$ satisfying the following conditions
(i)--(iv). Let $J=\{j\;|\; t_j=0\}\subset \{1,\dots,n\}$. 

(i) Let $m\in \Z^{\Phi}$. Then: $g_m=0$ unless $m(j)=0$ for all $j\in J$. $f_m=0$ unless $m(j)\leq 0$ for all $j\in J$. $h_m=0$ unless $m(j)\geq 0$ for all $j\in J$. 

(ii) Let $t'$  be any element of $\R^n_{>0}$ such that $t'_j =t_j$ for all $j\in \{1,\dots,n\}\smallsetminus J$. If $m\in \Z^n$ and $m(j)=0$ for all $j\in J$,   then $f_m=(\prod_{j=1}^n (t'_j)^{m(j)})g_m$ and $h_m =(\prod_{j=1}^n (t'_j)^{m(j)})^{-1}g_m$. 

(iii) We have $g\in R$ and we have $f_m+h_{-m}\in R$ for all $m\in \Z^n$. 

(iv) $\exp(k)\br\in G'_J(\R)\cdot\br$. 

\medskip

Let $Y_0= \{(t,f,g,h, k)\in Y\;|\; t\in \R^n_{>0}\}$. Then $$Y_0 \overset{\sim}\to \R^n_{>0}\times R\times S\;;\; (t,f,g,h, k) \mapsto (t,g,k).$$ Hence by Lemma \ref{lemopen1}, the map
$$Y_0 \to D\;;\; (t,f,g,h,k)\mapsto \tau_p(t)\exp(g)\exp(k)\br = \exp(f)\tau_p(t)\exp(k)\br$$
is an open map. 
\end{sbpara}

The following Theorems \ref{SL2loc} and \ref{SL2loc2} are variants of \cite{KNU2} Part II Theorem 3.4.4 and \cite{KNU2} Part IV Theorem 2.3.14.

  \begin{sbthm}\label{SL2loc} Assume that $G$ is reductive. Let the notation be as in $\ref{preloc}$. Then there exists an open neighborhood $U$ of $(0, \dots, 0)$ in $Y$ and an open immersion $U\to D_{\SL(2)}$ which sends $(0,\dots, 0)$ to $p$ such that 
  if $U_0$ denotes the subset of $U$ consisting of $(t, f,g,h,k)$ such that $t\in \R^n_{>0}$, then it sends $(t,f,g,h,k)\in U_0$ to $\tau_p(t)\exp(g)\exp(k)\br$. 
  \end{sbthm}  

  We return to the general $G$ not necessarily reductive, that is, in the next theorem, $G$ is a (general) algebraic linear group over $\bQ$.   

  \begin{sbthm}\label{SL2loc2} 
 Let $x\in D_{\SL(2)}$ (resp.\ $D^{\star}_{\SL(2)}$) and let $p$ be the image of $x$ in $D_{\red,\SL(2)}$. Take the spaces $Y$, $U$, $U_0$  in 
Theorem $\ref{SL2loc}$ for $G_{\red}$ and $p$,  and denote them by $Y_{\red}$, $U_{\red}$, and $U_{\red,0}$, respectively. Let $V_{\red}$ be the image of the open 
  immersion $U_{\red}\to D_{\red,\SL(2)}$ in Theorem $\ref{SL2loc}$. Let $L=\cL(\br)$. Let $V$ be the inverse image of $V_{\red}$ in $D^{II}_{\SL(2)}$ 
  (resp.\ $D^{\star}_{\SL(2)}$). Then there is an open immersion $U:=U_{\red} \times \spl(W) \times \bar L \to V$ which sends 
  $(t, f, g, h, k, s, \delta)\in U_{\red,0} \times \spl(W) \times L$ to the element of $D$ whose image 
  in $D_{\red}\times \spl(W) \times \cL$ is $$(\tau_p(t)\exp(g)\exp(k)\br,\  s, \ \Ad(\tau_p(t)\exp(g)\exp(k))\delta)$$ 
 $$ \text{(resp}. \;(\tau^{\star}_p(t)\exp(g)\exp(k)\br, \ s, \ \Ad(\tau^{\star}_p(t)\exp(g)\exp(k))\delta)\text{)}$$
 such that the diagram
  $$\begin{matrix}  U & \to & V\\
  \downarrow &&\downarrow \\
  U_{\red}&\to & V_{\red}\end{matrix}$$
  is cartesian.

  \end{sbthm}

\begin{sbpara}  We next consider the local property of $D^I_{\SL(2)}$. 
Let $x=(p, Z)\in D_{\SL(2)}$. 
 Take $z\in Z$ and let $s_0=\spl_W(z)\in \spl(W)$. Then $s_0$ is independent of the choice of $z$.. Let $\br=z_{\red}\in D_{\red}$. Let $\delta_0= \delta_W(z)\in L:=\cL(\br)$. Let $n$ be the rank of $p$.
 
Let $E=L$ in the case where $x$ is an $A$-orbit, and let $E=\bar L\smallsetminus \{0\}$ in the case where $x$ is a $B$-orbit. 

 Take the space $Y$  in 
Theorem $\ref{preloc}$ for $G_{\red}$ and $p$,  and denote it by $Y_{\red}$ here. In the case where $x$ is a $B$-orbit, take a real analytic map $\beta_0: L\smallsetminus \{0\}\to \R_{>0}$ such that $\beta_0(c\circ \delta)= c\beta_0(\delta)$ for all $c\in \R_{>0}$ and $\delta\in L\smallsetminus \{0\}$. We define a subset $X$ of $Y_{\red}\times  \Lie(G_{\R, u})\times \Lie(G_{\R.u})\times E$ as the set of 
 $(t,f,g,h,k,  u, v, \delta)$   ($(t,f,g,h,k)\in Y_{\red}$, $u,v\in \Lie(G_{\R,u})$, $\delta\in E$) satisfying the following conditions (i)--(iv). Let $n$ be the rank of $p$ and let $J:=\{j\;|\;1\leq j\leq n, t_j=0\}$, $m\in \Z^n$. 

(i)  $u_m=0$ unless $m(j)\leq 0$ for all $j\in J$.

(ii) $v_m=0$ unless $m(j)=0$ for all $j\in J$.

(iii) If $\delta\in \bar L\smallsetminus L$ (this happens only in the case where $x$ is a $B$-orbit), then $v=0$. 

(iv) Assume  $m(j)=0$ for all $j\in J$. Let $t'$ be an element of $\R^n_{>0}$ such that $t'_j=t_j$ if $j\in \{1,\dots, n\}\smallsetminus J$. Then if $x$ is an $A$-orbit, we have
$v_m = \sum_{j=1}^n (t'_j)^{-m(j)}u_m$. If $x$ is a $B$-orbit and $\delta\in L\smallsetminus \{0\}$, we have  $v_m =\Ad(\tau_{x,0}(\beta_0(\delta)))^{-1}\prod_{j=1}^n (t'_j)^{-m(j)}u_m$.

The following Theorem \ref{SL2loc3} is a variant of \cite{KNU2} Part II Theorem 3.4.6. 

\end{sbpara}

\begin{sbthm}\label{SL2loc3}

There are an open neighborhood $V$ of $(0, \dots, 0)$ in $Y \times \Lie(G_{\R,u})\times \Lie(G_{\bR, u})$ and an open immersion $U\to D^I_{\SL(2)}$, where $U$ is the inverse image of $V$ in $X$,  which sends  $(0, \dots, \delta_0)$ (resp.\ $(0, \dots, 0, 0\circ \delta)$) to $x$ 
 such that if $U_0$ denotes the subset of $U$ consisting of all $(t,f,g,h,k, u,v,\delta)$ such that $t\in \R^{\Phi}_{>0}$ and $\delta\in L$, it sends 
$(t,f,g,h,k,u,v,\delta)\in U_0$ to the element of $D$ whose image in $D_{\red}\times \spl(W) \times \cL$ is 
$$(\tau_p(t)\exp(g)\exp(k)\br,\  \exp(u)s_0, \ \Ad(\tau_p(t) \exp(g)\exp(k))\delta).$$
\end{sbthm}

\begin{sbrem}

In \cite{KNU2} Part II 3.4.5, in the definition of the space $Y^I(p, \br, R, S)$,   the last sentence \lq\lq If $t_0=0$, ...''   in the condition  (6${}'$) should be deleted. 
After this modification, Part II Theorem 3.4.6 becomes correct.
In our series of papers \cite{KNU2}, Part II Theorem 3.4.6 is used in Part IV 2.7.16, but Part IV
2.7.16 is correct with the modified Part II Theorem 3.4.6.

Another remark concerning Part II 3.4.5 is that 
the presentation of (5${}'$)  can be simplified: The condition \lq\lq$\exp(v)s_{\br}=s_{\br}$'' there is equivalent to a simpler one \lq\lq$v=0$''.

\end{sbrem}
    
\begin{sbpara}
  We can show that our $\SL(2)$-spaces belong to the category $\cB'_\R(\log)^+$ of nice objects in $\cB'_\R(\log)$ as in \cite{KNU2} Part IV Section 2.7, though the details are omitted. 
\end{sbpara}

\subsection{Proofs for Section \ref{SL2an}}\label{pfan}

We give proofs of the statements in Section \ref{SL2an} that have not yet proved.

\begin{sbpara}\label{plot1} 
  The proofs of Proposition \ref{emb1} for $\Phi\in\cW(G)$, Proposition \ref{emb2} for $\Phi\in\cW(G_{\red})$, and Proposition \ref{emb3} for $\Phi\in\cW(G)$ are similar to the corresponding parts in \cite{KNU2} Part II and Part IV, which are indicated after each proposition. 

\end{sbpara}

\begin{sbpara}\label{plot2} The proofs for the remaining results in Section \ref{SL2an} are also parallel to the corresponding parts in \cite{KNU2} Part II and Part IV. 

  That is, by Propositions \ref{emb1}, \ref{emb2}, and \ref{emb3} for respective $\Phi$ already proved, we first endow the $\Phi$-parts of each sets $D_{\SL(2)}$ and $D_{\SL(2)}^{\star}$ with the space structures. Then we %
consider the following version of Proposition \ref{U(p)} for $\Phi$-part: 
\medskip
  
  (1)
$U(p)$ in \ref{U(p)def} is open in $D_{\SL(2)}(\Phi)$ if $\Phi(p)\subset \Phi$.
\medskip
  
We prove this (1) in \ref{U(p)pf} after preparations \ref{U(p)pr} and Lemma \ref{U(p)pr2}.

\end{sbpara}

\begin{sbpara}\label{U(p)pr} %
Let $\Phi\in \cW(G)$. Fix a splitting  $\alpha: \bG_{m,\R}^{\Phi}\to G_{\Phi,\R}$ of $\Phi$ and  a distance to 
$\Phi$-boundary $\beta: D\to \R^{\Phi}_{>0}$. Let $J$ be a subset of $\Phi$. Let $n$ be the order of $\Phi$ and let $m$ be the order of $J$. Let $\alpha_J:\bG_{m,\R}^J\to G_{J,\R}$ be the splitting of $J$ defined as follows. If we write the inclusion map $J\to \Phi$ between totally ordered sets as an injective increasing map $\theta:\{1, \dots, m\}\to \{1,\dots, n\}$, then  $\alpha_J(t_1, \dots, t_m)= \alpha(t'_1, \dots, t'_n)$, where:
If $j \leq \theta(k)$  for some $k$, define $t'_j :=t_{\theta(k)}$  for the smallest integer $k$ such that $j \leq  \theta(k)$. Otherwise, $t'_j := 1$.

Let $P$ be the subset of $D_{\SL(2)}$ consisting of all elements $p$ such that $\Phi(p)=J$ and  $\tau_p=\alpha_J$.
Let $Q$ be the set of $\SL(2)$-orbits in $m$ variables whose classes belong to $P$.
We define a map 
$$\theta_{\alpha, \beta}: P\to Q$$ 
as follows. (There is an evident projection $Q\to P$, but the composition $P\overset{\theta_{\alpha,\beta}}\longrightarrow  Q\to P$ need not be the identity map.) 

Let $p\in P$. Take an $\SL(2)$-orbit  $(\rho, \varphi)$ with class $p$. Then $\theta_{\alpha, \beta}(p)$ is the $\SL(2)$-orbit $(\rho', \varphi')$ whose class $p'$ satisfies $\tau_{p'}=\tau_p=\alpha_J$ such that the $N_j$ of $\rho'$  (denote it  by $N'_j$) for $1\leq j\leq m$ is defined by using $N_j:= (N_j$ of $\rho)$ as $N_j':=\Ad(\alpha(\beta(\br)))^{-1}(N_j)$, where $\br=\varphi(\bi)$, and such that $\varphi'$ is defined as $\varphi'(z)= \alpha(\beta(\br))^{-1}\varphi(z)$.   Then $(\rho', \varphi')$  depends only on $p$ and is independent of the choice of $(\rho, \varphi)$. 

\end{sbpara}

\begin{sblem}\label{U(p)pr2}

Let $E$ be a field of characteristic $0$ and let $V$ be a finite-dimensional vector space over $E$ on which the Lie algebra $\frak{sl}(2,E)^m$ acts. Denote the action of $X\in \frak{sl}(2,E)^m$ on $V$ by $[X, \cdot]$. 
For $1\leq j\leq m$, let $N_j\in \frak{sl}(2, E)^m$ be the element whose $j$-th component is $\begin{pmatrix}0&1\\0&0\end{pmatrix}$ and whose $k$-th components are $0$ for all $k\neq  j$, and  let $Y_j\in \frak{sl}(2,E)^m$ be the element whose $j$-th component is $\begin{pmatrix} -1&0\\0&1\end{pmatrix}$ and whose $k$-th components are  $0$ for all $k\neq j$. Let $B_j$ be the set of all elements $v$ of $V$ such that $[Y_j, v]=-2v$, $[Y_k, v]=0$ for all $k\neq j$, and $[N_k, v]=0$ for all $k\neq j$. 
 Let $S_j\in B_j$ for $1\leq j\leq m$. Then there is an element $v\in V$ such that $[Y_j, v]=0$ for all $j$ and such that $[N_j, v]=S_j$ for all $j$.

\end{sblem}

\begin{pf} Let $A_j$ be the set of all elements $v$ of $V$ such that $[Y_k, v]=0$ for all $k$ and $[N_k, v]=0$ for all $k\neq j$. For a finite-dimensional $E$-vector space $V'$ on which the Lie algebra $\frak{sl}(2,E)$ acts, if we define $N, Y\in \frak{sl}(2,E)$ in the similar way to the above and we define $V'_a=\{v\in V'\;|\; [Y, v]=av\}$ for $a\in \Z$, the map $V'_0\to V'_ {-2}\;;\;v\mapsto [N, v]$ is surjective as is well known. 
By applying this to  $V'=\{v\in V\; |\; [Y_k, v]= [N_k,v]=0 \; \text{for all} \; k\neq j\}$, we have that the map $A_j \to B_j \;;\; v \mapsto [N_j, v]$ is surjective. Take $v_j\in A_j$ such that $[N_j, v_j]=S_j$. Let $v=\sum_j v_j$. Then $[N_j, v]=S_j$ for all $j$. 
\end{pf}

\begin{sbpara}\label{U(p)pf}
  Assume that $G$ is reductive. 
  We prove %
  \ref{plot2} (1).
  To do so, it is enough to show the following (1). 

(1) If $p' \in U(p)$ and if $p'_{\la}\in D_{\SL(2)}(\Phi)$ converges to $p'$, then $p'_{\la}\in  U(p)$ for all sufficiently large $\la$. (Here $(p'_{\la})_{\la}$ denotes a directed family.)

Since $U(p') \subset  U(p)$, by replacing $p'$  by $p$, we can reduce (1) to the following: 

(2) If $p_{\la} \in  D_{\SL(2)}(\Phi)$ converges to $p$, then $p_{\la} \in  U(p)$ for all sufficiently large $\la$.

Dividing the sequence $(p_{\la})_{\la}$ into subsequences, we may assume that, for a fixed $J \subset \Phi=\{1, . . . , n\}$, the family of weight filtrations associated to $p_{\la}$ is $J$ for every $\la$. Since $D_{\SL(2)}(\Phi(p))$ is open in $D_{\SL(2)}(\Phi)$ (in fact, the former is the inverse image of the open set of $\R^{\Phi}_{\geq 0}$ consisting of elements whose $\Phi\smallsetminus \Phi(p)$-components are non-zero), we may assume $J\subset \Phi(p)$. Assume $J\subset \Phi(p)$. Let $m$ be the order of $J$.

Fix a splitting $\alpha$ of $\Phi$. Let $\alpha_{\Phi(p)}: \bG_{m,\R}^{\Phi(p)} \to G_{\Phi,\R}$ (resp.\ $\alpha_J:\bG_{m,\R}^J\to G_{J,\R}$) be the splitting of $\Phi(p)$ (resp.\ $J$) induced by $\alpha$ defined as in \ref{U(p)pr}.

We prove the following. 

(3) We may assume $\tau_p= \alpha_{\Phi(p)}$  (\ref{U(p)pr}) and $\tau_{p_\la}=\tau_{p_J}=\alpha_J$ for all $\la$. 

We prove (3). 
  There exists a unique $u\in (G_{\Phi(p)})_u(\R)=(G'_{\Phi(p)})_u(\R)$ such that $\tau_p(t)= u\alpha_{\Phi(p)}(t)u^{-1}$ ($t\in \R^{\Phi(p)}_{>0}$) and there exists a 
unique $u_{\la}\in (G_J)_u(\R)=(G'_J)_u(\R)$ such that $\tau_{p_{\la}}(t) = u_{\la}\alpha_J(t)u^{-1}_{\la}$ ($t\in \R^J_{>0}$) . Then $\tau_{p_J}(t)= u\alpha_J(t)u^{-1}$ ($t\in \R^J_{>0}$). 
For $j\in J$, we have the convergence $\spl_j^{\BS}(p_{\la})\to \spl_j^{\BS}(p) =\spl_j^{\BS}(p_J)$. Since $\tau_{p_{\la}}$ gives $(\spl_j^{\BS}(p_{\la}))_{j\in J}$ and $\tau_{p_J}$ gives $(\spl_j^{\BS}(p_J))_{j\in J}$, 
$u_{\la}$ converges to $u$. Let $p'= u^{-1}p$ 
and $p'_{\la}:=u_{\la}^{-1}p_{\la}$. Then $\tau_{p'}=\alpha_{\Phi(p)}$ and $\tau_{p'_{\la}}=\tau_{p'_J}=\alpha_J$ and $p'_{\la}\to p'$. If we can prove that $p'_{\la}\in G'_J(\R)p'_J$ for $\la$ sufficiently large, we can obtain $p_{\la}\in G'_J(\R)p_J$ for $\la$ sufficiently large. This proves (3).

We now assume $\tau_p=\alpha_{\Phi(p)}$ and $\tau_{p_\la}=\tau_{p_J}= \alpha_J$ for all $\la$. 

Fix a distance $\beta$ to $\Phi$-boundary. 
Let $(\rho'_J, \varphi'_J):=\theta_{\alpha,\beta}(p_J)$ and $(\rho'_{\la}, \varphi'_{\la}):=\theta_{\alpha, \beta}(p_{\lam})$ (\ref{U(p)pr}). Let 
$(\rho'_J, \xi'_{1,J})$ and $(\rho_{\la}', \xi'_{1,\la})$ be the elements of the set 
(iv) in Lemma \ref{redp1} corresponding to   $(\rho'_J, \varphi'_J)$ and $(\rho_{\la}', \varphi'_{\la})$, respectively. 

  We prove the following. 

(4)  $(\rho'_\la, \xi'_{1,\la}) $ converges to $(\rho'_J, \xi'_{1, J})$ for  the compact-open topology of the space of continuous homomorphisms $\SL(2, \R)^m\times S^{(1)}_{\C/\R}(\R) \to G(\R)$. 

We prove (4). 
For $j\in J$, let $D^{(j)}$ be the subset of $D$ consisting of all $F$ such that for every $V\in \Rep(G)$, $(W^{(j)},F(V))$ is an $\R$-mixed Hodge structure. Then we have a continuous map $\delta_{W^{(j)}}: D^{(j)}\to \Lie(G_{\bR})$ (the $\delta$ in \ref{CKSdelta} for the weight filtration $W^{(j)}$). 
Let $(\rho_J, \varphi_J)$ and $(\rho_{\la}, \varphi_{\la})$ be  $\SL(2)$-orbits in $m$ variables whose classes in $D_{\SL(2)}$ are $p_J$ and $p_{\la}$, respectively. 
Then 
 $\nu_{\alpha, \beta}(p)=\nu_{\alpha, \beta}(p_J)= \alpha\beta(\br)^{-1}\br$, where  $\br=\varphi_J(\bi)$.  Since $\br= \exp(iN_1+\dots +iN_j)\varphi_J(\{0\}^j \times \{i\}^{m-j})$, where $N_j$ is the $N_j$ of $\rho_J$, we have $\nu_{\alpha, \beta}(p_J)\in \bigcap_{j\in J} D^{(j)}$ and $\delta_{W^{(j)}}(\nu_{\alpha,\beta}(p_J))=N_1'+\dots+N_j'$ for $1\leq j\leq m$, where $N'_j$ is the $N_j$ of $\rho'_J$. Similarly, $\nu_{\alpha, \beta}(p_{\la}) \in \bigcap_{j\in J} D^{(j)}$ and $\delta_{W^{(j)}}(\nu_{\alpha, \beta}(p_{\la}))=  N'_{1, \la}+\dots+ N'_{j,\la}$ for $1\leq j\leq m$, where $N'_{j,\la}$ is the $N_j$ of $\rho'_{\la}$. Since 
 $p_{\la}$ converges to $p$, $\nu_{\alpha, \beta}(p_{\la})$ converges to $\nu_{\alpha, \beta}(p)=\nu_{\alpha, \beta}(p_J)$. Hence $N'_{j, \la}$ converges to $N'_j$ for $1\leq j\leq m$. 
  Furthermore, $\varphi'_{\la}(\bi)=\nu_{\alpha,\beta}(p_{\la})$ converges to $\varphi'_J(\bi)= \nu_{\alpha, \beta}(p_J)$. 
  This proves (4). 

  We prove the following.

(5) If $\la$ is sufficiently large, there is $g_{\la}\in G'(\R)$ which commutes with $\alpha_J(t)$ for all $t\in (\R^\times)^m$ such that  $N'_{j, \la}=\Ad(g_{\la})N'_j$ for all $j\in J$ and such that $g_{\la}\to 1$. 

We apply Lemma \ref{U(p)pr2} to the case $E=\R$, $V=\R \otimes_\Q \Lie(G')$, and the action of $\frak{sl}(2, \R)^m$ on $V$ is induced by the adjoint action of $\SL(2)^m_\R$ on $V$ via $\rho_J$.  Let $H\subset G'_\R$ be the centralizer of the image of $\alpha_J: \bG_{m,\R}^m \to G'_\R$. Let $B_j$ ($1\leq j\leq m$) be as in Lemma \ref{U(p)pr2}, let $B= \prod_{j=1}^m B_j$, and let $b= (N_j')_{1\leq j\leq m}\in  B$. Then $\Lie(H)= \{v\in V\;|\; [Y_j, v]=0\; \text{for}\; 1\leq j\leq m\}$. The map $H \to B\;;\;g \mapsto \Ad(g)b$ induces the map $\Lie(H) \to T_b(B)=B$, where $T_b(B)$ denotes the tangent space of $B$ at $b$, and this last map is written as $v\mapsto ([v, N'_j])_{1\leq j \leq m}$. By Lemma \ref{U(p)pr2}, this last map is surjective. Hence the map $H \to B$ is smooth at $1\in H(\R)$ as a morphism of algebraic varieties. Hence there is $g_{\la}\in H(\R)$, $g_{\la}\to 1$ such that $\Ad(g_{\la})N_j= N_{j,\la}$. Thus (5) is proved. 

Now we complete the proof of  %
\ref{plot2} (1).
By (5), 
we may assume that $\rho'_{\la}=\rho'_J$. 
Consider $\xi'_{1,\la}, \xi'_{1,J}: S^{(1)}_{\C/\R}\to C_{G'}(\rho)$ (Lemma \ref{redp1}). 
  Here $C_{G'}(\rho)$ denotes the centralizer.
We use a result in \cite{LW} concerning homomorphisms from compact groups to locally compact groups, applied to homomorphisms from the compact  group $S^{(1)}_{\C/\R}(\R)$ to the locally compact group $C_{G'}(\rho)(\R)$ given by $\xi'_{1,\la},\xi'_{1,J}$. By the above-mentioned result, if $\la$ is sufficiently large, there is an $A\in C_{G'_\R}(\rho)(\R)$ 
 such that $\xi'_{1,\la}=A\xi'_{1,J}A^{-1}$. Hence $(\rho'_{\la}, \xi'_{1,\la})=A(\rho'_J, \xi'_{1,J})A^{-1}$ and hence, by Lemma \ref{redp1}, 
$(\rho'_{\la}, \varphi'_{\la})$ is the twist of $(\rho'_J, \varphi'_J)$ by $A$. 
Hence $p_{\la}$ is the twist of $p_J$ by an element of $G'_J(\R)$.
Thus \ref{plot2} (1) is proved.

\end{sbpara}

\begin{sbpara}\label{glwd} 
  As in \cite{KNU2} Part II Section 3.4, we can reformulate Theorems \ref{SL2loc}, \ref{SL2loc2}, and \ref{SL2loc3} in terms of $\Phi$-parts,  
which will be equivalent to the original statements as soon as the global space structures will be well-defined (Propositions \ref{str2} and \ref{str3}), and 
whose proofs are %
similar to those of \cite{KNU2} Part II Section 3.4 and Part IV Theorem 2.3.14. 
  Here, the openness of $U(p)$ in $D_{\SL(2)}(\Phi)$ in \ref{plot2} (1) 
is used in connection with the condition (iv) in \ref{preloc}.

The well-definedness of the global structures as objects of $\cB'_\R(\log)$  on the $\SL(2)$-spaces (Proposition \ref{str2}, Proposition \ref{str3}) follows from the versions of Theorems \ref{SL2loc}, \ref{SL2loc2} and \ref{SL2loc3} in terms of $\Phi$-parts explained above 
similarly as in the corresponding results in \cite{KNU2} Part II and Part IV indicated after each proposition. 
  That is, first we directly prove that the structure on each $\Phi$-part is independent on the choices of $(\alpha,\beta)$.
  Next, since the intersection of $\Phi$-part and $\Phi'$-part coincides with $\Phi \cap \Phi'$-part, it is enough to show that the localization map 
from $\Phi'$-part to $\Phi$-part is an open immersion whenever $\Phi'\subset \Phi$.
  The last statement is proved by using the version of Theorems \ref{SL2loc}, \ref{SL2loc2} and \ref{SL2loc3} in terms of $\Phi$-parts. %
\end{sbpara}

\subsection{The fan of parabolic subgroups}\label{Fan}

\begin{sbpara}\label{fan12} For a split torus $T$ over a field, let $X(T)=\Hom(T, {\bf G}_m)$ be the group of characters of $T$ and let $X_*(T)=\Hom({\bf G}_m, T)$ be the group of cocharacters of $T$.

In this Section \ref{Fan}, we give a variant  \ref{Wfan} of  the classical theory of Weyl fan (\ref{XTW}). Let $E$ be a field of characteristic $0$, let $\cG$ be a reductive algebraic group over $E$, let $T$ be a split torus over $E$, and let $a: T\to \cG$ be a homomorphism. Then in \ref{Wfan}, we will have a bijection between a certain set of parabolic subgroups of $\cG$ and the set of all cones of a certain cone decomposition of $\R \otimes X_*(T)$. In the case where $E$ is algebraically closed and $T$ is a maximal torus in $\cG$ with the inclusion map $a:T\to \cG$, this is the well-known  bijection (\ref{XTW}) between the set of all parabolic subgroups of $\cG$ which contain $T$ and the set of all cones of a cone decomposition called Weyl fan. This variant \ref{Wfan} should be also well-known, and %
is treated in our previous work \cite{KNU2} Part IV Section 2.6 in a certain situation.

 In the next Section \ref{s:star_to_BS}, we will use the results in \ref{Fan} to connect the space of $\SL(2)$-orbits and the space of Borel--Serre orbits

\end{sbpara}

\begin{sbpara}\label{Fan1}  
  Let $X^*$ and $X_*$ be finitely generated free abelian groups which are the $\Z$-duals of each other. We will denote the paring $X_*\times X^*\to \Z$ by $\langle\cdot,\cdot\rangle$.

Assume that we are given a finite subset $R$ of $X^*$ such that $R=-R$. 

In \ref{lemR}--\ref{SigR2}, we will show that we have a fan $\Sig(R)$  whose support is $\R \otimes X_*$. 
(Actually, the cones in this fan need not be sharp, and so $\Sig(R)$ should be called a quasi-fan. But we call it a fan for simplicity.)

For a finite subset $S$ 
of $X^*$, let $\langle S\rangle$ be the cone in $\R\otimes X^*$ generated by $S$ and let 
$$\sig(S)=\{y\in \R\otimes X_*\;|\; \langle y, x \rangle \geq 0\;\text{for all}\; x\in S\}.$$ 
We have $\langle S\rangle= \{x\in X^*\;|\; \langle y, x\rangle \geq 0\; \text{for all}\; y\in \sig(S)\}$. 

\end{sbpara}

\begin{sblem}\label{lemR} For a subset $R'$ of $R$, the following two conditions {\rm (i)} and {\rm (ii)} are equivalent.

{\rm (i)} There is $y\in X_*$ such that $R'=\{x\in R\;|\; \langle y, x\rangle \geq 0\}$.

{\rm (ii)} The following {\rm (ii-1)} and {\rm (ii-2)} are satisfied.

{\rm (ii-1)} $R=R'\cup {\rm (-R')}$.

{\rm (ii-2)} $\langle R'\rangle \cap R = R'$. 
\end{sblem}

\begin{pf}

The implication (i) $\Rightarrow$ (ii) is clear. 

We prove (ii) $\Rightarrow$ (i). Let $y$ be an interior point of $\sig(R')$. We prove that $R'$ satisfies (i) with this element $y$. Let $x\in R$ and assume $\langle y, x\rangle \geq 0$. We prove $x\in R'$. Assume $x\notin R'$. Then since $R=R'\cup (-R')$, we have $-x\in R'$. Hence $\langle y, -x\rangle \geq 0$ and hence $\langle y,-x\rangle =0$. Since $y$ is in the interior of $\sig(R')$ and $-x\in R'$, this shows that $-x=0$.  Hence $x=-x\in R'$.\end{pf}

\begin{sbpara}\label{SigR} Let $\Sig^*(R)$ be the set of all subsets $R'$ of $R$ satisfying the equivalent conditions in \ref{lemR}.

Note that $R'\in \Sig^*(R)$ is recovered from $\sig(R')$ as $R'=\{x\in R\;|\;\langle y, x\rangle\geq 0\; \text{for all}\; y\in \sig(R')\}$. 

Let $$\Sig(R)= \{\sig(R')\;|\; R'\in \Sig^*(R)\}.$$  We have a bijection $\Sig^*(R) \to \Sig(R)\;;\; R'\mapsto \sig(R')$. 

\end{sbpara}

\begin{sbprop}\label{SigR2} $\Sig(R)$ is a rational finite  fan whose support is $\R \otimes X_*$.
\end{sbprop}

\begin{pf} We need to prove the following (i), (ii) and (iii).

(i) If $\sig\in \Sig(R)$ and if $\tau$ is a face of $\sig$, then $\tau\in \Sig(R)$.

(ii) If $\sig, \tau\in \Sig(R)$, then $\sig\cap \tau\in \Sig(R)$.

(iii) If $\sig, \tau\in \Sig(R)$ and $\tau\subset \sig$, $\tau$  is a face of $\sig$.

Proof of (i).  Let $\sig=\sig(R')$ ($R'\in \Sig^*(R)$). Let $\tau$ be a face of $\sig$ and let $A=\{x\in \langle R'\rangle\;|\; \langle y, x\rangle =0\;\text{for all}\; y\in \tau\}$. Then $A$ is a face of $\langle R'\rangle$ and $\tau=\{y\in \sig\;|\; \langle y,x\rangle=0\; \text{for all}\; x\in A\}$. Let $S=R'\cap A$. Since the cone $\langle R'\rangle$ is generated by $R'$,  the cone $A$ is generated by $S$. Hence $\tau= \{y\in \sig\;|\; \langle y, x\rangle =0\; \text{for all}\; x\in S\}= \sig(R'\cup (-S))\in \Sig(R)$.

Proof of (ii). For $R'_1,R'_2\in \Sig^*(R)$, we have  $\sig(R'_1)\cap \sig(R'_2)=\sig(R'_3)$, where $R'_3=\langle R'_1\cup R'_2\rangle \cap R\in \Sig^*(R)$. 

Proof of (iii). If $R', R''\in \Sig^*(R)$ and if  $R'\subset R''$, then $R''=R'\cup (-S)$ for some subset $S$ of $R'$, and hence $\sig(R'')$ is a face of $\sig(R')$. 

We prove that the support of $\Sig(R)$ is $\R \otimes X_*$.   It is sufficient to prove that for each $y\in X_*$, there is $R'\in \Sig^*(R)$ such that $y\in \sig(R')$. In fact $R'=\{x\in  R\;|\; \langle y,x \rangle \geq 0\}$ has this property. 
\end{pf}

\begin{sbpara}\label{XT}
Now let $\cG$ be a reductive algebraic group over a field $E$ of characteristic $0$, let $T$ be an $E$-split torus, and let $a:T\to \cG$ be a homomorphism. 
Let $$X^*=X(T)=\Hom(T, \bG_m), \quad X_*=X_*(T)= \Hom(\bG_m, T)=\Hom(X(T), \Z).$$
Let $$R= \{\chi\in X^*\;|\; \Lie(\cG)_{\chi}\neq 0\},$$
where $\Lie(\cG)_{\chi}$ denotes the part of $\Lie(\cG)$ on which the adjoint action of $T$ via $a$ is given by $\chi$. 
\end{sbpara}

\begin{sbpara}\label{XTR} Let the notation be as in \ref{XT}. Then we have $R=-R$. This can be seen as follows. 

In the case where $E$ is algebraically closed and $T$ is a maximal torus in $\cG$ with $a:T \to \cG$  the inclusion map, $R=-R$  is well-known in the theory of root systems. 

The general case is reduced to this case by taking an algebraic closure %
$\bar E$ of $E$ and a maximal torus in %
$\cG \otimes_E \bar E$ which contains the image of %
$a:T\otimes_E \bar E \to \cG \otimes_E \bar E$.

\end{sbpara}

\begin{sblem}\label{XT2} Let the notation be as in $\ref{XT}$.  Then for a connected closed algebraic subgroup $P$ of $\cG$, the following two conditions {\rm (i)} and {\rm (ii)} are equivalent.  

 {\rm (i)} There is $y\in X_*$ such that $P$ is the parabolic subgroup of $\cG$ associated to the homomorphism $a\circ y: {\bf G}_m \to \cG$ in the sense of $\ref{P_is_para}$.

 {\rm (ii)} There is $R'\in \Sig^*(R)$ such that $\Lie(P)= \bigoplus_{\chi\in R'} \; \Lie(\cG)_{\chi^{-1}}$. 
Here we denote the group law of $X^*$ multiplicatively, and so $\chi^{-1}$ denotes the inverse of $\chi$.
\end{sblem}

\begin{pf} Assume that (i) is satisfied. Let $R'\in \Sig^*(R)$ be the set associated  to $y$ as in \ref{lemR} (i). Then (ii) is satisfied by this $R'$. 

Conversely assume that (ii) is satisfied. Take $y\in X_*$ which gives $R'$ as in \ref{lemR} (i), and let $P_1$ be the parabolic subgroup of $\cG$ associated to $a\circ y$. Then $\Lie(P)=\Lie(P_1)$. Since both $P$ and $P_1$ are connected, we have $P=P_1$. 
\end{pf}

\begin{sbpara}\label{Wfan} 
Let the notation be as in \ref{XT}.

Let $\cP$ be the set of all parabolic subgroups $P$ of $\cG$ satisfying the equivalent conditions in Lemma \ref{XT2}.

For $P\in \cP$, let $R(P):= \{\chi\in X^*\;|\; \Lie(P)_{\chi^{-1}}\neq 0\}.$

By Lemma \ref{XT2}, 
we have a bijection
$$\cP \quad\overset{1:1}\longleftrightarrow  \quad \Sig^*(R)$$
which sends $P\in \cP$ to $R(P)$ and conversely sends $R'\in \Sig^*(R)$ to the unique parabolic subgroup $P$ of $\cG$ such that $\Lie(P)= \bigoplus_{\chi\in R'}\;  \Lie(\cG)_{\chi^{-1}}$. 

Hence we have the composite bijection $$\cP\quad \overset{1:1}\longleftrightarrow \quad \Sig(R).$$
For $P\in \cP$, we denote the corresponding element of $\Sig(R)$ by $\sig(P)$. 

 \end{sbpara}

\begin{sbpara}\label{XTW}  Let the notation be as in \ref{XT}. Assume that $E$ is algebraically closed, $T$ is a maximal torus in $\cG$, and $a:T\to \cG$ is the inclusion map.

In this case, $\Sig(R)$ is called the Weyl fan and $\cP$ coincides with the set of 
  all parabolic subgroups of $\cG$ which contain $T$.

If $P\in \cP$ is a minimal parabolic subgroup, i.e., a Borel subgroup, the open cone of interior points of $\sig(P)$ is called the dominant Weyl chamber for $P$ (\cite{Mi2} Definition 21.35) 
and  $\Sig(R) =\bigcup_w\; (\text{faces of $w\sig(P)$})$,
where $w$ ranges over all elements of the Weyl group (while $P$ is fixed).

See \cite{Sp} Theorem 8.4.3.

\end{sbpara}

\begin{sbpara}\label{Wfan2}  Let the notation be as in \ref{XT}. Now fix an isomorphism $T\simeq \bG_{m,E}^n$. 

Then $X(T)$ is identified with $\Z^n$. Let $X(T)_+=\N^n\subset \Z^n=X(T)$. Let $R_+=R\cap X(T)_+$. Let $\Sig^*(R)_+$ be the subset of $\Sig^*(R)$ consisting of all  $R'$ such that $R_+\subset R'$. Let $\cP_+$ be the corresponding subset of $\cP$ and let $\Sig(R)_+\subset \Sig(R)$ be the corresponding subset. Then 
$\Sig(R)_+$ is a subfan of $\Sig(R)$. Its support is $X_*(R)_{\R,+}=\R^n_{\geq 0}=\sig(R_+)\subset \R^n=\R\otimes X_*$. (In fact, if $y\in X_*$ and $y\in \sig(R_+)$, then for $R':= \{x\in X^*\; |\; \langle y, x\rangle \geq 0\}$, we have $R'\in \Sig^*(R)_+$  and $y\in \sig(R')$.)

For $1\leq j\leq n$, let $W^{(j)}$ be the increasing filtration on the functor $V\mapsto V$ from $\Rep_E(\cG)$ to the category of $E$-vector spaces associated to $\bG_{m,E}\overset{\text{$j$-th}}\longrightarrow \bG_{m,E}^n \to \cG$ (\ref{Wlem}), where $j$-th means the $j$-th component. Let $\cG_\Phi\subset \cG$ be the stabilizer of $\Phi:=(W^{(j)})_{1\leq j\leq n}$, and let $\cG^{\circ}_\Phi$ be its connected component containing $1$. Then  $\cG^{\circ}_\Phi$ is the unique connected algebraic subgroup of $\cG$ such that $\Lie(\cG^{\circ}_\Phi)=\bigoplus_{\chi\in R_+} \Lie(\cG)_{\chi^{-1}}$.

For $P\in \cP$, $P\in \cP_+$ if and only if $\cG^{\circ}_{\Phi}\subset P$. 

\end{sbpara}

\begin{sbpara}\label{Wfan3} This is a complement to  \ref{Wfan2}.  Let $\cG$ be a reductive group and assume that we are given increasing filtrations 
$W^{(j)}$ ($1\leq j\leq n$)  on the functor $V\mapsto V$ from $\Rep_E(\cG)$ to the category of
 $E$-vector spaces such that there is a homomorphism $\bG_{m,E}^n \to \cG$  whose $j$-th $\bG_{m,E}\to \cG$ gives $W^{(j)}$ for $1\leq j\leq n$. Let  $\Phi=(W^{(j)})_{1\leq j\leq n}$ and let $\cG_\Phi$ be the stabilizer of  $\Phi$. Then we have a canonical homomorphism 
 $T:=\bG_{m,E}^n\to \cG_{\Phi, \red}$ whose every lifting $T\to \cG_\Phi$ gives $\Phi$. The  sets  $R_+$,  $\Sig^*(R)_+$, and  $\cP_+$, and the fan $\Sig(R)_+$ associated to $\bG_{m,E}^n\to \cG$ are independent of the choice of such lifting, for such liftings are conjugates in $\cG_\Phi$ of each other.

\end{sbpara}

\begin{sblem}\label{lemR+} Let the situation and the notation be as in $\ref{Wfan2}$. Then the following three conditions are equivalent.

{\rm (i)} $\cG^{\circ}_\Phi$ is a parabolic subgroup of $\cG$.

{\rm (ii)} $R_+\in \Sig^*(R)$, that is,  $R=R_+\cup (-R_+)$.

{\rm (iii)} $\Sig(R)_+$ coincides with the set of all faces of $\sig(R_+)$.

\end{sblem}

\begin{pf}
This is clear. 
\end{pf}

\subsection{Relation of $D^{\star}_{\SL(2)}$ and $D_{\BS}$}\label{s:star_to_BS}
\label{ss:star_to_BS} We  relate $D^{\star}_{\SL(2)}$ and $D_{\BS}$. To do this, the problem is that $D^{\star}_{\SL(2)}$ does not involve  parabolic subgroups though $D_{\BS}$ does. We define a modification $$D^{\star,\pa}_{\SL(2)}, %
\text{which is something like}\;\;%
\lq\lq D^{\star}_{\SL(2)}%
\;\text{plus parabolic subgroups,}
"%
$$
of $D^{\star}_{\SL(2)}$ and connect $D^{\star}_{\SL(2)}$ and $D_{\BS}$ via $D^{\star, \pa}_{\SL(2)}$. More precisely, we define a log modification $D^{\star,\pa}_{\SL(2)}\to D^{\star}_{\SL(2)}$ associated to cone decompositions related to parabolic subgroups, and define a morphism $D^{\star,\pa}_{\SL(2)}\to D_{\BS}$.

\begin{sbpara}\label{lgmdW}
Assume that $G$ is reductive.

For $\Phi\in \cW(G)$ and for a splitting $\alpha: \bG_{m,\R}^{\Phi}\to G_\R$ of $\Phi$, we apply \ref{Wfan3} to the case $E=\R$, $\cG=G_\R$, $T= \bG_{m,\R}^{\Phi}$, and $a=\alpha$.

Note that for $\Phi\in \cW(G)$, we have $X(\bG_m^{\Phi})_+=\bN^{\Phi} \to M_{>0}/\cO^\times_{>0}$ on  $D_{\SL(2)}(\Phi)$. The cone decomposition of $X_*(\bG_m^{\Phi})_{\R,+}$ in \ref{Wfan3}
defines a log modification $D^{\pa}_{\SL(2)}(\Phi)$ of $D_{\SL(2)}(\Phi)$. It is independent of the choice of the splitting $\alpha$. When $\Phi$ moves, these are glued to a log modification $D^{\pa}_{\SL(2)}\to D_{\SL(2)}$.
(Here the superscript W respects the Weyl fan.
For a log modification in the category $\cB'_\R(\log)$, see \cite{KNU2} Part IV 1.4.6.) %

In general, for a linear algebraic group $G$ over $\Q$, we define $D^{\star,\pa}_{\SL(2)}$ to be the fiber product of 
$D^{\star}_{\SL(2)}\to D_{\red,\SL(2)} \leftarrow  D^{\pa}_{\red,\SL(2)}$.

\end{sbpara}

\begin{sbpara}\label{lgmdW2} As a set, $D^{\star,\pa}_{\SL(2)}$ is identified with the set of triples $(x, P, Z)$, where $p\in D_{\red,\SL(2)}$, $x:=(p, Z')\in D^{\star}_{\SL(2)}$, $P$ is a parabolic subgroup of $G_{\red}$ satisfying the conditions in \ref{Wfan2} and the condition $G_{\red,\Phi,u} \subset P_u$ with $\Phi$ being the set of weight filtrations associated to $p$, and $Z$ is a subset of $Z'$ satisfying the following (i). Let $A_{p,P}$ be the inverse image of $A_P$ under $\tau_p: \R_{>0}^{\Phi} \to P_{\red}(\R)$. 

(i)  If $x$ is an $A$-orbit, $Z$ is a $\tau_x^{\star}(A_{p,P})$-orbit. If $x$ is a $B$-orbit, $Z$ is a $\tau^{\star}_x(\R_{>0}\times A_{p, P})$-orbit.

\medskip

The map $D^{\star,\pa}_{\SL(2)} \to D^{\star}_{\SL(2)}$ is understood as $(x, P, Z) \mapsto x$. 

\end{sbpara}

\begin{sbpara}\label{eta}  We have a map $$D^{\star,\pa}_{\SL(2)}\to D_{\BS}\;;\; (x, P, Z) \mapsto (P, A_P \circ Z).$$

The fact that this is a morphism is proved by using the local structure theorem (Theorem \ref{SL2loc2}). The proof is similar to the proof of \cite{KNU2}  Part IV Theorem 2.6.22 (1). 
\end{sbpara}

\begin{sbrem}
\label{r:par}
 In \cite{KNU2} Part IV Sections 2.5 and 2.6, we considered  
$$D^{\star,+}_{\SL(2)} \to D^{\star}_{\SL(2)} \to D^{\star,-}_{\SL(2)} \leftarrow D^{\star,\BS}_{\SL(2)} \to D_{\BS}.$$ 

In that situation, $D^{\pa}_{\red,\SL(2)}$ in \ref{lgmdW} coincides with $D_{\SL(2)}(\gr^W)^{\BS}$ in \cite{KNU2} Part IV 2.6.3, $D^{\star,\pa}_{\SL(2)}$ is the fiber product of $D^{\star}_{\SL(2)} \to D^{\star,-}_{\SL(2)} \leftarrow D^{\star,\BS}_{\SL(2)}$, and the map $D^{\star,\pa}_{\SL(2)} \to D_{\BS}$ in \ref{eta} coincides with the composition $D^{\star, \pa}_{\SL(2)}\to D^{\star,\BS}_{\SL(2)}\to D_{\BS}$. 
The fiber product property can be seen as 
 $$
D^{\star}_{\SL(2)}\times_{D^{\star,-}_{\SL(2)}}D^{\star,\BS}_{\SL(2)}
=D^{\star}_{\SL(2)}\times_{D^{\star,-}_{\SL(2)}}(D^{\star,-}_{\SL(2)}\times 
_{D_{\SL(2)}(\gr^W)}D_{\SL(2)}(\gr^W)^{\BS})
$$
$$
=D^{\star}_{\SL(2)}\times _{D_{\SL(2)}(\gr^W)}D_{\SL(2)}(\gr^W)^{\BS}
=D^{\star}_{\SL(2)}\times _{D_{\SL(2)}(\gr^W)}D_{\red,\SL(2)}^{\pa}
=D^{\star,\pa}_{\SL(2)},
$$
where the first equality is \cite{KNU2} Part IV Proposition 2.6.14.
\end{sbrem}

\begin{sbrem}
(1) In the second line of \cite{KNU2} Part IV 2.6.3, 
$D^{\star,-}_{\SL(2)}(\gr^W)$ should be 
$D_{\SL(2)}(\gr^W)$. 
In the last line of
 loc. cit.,  $D^{\star}_{\SL(2)}(\gr^W)^{\BS}$ should be 
$D_{\SL(2)}(\gr^W)^{\BS}$.

(2) In the proof of \cite{KNU2} Part IV Proposition 2.6.9, line 11 from the end of the proof, $L = \cS(\sig)\cup \cS(\sig)^{-1}$ must be corrected as $R(Q)\subset \cS(\sig)\cup \cS(\sig)^{-1}$. 

\end{sbrem}

\subsection{Case of Shimura varieties}\label{s:Shim}

\begin{sbpara}\label{Shi0} Assume that $G$ is reductive and that $h_0:S_{\C/\R}\to G_\R$ satisfies the condition that the Hodge type of $\Lie(G_{\bR})$ via $h_0$ is in $\{(1,-1), (0, 0), (-1,1)\}$ (as in Shimura data). Then $h_0$ is $\R$-polarizable
 by \cite{De} (Lemma \ref{pol2}).

\end{sbpara}
We prove 
\begin{sbthm}\label{Shim}  Let the assumption be as in $\ref{Shi0}$. 
Then we have an isomorphism  $$D^{\pa}_{\SL(2)}\overset{\sim}\to D_{\SL(2)}$$ in $\cB'_\R(\log)$. In particular, the identity map of $D$ extends uniquely to a morphism $$D_{\SL(2)}\to D_{\BS}$$ of locally ringed spaces with log structures with sign. 

\end{sbthm}

\begin{sbpara} Note that for a field $E$ of characteristic $0$ and for $n\geq 0$, a finite-dimensional representation of $\SL(2)^n_E$ over $E$ is semisimple and each irreducible representation is isomorphic to $\rho^{(r)}:= \Sym^{r(1)}(\rho_1)\otimes \dots \otimes \Sym^{r(n)}(\rho_n)$ for some $r\in \N^n$, where $\rho_j: \SL(2)^n_E\to \GL(2)_E$ is the composition of the $j$-th projection $\SL(2)_E^n 
\to \SL(2)_E$ and the inclusion homomorphism $\SL(2)_E\to \GL(2)_E$.  Consider the homomorphism $\bG_{m,E}^n \to \SL(2)^n_E$ whose restriction to the $k$-th $\bG_{m,E}$ ($1\leq k\leq n$) is
$$t \mapsto (g_1, \dots, g_n), \;\; g_j=\begin{pmatrix} t^{-1}&0 \\0 & t\end{pmatrix}\;\; \text{for}\;1\leq j\leq k, \;\; g_j=1\;\; \text{for}\;\; k<j\leq n.$$
Write the standard base of $E^2$ by $(e_1, e_2)$. Then the action of $\bG_{m,E}^n$ on $\bigotimes_{j=1}^n (e_1^{a(j)}e_2^{r(j)-a(j)})$ ($a\in \N^n$, $a\leq r$) via $\rho^{(r)}$ is given by the character 
$$(\sum_{j=1}^k r(j)-2 \sum_{j=1}^k a(j))_{1\leq k\leq n}\in \Z^n=X(\bG_{m,E}^n).$$

\end{sbpara}

We will use the following lemma later.

\begin{sblem}\label{rhor}

Let $r\in \N^n$ and let  $c\in \Z^n=X(\bG_{m,E}^n)$ be a character which appears in the representation of $\bG_{m,E}^n$  in $\rho^{(r)}$ via the above homomorphism $\bG^n_{m,E}\to \SL(2)^n_E$. 
Assume that either one of the following conditions (i)--(iii) is satisfied.
Then we have either $c\in \N^n$ or $-c\in \N^n$. 

(i) $r(k)=0$ for all $k$.

(ii)  There is $k$ such that $1\leq k\leq n$, $r(k)>0$, and $r(j)=0$ for all $j\neq k$.

(iii)  There are $k,l$ such that $1\leq k<l\leq n$, $r(k)>0$, $r(l)=1$, and $r(j)=0$ for $j\neq k, l$.
\end{sblem}

\begin{pf} We consider the character $c$ of the action of $\bG_{m,\R}^n$ on $\bigotimes_{j=1}^n (e_1^{a(j)}e_2^{r(j)-a(j)})$ ($0\leq a(j)\leq r(j)$ for $1\leq j\leq n$). 

In the case (i), $c=0$. 

In the case (ii), $c(j)=0$ if $j<k$ and $c(j)=r(k)-2a(k)$ if $j\geq k$.

In the case (iii), $c(j)= 0$ if $j<k$, $c(j)=c(k)=r(k)-2a(k)$ if $k\leq j<l$, and $c(j)=c(l)=r(k)+r(l)-2(a(k)+a(l))$  if $j\geq l$. Since  $r(l)=1$ and $a(l)\in \{0,1\}$, we have  $|c(k)-c(l)|=1$. Hence we have either $\{c(k), c(l)\}\subset \N$ or $\{-c(k), -c(l)\}\subset \N$. Since  $c(j)\in \{0, c(k), c(l)\}$ for all $j$, we have either   $c\in \N^n$ or $-c\in \N^n$. 
\end{pf}

\begin{sbpara}\label{rhor2}  Assume that $G$ is reductive. Let  $(\rho, \varphi)$ be an $\SL(2)$-orbit in $n$ variables of rank $n$ for $(G, h_0)$, and let $V\in \Rep(G)$. Let $r\in \N^n$ and assume that  $\rho^{(r)}$ appears in the action of $\SL(2)^n_\R$ on $V_\R$ induced by $\rho$. 
Then, by Claim in the proof of \ref{redp2},  for the Hodge structure of $V$ given by any element of $D$, 
there is $p\in \Z$ such that the $(p+b, p-b)$-Hodge component of $V$ is non-zero for $0\leq b \leq \sum_{j=1}^n r(j)$.

\end{sbpara}

\begin{sbpara} We prove Theorem \ref{Shim}.  Let $\Phi\in \cW(G)$. Let $(\rho, \varphi)$ be an $\SL(2)$-orbit in $n$ variables of rank $n$ whose associated family of weight filtrations is $\Phi$. Then 
by \ref{rhor2} and by the fact that only the Hodge type $(1,-1), (0, 0), (-1,1)$ appears in $\Lie(G)$,  we have that if $\rho^{(r)}$ appears  in the representation $\Lie(G_{\bR})$ of $\SL(2)^n_\R$, then we have $\sum_{j=1}^n r(j) \leq 2$. 
Hence if $\rho^{(r)}$ appears in $\Lie(G_{\bR})$, the assumption of \ref{rhor} for $r$ is satisfied. Hence 
by \ref{rhor}, 
 each character of $\bG_{m,\R}^n$ which appears in $\Lie(G_{\bR})$ is either in $R_+$ or $-R_+$, that is, $R=R_+\cup (-R_+)$. Hence the condition (ii) of \ref{lemR+} is satisfied for $\Phi$.  Hence the condition (iii) of \ref{lemR+} is satisfied. This proves that $D^{\pa}_{\SL(2)}(\Phi)\to D_{\SL(2)}(\Phi)$ is an isomorphism and hence proves Theorem \ref{Shim}.

\end{sbpara}

\begin{sbpara}
In the classical case of $h_0: S_{\C/\R}\to \text{GSp}(g)_\R$ which gives the Siegel upper half space
$\frak H_g$ of degree $g$, $D_{\SL(2)}\to D_{\BS}$ is a homeomorphism (\cite{KU1} Theorem 6.7). 

But even in the case $g=1$,  the real analytic structures of  $D_{\BS}$ and $D_{\SL(2)}$ are slightly different as is seen in \ref{Ex1}. 
For the case of some Shimura variety as in \ref{Ex2}, the map $D_{\SL(2)}\to D_{\BS}$ is not bijective. 

\end{sbpara}

\subsection{Relation of $D^{\star}_{\SL(2)}$ and $D^{II}_{\SL(2)}$}
\label{ss:star_to_II}

We define a log modification 
$D^{\star,+}_{\SL(2)}\to D^{\star}_{\SL(2)}$ associated to a cone decomposition, and define a morphism $D^{\star,+}_{\SL(2)}\to D^{II}_{\SL(2)}$. 

\begin{sbpara}\label{lgmd+}
  Let $\Phi$ be a finite set of weight filtrations which has a common splitting.

Note that we have $X(\bG_m \times \bG_m^{\Phi})_+=\bN \times \bN^{\Phi}\to M_{>0}/\cO^\times_{>0}$ on $D^{\star}_{\SL(2)}(\Phi)$. 
 Let $\beta_0^{\star} \in M/\cO^{\times}$ be the image of $(1,0,...,0)$ and 
let $\beta_{\mathrm{tot}} \in M/\cO^{\times}$ be the image of $(0,1,...,1)$.
These $\beta_0^{\star}$ and $\beta_{\mathrm{tot}}$ are glued to global sections of $M_{>0}/\cO^\times_{>0}$ on $D^{\star}_{\SL(2)}$ which we still denote by $\beta_0^{\star}$ and $\beta_{\mathrm{tot}}$, respectively.

Let $D^{\star,+}_{\SL(2)}(\Phi)$ be the log modification of $D^{\star}_{\SL(2)}(\Phi)$ associated to the cone decomposition of $X_*(\bG_m \times \bG_m^{\Phi})_+=\bN \times \bN^{\Phi}$ consisting of cones 
$$\sig_1:=\{(x_0, x_1, \dots, x_n)\;|\; x_0\leq {\textstyle{\sum}}_{j=1}^n x_j\}\text{ and }\sig_2:=\{(x_0, x_1, \dots, x_n)\;|\; x_0\geq {\textstyle{\sum}}_{j=1}^n x_j\}$$ and their faces. 

When $\Phi$ moves, these are glued to a log modification $D^{\star,+}_{\SL(2)}\to D^{\star}_{\SL(2)}$.

For $j=1,2$, let $D^{\star,+}_{\SL(2)}(\sig_j)$ be the open set of  $D^{\star,+}_{\SL(2)}$ whose intersection with $D_{\SL(2)}^{\star,+}(\Phi)$ coincides with its $\sig_j$-part. Then $D^{\star,+}_{\SL(2)}(\sig_1)$ (resp.\ $D^{\star, +}_{\SL(2)}(\sig_2)$) coincides with the set of all points $s$ of $D^{\star.+}_{\SL(2)}$ such that $\beta_{\mathrm{tot}}$ (resp.\ $\beta_0^{\star}$) is divided by $\beta_0^{\star}$ (resp.\ $\beta_{\mathrm{tot}}$) at $s$. 

\end{sbpara}

\begin{sbpara}\label{+II} We define a map $D^{\star,+}_{\SL(2)}\to D^{II}_{\SL(2)}$ as follows (cf.\ \cite{KNU2} Part IV 2.5.4).

   Let $x^+$ be a point of $D^{\star,+}_{\SL(2)}$ lying over $x \in D^{\star}_{\SL(2)}$. 
  We define the image $x^{II}$ of $x^+$ in $D^{II}_{\SL(2)}$. 
  There are four cases. 

Case 1. Both $\beta_0^{\star}$ and $\beta_{\mathrm{tot}}$ are trivial at $x^+$. That is, $x^+=x\in D$.

Case 2. $\beta_0^{\star}$ is strictly divided by $\beta_{\mathrm{tot}}$ at $x^+$. 

Case 3. $\beta_{\mathrm{tot}}$ is strictly divided by $\beta_0^{\star}$ at $x^+$.

Case 4. $\beta_0^{\star}$ and $\beta_{\mathrm{tot}}$ coincide at $x^+$ but are nontrivial. 

  In Case 1, $x^{II}=x^+=x\in D$.  
  
 In Cases 2--4, write $x=(p, Z)$.  
  
  In Case 2, $x^{II}$ is $x=(p,Z)$ regarded as an element (a $B$-orbit) of  $D_{\SL(2)}$. 

  In Case 3, $x^{II}$ is $(p,Z_{\spl})$, where $Z_{\spl}:=\{\spl_W(z)(z_{\red})\,|\,z \in Z\}$. 
  (See \ref{Dspl} for the notation.)
  
  In Case 4, $Z$ in  $x=(p, Z)$ is a $\tau_x^{\star}(\R_{>0}\times \R_{>0}^{\Phi(p)})= \tau_x(\R_{>0}\times \R_{>0}^{\Phi(p)})$-orbit, 
and $x^+$ is identified with a triple $(p, Z, Z')$ where $Z'$ is a $\tau_x(\{1\} \times \R_{>0}^{\Phi(p)})$-orbit contained in $Z$. %
We define $x^{II}= (p, Z')$.

The proof of the fact that the map just defined is a proper and surjective morphism is similar to that of \cite{KNU2} Part IV Theorem 2.5.5 (1). 
\end{sbpara}

  \begin{sbprop}\label{mtoII} There is a unique morphism $D^{\star,\mild}_{\SL(2)}\to D^{II}_{\SL(2)}$ of $\cB'_\R(\log)$ which extends the identity map of $D$. 
  
  \end{sbprop}
  
 This is a $G$-MHS version of a part of \cite{KNU2} Part IV Theorem 2.5.5 (1), and proved similarly as follows. The map $D^{\star,+}_{\SL(2)}\to D^{\star}_{\SL(2)}$ is an isomorphism over the open set $D^{\star,\mild}_{\SL(2)}$ of $D^{\star}_{\SL(2)}$ as is easily seen. Hence the morphism $D^{\star,+}_{\SL(2)}\to D^{II}_{\SL(2)}$ induces a morphism $D^{\star,\mild}_{\SL(2)}\to D^{II}_{\SL(2)}$.
 
\begin{sbrem} Let $\cB'_\R(\log)^+$ be the full subcategory of $\cB'_\R(\log)$ defined in Part IV 2.7.5 of \cite{KNU2}. 
  Then by the method in Section 2.7 of ibid. basing on the local structure theorems (Theorems \ref{SL2loc2} and \ref{SL2loc3}), we can prove the following $G$-MHS version of Theorem 2.7.14 of ibid. The spaces $D^I_{\SL(2)}$, $D^{II}_{\SL(2)}$, $D^{\star}_{\SL(2)}$, $D^{\star,\pa}_{\SL(2)}$, $D^{\star,+}_{\SL(2)}$ belong to $\cB'_\R(\log)^+$. We do not give the details of the proof.

\end{sbrem}

\subsection{Valuative spaces, I}
\label{ss:val1}

Recall that, for an abelian group $L$, a submonoid $V$ of $L$ is said to be {\it valuative} if $V\cup V^{-1}=L$.

  \begin{sbpara}\label{val1} 
  We review the associated valuative space. 
  
  For an object $S$ of $\cB'_\R(\log)$, we have a locally ringed space $S_{\val}$ endowed with a  log structure with sign defined as in \cite{KNU2} Part IV 3.1.13. As a set, $S_{\val}$ is the set of triples $(s, V, h)$, where $s\in S$,  $V$ is a valuative submonoid of $(M_S/\cO_S^\times)_s^{\gp}$ such that $V\supset (M_S/\cO^\times_S)_s$ and such that $V^\times \cap (M_S/\cO_S^\times)_s= \{1\}$, and, $\tilde V_{>0}$ being the inverse image of $V$ in $M_{S,>0,s}^{\gp}$, $h$ is a homomorphism $(\tilde V_{>0})^\times \to \R_{>0}^{\mult}$  extending  the evaluation homomorphism $f \mapsto  f(s)$ on $\cO_{S,>0,s}^\times$ at $s$. 

\cite{KNU2} Part IV Proposition 3.1.9 explicitly describes the projection $S_{\val}\to S$ as a projective limit of log modifications of $S$ (\cite{KNU2} Part IV Proposition 1.4.6).
It follows that the projection $S_{\val}\to S$ is proper and surjective (Corollary 3.1.10).
  \end{sbpara}
  
  \begin{sbpara}\label{val2} 
  By \ref{val1}, we have the following locally ringed spaces with a log structure with sign
  $$D_{\BS,\val}, \;D^{\star}_{\SL(2),\val}, \; D^I_{\SL(2),\val}, \; D^{II}_{\SL(2),\val}$$
  associated to the objects $D_{\BS}$, $D^{\star}_{\SL(2)}$,  $D^I_{\SL(2)}$, $D^{II}_{\SL(2)}$ of $\cB'_\R(\log)$, respectively. 
  
  The underlying sets of $D^I_{\SL(2),\val}$ and $D^{II}_{\SL(2),\val}$ are identified because the log structure with sign of 
 $D^I_{\SL(2)}$ is the inverse image of that of $D^{II}_{\SL(2)}$.
We will denote their common underlying set by $D_{\SL(2),\val}$.
  
  \end{sbpara}
  
  \begin{sbthm}\label{valthm} 
$(1)$ There is a unique morphism 
   $D^{\star}_{\SL(2),\val}\to D_{\BS,\val}$ which extends the identity map of $D$. This map is injective.
   
$(2)$ There is a unique morphism $D^{\star}_{\SL(2),\val}\to D^{II}_{\SL(2),\val}$ which extends the identity map of $D$. It is proper and surjective. 
  \end{sbthm}
  
\begin{pf}  %
(1) The morphism $D^{\star,\pa}_{\SL(2)} \to D_{\BS}$ in Section \ref{ss:star_to_BS} induces the morphism $D^{\star}_{\SL(2),\val}\to D_{\BS,\val}$ which extends the identity map of $D$ because $D^{\star,\pa}_{\SL(2)}$ is a log modification of $D^{\star}_{\SL(2)}$ (\ref{lgmdW2}). 
  The uniqueness is by the density of $D$ in $D^{\star}_{\SL(2),\val}$. The proof of the injectivity is similar to the proof of \cite{KNU2} Part IV Theorem 3.4.4 (1).

  (2) The morphism $D^{\star,+}_{\SL(2)} \to D^{II}_{\SL(2)}$ in Section \ref{ss:star_to_II} induces the morphism $D^{\star}_{\SL(2),\val}\to D^{II}_{\SL(2),\val}$ which extends the identity map of $D$ because $D^{\star,+}_{\SL(2)}$ is a log modification of $D^{\star}_{\SL(2)}$ (cf.\ \ref{+II}). 
  The uniqueness is by the same reason as in (1).

This morphism is proper, because both terms are proper over $D_{\red,\SL(2)}\times \spl(W)$ by Proposition \ref{SLbdl}. The surjectivity of this map follows from its properness and the fact that the image of this map contains $D$ and hence is dense in $D^{II}_{\SL(2),\val}$. 
\end{pf}

  \begin{sbpara}
   The inverse image of $D^{\mild}_{\BS,\val}$ under $D^{\star}_{\SL(2),\val}\to D_{\BS,\val}$ is $D^{\star,\mild}_{\SL(2),\val}$. 
  \end{sbpara}
  
  \begin{sbrem} Let $\cC_\R(\val)^+$ be the category defined in Part IV 3.2.5 of \cite{KNU2}. Then by the method in Section 3.2 of ibid., we can prove the following $G$-MHS version of a statement in 3.3.1 of ibid. The spaces $D_{\BS, \val}$, $D^I_{\SL(2),\val}$, $D^{II}_{\SL(2),\val}$, $D^{\star}_{\SL(2),\val}$ belong to $\cC_\R(\val)^+$. We do not give the details of the proof.

\end{sbrem}

\subsection{Global properties of $D^I_{\SL(2)}$, $D^{II}_{\SL(2)}$, $D^{\star}_{\SL(2)}$ etc.} \label{ss:SL2gl}

  \begin{sbthm}\label{SL2gl}  Let $X$ be one of $D^I_{\SL(2)}$, $D^{II}_{\SL(2)}$, $D^{\star}_{\SL(2)}$, $D^I_{\SL(2),\val}$, $D^{II}_{\SL(2),\val}$, $D^{\star}_{\SL(2),\val}$, $D_{\BS,\val}$, $D^{\star,+}_{\SL(2)}$, $D^{\star,\pa}_{\SL(2)}$.
Let $\Gamma$ be a semi-arithmetic subgroup of $G'(\Q)$ ($\ref{Gamma}$).

$(1)$ The action of $\Gamma$ on $X$ is proper and the quotient space $\Gamma\bs X$ is Hausdorff. In particular, $X$ is Hausdorff.

$(2)$ If $\Gamma$ is torsion-free, the action of $\Gamma$ on $X$ is free and  the map $X\to \Gamma\bs X$ is a local homeomorphism. 
\end{sbthm}

  The proof, given in \ref{SL2H} and \ref{SL2H2} below, is similar to that of \cite{KNU2} Part II Theorem 3.5.17 and that of \cite{KNU2} Part IV Theorem 6.1.1.
  Starting from $D_{\BS}$, we transport various properties along the fundamental diagram in Introduction.

\begin{sbpara}\label{SL2H} We first prove that the case $\Gamma=\{1\}$ of (1), that is, the space $X$ is Hausdorff. 

  We have an injective continuous map $D^{\star}_{\SL(2),\val}\to D_{\BS,\val}$ (Theorem \ref{valthm} (1)).  
 Since $D_{\BS}$ is Hausdorff (Proposition \ref{p:DBSHaus}), $D_{\BS,\val}$ is Hausdorff.
  Hence, $D^{\star}_{\SL(2),\val}$ is Hausdorff. 
 Since $D^{\star}_{\SL(2),\val} \to D^{\star,\pa}_{\SL(2)}$, $D^{\star}_{\SL(2),\val} \to D^{\star,+}_{\SL(2)}$ and $D^{\star}_{\SL(2),\val} \to D^{\star}_{\SL(2)}$ are proper and surjective, $D^{\star,\pa}_{\SL(2)}$, $D^{\star,+}_{\SL(2)}$ and $D^{\star}_{\SL(2)}$ are Hausdorff. 
  Since the maps $D^{\star}_{\SL(2),\val}\to D^{II}_{\SL(2),\val} \to D^{II}_{\SL(2)}$ are proper and surjective (Theorem \ref{valthm} (2)), $D^{II}_{\SL(2),\val}$ and $D^{II}_{\SL(2)}$ are Hausdorff. 
  Since we have a bijective continuous map $D^I_{\SL(2)}\to D^{II}_{\SL(2)}$, $D^I_{\SL(2)}$ is Hausdorff. 
   Hence $D^I_{\SL(2),\val}$ is also Hausdorff. 
  \end{sbpara}
  
\begin{sbpara}\label{SL2H2} We prove Theorem \ref{SL2gl}. 

  (1) We prove the former part, that is, that the action is proper. 
  Since it is valid for $X=D_{\BS}$ (Theorem \ref{BSgl} (1)), 
we see that it is valid for $D_{\BS, \val}$ and $D^{\star}_{\SL(2),\val}$ by using continuous maps
  $D^{\star}_{\SL(2),\val}\to D_{\BS,\val}\to D_{\BS}$,  \ref{SL2H},  and the fact \ref{proper} (3.1). 
  Then we see that it is valid also for $D^{\star,\pa}_{\SL(2)}$, $D^{\star,+}_{\SL(2)}$, $D^{\star}_{\SL(2)}$, $D^{II}_{\SL(2),\val}$ and $D^{II}_{\SL(2)}$ by using the proper and surjective maps $D_{\SL(2),\val}^{\star}\to D^{\star,\pa}_{\SL(2)}$, $D_{\SL(2),\val}^{\star}\to D^{\star,+}_{\SL(2)}$, $D_{\SL(2),\val}^{\star}\to D^{\star}_{\SL(2)}$, $D^{\star}_{\SL(2),\val} \to D_{\SL(2),\val} \to D_{\SL(2)}$ and the fact \ref{proper} (3.2). 
  Then we see it also for $D^I_{\SL(2)}$ and $D^I_{\SL(2),\val}$ by using the continuous maps $D^I_{\SL(2)}\to D^{II}_{\SL(2)}$ and $D^I_{\SL(2),\val}\to D^{II}_{\SL(2),\val}$. 

The latter part of (1) follows from the former part by \ref{proper} (1).

  (2) We prove the former part, that is, the action is free, by using a similar argument as in \ref{pfBSgl}.  
We apply \ref{proper} (4) to $H=\Gamma$, $H_1=\Gamma_u$, $X_1=\spl(W)$, $X_2=D_{\red,\SL(2)}$.
The action of $\Gamma/\Gamma_u$ on $D_{\red,\SL(2)}$ is free, by a similar argument in the proof of \cite{KU1} Lemma 5.7, and the action of $\Gamma_u$ on $\spl(W)\simeq G_u(\R)$ is free.
Hence the action of $\Gamma$ on $\spl(W)\times D_{\red,\SL(2)}$ is free.
  By using the canonical maps from $D_{SL(2)}^{I}$, $D_{SL(2)}^{II}$, $D^{\star}_{\SL(2)}$ to $\spl(W)\times D_{\red,\SL(2)}$ together with the related maps in the fundamental diagram, we see that the action of $\Gamma$ on $X$ is free.

The latter part of (2) follows from the former by \ref{proper} (2), the properness of the action proved in (1), and \ref{SL2H}. 

\end{sbpara}

\begin{sbpara}\label{SL2SA1} 
  
 In Theorem \ref{SL2gl}, we can use  a semi-arithmetic subgroup of $G(\Q)$ (not of $G'(\Q)$)  in the following two situations. 
 
 First, 
  if either $G$ is semisimple or  the condition (1) in Lemma \ref{pol2} is satisfied,   \ref{SL2gl} holds for a semi-arithmetic subgroup $\Gamma$ of $G(\Q)$. In fact, $\Gamma \cap G'(\Q)$ is of finite index (cf. Proposition \ref{A=A*} for the latter case). Hence by \ref{proper} (5), we can replace $\Gamma$ by the semi-arithmetic subgroup $\Gamma\cap G'(\Q)$ of $G'(\Q)$.

 Next

\end{sbpara}

\begin{sbprop}\label{SL2SA2} Assume that $G$ is reductive.
Let $X$ be one of $D_{\SL(2)}$, $D_{\SL(2),\val}$,  $D_{\BS,\val}$. 
Let $\Gamma$ be a semi-arithmetic subgroup of $G(\Q)$.

$(1)$ The quotient space $\Gamma\bs X$ is Hausdorff.

$(2)$ Let $Z$ be the center of $G$. If the image of $\Gamma$ in $(G/Z)(\Q)$ is torsion-free, the map $X\to \Gamma\bs X$ is a local homeomorphism. 
\end{sbprop}

See \ref{redpf} for the proof. 

%


\section{The space of nilpotent orbits}\label{s:DSig}

Let $D=D(G,h_0)$ be as in \ref{D}.
We assume that $h_0$ is $\R$-polarizable (\ref{pol}).

  In this section, we define and study  the toroidal partial compactification $\Gamma \bs D_{\Sig}$ of $\Gamma \bs D$, an extended period domain consisting of nilpotent orbits. 
  We consider $\Gamma \bs D_{\Sig}$ as the moduli of $G$-log mixed Hodge structures. 
  It is  the $G$-MHS version of the toroidal partial compactification  $\Gamma \bs D(\La)_\Sig$ of $\Gamma\bs D(\Lambda)$ for the classical period domain $D(\Lambda)$ (\ref{ss:Dold}). This $\Gamma \bs D(\La)_{\Sig}$  is  
   defined and studied in \cite{KU2} (in the pure case) and in \cite{KNU2} Part III and Part IV (in the mixed case)  and is the moduli space of LMH.

\subsection{The sets $D_{\Sig}$ and $D_{\Sig}^{\sharp}$}
\label{ss:setDSig}

\begin{sbpara}\label{nilp1}
  A {\it nilpotent cone} is a subset $\sig$ of $\Lie(G'_\R)$ 
  satisfying the following (i)--(iii).
  
  (i) $\sig=\R_{\geq 0}N_1+\dots +\R_{\geq 0}N_n$ for some $N_1,\dots,N_n\in \Lie(G'_\R)$. 
  
  (ii) For every $V\in \Rep(G)$, the image of $\sig$ under the induced map $\Lie(G_\R)\to \End_\R(V_{\bR})$ consists of nilpotent operators. 

  (iii) We have $[N,N']=0$ for $N, N'\in \sig$. 
\end{sbpara}

\begin{sbpara}
Recall that for a vector space $V$ over a field $E$, for an increasing filtration $W$ on $V$, and for a nilpotent linear map $N: V\to V$ such that $NW_w\subset W_w$ for all $w$, an increasing filtration $M$ on $V$ is called the {\it relative monodromy filtration} of $N$ relative to $W$ if $NM_w\subset M_{w-2}$ for all $w$ and $N^m: \gr^M_{w+m}\gr^W_w\overset{\sim}\to \gr^M_{w-m}\gr^W_w$ for all $w$ and all $m\geq 0$. 

The relative monodromy filtration $M$ need not exist, but it is unique if it exists.

\end{sbpara}

\begin{sbpara}\label{nilp2}
    Let $F\in \Dc$ and let $\sig$ be a nilpotent cone.
     We say that the pair $(\sig, F)$ {\it generates a nilpotent orbit}  if the following conditions (i)--(iii) are satisfied.
\medskip
     
 {\rm (i)} Let $N_1, \dots, N_n$ be as in (i) in \ref{nilp1}. Then $\exp(\sum_{j=1}^n z_jN_j)F\in D$ if $z_j\in \C$ and $\text{Im}(z_j)\gg 0$. 

{\rm (ii)}  $NF^p\subset F^{p-1}$  for all $N\in \sig$ and $p\in \Z$.

{\rm (iii)} For every $N\in \sig$ and every $V\in \Rep(G)$, the relative monodromy filtration of $N:V_\R\to V_\R$ with respect to $W(V)_\R$ exists. 
           \medskip

In this case, we also say that $(N_1,\dots,N_n,F)$ {\it generates a nilpotent orbit}.

Note that the above condition (i) is independent of the choice of $(N_1, \dots, N_n)$ as in (i) in \ref{nilp1}.  Note also that it is equivalent to the condition that $\exp(\sum_{j=1}^n iy_jN_j)F\in D$ if $y_j\in \R$ and $y_j\gg 0$.

\end{sbpara}

\begin{sbprop}\label{nilp3} Let $F\in \Dc$ and let $\sig$ be a nilpotent cone. Then the following conditions {\rm (i)}, {\rm (ii)}, and {\rm (iii)} are equivalent.

{\rm (i)} $(\sig,F)$ generates a nilpotent orbit in the sense of $\ref{nilp2}$.

{\rm (ii)} For every $V\in \Rep(G)$, we have the following 
{\rm (ii-1)}, {\rm (ii-2)}, {\rm (ii-3)}. 

\quad {\rm (ii-1)} 
For each $w\in \Z$, there  is a $G'_\R$-invariant $\R$-bilinear form $\langle\cdot, \cdot\rangle_w: \gr^W_wV_\R\times \gr^W_wV_\R \to \R$ such that if $z_j\in \C$ and $\text{Im}(z_j)\gg 0$ ($1\leq j\leq n$), $(\gr^W_wV_\R, \exp(\sum_{j=1}^n z_jN_j)F(\gr^W_wV))$ is a Hodge structure of weight $w$ polarized by $\langle \cdot,\cdot\rangle_w$. 

\quad{\rm (ii-2)} 
$NF^p(V)\subset F^{p-1}(V)$ for all $N\in \sig$ and $p\in \Z$. 

\quad{\rm (ii-3)} For every $N\in \sig$, the relative monodromy filtration of $N:V_\R\to V_\R$ with respect to $W$ exists.

{\rm (iii)} For some faithful $V\in \Rep(G)$, the above conditions {\rm (ii-1)}--{\rm (ii-3)} are satisfied. 

\end{sbprop}

\begin{pf} Assume (i). We prove (ii). 
Consider the continuous  map $\C^n \to \Dc\;;\; z\mapsto \exp(\sum_{j=1}^n z_jN_j)F$. For some $c\in \R$, the image of 
$S:=\{z\in \C^n\;|\; \text{Im}(z_j)\geq c\}\subset \C^n$ under this map is contained in $D$. Hence it induces a continuous map $S\to D$.
Since $D$ is a disjoint union of $G'(\R)G_u(\C)$-orbits which are open and closed and since $S$ is connected, the image of $S$ in $D$ is contained in one $G'(\R)G_u(\C)$-orbit $D'$. 
Take $F'\in D'$. Let $V\in \Rep(G)$. By the $\R$-polarizability and Lemma \ref{pol2}, for each $w\in \Z$, there is a $G'_\R$-invariant $\R$-bilinear form $\langle\cdot,\cdot\rangle_w: \gr^W_wV_\R\times \gr^W_wV_\R\to \R$ which polarizes $F'(\gr^W_w)$. If $z\in S$,  $\exp(\sum_{j=1}^n z_jN_j)F=gF'$ for some $g\in G'(\R)$. Since  $\langle\cdot, \cdot\rangle_w$ is fixed by $g$, it  polarizes $\exp(\sum_{j=1}^n z_jN_j)F$. 

The implication (ii) $\Rightarrow$ (iii) is clear. 

By \ref{faith} and by \cite{KNU2} Part III 1.2.2.1, 1.2.2.2, we have the implication (iii) $\Rightarrow $ (i). 
\end{pf}

\begin{sbprop}\label{nilp4} Assume that $(\sig, Z)$ generates a nilpotent orbit. 
Then we have a filtration  $M(\tau)\in \frak W(G_{\bR})$ ($\ref{Wlem}$) for each face $\tau$ of $\sig$ satisfying the following conditions {\rm (i)}--{\rm (iii)}.
  If $\sig$ is rational, then all $M(\tau)$ belong to $\frak W(G)$. 

{\rm (i)}  $NM(\tau)_w\subset M(\tau)_w$ for all $N\in \sig$ and $w\in \Z$. 

{\rm (ii}) $M(0)=W$. 

{\rm (iii)} If $\tau$ and $\tau'$ are faces of $\sig$ and if $N\in \sig$, and if $\tau'$ is the smallest face of $\sig$ containing $\tau$ and $N$, then $M(\tau')$ is the relative monodromy filtration of $N$ with respect to $M(\tau)$. 

\end{sbprop}

\begin{pf} This follows from  Kashiwara \cite{Kas}  4.4.1 and 5.2.5.
\end{pf}

\begin{sbpara}
  A {\it nilpotent orbit} (resp.\ {\it nilpotent $i$-orbit})
is a pair %
$(\sig, Z)$ of a nilpotent cone $\sig$ and an $\exp(\sig_\C)$ (resp.\ $\exp(i\cdot \sig_\R)$)-orbit in $\Dc$
satisfying that, for some $F\in Z$, $(\sig, F)$ generates a nilpotent orbit in the sense of \ref{nilp2}. Here $\sig_\C$ (resp.\ $\sig_\R$) denotes the $\C$ (resp.\ $\R$)-linear subspace spanned by $\sig$ in $\Lie(G'_\C)$ (resp.\ $\Lie(G'_\R)$). 

\end{sbpara}

\begin{sbpara}\label{fan}

  A {\it weak fan $\Sig$ in $\Lie(G')$} is a nonempty set of sharp rational nilpotent cones in $\Lie(G'_\R)$ 
satisfying the following conditions (i) and (ii).

(i) If $\sig\in \Sig$ and if $\sig'$ is a face of $\sig$, then $\sig'\in \Sig$.

(ii) Let $\sig, \sig' \in \Sig$, and assume that $\sig$ and $\sig'$ have a common interior 
point and that there is an $F\in \Dc$ such that both $(\sig, F)$ and $(\sig',F)$ generate nilpotent 
orbits. Then $\sig=\sig'$. 

\end{sbpara}

\begin{sbpara}\label{DSigdef}
  Let $D_{\Sig}$ be the set of all nilpotent orbits $(\sig, Z)$ such that $\sig\in \Sig$.
  Then $D$ is naturally embedded in $D_{\Sig}$ via $F \mapsto (\{0\}, \{F\})$. 

  Let $D^{\sharp}_{\Sig}$ be the set of all nilpotent $i$-orbits $(\sig, Z)$ such that $\sig\in \Sig$.
  Then $D$ is also naturally embedded in $D^{\sharp}_{\Sig}$
via $F \mapsto (\{0\}, \{F\})$. 

We have a canonical map $D^{\sharp}_{\Sig}\to D_{\Sig}\;;\; (\sig, Z) \mapsto (\sig, \exp(\sig_\C)Z)$. 

  For a rational nilpotent cone $\sig$, we define 
$D_{\sig}:=D_{\{\text{face of $\sig$}\}}, 
D^{\sharp}_{\sig}:=D^{\sharp}_{\{\text{face of $\sig$}\}}.$
\end{sbpara}

\begin{sbpara}\label{scomp}
  Let $\Gamma$ be a subgroup of $G(\Q)$ satisfying $(SA)$ (\ref{Gamma}). 
  
We say that $\Sigma$ and $\Gamma$ are {\it compatible} if $\Sig$ is stable under 
the adjoint action of $\Gamma$.   If this is the case, $\Gamma$ naturally acts on $D_{\Sig}$. 

We say that $\Sig$ and $\Gamma$ are {\it strongly compatible} if they are compatible and if every $\sig \in \Sig$ is generated by elements whose $\exp$ in $G(\R)$ belong to $\Gamma$. 
 \end{sbpara}

\subsection{$E_{\sig}$ and the spaces of nilpotent orbits}
\label{ss:Esig}

For $\Sig$ and $\Gamma$ which are strongly compatible, we endow $\Gamma\bs D_{\Sig}$ with a structure of a locally ringed space over $\C$ and with a log structure. We  endow $D^{\sharp}_{\Sig}$ with a topology.

\begin{sbpara}\label{toric}
 Let $\Sig$ and $\Gamma$ be as in \ref{fan} and in \ref{scomp}. 
  Assume that they are strongly compatible (\ref{scomp}). 
  Let $\sig \in \Sig$. 
  
   Let $\Gamma(\sig):= \Gamma \cap \exp(\sig)$ in $G(\R)$. Then $\Gamma(\sig)$ is an fs monoid and $\Gamma(\sig)^{\gp}= \Gamma \cap \exp(\sig_\R)$ is a finitely generated free abelian group. Let 
   $$\torus_{\sig}=\torus_{\sig,\Gamma}=\C^\times \otimes \Gamma(\sig)^{\gp}.$$
 
Let  $P(\sig)=\Hom(\Gamma(\sig),\bN)$. Let
$$\toric_{\sig}=\toric_{\sig,\Gamma}=\Hom(P(\sig),\bC^{\mult})=\Spec(\C[P(\sig)])^{\an}.$$
\noindent
 Here $\C^{\mult}=\C$ regarded as a multiplicative monoid.  The standard log structure of the toric variety $\Spec(\C[P(\sig)])$ induces the log structure of the analytic toric space $\toric_{\sig}$.

We regard $\torus_\sig$ as an open set of $\toric_\sig$  via the embedding
$$
\torus_\sig=\Hom(P(\sig)^{\gp}, \bC^\times)\subset \toric_\sig.
$$ 
We have a natural action of $\torus_\sig$ on $\toric_\sig$. 
We have an exact sequence 
$$
0\to \G(\sig)^{\gp}\overset{\log}\longrightarrow  \sig_\bC\overset{\bold e}\longrightarrow  \torus_\sig\to 0, 
$$ 
where 
$$
{\bold e}(z \otimes \log(\g))=e^{2\pi iz}\otimes \g\quad \text{for} \;\;
z\in \bC, \g\in \G(\sig)^{\gp}.
$$

For a face $\tau$ of $\sig$, 
let $0_\tau\in \toric_\sig$ be  the homomorphism $P(\sig)=\Hom(\Gamma(\sig), \bN) \to \bC^{\mult}$ which sends $h \in \Hom(\Gamma(\sig), \bN)$ to 
$1$ if $h(\Gamma(\tau))=0$ and to $0$ otherwise.

Each element $q$ of $\toric_\sig$ is written in the form $q={\bold e}(a)\cdot 0_\tau$ for $a\in \sig_\bC$ and for a face $\tau$ of $\sig$. 
The face $\tau$ of $\sig$ is determined by $q$ and called the {\it face associated to $q$}, and $a$ modulo $\tau_\bC+ \log(\G(\sig)^{\gp})$ is determined by $q$. The stalk of $M/\cO^\times$ of $\toric_{\sig}$ at $q$ is identified with $\Hom(\Gamma(\tau), \N)$. 
\medskip

\end{sbpara}

\begin{sbpara}\label{|toric|} Let the notation be as in \ref{toric}. Define  $$|\toric|_{\sig}:=\Hom(P(\sig), \R_{\geq 0}^{\mult})\supset |\torus|_{\sig}:= \Hom(P(\sig)^{\gp}, \R_{>0})= \R_{>0}\otimes \Gamma(\sig)^{\gp}.$$
Here $\R_{\geq 0}^{\mult}$ denotes the multiplicative monoid $\R_{\geq 0}$. Thus we have $|\toric|_{\sig}\subset \toric_{\sig}$ and $|\torus|_{\sig}\subset \torus_{\sig}$.

We have projections $q\mapsto |q|\;;\; \toric_{\sig}\to |\toric|_{\sig}$ and $\torus_{\sig}\to |\torus|_{\sig}$ induced by taking the absolute value $|\cdot|\;:\; \C^{\mult} \to \R_{\geq 0}^{\mult}$.

\end{sbpara}

\begin{sbpara}\label{add} We give additive  presentations of $|\torus|_{\sig}$ and $|\toric|_{\sig}$. 

We have an isomorphism of topological groups 
$$\sig_\R\simeq |\torus|_{\sig}\;;\; b \mapsto {\bf e}(ib).$$ We will often identify $|\torus|_{\sig}$ with $\sig_{\bR}$ via this isomorphism. 
This identification is extended to an identification of  
$|\toric|_{\sig}$ and the set of equivalence classes of pairs $(\tau, b)$, where $b\in \sig_\R$ and  $\tau$ is a face of $\sig$.
  Here $(\tau, b)$ and $(\tau', b')$ are equivalent if and only if $\tau'=\tau$ and $b'\equiv b \bmod \tau_\R$. We identify ${\bf e}(ib)0_{\tau}\in |\toric|_\sig$ with the class of $(\tau, b)$. 
The topology of $|\toric|_{\sig}$ is understood as follows. Let $x=\text{class}\,(\tau, b)\in |\toric|_{\sig}$. 
Take a finite set $(N_j)_j$ of generators of the cone $\tau$. 
Then the following sets $V(U, c)$ form a base of neighborhoods of  $x$.   
  Let $U$ be a neighborhood of $b$ in $\sig_\R$ and let $c\in \R_{>0}$. 
  Then, $V(U, c)$ is the set of $\text{class}\,(\tau', b')$, where 
 $\tau'$ is a face of $\tau$ and  $b' =b''+\sum_j y_j N_j$ for some $b''\in U$ and some real numbers  $y_j\geq c$.

\end{sbpara}

\begin{sbpara}\label{ind1}  Note that $\toric_{\sig}$ depends on the choice of $\Gamma$ (actually it depends on $\Gamma(\sig)$) though the notation $\toric_{\sig}$ does not tell this dependence. 

However, the topological space $|\toric|_{\sig}$ does not depend on $\Gamma$ as is seen in \ref{add}. 

\end{sbpara}

\begin{sbpara}\label{Esig} 

Let  $$\Ec_{\sig}:=\toric_{\sig} \times \Dc, \quad \Ec^{\sharp}_{\sig}:=|\toric|_{\sig} \times \Dc.$$
  
Let $E_{\sig}$ (resp.\ $\tilde E_{\sig}$)  be the set of all $(q, F)\in \Ec_{\sig}$ satisfying the following condition (i) (resp.\ (i)${}'$). 
  Write  $q={\bold e}(a)0_{\tau}$, where $a \in \sig_{\bC}$ and $\tau$ is a face of $\sig$.

(i) The pair $(\tau, \exp(a)F)$ generates a nilpotent orbit.

(i)${}'$ We have  $NF^p\subset F^{p-1}$ for $N\in \tau$ and $p\in \Z$.

We have $E_{\sig}\subset \tilde E_{\sig} \subset \Ec_{\sig}$.

Endow the complex analytic space $\Ec_{\sig}$ with the pullback log structure from $\toric_{\sig}$. 
Endow $E_{\sig}$ and $\tilde E_{\sig}$ with the strong topologies (\cite{KU2},  3.1.1) in $\Ec_{\sig}$ as reviewed below, 
with the inverse images of the sheaf of holomorphic functions on $\Ec_{\sig}$, and with the inverse images of the log structure of $\Ec_{\sig}$. 

Recall that for a complex analytic space $X$ and a subset $S$ of $X$, the {\it strong topology of $S$ in $X$} is the strongest topology on $S$ such that for every complex analytic space $Y$ and for every morphism $\lam:Y\to X$ of complex analytic spaces such that $\lam(Y)\subset S$, the map $\la: Y\to S$ 
is continuous. It is stronger than the topology on $S$ as a subspace of $X$.

 Let
$$E^{\sharp}_{\sig}:= E_{\sig}\cap \Ec^{\sharp}_{\sig},  \quad \tilde E^{\sharp}_{\sig}:= \tilde E_{\sig}\cap \Ec^{\sharp}_{\sig},$$
endowed with the topologies as the subspaces of $E_{\sig}$ and $\tilde E_{\sig}$, respectively. 
Using the additive presentation of $|\toric|_{\sig}$ in \ref{add}, we can identify
$E^{\sharp}_{\sig}$ (resp.\  $\tilde E^{\sharp}_{\sig}$) with the set of $(\text{class}\,(\tau, b), F)\in |\toric|_{\sig} \times \Dc=\Ec^{\sharp}_{\sig}$ such that $(\tau, \exp(ib)F)$ generates a nilpotent orbit (resp.\  
such that $NF^p\subset F^{p-1}$ for all $N\in \tau$ and $p\in \Z$). 
\end{sbpara}

\begin{sbpara}\label{indepGam1}
  We show that as a topological space, $E^{\sharp}_{\sig}$ does not depend on  $\Gamma$.  Write $E_{\sig}$, $E^{\sharp}_{\sig}$ and $P(\sig)$ (\ref{toric}) for $\Gamma$ as $E_{\Gamma(\sig)}$, $E^{\sharp}_{\Gamma(\sig)}$ and $P(\Gamma(\sig))$, respectively, to express the dependence on $\Gamma(\sig)$. 

Assume that $\Gamma_1$ and $\Gamma_2$ are strongly compatible with $\Sig$. Then $\Gamma_1\cap \Gamma_2$ is strongly compatible with $\Sig$. We show that $E_{(\Gamma_1\cap \Gamma_2)(\sig)}\to E_{\Gamma_j(\sig)}$ are homeomorphisms for $j=1,2$. Replacing $\Gamma_1\cap \Gamma_2$ by $\Gamma_1$, we may assume $\Gamma_1\subset \Gamma_2$. Write $\check{E}_{\Gamma_j(\sig)}$ as $Z_j$, $E_{\Gamma_j(\sig)}$ as $Y_j$, and $E^{\sharp}_{\Gamma_j(\sig)}$ as $Y^{\sharp}_j$ ($j=1,2$). 
 Since $P(\Gamma_2(\sig))\subset P(\Gamma_1(\sig))$, and $\C[P(\Gamma_1(\sig))]$ is a finitely generated $\C[P(\Gamma_2(\sig)])$-module,  the map $\toric_{\Gamma_1(\sig)}\to \toric_{\Gamma_2(\sig)}$ is proper and surjective  and hence 
the map $Z_1\to Z_2$ is proper and surjective. 

We prove that for the map  $Y_1\to Y_2$, the topology of the latter is the image of the topology of the latter. Let $U$ be a subset of $Y_2$ and assume that the inverse image $U'$  of $U$ in $Y_1$ is open. We prove that $U$ is open. Let $S$ be an analytic space over $\C$ and assume that we have a morphism $S\to Z_2$ whose image is contained in $Y_2$. Our task is to prove that the inverse image of $U$ in $S$ is open. Let $S'$ be the fiber product of $S\to Z_2\leftarrow Z_1$ in the category of analytic spaces over $\C$. Since the set $Y_1$ is the fiber product 
of $Y_2 \to Z_2 \leftarrow Z_1$, the image of $S'\to Z_1$ is contained in $Y_1$. Hence the inverse image of $U'$ in $S'$, which is the inverse image of $U$ in $S'$, is open. Since $S'\to S$ is proper surjective, this proves that the inverse image of $U$ in $S$ is open. 

Since $Y_j^{\sharp}$ is closed in $Y_j$, this tells that for the bijection $Y_1^{\sharp}\to Y^{\sharp}_2$, the topology of the latter is the image of the topology of the former. Hence the last map is a homeomorphism.

\end{sbpara}

\begin{sbpara}\label{phisig}

 For $\sig\in \Sig$, consider the  map 
$$\varphi_{\sig}: E_{\sig}\to \Gamma \bs D_{\Sig}\;;\; (q,F) \mapsto \text{class}\,(\tau, Z)$$ with $Z=\exp(\tau_\C)\exp(a)F$, where $\tau$ and $a\in \sig_\C$ are such that $q= {\bold e}(a) 0_{\tau}$.
We endow $\Gamma\bs D_{\Sig}$ with a structure of a locally ringed space over $\C$ and with a log structure as follows. The topology of $\Gamma \bs D_{\Sig}$ is the strongest topology for which the maps 
$\varphi_{\sig}$ are continuous for all $\sig\in \Sig$. The structure sheaf of $\Gamma \bs D_{\Sig}$ consists of functions whose pullbacks to $E_{\sig}$ belong to the structure sheaf of $E_{\sig}$ for all $\sig\in \Sig$. The log structure of $\Gamma \bs D_{\Sig}$ is  the subsheaf of the structure sheaf of $\Gamma \bs D_{\Sig}$ consisting of functions whose pullbacks to $E_{\sig}$ belong to the log structure of $E_{\sig}$ for all $\sig\in \Sig$. 

\end{sbpara}

\begin{sbpara}

In this topology of $\Gamma \bs D_{\Sig}$, if $(\sig, Z)\in D_\Sig$ and if $F\in Z$ and $N_1, \dots, N_n\in \sig$ generate the cone $\sig$, $\text{class}\,(\sig, Z)\in \Gamma \bs D_{\Sig}$ is the limit of $\text{class}\,(\exp(\sum_{j=1}^n z_jN_j)F)\in \Gamma \bs D$ where $z_j\in \C$ and $\text{Im}(z_j)$ tends to $\infty$ for $1\leq j\leq n$.

\end{sbpara}

\begin{sbpara} Here we do not need $\Gamma$.

For $\sig\in \Sig$, consider the surjective map $$\varphi^{\sharp}_{\sig}: E^{\sharp}_{\sig}\to D^{\sharp}_{\sig}\;;\; (\text{class}\,(\tau, b), F) \mapsto \text{class}\,(\tau, Z)$$ with $Z=\exp(i\tau_\R)\exp(ib)F$. 

We define the topology of $D^{\sharp}_{\Sig}$ as the strongest topology such that the composite $E^{\sharp}_{\sig}\to D^{\sharp}_{\sig}\to D^{\sharp}_{\Sig}$ also denoted by $\varphi^{\sharp}_{\sig}$ are continuous for all $\sig\in \Sig$. 

\end{sbpara}

\begin{sbpara}

In this topology, if $(\sig, Z)\in D^{\sharp}_\Sig$ and if $F\in Z$ and $N_1, \dots, N_n\in \sig$ generate the cone $\sig$, $(\sig, Z)$ is the limit of $\exp(\sum_{j=1}^n iy_jN_j)F\in D$ where $y_j\in \R_{>0}$ and $y_j\to \infty$  for $1\leq j\leq n$.

\end{sbpara}

\begin{sbpara} Assume that $(\Sig, \Gamma)$ is strongly compatible. 

 Since $$\begin{matrix} E^{\sharp}_{\sig}& \overset{\varphi^{\sharp}_{\sig}}\longrightarrow &  D_{\Sig}^{\sharp}\\ \downarrow &&\downarrow \\
E_{\sig}&\overset{\varphi_{\sig}}\longrightarrow &  \Gamma\bs D_{\Sig}\end{matrix}$$ is commutative, we have that the map $D_{\Sig}^{\sharp}\to \Gamma\bs D_{\Sig}$ is continuous.

\end{sbpara}

\begin{sbpara}
\label{logmfd} For an fs log analytic space $X$, by a {\it strong subspace} of $X$, we mean a subset of $X$ endowed with the strong topology in $X$ (\ref{Esig}), with the inverse image of the sheaf of holomorphic functions on $X$, and with the inverse image of the log structure of $X$. 

For example, $E_{\sig}$ is a strong subspace of $\check{E}_{\sig}$.

  Let $\cB(\log)$ be the category of locally ringed spaces endowed with a log structure  which are locally isomorphic to a strong subspace of an fs log analytic space. 
  See \cite{KU2} 3.2.4 and \cite{KNU2} Part III 1.1.4. 

  Both $E_{\sig}$ and $\tilde E_{\sig}$ are objects of $\cB(\log)$. 

  An object of $\cB(\log)$ is a {\it log manifold} if it is locally isomorphic to an open set of a strong subspace $S$ of a log smooth fs log analytic space $X$ satisfying the following condition: There is a finite family of log differential forms $(\omega_j)_j$ on $X$ such that $S=\{x\in X\;|\; \omega_j(x)=0$ for all $j%
\}$. Here $\omega_j(x)$ denotes the pullback of $\omega_j$ to the log point $x$ (it is not the germ of $\omega_j$ at $x$).
  See \cite{KU2} 3.5.7 and \cite{KNU2} Part III 1.1.5. 

  Later we will show that $E_{\sig}$, $\tilde E_{\sig}$, and $\Gamma \bs D_{\Sig}$ for a strongly compatible $(\Sig, \Gamma)$ with $\Gamma$ being a neat semi-arithmetic subgroup of $G'(\Q)$ (\ref{Gamma}), are log manifolds. 
  See Theorem \ref{Elm}  and Theorem \ref{t:property} (4). We will also show that for such $\Gamma$, the quotient space $\Gamma \bs D^{\sharp}_{\Sig}$ is identified with the space $(\Gamma\bs D_{\Sig})^{\log}$ (Theorem \ref{t:property} (5)). 
\end{sbpara}

\subsection{The space of ratios}
\label{ss:ratio}

We review the space of ratios $S_{[:]}$, and discuss $D_{\Sig, [:]}$ and $D^{\sharp}_{\Sig, [:]}$ which are modified versions of $D_{\Sig}$ and $D^{\sharp}_{\Sig}$, respectively.

\begin{sbpara}\label{rat1}  
  This is a review of \cite{KNU2} Part IV 4.2. 
  Let
$S$ be a locally ringed space over $E=\R$ or $\C$ with an fs log structure $M_S$ satisfying the following conditions (i) and (ii).

(i) For every $s\in S$, the natural homomorphism 
$E\to \cO_{S,s}/m_s$ is an isomorphism. 
Here, we denote by $m_s$ the maximal ideal of $\cO_{S,s}$.

(ii) For every open set $U$ of $S$ and for every $f\in \cO(U)$, the map $U\to E\;;\; s \mapsto f(s)$ is continuous. Here $f(s)$ is the image of $f$ in $\cO_{S,s}/m_s=E$. 

Then we have a topological space $S_{[:]}$ over $S$, called  the {\it space of ratios}, defined as follows.

 \end{sbpara}
 
 \begin{sbpara}\label{rat2a}

 For a sharp fs monoid $\cS$, let $R(\cS)$ be the set of all maps $r:(\cS\times \cS)\smallsetminus \{(1,1)\}\to [0,\infty]$ satisfying the following conditions (i)--(iii).  
  
(i) $r(g,f)=r(f,g)^{-1}$.

(ii) $r(f, g)r(g,h)=r(f,h)$ if $\{r(f,g), r(g,h)\}\neq \{0,\infty\}$. 

(iii) $r(fg, h)=r(f, h)+r(g,h)$. 

\end{sbpara}

\begin{sbpara}\label{rat2} Let $\sig$ be the cone $\Hom(\cS, \R_{\geq 0}^{\add})$, where $\R_{\geq 0}^{\add}$ denotes the additive monoid $\R_{\geq 0}$. Then 
$R(\cS)$ is identified with the set of equivalence classes of 
$(\sig_j, N_j)_{1\leq j\leq n}$, where $n\geq 0$,  $\sig_j$ are faces of $\sig$ 
such that $\{0\}=\sig_0\subsetneq \sig_1\subsetneq \dots \subsetneq \sig_n=\sig$, and $N_j$ is an element of the interior of $\sig_j$. The equivalence relation is that $(\sig_j, N_j)_{1\leq j\leq n} \sim (\sig'_j, N'_j)_{1\leq j\leq n'}$ if and only if $n=n'$, $\sig_j'=\sig_j$ for all $j$,  and for each $j$, there is a $c_j>0$  such that $N_j'\equiv c_jN_j\bmod \sig_{j-1,\R}$.

Such a $\text{class}\,((\sig_j, N_j)_j)$ is identified with $r\in R(\cS)$ defined as follows. For $(f,g)\in \cS\times \cS\smallsetminus \{(1,1)\}$, take the biggest $j$ such that $\sig_{j-1}$ kills $f$ and $g$. Then $r(f,g)=N_j(f)/N_j(g)$. 
\end{sbpara}

\begin{sbrem}
  There is an error in a related part 4.1.6 in \cite{KNU2} Part IV. 
  In the definition of the map $R'(\cS) \to R(\cS)$, \lq\lq$\infty$'' in (2) and \lq\lq $0$'' in (3) should be interchanged. 
\end{sbrem}

\begin{sbpara}
\label{ratSs}
 Let $S$ be as in \ref{rat1}.
We define the set 
$S_{[:]}$ as 
the set of $(s,r)$ with $s\in S$ and  $r\in R((M_S/\cO_S^\times)_s)$. 
\end{sbpara}

\begin{sbpara}
\label{ratSt}
 We define the topology of $S_{[:]}$ as the weakest topology for which the following conditions (1) and (2) are satisfied.

(1) The map $S_{[:]}\to S\;;\;(s,r)\mapsto s$ is continuous. 

(2) Let $U$ be an open set of $S$, let $f,g\in M_S(U)$, and assume that $|f(s)|<1$ and $|g(s)|<1$ for all $s\in U$.  Here $f(s)$ is the value at $s$ of the image of $f$ in $\cO_U$ and $g(s)$ is defined similarly. Let $\tilde U$ be the inverse image of $U$ in $S_{[:]}$. 
  Then the following map  $r_{f,g}: \tilde U \to [0, \infty]$ is continuous.

Let $s\in U$.  If $f,g\in \cO_{S,s}^\times$, then $r_{f,g}(s,r):= \log(|f(s)|)/\log(|g(s)|)$. Otherwise, $r_{f,g}(s,r)= r({\bar f}(s),{\bar g}(s))$, where ${\bar f}(s)$ is the image of $f$ in $(M_S/\cO_S^\times)_s$ and ${\bar g}(s)$ is defined similarly. 
\end{sbpara}

\begin{sbpara}\label{rat4} 
Assume that we are given a chart  $\cS\to M_S$.
Let $\Phi$ be a set of faces of $\cS$ which is totally ordered for the inclusion relation and which contains $\cS$. 
Let
$S_{[:]}(\Phi)$ be the subset of $S_{[:]}$ consisting of all $(s, r)$ such that if $\text{class}((\sig_j, N_j)_{1\leq j\leq n})$ corresponds to $r$ as in \ref{rat2}, the annihilator of $\sig_j$ in $\cS$, which is a face of $\cS$, belongs to $\Phi$ for every $1\leq j\leq n$. 
Here the annihilator of $\sig_j$ in $\cS$ means the subset of $\cS$ consisting of all elements which are sent to $0$ by $\cS\to (M_S/\cO_S^\times)_s\overset{h}\to \R_{\geq 0}$ for all $h\in \sig_j$. 

Then $S_{[:]}(\Phi)$ for all above $\Phi$ forms an open covering of $S_{[:]}$ (ibid.\ Part IV 4.2.11).

For a geometric meaning of $S_{[:]}(\Phi)$, see ibid.\ Part IV 4.2.11--4.2.19 in the general case, and ibid.\ Part IV 4.2.20--4.2.22 in the case $S=|\Delta|^n$ with $|\Delta|:=\{t\in\R\,|\,0\le t<1\}$.
\end{sbpara}

\begin{sbpara} Let $|\toric|_{\sig,[:]}\subset \toric_{\sig,[:]}$ be the inverse image of $|\toric|_{\sig}$ under the map $\toric_{\sig,[:]}\to \toric_{\sig}$. We have a projection  $\toric_{\sig,[:]}\to |\toric|_{\sig,[:]}$.
  This is a unique continuous map which is compatible with the projection $\toric_{\sig} \to |\toric|_{\sig}$. 
\end{sbpara}
 
\begin{sbpara}\label{add2} We give an additive presentation of $|\toric|_{\sig,[:]}$. 

We can identify
$|\toric|_{\sig,[:]}$ with the set of equivalence classes of $((\sig_j, N_j)_{1\leq j\leq n},b)$, where 
$n\geq 0$,  $\sig_j$ are faces of $\sigma$ such that $\{0\}\subsetneq\sig_1\subsetneq \dots \subsetneq \sig_n$ (here  $\sig_n$ need not coincide with $\sig$),
 $N_j$ is an element of the interior of $\sig_j$, and $b\in 
\sig_\R$. The equivalence relation is that $((\sig_j, N_j)_{1\leq j\leq n}, b)$ and $((\sig'_j, N'_j)_{1\leq j\leq n'}, b')$ are equivalent if and only if 
$n'=n$, $\sig'_j=\sig_j$, $N_j'\equiv c_jN_j \bmod \sig_{j-1,\R}$  for some $c_j\in \R_{>0}$ ($1\leq j\leq n$), and $b'\equiv b \bmod \sig_{n,\R}$.
  Here $\sig_0$ denotes $\{0\}$. 

The projection $|\toric|_{\sig, [:]}\to |\toric|_{\sig}$ is understood as 
 $\text{class}\,((\sig_j,N_j), b)\mapsto \text{class}\,(\sig_n, b)$.

\end{sbpara}

\begin{sbpara}\label{excel1} 
Let $x=\text{class}\,((\sig_j, N_j)_{1\leq j \leq n}, b)\in |\toric|_{\sig, [:]}$. 
By 
a {\it good base} for $x$, we mean a family of elements $(N_s)_{s\in S}$ of $\sig_n$  satisfying the following conditions (i) and (ii).

(i)   The index set $S$ is the disjoint union of some subsets $S_j$ ($1\leq j\leq n$), and the following holds for each $j$ ($1\leq j\leq n$).  
  Let $S_{\leq j}:=\bigsqcup_{k\leq j} S_k\subset S$. Then $N_s\in \sig_j$ if $s\in S_{\leq j}$, and $(N_s)_{s\in S_{\leq j}}$ is a base of the $\R$-vector space $\sig_{j, \R}$.

(ii) There are $a_s\in \R_{>0}$ ($s\in S_j$) such that $N_j \equiv  \sum_{s\in S_j} a_sN_s\bmod \sig_{j-1, \R}$ for $1\leq j\leq n$.

  This is a [:]-version of 
  good base discussed in \cite{KU2} 6.3.
\end{sbpara} 

\begin{sbprop}\label{excel2} 
A good base for $x$ exists. 

\end{sbprop}

\begin{pf}

For each $j$, 
take a simplicial subcone of the cone $(\sig_j+\sig_{j-1,\bR})/\sig_{j-1,\bR}$ 
consisting of the interior points except the origin and including the class of $N_j$. 
  Take a base $(\overline {N}_s)_{s \in S_{j}}$ of this cone and lift $\overline {N}_s$ ($s \in S_{j}$) to an element $N_s$ of $\sig_j$ for each $j$. 
  Then these $N_s$'s form 
  a good base.
\end{pf}

\begin{sbpara}\label{add3}  We describe the topology of $|\toric|_{\sig, [:]}$ by using the additive presentation \ref{add2}.

Let $n\geq 0$, let $\{0\}\subsetneq \sig_1\subsetneq \dots \subsetneq \sig_n$ be faces of $\sig$, and let  $(N_s)_{s\in S}$ be a finite family of elements of $\sig_n$ satisfying the condition (i) in \ref{excel1}.

Fix an $\R$-subspace $B$ of $\sig_\R$ such that $\sig_\R=\sig_{n, \R}\oplus B$. 

Let $U$ be the subset of $|\toric|_{\sig,[:]}$ consisting of classes of  $((\sig'_k, N _k)_{1\leq k \leq n'}, b)$ 
satisfying the following conditions. There is an injective increasing map $\theta: \{1, \dots, n'\}\to \{1, \dots, n\}$  such that $\sig'_k=\sig_{\theta(k)}$  
for $1\leq k \leq n'$. For $1\leq k\leq n'$, if we write $N_k \equiv \sum_s y_s N_s\bmod \sig'_{k-1, \R}$, where $s$ ranges over all elements of $S$ 
belonging to $S_j$ for some $j$ such that $\theta(k-1)<j\leq \theta(k)$, then $y_s\in \R_{>0}$. (In the case $k=1$, $\theta(k-1)$ means $0$.) If we write $b\equiv b' +\sum_s y_sN_s \bmod \sig'_{n',\R}$, 
where $b'\in B$ and  $s$ ranges over all elements of $S$ belonging to $S_j$ for some $j$ such that $\theta(n')<j$, then $y_s\in \R_{>0}$. 

Then $U$ is an open set of $|\toric|_{\sig, [:]}$. 

For each $j$, choose an element $c_j$ of $S_j$. 
Let $T$ be the complement of $\{c_j\;|\; 1\leq j\leq n\}$ in $S$. 
Then we have a homeomorphism 
$$U\simeq \R^n_{\geq 0} \times \R_{>0}^T \times B$$
defined as follows. Let $x' =\text{class}\,((\sig'_k, N_k)_{1\leq k \leq n'}, b)\in U$, and let $\theta$, $y_s$ ($s\in S$), $b'\in B$ be as above. Then the image of $x'$ in $\R^n_{\geq 0}\times \R_{>0}^T \times B$ is $((t_j)_{1\leq j\leq n}, (y_s/y_{c_j})_{s\in T}, b')$, where $t_j$  are as follows. 

If $j$ is in the image of $\theta$, then $t_j=0$.

If $j$ is not in the image of $\theta$ and $j<\theta(n')$, then $t_j= y_{c_{j+1}}/y_{c_j}$.

If $\theta(n')<j$, then $t_j= 1/y_{c_j}$.

Note that if $(N_s)_{s\in S}$ is 
a good base for $x=\text{class}\,((\sig_j, N_j)_{1\leq j \leq n}, b)\in |\toric|_{\sig,[:]}$, $x$ belongs to the above $U$ and hence $U$ is an open neighborhood of $x$.

\end{sbpara}

\begin{sbpara}\label{Esig:} 
  Let the notation be as in \ref{toric}. 
  Let
$$\Ec^{\sharp}_{\sig,[:]}:=|\toric|_{\sig,[:]} \times \Dc \subset \Ec_{\sig,[:]}=\toric_{\sig,[:]}\times \Dc.$$
Let
$$E^{\sharp}_{\sig,[:]}:=E_{\sig,[:]}\cap \Ec^{\sharp}_{\sig,[:]}, \quad 
\tilde E^{\sharp}_{\sig,[:]}:=\tilde E_{\sig,[:]}\cap \Ec^{\sharp}_{\sig,[:]}.$$
  We endow $E^{\sharp}_{\sig,[:]}$ and $\tilde E^{\sharp}_{\sig,[:]}$ with the topologies as the subspaces of $E_{\sig,[:]}$ and $\tilde E_{\sig,[:]}$, respectively. 
  These coincide with the topologies defined via the projections $E_{\sig,[:]}\to E^{\sharp}_{\sig,[:]}$ and $\tilde E_{\sig,[:]}\to \tilde E^{\sharp}_{\sig,[:]}$, respectively. 
\end{sbpara}

\begin{sbpara}\label{indepGam2}  We show that as topological spaces, $|\toric|_{\sig, [:]}$ and 
$E^{\sharp}_{\sig, [:]}$ are independent of the choice of $\Gamma$. 

We describe the proof for $E^{\sharp}_{\sig, [:]}$. The proof for $|\toric|_{\sig, [:]}$ is similar. 
In fact, assume that $\Gamma_j$ ($j=1, 2$) are strongly compatible with $\Sig$ and assume $\Gamma_1\subset \Gamma_2$. Since $E^{\sharp}_{\Gamma_j(\sig),[:]}$ is identifies with $(E^{\sharp}_{\Gamma_j(\sig)})_{[:]}$ (the $[:]$-space of $E^{\sharp}_{\Gamma_j(\sig)}$ which is endowed with the inverse image of the log structure of $E_{\Gamma_j(\sig)}$, $E^{\sharp}_{\Gamma_j(\sig), [:]}$ is proper over $E^{\sharp}_{\Gamma_j(\sig)}$. Since $E^{\sharp}_{\Gamma_1(\sig)}\to E^{\sharp}_{\Gamma_2(\sig)}$ is a homeomorphism (\ref{indepGam1}), the map 
$E^{\sharp}_{\Gamma_1(\sig),[:]}\to E^{\sharp}_{\Gamma_2(\sig),[:]}$ is  proper continuous and it is bijective. Hence it is a homeomorphism. 

\end{sbpara}

\begin{sbpara}\label{rat3}

 Let $D_{\Sig, [:]}$ (resp.\ $D^{\sharp}_{\Sig,[:]}$) be the set of $(\sig, Z, \text{class}\,((\sig_j, N_j)_{1\leq j\leq n}))$, where $\sig\in \Sig$, $(\sig,Z)\in D_{\Sig}$ (resp.\ $(\sig,Z)\in D^{\sharp}_{\Sig}$), $n\geq 0$, $\sig_j$ are faces of $\sig$ such that $\{0\}\subsetneq \sig_1\subsetneq \dots \subsetneq \sig_n=\sig$, and $N_j$ is an element of the interior of $\sig_j$. The notation class is for the equivalence relation  that
 $(\sig_j, N_j)_{1\leq j\leq n}\sim (\sig'_j, N'_j)_{1\leq j\leq n'}$ if and only if  $n'=n$, $\sig'_j=\sig_j$, and $N'_j\equiv a_jN_j\bmod \sig_{j-1,\R}$
 for some $a_j\in \R_{>0}$ ($1\leq j\leq n$). 
 
 We have canonical maps $$D_{\Sig, [:]}\to D_{\Sig}, \quad D^{\sharp}_{\Sig, [:]}\to D^{\sharp}_{\Sig}$$
 by sending the class of $(\sig, Z, \text{class}\,((\sig_j,N_j)_{1\leq j\leq n}))$ to $(\sig, Z)$. 
 \end{sbpara}

\begin{sbpara}\label{rat8} Assume that $(\Sig, \Gamma)$ strongly compatible. We define a topology of $\Gamma \bs D_{\Sig, [:]}$. 

For $\sig\in \Sig$, by \ref{rat2}, $E_{\sig}$ is identified with the set of triples $(q, \text{class}\,((\sig_j, N_j)_{1\leq j\leq n}), F)$, where $q\in \toric_{\sig}$ and writing $q$ as ${\bold e}(a)0_{\tau}$, $n\geq 0$,  $\sig_j$ are faces of $\tau$ such that 
 $\{0\}\subsetneq \sig_1\subsetneq \dots \subsetneq \sig_n=\tau$, and $N_j$ is an element of the interior of $\sig_j$. The equivalence relation is defined in the same way as in \ref{rat2}. 
 
 For $\sig\in \Sig$, consider the  map 
$$\varphi_{\sig,[:]}: E_{\sig,[:]}\to \Gamma \bs D_{\Sig,[:]}\;;\; (q,\text{class}\,(\sig_j, N_j)_{1\leq j\leq n}, F) \mapsto \text{class}\,(\tau, Z, \text{class}\,(\sig_j, N_j)_{1\leq j\leq n})$$ with $q= {\bold e}(a) 0_{\tau}$ and  $Z=\exp(\tau_\C)\exp(a)F$.

We endow $\Gamma\bs D_{\Sig,[:]}$ with the strongest topology for which the maps 
$\varphi_{\sig,[:]}$ are continuous for all $\sig\in \Sig$.

In this topology, if $(\sig, Z, \text{class}\,((\sig_j, N_j)_{1\leq j\leq n}))\in D_{\Sig,[:]}$ and if $F\in Z$, then 
 $\text{class}\,(\sig, Z, (\sig_j, N_j)_{1\leq j\leq n})\in \Gamma \bs D_{\Sig}$ is the limit of $\text{class}\,(\exp(\sum_{j=1}^n z_j\tilde N_j)F)\in \Gamma \bs D$ where $z_j\in \C$, $\text{Im}(z_j)\to \infty$ for $1\leq j\leq n$ and $\text{Im}(z_j)/\text{Im}(z_{j+1})\to \infty$ for $1\leq j\leq n-1$.

\end{sbpara}
\begin{sbpara}\label{rat9}  Here we do not need $\Gamma$.

 For $\sig\in \Sig$, we have a map $$\varphi^{\sharp}_{\sig, [:]}: E^{\sharp}_{\sig,[:]}\to D^{\sharp}_{\Sig,[:]}\;;\; (\text{class}\,((\sig_j, N_j)_{1\leq j\leq n}), b, F) \mapsto (\sig_n,\exp(i\sig_{n,\R}+ ib)F,\text{class}\,((\sig_j, N_j)_{1\leq j\leq n})).$$

 We define the topology of 
$D^{\sharp}_{\Sig,[:]}$ as the strongest topology for which the maps $\varphi^{\sharp}_{\sig,[:]}: E^{\sharp}_{\sig,[:]}\to D^{\sharp}_{\Sig,[:]}$ are continuous for all $\sig \in \Sig$.

In this topology, if $x=(\sig, Z, \text{class}\,((\sig_j, N_j)_{1\leq j\leq n}))\in D^{\sharp}_{\Sig,[:]}$ and if $F\in Z$ and $\tilde N_j$, then 
 $x$ is the limit of $\exp(\sum_{j=1}^n iy_jN_j)F\in D$ where $y_j\in \R_{>0}$, $y_j/y_{j+1}\to \infty$ for $1\leq j\leq n$ ($y_{n+1}:=1$). 

\end{sbpara}

\begin{sbpara} The maps  $\Gamma\bs D_{\Sig,[:]}\to \Gamma \bs D_{\Sig}$ and $D^{\sharp}_{\Sig,[:]}\to D^{\sharp}_{\Sig}$ are continuous.

\end{sbpara}

\subsection{Valuative spaces, II}
\label{ss:val2}

\begin{sbpara} Let $E$ and $S$ be as in \ref{rat1}. 

We recall that there are two kinds of valuative spaces $S_{\val(E)}$ and $S_{\val(|\cdot|)}$ associated to $S$ (\cite{KNU2} Part IV 3.1.3). 

In the case where $S$ is an object of $\cB_\R'(\log)$, $S_{\val}$ in Section 
\ref{ss:val1} 
is identified, as a topological space, with the topological space $S_{\val(|\cdot|)}$. On the other hand, $S_{\val}$ for an object $S$ of $\cB(\log)$ (\ref{logmfd}) defined in \cite{KU2} 3.6.18, 3.6.23 is a locally ringed space over $\C$ with log structure, and it is $S_{\val(\C)}$.

These two kinds of $S_{\val}$ played important roles in our previous works \cite{KU2} and \cite{KNU2}. 

Our rule of the notation $S_{\val}$ in \cite{KU2} and \cite{KNU2} and in this paper is that  if $S$ is an object of $\cB(\log)$, then  $S_{\val}$ means $S_{\val(\C)}$, but otherwise, $S_{\val}$ means $S_{\val(|\cdot |)}$.

We also consider, for a weak fan $\Sig$, the maps $D_{\Sig, \val}\to D_{\Sig, [:]}$ and $D^{\sharp}_{\Sig, \val}\to D^{\sharp}_{\Sig,[:]}$.
\end{sbpara}

\begin{sbpara} 

Let $E$ and $S$ be as in \ref{rat1}. 

As a set, let $S_{\val(|\cdot |)}$ (resp.\ $S_{\val(E)}$) be the set of triples $(s, V, h)$, where $s\in S$, $V$ is a submonoid of $(M_S/\cO_S^\times)_s^{\gp}$ such that $(M_S^{\gp}/\cO_S^\times)_s=V\cup V^{-1}$, $(M_S/\cO_S^\times)_s\subset V$ and $V^\times\cap (M_S/\cO_S^\times)_s=\{1\}$ (here $V^\times$ denotes the subgroup $V\cap V^{-1}$ of  $(M_S^{\gp}/\cO_S^\times)_s$), and $h$ is a homomorphism from 
$\{f\in M^{\gp}_{S,s}\;|\; f\bmod \cO_{S,s}^{\times} \in V^{\times}\}$ 
to $\R_{>0}$ (resp.\ $E^\times$) such that $h(u)=|u(s)|$ (resp.\ $h(u)=u(s)$) for all $u\in \cO_{S,s}^\times$. 

The topologies of $S_{\val(|\cdot|)}$ and $S_{\val(E)}$ are defined 
to be the weakest topology having the following properties (i) and (ii).

(i) The projection $(s,V,h)\mapsto s$ to $S$ is continuous.

(ii) Let $U$ be an open set of $S$ and let $f\in M_S(U)$. Then the subset $W=\{(s, V, h)\in S_{\val(|\cdot|)}$ (resp.\ $S_{\val(E)})
\;|\; s\in U,\ f\bmod \cO_{S,s}^\times\in V\}$ of $S_{\val(|\cdot|)}$ (resp.\ $S_{\val(E)}$) is open, and the map $W\to \R_{\geq 0}$ (resp.\ $W\to E$), which sends $(s,V, h)$ to $h(f)$ if $f\bmod \cO_{S,s}^\times\in V^\times$ 
and to $0$ otherwise, is continuous. 

For an object $S$ of $\cB(\log)$, $S_{\val(E)}$ has a structure of a locally ringed space with log structure given in \cite{KU2} 3.6.23.  

For the understandings of these val spaces using inverse limits for blowing ups along the log structure, see \cite{KNU2} Part IV Section 3.1.

\end{sbpara}

\begin{sbpara}\label{VtoR} 
For a sharp fs monoid $\cS$, we recall the set $V(\cS)$ and the canonical map $V(\cS)\to R(\cS)$ from \cite{KNU2} Part IV 4.1.7, 4.1.8.
Let $R(\cS)$ be as in \ref{rat2a} and let $V(\cS)$ be the set of all valuative submonoids  $V$ 
of $\cS^{\gp}$ such that 
$V \supset \cS$ and  $V^\times \cap \cS=\{1\}$.
We have a surjective map
$$V(\cS)\to R(\cS)$$ 
sending $V\in V(\cS)$ to the element $r_V$ of $R(\cS)$ which is the map $\cS\times \cS\smallsetminus \{(1,1)\}\to[0,\infty]$ defined by
$$r_V(f, g)=\sup\{a/b\;|\; (a,b)\in \N^2\smallsetminus  \{(0,0)\}, f^b/g^a\in V\}$$
$$= \inf\{a/b\;|\; (a, b)\in \N^2\smallsetminus \{(0,0)\}, g^a/f^b \in V\}.$$
\end{sbpara}

\begin{sbpara} We have proper surjective continuous maps
$S_{\val(E)}\to S_{\val(|\cdot |)}\to S_{[:]}\to S$. 
Here the second arrow is defined by  \ref{VtoR}. 
\end{sbpara}

\begin{sbpara}  Define 
$$\toric_{\sig,\val}:= (\toric_{\sig})_{\val(\C)}\supset  |\toric|_{\sig, \val}:=(|\toric|_{\sig})_{\val(|\cdot |)},$$
$$E_{\sig,\val}:=(E_{\sig})_{\val(\C)}\supset E^{\sharp}_{\sig, \val}:= (E^{\sharp}_{\sig})_{\val(|\cdot |)}.$$
Then $E^{\sharp}_{\sig \val}$ is identified with the inverse image of  $|\toric|_{\sig,\val}$ under $E_{\sig, \val}\to \toric_{\sig,\val}$. We have the projections 
$$\toric_{\sig,\val}\to |\toric|_{\sig,\val}, \quad E_{\sig, \val}\to E^{\sharp}_{\sig, \val}.$$

\end{sbpara}

\begin{sbpara}\label{indepGam3}

The topological spaces $|\toric|_{\sig, \val}$ and   $E^{\sharp}_{\sig, \val}$ are independent of $\Gamma$. The proofs are similar to the proof for $E^{\sharp}_{\sig, [:]}$ given in \ref{indepGam2}.

\end{sbpara}

\begin{sbpara}\label{Dval} Let $\Sig$ be a weak fan. We define sets $D_{\Sig, \val}$ and $D^{\sharp}_{\Sig,\val}$. 

For $\sig\in \Sig$, let $Q(\sig)$ be the set of all  rational linear maps $\sig_\R\to \R$, and let $\tilde P(\sig)$ be the set of all elements $h$ of $Q(\sig)$ such that $h(\sig)\subset \R_{\geq 0}$. 

Let $D_{\Sig, \val}$ (resp.\ $D^{\sharp}_{\Sig, \val}$) be the set of quadruple $(\sig, Z, V, Z')$, where $(\sig, Z)\in D_{\Sig}$ (resp.\ $D^{\sharp}_{\Sig}$), $V$ is a submonoid of $Q(\sig)$ such that $Q(\sig)= V\cup (-V)$, $\tilde P(\sig)\subset V$, and $V\cap (-V)\cap \tilde P(\sig)= \{0\}$, and if $A\subset \sig_\R\cap \Lie(G')$ denotes the annihilator of $V\cap (-V)$ in the perfect pairing $\sig_\R\cap \Lie(G')\times Q(\sig) \to \Q$ of $\Q$-vector spaces, $Z'$ is an $\exp(A_\C)$ (resp.\ $\exp(iA_\R)$)-orbit in $Z$.

The canonical maps $D_{\Sig, \val}\to D_{\Sig}$, $D^{\sharp}_{\Sig,\val}\to D^{\sharp}_{\Sig}$ are given by $(\sig, Z, V, Z') \mapsto (\sig, Z)$. 

The shapes of the definitions of $D_{\Sig, \val}$ and $D^{\sharp}_{\Sig,\val}$ seem to be slightly different from those in \cite{KU2}, \cite{KNU2} Part III, but are the same.  We hope the presentations of the present definitions are better.

\end{sbpara}

\begin{sbpara}
We have canonical maps of sets
$$D_{\Sig, \val} \to D_{\Sig, [:]}, \quad D^{\sharp}_{\Sig, \val} \to D^{\sharp}_{\Sig, [:]}$$
 sending $(\sig, Z, V, Z')$ to $(\sig, Z, \text{class}\,((\sig_j, N_j)_j))$, where $V\mapsto r_V\leftrightarrow\text{class}\,((\sig_j, N_j)_j)$ by \ref{VtoR} and \ref{rat2}.

\end{sbpara}

\begin{sbpara}
\label{stval1}
 Assume that $(\Sig, \Gamma)$ is strongly compatible. 

The structure of $\Gamma\bs D_{\Sig, \val}$ as a locally ringed space with log structure is defined by using $\varphi_{\sig, \val}: E_{\sig, \val}\to \Gamma \bs D_{\Sig, \val}$ ($\sig\in\Sig$) analogously as \ref{phisig}.

\end{sbpara}

\begin{sbpara}
\label{stval2}
 Here we do not need $\Gamma$.

The topology of $D^{\sharp}_{\Sig, \val}$ is defined by using $\varphi^{\sharp}_{\sig, \val}: E^{\sharp}_{\sig, \val}\to D^{\sharp}_{\Sig, \val}$ ($\sig\in\Sig$) analogously as \ref{rat9}. 
\end{sbpara}

\begin{sbpara}
The  continuity of $\Gamma\bs D_{\Sig, \val}\to \Gamma \bs D_{\Sig,[:]}$ and the continuity of $D^{\sharp}_{\Sig, \val}\to D^{\sharp}_{\Sig, [:]}$ are  clear because the topologies of the val spaces are   defined by using $E_{\sig,\val}$, $E^{\sharp}_{\sig, \val}$, and the topologies of  $[:]$-spaces are defined by $E_{\sig,[:]}$ and $E^{\sharp}_{\sig, [:]}$, and $E_{\sig,\val}\to E_{\sig,[:]}$ and $E^{\sharp}_{\sig, \val}\to E^{\sharp}_{\sig, [:]}$ are continuous. 

\end{sbpara}

\subsection{Nilpotent orbits and $\SL(2)$-orbits}\label{ss:CKS}

  In this section we prove the continuity of the CKS map (Theorem \ref{CKS1}), which is the most important bridge in the fundamental diagram. 
  We prove this together with that $E_{\sig}$ is open in $\tilde E_{\sig}$ (Theorem \ref{Eopen}), which implies that $E_{\sig}$ is a log manifold (Theorem \ref{Elm}). 

  Note that there are mistakes in the corresponding parts \cite{KU2} and \cite{KNU2} Part III, which are corrected in \cite{KNU2} Part IV Appendix in a rather complicated way and, unfortunately, the correction itself contains several errors.
So in this section \ref{ss:CKS}, we present the whole structure of the corrected (and improved) proof in the present context.
For more precise explanations about the mistakes and errors in the previous works, see Remark \ref{c_to_IVApp} at the end of this section.

\begin{sbpara}
\label{assSL2orb}
In order to define the map $\psi$ in Theorem \ref{CKS1} below as a map of sets, we explain how to associate an $\SL(2)$-orbit to a nilpotent orbit in \ref{assSL2orb}--\ref{assSL2orb3}.
  This construction is a generalization of \cite{KNU2} Part II 2.4.6 based on the SL(2)-orbit theorem in many variables for MHS (\cite{KNU1}).
  The pure case of the construction is essentially due to Cattani, Kaplan and Schmid (\cite{CKS}).
  
  Let $(N_1,\ldots,N_n,F)$ generate a nilpotent orbit (\ref{nilp2}). We give the associated triple $(\rho,\varphi, Z)$, where $(\rho,\varphi)$ is an $\SL(2)$-orbit in $n$ variables for $(G_{\red}, h_0)$ and $Z\subset D$ such that $(p, Z)\in D_{\SL(2)}$ for $p=\text{class}\,(\rho, \varphi)\in D_{\red,\SL(2)}$. 

  For $0 \le j \le n$, we denote by $W^{(j)}$ the functor 
$$V \mapsto (V, M(N_{1,V}+\cdots +N_{j,V}, W(V)))$$
from $\Rep(G)$ to the category of finite-dimensional $\bQ$-vector spaces endowed with an increasing filtration on its realification.
  Here $N_{j,V}$ denotes the image of $N_j$ in $\End_{\bR}(V_{\bR})$ for each $j$ and $M(N_{1,V}+\cdots +N_{j,V}, W(V))$ is the relative monodromy filtration of $N_{1,V}+\dots+N_{j, V}$ with respect to $W(V)$. In particular, $W^{(0)}=W$.
  
  The functor $(W^{(n)}, F):V\mapsto (V, W^{(n)}(V), F(V))$ is a $G$-MHS. Let $(W^{(n)}, \hat F_{(n)})$ be the $\R$-split $G$-MHS defined by the canonical splitting of $W^{(n)}$ associated to the $G$-MHS  $(W^{(n)}, F)$ (\ref{splW}). Then $(W^{(n-1)}, \exp(iN_n)\hat F_{(n)})$ is a $G$-MHS. Let $(W^{(n-1)}, \hat F_{(n-1)})$ be the associated $\R$-split $G$-MHS defined by the canonical splitting of $W^{(n-1)}$ associated to the $G$-MHS $(W^{(n-1)}, \exp(iN_n)\hat F_{(n)})$. 
  Then $(W^{(n-2)}, \exp(iN_{n-1})\hat F_{(n-1)})$ is a $G$-MHS. $\dots$. In this way, we obtain inductively an $\R$-split $G$-MHS  $(W^{(j)}, \hat F_{(j)})$ for $0\leq j\leq n$. We have $\hat F_{(j)}\in \check{D}$. 
  
  We define $\br \in D$ as follows. If all $N_j$ are $0$, let $\br
  =F$. 
  If $N_j\neq 0$ for some $j$, let $k$ be the smallest such $j$ ($1\leq j\leq n$), and let 
  $\br
  =\exp(iN_k)\hat F_{(k)}$.   
  (Note that we have used the symbol $\br_1$ in \cite{KNU2} Part II 2.4.6 instead of $\br$ here.
  Thus the notation here is not compatible with that in \cite{KNU1} Theorem 0.5 and in \cite{KNU2} Part II 2.4. 
  See the remark after \cite{KNU2} Part II Theorem 2.4.2.)
\end{sbpara}

\begin{sbpara}
\label{assSL2orb2}
  Now assume that $G$ is reductive. Then there is a unique homomorphism 
  $\tau: \bG_{m,\R}^n\to G_\R$ whose $j$-th component $\tau_j$ gives 
the $\R$-splitting of $W^{(j)}$ associated to the $\R$-split  $G$-MHS $(W^{(j)}, \hat F_{(j)})$. This $\tau_j$ also gives  $\spl^{\BS}_{W^{(j)}}(\br)$ for $1\leq j\leq n$.  Let $\tau^{\star}: \bG_{m,\R}^n\to G'_\R\subset G_\R$ be the homomorphism defined by $\tau^{\star}(t)= \tau(t) k_0(\prod_{j=1}^n t_j)^{-1}$.

Consider the direct sum decomposition of $\Lie(G'_\R)$ by the adjoint action of $\bG_{m,\R}^n$ on it via $\tau^{\star}$. For $1\leq j\leq n$, let $\hat N_j$ be the component of $N_j$ in this direct sum decomposition on which the $s$-th factors of $\bG_m^n$ for $1\leq s<j$ act trivially. Then there is a unique homomorphism $\rho: \SL(2)^n_\R\to G_\R$ such that $\tau^{\star}_j(t) =\rho(g_1, \dots, g_n)$, where $g_s= \begin{pmatrix} 1/t& 0\\ 0&t\end{pmatrix}$ for $1\leq s\leq j$ and $g_s=1$ for $j<s\leq n$
and such that the induced map $\Lie(\rho):{\frak {sl}}(2,\R)^n  \to \Lie(G_\R)$  sends the $j$-th $\begin{pmatrix} 0&1\\0&0\end{pmatrix}$ to $\hat N_j$ for $1\leq j\leq n$. There is a unique holomorphic map $\varphi: \bP^1(\C)^n \to D$ such that $\varphi(gz)=\rho(g)\varphi(z)$ for all $g\in \SL(2,\C)^n$ and for all $z\in \bP^1(\C)^n$ and such that $\varphi(\bi)=\br$. This holomorphic map $\varphi$ is also characterized by the properties $\varphi(gz)=\rho(g)\varphi(z)$ for all $g\in \SL(2,\C)^n$ and for all $z\in \bP^1(\C)^n$ and $\varphi(\bold0)=\hat F_{(n)}$.
It satisfies $\varphi(\{0\}^j\times \{i\}^{n-j})= \hat F_{(j)}$ for $0\leq j\leq n$. 

We call this  $(\rho, \varphi)$  the $\SL(2)$-orbit associated to $(N_1, \dots, N_n, F)$. 

See \cite{CKS}  4.20 for the above facts. 
\end{sbpara}

\begin{sbpara}
\label{assSL2orb3}

Now $G$ is not necessarily assumed to be reductive. 
Let $p= \text{class}\,(\rho, \varphi)\in D_{\red, \SL(2)}$. We define $x=(p, Z)\in D_{\SL(2)}$ as follows. 
  If $N_j \neq 0$ for some $j$ and if $\gr^W(N_k)=0$ for $k=\min\{j\;|\; N_j\neq 0\}$, let $x$ be the unique $B$-orbit lying over $p$ such that $\br \in Z$.
    Otherwise, let  $x$ be the unique $A$-orbit  lying over $p$ such that $\br\in Z$. 
    
    This  $x=(p, Z)$ is called the element of $D_{\SL(2)}$ associated to $(N_1, \dots, N_n, F)$. 
    
    The rank of $p$, the family of weight filtrations $(W_x^{(k)})_{1\leq k\leq r}$ (resp.\ $(W_x^{(k)})_{0\leq k\leq r}$) associated $x$,  and the homomorphism $\tau_x=(\tau_{x,j})_{1\leq j\leq r}$ (resp.\ $(\tau_{x,k})_{0\leq k\leq r}$) associated to $x$ in the case where $x$ is an $A$ (resp.\ $B$)-orbit (\ref{taux}) are described as follows.

 (1)  Let $W^{(j)}$ ($0\leq j\leq n$) be as in \ref{assSL2orb} and let $\{j\;|\; 1\leq j\leq n, W^{(j)}\neq W^{(j-1)}\}=\{s(1), \dots, s(r)\}$ with $1\leq s(1)<\dots<s(r)\leq n$.
 Then $r$ is the rank of $p$. We have $W_x^{(k)}=W^{(s(k))}$ for $1\leq k\leq r$. In the case where $x$ is a $B$-orbit, $W_x^{(0)}=W$.

 (2) Let $\tau_j: \bG_{m,\R} \to G_\R$ ($0\leq j \leq n$)  be the homomorphism corresponding to the splitting of $W^{(j)}$ given by $\hat F_{(j)}$ in \ref{assSL2orb}. Then 
 $\tau_{x,k}=\tau_{s(k)}$ for $1\leq k\leq r$. In the case where $x$ is a $B$-orbit, $\tau_{x,0}=\tau_0$.

 Let $s(r+1):=n+1$. Then $\tau_j=\tau_{x,k}$ if $1\leq k \leq r$ and $s(k) \leq j<s(k+1)$, and $\tau_j=\tau_0$ if  $0\leq j<s(1)$.  Hence we have: 
 
 Assume $N_1\neq 0$ in the case $n\geq 1$. 
 Then
  for $(t_j)_{1\leq j\leq n} \in \bG_{m,\R}^n$, 
we have $$\tau(t)= \tau_x(t')$$ ($\tau(t):=\prod_{j=1}^n \tau_j(t_j)$, $t'\in \bG_{m,\R}^{\{1,\dots, r\}}$ if $x$ is an $A$-orbit, $t'\in \bG_{m,\R}^{\{0, \dots, r\}}$ if $x$ is a $B$-orbit), where $t'$ is defined as follows. For $1\leq k\leq r$, $t'_k= \prod_{s(k)\leq j<s(k+1)} t_j$. In the case where $x$ is a $B$-orbit, 
 $t'_0=\prod_{1\leq j<s(1)} t_j$.  
 
\end{sbpara}

\begin{sblem}\label{assoc2}
  There is a unique map $\psi: D^{\sharp}_{\Sig,[:]}\to D_{\SL(2)}$ which sends $(\sig, Z, \mathrm{class}\,((\sig_j, N_j)_{1\leq j\leq n}))$ to the element of $D_{\SL(2)}$ associated  to $(N_1, \dots, N_n, F)$ for $F\in Z$ in $\ref{assSL2orb3}$. 
\end{sblem}

This is proved similarly to \cite{KNU2} Part IV 4.5.9. %

\begin{sbthm}\label{CKS1} The map $\psi: D^{\sharp}_{\Sig,[:]}\to D^I_{\SL(2)}$  is continuous.  

\end{sbthm}

This map is the unique continuous extension of the identity map of $D$. 
     We call this map $\psi$ the {\it CKS map} respecting  the work of Cattani, Kaplan and Schmid on their $\SL(2)$-orbit theorem in many variables in pure case in \cite{CKS}. 
     This is the most important map in the fundamental diagram.
 
For the proof of Theorem \ref{CKS1}, we use the following theorem.

\begin{sbthm}\label{Eopen} Let $\sig\in\Sig$.
For the topologies in $\ref{Esig}$, $E_{\sig}$ is open in $\tilde E_{\sig}$.

\end{sbthm}

We will prove Theorem \ref{Eopen} in Proposition \ref{KU7.1.1}--\ref{8.1.b}.

From %
Theorem \ref{Eopen}, we obtain
\begin{sbthm}\label{Elm} 
  Let $\sig\in \Sig$. 
  Then 
$\tilde E_{\sig}$ ($\ref{Esig}$) and  $E_{\sig}$ are  log manifolds. 
\end{sbthm}

\begin{pf} For $\tilde E_{\sig}$, 
the argument of the proof of \cite{KU2} Proposition 3.5.10 for the pure standard case also works %
for the present case (see also \cite{KNU2} Part III 4.1.1).  

The result for $E_{\sig}$ follows from this by Theorem \ref{Eopen}. 
\end{pf}

  We start to prove Theorem \ref{Eopen}.
  The following Proposition \ref{KU7.1.1} %
  is a $G$-MHS version of \cite{KU2} Proposition 7.1.1.

 \begin{sbprop}\label{KU7.1.1}

 Let $\Sig$ be a weak fan and let 
   $\sig\in \Sig$. Let $\tilde A(\sig)$ be the closed analytic subset of $\Dc$ defined by
$\tilde A(\sig)=  \{ F \in \Dc\;|\;  N (F^pV_\C ) \subset  F^{p-1}V_\C\; \text{ for all}\;  N\in \sig, V\in \Rep(G), p\in \Z\}$,  
and let $A(\sig)$  be the subset of $\Dc$  consisting of all elements $F $ such that $(\sig, \exp(\sig_\C)F)$  is
a nilpotent orbit. Then $A(\sig)$  is an open set of $\tilde A(\sig)$ in the strong topology of $\tilde A(\sig)$ in $\Dc$ ($\ref{Esig}$).

  \end{sbprop}
  
  \begin{pf}   Assume $F\in A(\sig)$ and assume that a directed family $F_{\la}\in \tilde A(\sig)$ converges to $F$. We prove that $F_{\la}\in A(\sig)$ for $\la$ sufficiently large. Let $H\in D$. Take a faithful representation $V\in \Rep(G)$ and take  $\R$-polarizations $\langle\cdot, \cdot\rangle_w$  of $H(\gr^W_wV)$ which are stable under $G'$. 
  By \ref{lemopen1} (2),  for $\la$ sufficiently large, the annihilator of $F^p_{\la}(\gr^W_wV)$ with respect to $\langle \cdot, \cdot\rangle_w$ is $F^{w+1-p}(\gr^W_wV)$. 
  Take a finite set $N_1, \dots, N_n$ of generators of the cone $\sig$. The proof of \cite{KU2} Proposition 7.1.1 shows that if $\la$ is sufficiently large, then $\exp(\sum_{j=1}^n z_jN_j)F_{\la}(\gr^W_wV)\in G(\R)H(\gr^W_wV)$ for every $w$ if $\text{Im}(z_j)$ are sufficiently large. 
  By \ref{faith}, $F_{\la}\in A(\sig)$ for $\la$ sufficiently large. 
    \end{pf}

 \begin{sbpara}\label{oneN}  The proof of Theorem \ref{Eopen}  consists of complicated inductions. 
  To clarify the idea, we explain here the proof when the cone $\sig$ is of rank one. 
  
    By the canonical isomorphisms $\Gamma(\sig)\simeq \N$ and $P(\sig)\simeq \N$, we have  canonical isomorphisms $\toric_{\sig}\simeq \C$, $|\toric|_{\sig}\simeq \R_{\geq 0}$.  
    We identify $\check{E}_{\sig}$ with $\C\times \check{D}$ and $\check{E}^{\sharp}_{\sig}$ with $\R_{\geq 0} \times \Dc$. Since $E_{\sig}$ is the inverse image of $E^{\sharp}_{\sig}$ under $|\cdot|\,:\, \check{E}_{\sig}\to \check{E}^{\sharp}_{\sig}$, it is sufficient to prove that $E^{\sharp}_{\sig}$ is open in $\tilde E^{\sharp}_{\sig}$.

  Assume that $(q_{\lam}, F_{\lam}) \in \tilde E^{\sharp}_{\sig}$ converges to $(0, F)\in E^{\sharp}_{\sig}$. 
  We prove that $(q_{\lam}, F_{\lam}) \in E^{\sharp}_{\sig}$ for sufficiently large $\lam$.  
  Taking subsequences, we divide it into two cases: the case where $q_{\lam}=0$ for all $\lam$ and the case $q_\lam\not=0$ for all $\lam$.

Assume $q_{\lam}=0$.  
This is the case Proposition \ref{KU7.1.1} for $\rank(\sig)=1$.

Assume $q_{\lam} \neq 0$.
  Let $q_\lam=e^{-2\pi y_\lam}$.   We have to see that, for any sufficiently large $\lam$, 
$\exp(iy_{\lam}N)F_{\lam}\in D$  (which means $(q_{\lam}, F_{\lam})\in E_{\sig}$ by definition \ref{Esig}). 
Here $N\in\sig$ is such that $\exp(N)$ corresponds $1\in \N$ via $\Gamma(\sig)\simeq\N$.

  First, applying \cite{KU2} Proposition 3.1.6 to $S:=\tilde E_{\sig}\subset X:=\Ec_\sig=\C\times\Dc$, $s=(0,F)$, $s_{\lam}=(q_{\lam}, F_{\lam})$, and $A:=\{0\}\times \Dc$, we find $F^*_{\lam}\in\Dc$ which is very near to 
$F_{\lam}$ such that $(0, F^*_{\lam})\in \tilde E_{\sig}$ and $F^*_{\lam}$ converges to $F$. 
  For the precise meaning of \lq\lq very near'' here, see the condition (i) in Lemma \ref{lem712}.

  By Proposition \ref{KU7.1.1} again, we may assume that 

\smallskip

(1) $(N, F^*_{\lam})$ generates a nilpotent orbit.

Next, let $\tau:\bG_{m,\R}\to G_\R$ be the homomorphism corresponding to the splitting of $W^{(1)}:=M(N,W)$ given by $\hat F_{(1)}$ in \ref{assSL2orb3}, and let $\tau_{\lam}=\tau(1/\sqrt{y_{\lam}})$.
  By $(1)$, %
$(W^{(1)},F_{\lam}^*)$ is a $G$-MHS.
  Together with the fact that $\tau_{\lam}$ splits %
  $W^{(1)}$, we have 
  
(2) $\tau_{\lam}^{-1}F^*_{\lam} \to \hat F_{(1)}$.  

  From this, we have 
  
(3)  $\tau_{\lam}^{-1}\exp(iy_{\lam}N)F^*_{\lam}=\exp(iN)\tau_{\lam}^{-1}F^*_{\lam} \to \exp(iN)\hat F_{(1)}$. 
  
  Since $F_{\lam}$ and $F^*_{\lam}$ are very near, we have $\tau_{\lam}^{-1}\exp(iy_{\lam}N)F_{\lam} \to \exp(iN)\hat F_{(1)} \in D$.
  Since $D$ is open in $\Dc$, we conclude 
  
(4) $\exp(iy_{\lam}N)F_{\lam}\in D$ for sufficiently large $\lam$.

\end{sbpara}

\begin{sbpara}\label{8.1.a}
  Note that a large part of the following %
\ref{8.1.a}--\ref{8.1.b} is a copy from \cite{KNU2} Part IV Appendix (after some modifications explained in Remark \ref{c_to_IVApp}). 
  For example, the statement of Proposition \ref{prop712} is almost identical to that of \cite{KNU2} Part IV Proposition A.1.7. 

  To prove Theorem \ref{Eopen}, 
since $\tilde E_{\sig,\val}\to\tilde E_\sig$ 
is proper and surjective, and $E_{\sig,\val}\subset\tilde E_{\sig,\val}$ 
is the inverse image of $E_\sig\subset\tilde E_\sig$, 
it is sufficient to prove that $E_{\sig,\val}$ is open in $\tilde E_{\sig,\val}$. 
  Since we have the continuous  projection $\tilde E_{\sig,\val}\to \tilde E^{\sharp}_{\sig,\val}$ for which $E_{\sig,\val}$ is the inverse image of $E^{\sharp}_{\sig,\val}$, it is sufficient to prove that $E^{\sharp}_{\sig,\val}$ is open in $\tilde E^{\sharp}_{\sig, \val}$.

  Assume that a directed system $w_{\lam}=(q_\lam, F'_\lam)\in \tilde E^{\sharp}_{\sig,\val}$ ($\lam\in L$, where  $L$ is a directed ordered set) converges in $\tilde E^{\sharp}_{\sig,\val}$ to $w=(q, F')\in E^{\sharp}_{\sig,\val}$.
  We prove that $w_{\lam}\in E^{\sharp}_{\sig,\val}$ for any sufficiently large $\lam$. 
\end{sbpara}

\begin{sbpara}\label{N,F}
We fix notation.

Let $w_{[:]}=(t, F')$ (resp.\ $w_{\lam, [:]}=(t_{\lam}, F'_{\lam})$) be the image $w$ (resp.\ $w_{\lam}$) in $\check E^{\sharp}_{\sig,[:]}=|\toric|_{\sig,[:]}\times \check{D}$ ($t, t_{\lam}\in |\toric|_{\sig,[:]}$, $F', F'_{\lam}\in \Dc$). 

We use the additive presentation of $|\toric|_{\sig, [:]}$ in \ref{add2}.
Let $n$ and $\sig_n$ be as those in there for the point $t$.
Take an $\bR$-subspace $B$ of $\sig_\bR$ such that
$\sig_{\bR}=\sig_{n,\bR}\oplus B$. 
Take a good base $(N_s)_{s\in S}$ for $t$ and let $S_j\subset S$ be as in the definition of a good base (\ref{excel1}). Write 
$t=(\text{class}\,((\sig_j, N_j)_{1\le j \le n}), b)$, where $N_j$ is an interior point of $\sig_j$ and $b\in B$. 
Let $a_s\in\R_{>0}$ such that $N_j\equiv\sum_{s\in S_j} a_sN_s \mod \sig_{j-1,\R}$ for all $1\le j\le n$. 

Let $F=\exp(ib)F'$.  
Then $((N_s)_{s\in S}, F)$ generates a nilpotent orbit.
We have also the $B$-component $b_{\la}\in B$ of $t_{\la}$. Let $F_{\lam}=\exp(ib_{\la})F'_{\la}$. Then $F_{\lam}$ converges to $F$ in $\Dc$. 

Take an open neighborhood $V$ of $b$ in $B_\C$ such that the map $\toric_{\sig_n} \times V \to \toric_{\sig}$ induced by the canonical map $\toric_{\sig_n} \times B_\C \to \toric_{\sig}$ is an open immersion. We may assume $b_\la\in V\cap B$. 

Let $\sig'\subset \sig$ be the cone generated by $N_s$ ($s\in S$). Replacing $N_s$ by $N_s/r$ for some integer $r>0$, we may assume that $\log\Gamma(\sig')\subset \bigoplus_{s\in S} \N N_s\subset \bigoplus_{s\in S} \Z N_s$. 
  We have the injective homomorphism $\Gamma(\sig')\overset{\subset}\to \N^S$. For $s\in S$ we have the $s$-component $\Gamma(\sig')\to \N$, which is an element $q_s$ of $P(\sig')$,  regarded as a holomorphic function on $\Spec(\C[P(\sig')])^{\an}\times V=\toric_{\sig'}\times V$.

We use the description of the topology of $|\toric|_{\sig,[:]}$ in \ref{add3} associated to the good base $(N_s)_{s\in S}$ for $t$.
  We use the notation there. 
Let $U$  be the open neighborhood of $t$. %
We may assume that all $t_{\lam}$ belong to $U$.

Since $w_{\lam}\to w$ and since $|\toric|_{\sig', \val} \times B$ 
is open in $|\toric|_{\sig, \val}$, we may assume that all $w_{\lam}$ belong to the open set $\tilde E_{\sig', \val} \times V$ of $\tilde E_{\sig,\val}$.

  Note that we have  $c_j\in S_j$ and $T$ as in \ref{add3}.  
In the case $q_s(w_{\lam})\neq 0$, $y_{\lam,s}\in \R_{>0}$ is defined by $|q_s(w_{\lam})|= \exp(- 2\pi y_{\lam,s})$.  
If $s\in S_j$ and $q_s(w_{\la})\neq 0$, then $q_{s'}(w_{\la})\neq 0$ for
any $s'\in S_{j'}$ with $j'\geq j$,  and in the case $s\neq c_j$,
$y_{\la,s}/y_{\la,c_j}$ gives the $s$-component of $t_{\la}$ in $U\simeq
\R^n_{\geq 0} \times \R_{\geq 0}^T \times B \to \R_{\geq 0}^T$. If $s\in
S_j$ and if the set $I=\{\la\;|\; q_s(w_{\la})\neq 0\}$ is cofinal in the
directed index set, then for $s'\in S_{j'}$ with $j'> j$ (resp.\ for $s'\in
S_j$) and for $\la\in I$, we have $y_{\la, s'}/y_{\la,s}\to 0$ (resp.\
$y_{\la,s'}/y_{\la,s}\to a_{s'}/a_s$).

Let $x_{\SL(2)}\in D_{\SL(2)}$ be the element associated to $(N_1, \dots, N_n, F)$ in \ref{assSL2orb3}.  
We have a homomorphism $\tau: \bG_m^n \to G_\R$ as in (2) in \ref{assSL2orb3}. In the following, $\tau_j$ denotes its $j$-th component. 
\end{sbpara}

\begin{sbpara}\label{e,tau}
We may assume that there exists an integer $m$ such that $1\leq m\leq n+1$ and such that $q_s(x_\lam)=0$ for any $\lam$ and $s\in S_{\leq m-1}$ and $q_s(x_\lam)\neq 0$ for any $\lam$ and $s\in S_{\geq m}:=\bigsqcup_{k \ge m}S_k$ ($S_{\leq 0}$ and $S_{\geq n+1}$ are defined as the empty set).

For $m \leq j\leq n$, let 
$$
e_{\lam,\geq j}=\exp(\sum_{s\in S_{\geq j}} iy_{\lam,s}N_s)\in G_\bC,
$$ 
$$
\tau_{\lam,j}= 
\tau_j\left(\sqrt{y_{\lam,c_{j+1}}/y_{\lam,c_j}}\right)\in G_\bR, 
\quad  \tau_{\lam, \geq j}
 = \prod_{k= j}^n \tau_{\lam,k}\in G_\bR
$$
$(y_{\lam,c_{n+1}}$ denotes $1$). 
Let $\hat F_{(j)}$ ($1\leq j\leq n$) be as in \ref{assSL2orb} associated to $(N_1,\dots, N_n, F)$.

\end{sbpara}

\begin{sblem}\label{lem712} Let the situation and the notation be as above. 
  Let $d$ be a local metric on $\Dc$ that is compatible with the analytic structure. 
Let $m\leq j \leq n$ and let $e\geq 0$. Then for any sufficiently large $\lam$, there exist $F^*_{\lam}\in \Dc$ satisfying the following {\rm (i)} and {\rm (ii)}.

\medskip

{\rm (i)} $y_{\lam,s}^ed(F_{\lam}, F^*_{\lam})\to 0 \;\;(\forall s\in S_j)$.

\smallskip

{\rm (ii)} $(N_s, F^*_{\lam})$ {\it satisfies Griffiths transversality for any}\; $s\in S_{\leq j}.$

\end{sblem}

\begin{pf}
  We apply \cite{KU2} Proposition 3.1.6 to the following situation:
$S:=\tilde E_{\sig'}\times V\subset X:=\Ec_{\sig'}\times V=\toric_{\sig'}\times V\times \Dc$ (\ref{Esig}), where $\sig'$ and $V$ are as in \ref{N,F}, and $A$ is the closed analytic subspace of $X$ defined by $q_s=0$ for $s\in S_{\leq j}$.

Consider the images of $w_{\lam}=(q_\lam,F'_\lam)$ and $w=(q,F')$ of \ref{8.1.a} in $X$ and $S$, respectively, and let $F_\lam$, $F$, and $y_{\lam,s}$ be as in \ref{N,F}.

Take a metric on a neighborhood of the image of $w$ in $X$ as the direct sum of a local metric on $\toric_{\sig'}\times V$ and the local metric $d$ on $\Dc$.
Then we get Lemma \ref{lem712}.
\end{pf}

  The next is a key proposition of this section.

\begin{sbprop}\label{prop712}
Let the situation and the assumption be as above.
 Then the following assertions $(A_j)$ $(m-1\leq j\leq n)$, $(B_j)$ $(m\leq j\leq n)$, $(C_j)$ $(m\leq j\leq n)$ are true.

$(A_j)$ (resp.\ $(B_j)$, resp.\ $(C_j)$) for $m\leq j\leq n${\rm:}
Let $e\geq 1$. Then for any sufficiently large $\lam$, there are $F^{(j)}_{\lam}\in \Dc$ satisfying the following $(1)$--$(3)$.

$(1)$ $y^e_{\lam,c_j}d(F_{\lam}, F^{(j)}_{\lam}) \to 0$.

$(2)$ $((N_s)_{s\in S_{\leq j}}, e_{\lam,\geq j+1}F^{(j)}_{\lam})$ generates a nilpotent orbit.  

$(3)$  $\tau_{\lam,\geq j+1}^{-1}e_{\lam,\geq j+1} F^{(j)}_{\lam}\to \exp(iN_{j+1})\hat F_{(j+1)}$. 

(resp.\ $(3)$ $\tau_{\lam,\geq j}^{-1}e_{\lam,\geq j+1} F^{(j)}_{\lam}\to \hat F_{(j)}$. 

resp.\ $(3)$ $\tau_{\lam,\geq j}^{-1}e_{\lam,\geq j}F^{(j)}_{\lam}\to \exp(iN_j)\hat F_{(j)}$.)

\noindent
Here $(A_n)$ is formulated by understanding $N_{n+1}=0$ and $\hat F_{(n+1)}=F$. 

\medskip

$(A_{m-1})${\rm :} For any sufficiently large $\lam$, we have the following $(2)$ and $(3)$. 

$(2)$ $((N_s)_{s\in S_{\leq m-1}}, e_{\lam,\geq m}F_{\lam})$ generates a nilpotent orbit.  

$(3)$ $\tau_{\lam,\geq m}^{-1}e_{\lam,\geq m} F_{\lam}\to \exp(iN_m)\hat F_{(m)}$. 
\end{sbprop}

\begin{sbrem}\label{oneNrem} The case $n=1$ was treated in \ref{oneN}. The case $n=m=1$ is the case $q_{\lam}\neq 0$ for all $\la$ in \ref{oneN}. The case $n=1$ and $m=2$ is the case $q_{\lam}=0$ for all $\la$ in \ref{oneN}. In the case $n=m=1$, the course of the arguments about (1), (2), (3), (4) in \ref{oneN} is understood as 
$(1)=(A_1)\Rightarrow(2)=(B_1)\Rightarrow(3)=(C_1)\Rightarrow(4)=(A_0)$.

\end{sbrem}
\begin{sbpara}
We prove Proposition \ref{prop712} by using the downward induction of the form

$(A_j)$ $\Rightarrow$ $(B_j)$ $\Rightarrow$ $(C_j)$ $\Rightarrow$ $(A_{j-1})$. (Here $m\leq j\leq n$.)

$(B_j)$ $\Rightarrow$ $(C_j)$ is clear. 

We prove $(A_j)$ $\Rightarrow$ $(B_j)$. 
By (2), for any sufficiently large $\lam$, $(W^{(j)}, e_{\lam,\geq j+1} F^{(j)}_{\lam})$ is a mixed Hodge structure, where $W^{(j)}$ denotes the relative monodromy filtration of $N_1+\dots+N_j$ with respect to $W$. 
  Since $\tau_{\lam,j}$ splits $W^{(j)}$, we have $(B_j)$. 

We prove $(C_{j+1})$ $\Rightarrow$ $(A_j)$ if $m-1 \leq j < n$ and $(A_n)$. 

  Let $e \ge1$. 
If $m\leq j \le n$ (resp.\ $j=m-1$), then, by Lemma \ref{lem712}, there are $F^{(j)}_{\lam}\in \Dc$ satisfying (1) and (resp.\ $F^{(j)}_{\lam}:=F_{\lam}$ satisfies) the condition 

\medskip

$(2')$ $(N_s, F^{(j)}_{\lam})$ satisfies Griffiths transversality for any $s\in S_{\le j}$.

\medskip

  When $m-1\le j<n$, by $(C_{j+1})$, there are $F^{(j+1)}_{\lam}\in \Dc$ satisfying 
(resp.\ when $j=n$, $F^{(j+1)}_{\lam}=F_{\lam}$ satisfies) the following. 

\medskip

$(1'')$ $y^e_{\lam,c_{j+1}}d(F_{\lam}, F^{(j+1)}_{\lam}) \to 0$.

$(3'')$  $\tau_{\lam,\geq j+1}^{-1}e_{\lam,\geq j+1} F^{(j+1)}_{\lam}\to \exp(iN_{j+1})\hat F_{(j+1)}$.  (Recall $\hat F_{(n+1)}$ means $F$.)

\medskip

By $(1'')$ and $(3'')$, we have

\medskip

(4) $\tau_{\lam,\geq j+1}^{-1}e_{\lam,\geq j+1} F_{\lam}\to \exp(iN_{j+1})\hat F_{(j+1)}$.

\medskip

By (4) and by (1), we have 

\medskip

(5) $\tau_{\lam,\geq j+1}^{-1}e_{\lam,\geq j+1} F^{(j)}_{\lam}\to \exp(iN_{j+1})\hat F_{(j+1)}$.

\medskip

For the left-hand side of (5), by $(2')$, $(N_s, \tau_{\lam,\geq j+1}^{-1}e_{\lam,\geq j+1} F^{(j)}_{\lam})$ satisfies Griffiths transversality for any $s\in S_{\leq j}$. On the other hand, concerning the right-hand side of (5), 
$((N_s)_{s\in S_{\leq j}}, \exp(iN_{j+1})\hat F_{(j+1)})$ generates a nilpotent orbit. Hence (5) and Proposition \ref{KU7.1.1} show that
$((N_s)_{s\in S_{\leq j}}, \tau_{\lam,\geq j+1}^{-1}e_{\lam,\geq j+1} F^{(j)}_{\lam})$ generates a nilpotent orbit for any sufficiently large $\lam$. This proves that $((N_s)_{s\in S_{\leq j}}, e_{\lam,\geq j+1} F^{(j)}_{\lam})$ generates a nilpotent orbit for any sufficiently large $\lam$. 
By this and by (1) and (5), we have 
$(A_j)$. 
\end{sbpara}

\begin{sbpara}\label{8.1.b}

By $(A_{m-1})$ (2) of Proposition \ref{prop712}, $w_{\lam}$ belongs to $E^{\sharp}_{\sig,\val}$ if $\lam$ is sufficiently large. 

  This proves that $E^{\sharp}_{\sig,\val}$ is open in $\tilde E^{\sharp}_{\sig,\val}$, and hence proves that $E_{\sig}$ is open in $\tilde E_{\sig}$ (\ref{8.1.a}). 
  The proof of Theorem \ref{Eopen} is completed.
  
\end{sbpara}

\begin{sbpara}\label{CKS2}   We next prove Theorem \ref{CKS1} (the continuity of the CKS map) in \ref{regular}--\ref{8.3.d}.
   
  The outline of the proof of Theorem \ref{CKS1} is as follows. By using the properties of regular spaces reviewed in \ref{regular} applied to the spaces $D^I_{\SL(2)}(\Phi)$, it becomes enough to discuss the convergence of  ordinary points to a boundary point (we do not need to discuss the convergence of boundary points to a boundary point). We reduce the proof of Theorem \ref{CKS1} to the three  convergences in Proposition \ref{III3.3.4}. 
These three convergences are proved in \ref{III3.3.3}, 
\ref{8.3.c}, and \ref{KU6.2.1}--\ref{8.3.d}, respectively. 
  The structure of these proofs is essentially the same as that in \cite{KU2}, \cite{KNU2} Parts III and IV.
\end{sbpara}

\begin{sbpara}\label{regular} Concerning regular spaces, we first recall two facts.

\cite{KU2} Definition 6.4.6 (\cite{Bou}, Ch.\ 1, \S 8, no.\ 4, Definition 2). A topological space $X$ is called {\it regular} if it is Hausdorff and satisfies the following axiom: Given any closed subset $F$ of $X$ and any point $x \notin F$, there is a neighborhood of $x$ and a neighborhood of $F$ which are disjoint.

(\cite{Bou}, Ch.\ 1, \S 8, no.\ 5, Theorem 1.) Let $X$ be a topological space, $A$ a dense subset of $X$, $f:A\to Y$ a map from $A$ into a regular space $Y$. A necessary and sufficient condition for $f$ to extend to a continuous map $f : X \to Y$ is that, for each $x \in X$, $f(y)$ tends to a limit in $Y$ when $y$ tends to $x$ while remaining in $A$. The continuous extension ${\bar f}$ of $f$ to $X$ is then unique.
\end{sbpara}

\begin{sbpara}

   To prove that $\psi: D^{\sharp}_{\Sig, [:]}\to D^I_{\SL(2)}$ is continuous, it is sufficient to prove that the composition  
$$E^{\sharp}_{\sig, \val}\to E^{\sharp}_{\sig, [:]} 
\to D^{\sharp}_{\Sig, [:]}\overset{\psi}\to D^I_{\SL(2)}$$
is continuous 
for each $\sig \in \Sig$.

  Let the notation and the assumption be as in \ref{8.1.a} and \ref{N,F}.
  Assume further that all $w_{\lambda}$ belong to $E^{\sharp}_{\sig, \val}$. 
  Then it is enough to show that the image of $w_{\la}$ in $D^I_{\SL(2)}$ converges to the image of $w$ in $D^I_{\SL(2)}$.

If $(\sig, Z)\in D_{\Sig}$, the relative monodromy filtration $M(N', W)$ exists for any $N'\in \sig$, and $M(N', W)$ depends only on the face $\sig'$ of $\sig$ such that $N'$ is in an interior point of $\sig'$. We will denote $M(N', W)$ as $W(\sig')$. 

Recall that the image of $w\in E^\sharp_{\sig,\val}$ in $D_{\Sig, [:]}$ is $(\sig, Z, \text{class}\,((\sig_j, N_j)_{1\leq j\leq n}))$. 
Since $N_j$ is in the interior of $\sig_j$ for each $j$,  $N_1+\dots+ N_j$ is in the interior of $\sig_j$, and hence $W^{(j)}$ in \ref{assSL2orb} for $(N_1, \dots, N_n, F)$ coincides with $ W(\sig_j)$. 

Hence if $w_{\lam}\in E^\sharp_{\sig,\val}$ is near to $w$, the $W(\sig_j)$ of $w_{\la}$ are contained in the set $\Phi$ of $W(\sig_j)$
($1\leq j\leq n$), 
and hence the image of $w_{\la}$ 
in $D_{\SL(2)}$ is contained in $D_{\SL(2)}(\Phi)$. 
  Since $D_{\SL(2)}(\Phi)$ is regular by Proposition \ref{str3}, 
we can use the above regular point method, and hence in the proof of the fact that the image of $w_{\la}$ in $D^I_{\SL(2)}$ converges to the image of $w$ in $D^I_{\SL(2)}$, we may assume that $w_{\la}\in |\torus|_{\sig}\times D$. 
\end{sbpara}

\begin{sbpara}\label{4.5.21} %
  To state Proposition \ref{III3.3.4} on three kinds of convergences, which are used to prove the continuity of the CKS map, we recall that we are in  the following situation. We keep the notation in \ref{N,F} and assume $w_{\la}\in |\torus|_{\sig}\times D$. So  the $|\torus|_{\sig}$-component of $w_{\la}$ is written as $(\sum_{s\in S} y_{\la,s}N_s) +b_{\la}$. 

For each $1\le j\le n$, we have $c_j\in S_j$ and we have $N_j\equiv\sum_{s\in S_j}a_sN_s \mod \sig_{j-1,\R}$. 

Note that $(N_1, \cdots, N_n, F)$ generates a nilpotent orbit.
  Let $\br \in D$ be the point associated to $(N_1, \dots, N_n, F)$ in \ref{assSL2orb}.

\end{sbpara}

By \ref{regular}, for the proof of Theorem \ref{CKS1}, it is sufficient to prove that 
$\exp(\sum_{s\in S}iy_{\lam,s}N_s)F_\lam $ converges to the class in $D^I_{\SL(2)}$ of the $\SL(2)$-orbit associated to $(N_1, \dots, N_n,F)$.
As %
will be explained in \ref{conv} below, this is reduced to the following proposition.  

\begin{sbprop} \label{III3.3.4}

We have the following convergences {\rm(1)}--{\rm(3)} as $\lam\to\infty$.
\smallskip

{\rm(1)}
$\tau\left(\sqrt{\frac{y_{\lam,c_1}}{y_{\lam,c_2}}},\dots,
      \sqrt{\frac{y_{\lam,c_n}}{y_{\lam,c_{n+1}}}}\right)
\exp(\sum_{s\in S}iy_{\lam,s}N_s)F_\lam \to \br$ in $D$, 
where $y_{\lam,c_{n+1}}=1$.

\smallskip

{\rm(2)} $\spl_W\bigl(\exp(\sum_{s\in S}iy_{\lam,s}N_s)F_\lam\bigr) \to \spl_W(\br)$ in $\spl(W)$.

\smallskip

{\rm (3)} $(\exp(\sum_{s\in S}iy_{\lam,s}N_s)F_\lam)_{\red}$ converges to the class in $D_{\red,\SL(2)}$ of the $\SL(2)$-orbit associated to $(N_1, \dots, N_n, F)$. 
\end{sbprop}

\begin{sbpara}\label{conv}  
  We explain that this Proposition \ref{III3.3.4} shows that 
$x_{\la}:=\exp(\sum_{s\in S}iy_{\lam,s}N_s)F_\lam $ converges to the class $x_{\SL(2)}$ in $D^I_{\SL(2)}$ of the $\SL(2)$-orbit associated to $(N_1, \dots, N_n,F)$ and hence proves Theorem \ref{CKS1}.

 By (2) and (3) of Proposition \ref{III3.3.4} and by Proposition \ref{emb2} and  Proposition \ref{emb3}, for the proof of $x_{\la} \to x_{\SL(2)}$, it is sufficient to prove
the following (1) and (2). Let $p=x_{\SL(2), \red}\in D_{\red,\SL(2)}$.

\smallskip

(1) 
 Let $\beta$ be a  distance to $\Phi(p)$-boundary. Then if $x_{\SL(2)}$ is an $A$-orbit (resp.\ $B$-orbit),
$\Ad(\tau_p(\beta(x_{\la})))^{-1}\delta(x_{\la})$ converges to $\Ad(\tau_p(\beta(\br)))^{-1}\delta(\br)$ (resp.\ $0\circ \Ad(\tau_p(\beta(\br)))^{-1}\delta(\br)$) in $\overline{\cL}$. 

\smallskip

(2)  
Let $\beta$ be a distance to $\Phi(x_{\SL(2)})$-boundary. 
Then $\tau_{x_{\SL(2)}}(\beta(x_{\la}))^{-1}x_{\la}$ converges to $\tau_{x_{\SL(2)}}(\beta(\br))^{-1}\br$ in $D$.

\smallskip

We prove that these (1) and (2) follow from (1) of Proposition \ref{III3.3.4}.

Let  $t_{\la}=  \left(\sqrt{\frac{y_{\lam,c_2}}{y_{\lam,c_1}}},\dots,
      \sqrt{\frac{y_{\lam,c_{n+1}}}{y_{\lam,c_n}}}\right)\in \R_{>0}^n$.  We have $\tau(t_{\la})= \tau_{x_{\SL(2)}}(t'_{\la})$, where $t'_{\la}\in \R_{>0}^{\Phi(x_{\SL(2)})}$ is as in \ref{assSL2orb3} (2). 
      Since  $t_{\la}$ converges to $ \bold 0$ in $\R^n_{\geq 0}$, $t'_{\la}$ converges to $\bold 0$ in $\R_{\geq 0}^{\Phi(x_{\SL(2)})}$.

       By (1) of Proposition \ref{III3.3.4}, we have
       
       \smallskip

\noindent       
       $(*)$ \; $x_{\la}= \tau_{x_{\SL(2)}}(t'_{\la})u_{\la}\br$, $u_{\la}\in G(\R)G_u(\C)$, $u_{\la}\to 1$. 
      
    \smallskip
  
  If $x_{\SL(2)}$ is a $B$-orbit,  since $\br\in D_{\nspl}$, $u_{\la}\br\in D_{\nspl}$ if $\la$ is sufficiently large and hence by $(*)$, $x_{\la} \in D_{\nspl}$ if $\la$ is sufficiently large.

     We prove (1). By applying $\beta$ to $(*)_{\red}$, we obtain that if $x_{\SL(2)}$ is an $A$-orbit (resp.\ $B$-orbit),  
     $\Ad(\tau_p(\beta(x_{\la,\red})))^{-1}\delta(x_{\la})$ is equal to $\tau_p(\beta(u_{\la}\br_{\red}))^{-1}\delta(u_{\la}\br)$ (resp.\ $t'_{\la,0}\circ \Ad(\tau_p(\beta(u_{\la}\br_{\red})))^{-1}\delta(u_{\la}\br)$) and this converges to $\Ad(\tau_p(\beta(\br_{\red})))^{-1}\delta(\br)$ (resp.\ $0\circ \Ad(\tau_p(\beta(\br_{\red})))^{-1}\delta(\br)$).

     We prove (2).  By applying $\beta$ to $(*)$, we have  $\tau_{x_{\SL(2)}}(\beta(x_{\la}))^{-1}x_{\la} = \tau_{x_{\SL(2)}}(\beta(u_{\la}\br))^{-1}u_{\la}\br$, and 
     this converges to $\tau_{x_{\SL(2)}}(\beta(\br))^{-1}\br$.

\end{sbpara}
   
\begin{sbpara}\label{III3.3.3}
  We start to prove Proposition \ref{III3.3.4}.

  Proposition \ref{III3.3.4} (1) follows from $(A_0)$ of the case $m=1$ of Proposition \ref{prop712}. 

  We prepare the proofs of Proposition \ref{III3.3.4} (2) and (3).

For $1\le j\le n$, let $\Dc_j$ be the subset of $\Dc$ consisting
of all $F'\in\Dc$ such that $((N_s)_{s\in S_{\le j}},F')$ generates
a nilpotent orbit. 
We have 
 $F\in \Dc_n$,  $F_\lam\in\Dc$,
and we know from \ref{N,F}, Proposition \ref{prop712} and Theorem \ref{Eopen} 
that the following five conditions are satisfied.

\medskip

\noindent
{\rm(1)} $F_\lam$ converges to $F$ in $\Dc$.
\smallskip

\noindent
{\rm(2)} $y_{\lam,s}\to\infty$ for any $s\in S$.
\smallskip

\noindent
{\rm(3)} If $1\le j<n$, $s\in S_{\le j}$ and
$t\in S_{\ge j+1}$, then
$\tfrac{y_{\lam,s}}{y_{\lam,t}}\to\infty$.
\smallskip

\noindent
{\rm(4)} If $1\le j\le n$ and $s,t\in S_j$, then
$\tfrac{y_{\lam,s}}{y_{\lam,t}}\to
\tfrac{a_s}{a_t}$.
\smallskip

\noindent
{\rm(5)} For $1\le j\le n$ and $e\ge0$,
there exist
$F^*_\lam\in\Dc$ $(\lam\in L)$ satisfying the following conditions
(5.1) and (5.2).

(5.1) $\exp\big(\sum_{t\in S_{\ge j+1}}
iy_{\lam,t}N_t\big)F_\lam^*\in\Dc_j\quad
(\lam:\text{sufficiently large}),$

(5.2) $y_{\lam,s}^e d(F_\lam,F^*_\lam)\to0\quad
(\forall s\in S_j).$

\noindent
Here $d$ is a metric on a neighborhood of $F$ in $\Dc$ which is compatible with the analytic topology of $\Dc$. 
\medskip

\end{sbpara}

\begin{sbpara}\label{8.3.c}  
  Proof of Proposition \ref{III3.3.4} (2).
(This is almost a copy from \cite{KNU2} Part III 3.3.7.)

We prove the following assertion $(A_j)$ by downward induction on $j$.
(Note that $(A_0)$ is what we want to prove.) 

\smallskip

$(A_j)$ : Proposition \ref{III3.3.4} (2) is true in the case where $\exp\big(\sum_{t\in S_{\ge j+1}} iy_{\lam,t}N_t\big)F_\lam\in\Dc_j$ for all $\lam$.   
Here $\Dc_j$ for $1\le j \le n$ are defined in \ref{III3.3.3}, and $\check D_0:=D$. 

\medskip

{\it Proof.} 
Let $0\leq j\leq n$.

Let $x_{\la} = \exp(\sum_{s\in S}iy_{\lam,s}N_s)F_\lam$, 
$\tau_{\lam,>j}= \prod_{j< k\le n}\tau_k\Big(
\sqrt{\tfrac{y_{\lam,c_k}}{y_{\lam,c_{k+1}}}}\Big)$, and $x_{\la,j} = \tau_{\lam,>j}x_{\la}$. 
Then, by \cite{KNU1} 10.3, we have  
$$
x_{\la,j} = \tau_{\lam,>j}x_{\la} =\exp\Bigl(\sum_{s\in S_{\leq j}}i
\tfrac{y_{\lam,s}}{y_{\lam,c_{j+1}}}N_s\Bigr)U_{\lam,j},
$$
$$
\text{where}\quad 
U_{\lam,j}:=\tau_{\lam,>j}\exp\Bigl(\sum_{t\in S_{\ge j+1}}iy_{\lam,t}N_t\Bigr)F_\lam.
$$
  Assume $\exp\big(\sum_{t\in S_{\ge j+1}} iy_{\lam,t}N_t\big)F_\lam\in\Dc_j$.
  Then 
$(N'_1,\ldots,N'_j, U_{\lam,j})$ generates a nilpotent orbit, 
where $N'_k=\sum_{s \in S_k}\frac{y_{\lam,s}}{y_{\lam,c_k}}N_s$ 
($1 \le k \le j$).
  Let $s_{\lam}$ be the associated limit splitting. %
  By \cite{KNU1} Theorem 0.5 (2) and ibid.\ Proposition 10.8 (1), there is a convergent power series 
$u_{\lam}$ whose values are in $G(\bR)$,
whose constant term is $1$ and 
whose coefficients depend on $U_{\lam,j}$ and $y_{\lam,s}/y_{\lam,c_k}$ ($1\leq k\leq j$, $s\in S_k$) real analytically 
such that 
$\spl_W(x_{\la,j}) =u_{\lam}\Bigl(\frac{y_{\lam,c_2}}{y_{\lam,c_1}}, \dots, 
\frac{y_{\lam,c_{j+1}}}{y_{\lam,c_j}}\Bigr)s_{\lam}$.
  Since $s_{\lam}$ also depends real analytically on $U_{\lam,j}$, 
we have 

\smallskip

\noindent
(1) $\spl_W(x_{\la,j})$ converges to $\spl_W(\br)$.

\medskip

  This already showed ($A_n$). 

  Next, assume $j < n$ 
 and assume that $(A_{j+1})$ is true. We prove that $(A_j)$ is true. 

  Choose a sufficiently big $e>0$ depending on $\tau_{j+1},\ldots,\tau_n$.

  Take $F^*_{\lam}$ %
as in \ref{III3.3.3} (5). 
  Define $x_{\la}^*$ and $x_{\la,j}^*$ similarly 
as $x_{\la}$ and $x_{\la,j}$, respectively. 

Then we have 

\smallskip

\noindent
(2) $\spl_W(x^*_{\la,j})$ converges to $\spl_W(\br)$, and
$y_{\lam, c_{j+1}}^ed(\spl_W(x_{\la,j}), \spl_W(x^*_{\la,j})) \to 0$.

\smallskip

By downward induction hypothesis on $j$, we have

\smallskip

\noindent
(3) $\spl_W(x^*_{\lam})= \tau_{\lam,>j}^{-1} \spl_W(x^*_{\lam, j})
\tau_{\lam,>j}(\gr^W)$ converges to $\spl_W(\br)$.

\smallskip

By (1)--(3), we have $\spl_W(x_\lam)=\tau_{\lam,>j}^{-1} \spl_W(x_{\lam,j}) %
\tau_{\lam,>j}(\gr^W)$ also converges to $\spl_W(\br)$. %

\end{sbpara}

We next prove (3) of Proposition \ref{III3.3.4}.
We assume that $G$ is reductive in \ref{KU6.2.1}--\ref{8.3.d}.

\begin{sbpara}\label{KU6.2.1} 
(Cf.\ \cite{KU2} 6.2.1.) Let $N_j \in  \frak g_\R$ ($1 \leq  j \leq n$) be mutually commuting nilpotent elements, let
$F \in \Dc$, and assume that $(N_1,\dots N_n,F)$ generates a nilpotent orbit. Let $(\rho,\varphi)$ be
the associated $\SL(2)$-orbit in $n$ variables and let $\br = \varphi({\bf i})$ as in \ref{4.5.21}.  Then, from \cite{KU2} 6.1.5 (recall that $G$ is reductive now),
we
see that there are $c_h \in \frak g_\R^{-}, k_h \in \frak g_\R^{+}$  ($h\in \N^n$), where $\frak g^{\pm} =\frak g^{\pm,\br}$ are the 
$(\pm 1)$-eigen subspaces of $\frak g_\R$ under the Cartan involution associated to $K_{\br}$, respectively, such that  $\sum_{h\in \N^n} c_h  \prod_{j=1}^n t_j^{h(j)}$ and   $\sum_{h\in \N^n} k_h  \prod_{j=1}^n t_j^{h(j)}$ ($t_j \in \C$) converge when $|t_j|$
are sufficiently small, that $c_0 =k_0 =0$ and that, for $y_j/y_{j+1} \gg %
0$ ($1\leq j\leq n,y_{n+1} =1$), 
\begin{align*}
&\tau\left(\sqrt{\frac{y_1}{y_2}},\dots,
      \sqrt{\frac{y_n}{y_{n+1}}}\right)\exp\Bigl(\sum_{1\leq j\leq n} iy_j N_j\Bigr)F \tag 1\\
=\; &\exp\Bigl(\sum_{h\in \N^n}  c_h  \prod_{j=1}^n (y_j/y_{j+1})^{-h(j)/2}\Bigr) 
\cdot  \exp\Bigl(\sum_{h\in \N^n}  k_h  \prod_{j=1}^n (y_j/y_{j+1})^{-h(j)/2}\Bigr)\cdot \br.
\end{align*}
Note that these $c_h$ ($h \in  \N^n$) are uniquely determined (whereas the $k_h$ are not).
\end{sbpara}

\begin{sbprop}\label{KU6.2.2} (Cf.\ {\rm \cite{KU2} Proposition 6.2.2.}) 
We use the notation in $\ref{KU6.2.1}$. For $v \in \frak  g_\R$ and $e \in  \Z^n$, let $v(e)$ be the component of $v$ on which $\Ad(\tau(t))$  ($t = (t_j)_{1\leq j\leq n} \in  \bG^n_{m,\R} $) acts as $\prod_{j=1}^n  t_j^{e(j)}$. 
  We denote $|e| = (|e(j)|)_{1\leq j\leq n}$. Let $h \in  \N^n \smallsetminus \{0\}$, $e \in \Z^n$. 

  Then, $c_h(e)=0$ unless $|e|<h$ for the product order in $\N^n$, i.e., $|e(j)|\leq h(j)$ ($1\leq j\leq n$)
and $|e| \neq  h$.

\end{sbprop}

\begin{pf}
This proposition is proved in the same way as \cite{KU2} 6.2.4--6.2.6.
\end{pf}

\begin{sbpara}
The following arguments are  almost copied from the proof \cite{KU2} 6.4.4 of \cite{KU2} Proposition 6.4.1.

We prove the following assertion (C$_j$) by a
downward induction on $j$.
Let $0\le j\le n$, and let %
$\Dc_0:=D$ (this is the correction of a typo \lq\lq$\Dc_0:=\Dc$'' in \cite{KU2} 6.4.4).
\medskip

\noindent
(C$_j$)\quad{\it
Assume that $\exp\big(\ts_{t\in S_{\ge j+1}}
iy_{\lam,t}N_t\big)F_\lam\in\Dc_j$ for any
$\lam$.
Then, for a sufficiently large $\lam$, we have
\begin{align*}\
&\tau\left(\sqrt{\frac{y_{\lam,c_1}}{y_{\lam,c_2}}},\dots,
      \sqrt{\frac{y_{\lam,c_n}}{y_{\lam,c_{n+1}}}}\right)\exp\big(\ts_{s\in S}
iy_{\lam,s}N_s\big)F_\lam\\
=\;&\exp\Big(\ts_{h\in\bN^j}b_{\lam,h}
\tp_{1\le k\le j}\sqrt{\tfrac{y_{\lam,c_k}}
{y_{\lam,c_{k+1}}}}^{\,-h(k)}\Big)k_\lam\cdot\br,
\end{align*}
where $b_{\lam,h}\in\fg_\bR^{-}$ $(\fg_\bR^{-}$
denotes the $(-1)$-eigenspace of $\fg_\bR$ under
the Cartan involution associated to $K_\br)$,
$k_\lam\in K_\br$, and $\ts_{h\in\bN^j}
b_{\lam,h}\tp_{1\le k\le j}x_k^{h(k)}$
$(x_k\in\bC)$ absolutely converges when $|x_k|$
$(1\le k\le j)$ are sufficiently small, which
satisfy the following three conditions.
Here $y_{c_{n+1}}:=1$.
\medskip

\noindent
{\rm(1)} $k_\lam\to1$.
\smallskip

\noindent
{\rm(2)} $|(\Ad\circ\tau_k)_{1\le k\le j}
\text{-weight of $b_{\lam,h}$}|\le h$ for the product order in $\bN^j$.
\smallskip

\noindent
{\rm(3)} For each $h\in\bN^j$,
$\Ad\Big(\tp_{k\ge j+1}\tau_k\Big(
\tsqrt{\tfrac{y_{\lam,c_k}}{y_{\lam,c_{k+1}}}}
\Big)\Big)^\nu(b_{\lam,h})$ converges for
$\nu=0,\pm1$.
Moreover, if $h=0$ then it converges to $0$.
}

\medskip

 Note that by the local structure theorem of $D_{\SL(2)}$ in the reductive group case (Theorem \ref{SL2loc}), Proposition \ref{III3.3.4} (3) is 
equivalent to  (C$_0$). 
\end{sbpara}

\begin{sblem}
\label{l:KU6.4.3}
(Cf.\ {\rm \cite{KU2}} Lemma $6.4.3.$)
Let the notation be as above.
Fix $j$ such that $1\le j\le n$, and assume that
$\exp\big(\ts_{t\in S_{\ge j+1}}iy_{\lam,t}
N_t\big)F_\lam\in\Dc_j$ for any $\lam$.
Put
\begin{align*}
&N_{\lam,k}:=\ts_{s\in S_k}\tfrac{y_{\lam,s}}{y_{\lam,c_k}}N_s\quad(1\le k\le j),\\
&U_\lam:=\Big(\tp_{l\ge j+1}\tau_l\Big(\tsqrt{\tfrac{y_{\lam,c_l}}{y_{\lam,c_{l+1}}}}\Big)\Big)\exp\big(\ts_{t\in S_{\ge j+1}}iy_{\lam,t}N_t\big)F_\lam.
\end{align*}
For each $\lam$, let $(\rho_\lam,\vf_\lam)$ be
the $\SL(2)$-orbit in $j$ variables associated to $((N_{\lam,k})_{1\le k\le j},U_{\lam})$, and let $\br_\lam:=\vf_\lam(\bi)$.
Then, $\br_\lam$ converges to $\br$.
\end{sblem}

  The proof of \cite{KU2} Lemma 6.4.3 also works for this Lemma \ref{l:KU6.4.3}, which uses the important fact that an $\SL(2)$-orbit moves real analytically if the monodromy weight filtration is constant (\cite{CKS}, see also \cite{KU2} 6.1.6).

\begin{sbpara}\label{8.3.d}
  We prove (C$_j$).
By \cite{KNU1} 10.3, %
we have
\begin{align*}
&\tau\left(\sqrt{\frac{y_{\lam,c_1}}{y_{\lam,c_2}}},\dots,
      \sqrt{\frac{y_{\lam,c_n}}{y_{\lam,c_{n+1}}}}\right)\exp\big(
\ts_{s\in S}iy_{\lam,s}N_s\big)F_\lam\\
=\;&\Big(\tp_{k\le j}\tau_k\Big(
\tsqrt{\tfrac{y_{\lam,c_k}}{y_{\lam,c_{k+1}}}}
\Big)\Big)
\exp\big(\ts_{s\in S_{\le j}}
i\tfrac{y_{\lam,s}}{y_{\lam,c_{j+1}}}N_s\big)
U_\lam\\
=\;&\Big(\tp_{k\le j}\tau_k\Big(
\tsqrt{\tfrac{y_{\lam,c_k}}{y_{\lam,c_{k+1}}}}
\Big)\Big)
\exp\big(\ts_{k\le j}i
\tfrac{y_{\lam,c_k}}{y_{\lam,c_{j+1}}}N_{\lam,k}
\big)U_\lam.
\end{align*}
Then, by \ref{KU6.2.1} at $\br_\lam$ (see Lemma \ref{l:KU6.4.3} for $\br_\lam$),
\begin{align*}
&\tau\left(\sqrt{\frac{y_{\lam,c_1}}{y_{\lam,c_2}}},\dots,
      \sqrt{\frac{y_{\lam,c_n}}{y_{\lam,c_{n+1}}}}\right)
\exp(\ts_{s\in S}iy_{\lam,s}N_s)F_\lam\\
=\;&\Big(\tp_{k\le j}\tau_k\Big(
\tsqrt{\tfrac{y_{\lam,c_k}}{y_{\lam,c_{k+1}}}}
\Big)\Big)
\Big(\tp_{k\le j}\tau_{\lam,k}\Big(
\tsqrt{\tfrac{y_{\lam,c_k}}{y_{\lam,c_{k+1}}}}
\Big)\Big)^{-1}f_\lam k_{1,\lam}\cdot\br_\lam,
\quad\text{where}\\
&f_\lam:=\exp\Big(\ts_{h\in\bN^j}a_{\lam,h}
\tp_{k\le j}
\sqrt{\tfrac{y_{\lam,c_k}}{y_{\lam,c_{k+1}}}}
^{\,-h(k)}\Big),\;\;
a_{\lam, h}\in\fg_\bR^{-,\br_\lam},\\
&k_{1,\lam}\in K_{\br_\lam},\;\;
k_{1,\lam}\to1.
\end{align*}
Here $\fg_\bR^{-,\br_\lam}$ denotes the
$(-1)$-eigenspace of $\fg_\bR$ under the Cartan
involution associated to the maximal compact
subgroup $K_{\br_\lam}$ of $G_\bR$.

\medskip

\noindent {\bf Claim 1.} {\it
We can write
\begin{align*}
&\br_\lam=g_\lam k_{2,\lam}\cdot\br,\;\;
\tau_{\lam,k}=\Int(g_\lam)\circ\tau_k,\\
&g_\lam\in (G^\circ)_{W^{(1)},\dots,W^{(j)},\bR},
\;\;k_{2,\lam}\in K_\br,\;\;g_\lam\to1,\;\;
k_{2,\lam}\to1.
\end{align*}
}
\smallskip

  This claim is proved as exactly in the same way as in \cite{KU2}. 
  We obtain a proof of this claim from the corresponding proof in \cite{KU2} by replacing ``$\tilde \rho_{\lambda,k}$, $\tilde \rho_k$'' by $\tau_{\lambda,k}$, $\tau_k$, respectively, and by 
replacing \lq\lq Lemma 6.4.3'' by Lemma \ref{l:KU6.4.3}. 
The key point of the proof is the fact that $\tau_{\lam,k}$ and $\tau_\lam$ are the Borel--Serre liftings at $\br_\lam$ and $\br$, respectively, of the common homomorphism $\bold G_m\to P/P_u$, where $P$ is a $\Q$-parabolic subgroup of $G$ containing $(G^\circ)_{W^{(1)},\dots,W^{(j)},\R}$.

\medskip

By Claim 1, we have
\begin{align*}
&\tau\left(\sqrt{\frac{y_{\lam,c_1}}{y_{\lam,c_2}}},\dots,
      \sqrt{\frac{y_{\lam,c_n}}{y_{\lam,c_{n+1}}}}\right)
\exp\big(\ts_{s\in S}
iy_{\lam,s}N_s\big)F_\lam\tag4\\
=\;&\Int\Big(\tp_{k\le j}\tau_k\Big(
\tsqrt{\tfrac{y_{\lam,c_k}}{y_{\lam,c_{k+1}}}}
\Big)\Big)(g_\lam)\Int(g_\lam)^{-1}(f_\lam)
\Int(g_\lam)^{-1}(k_{1,\lam})k_{2,\lam}\cdot\br.
\end{align*}
Here $\Int(g_\lam)^{-1}(k_{1,\lam})\in K_\br$
and, concerning
$a_{\lam,h}\in\fg_\bR^{-,\br_\lam}$ in the
definition of $f_\lam$, we have
$\Ad(g_\lam)^{-1}(a_{\lam,h})\in\fg_\bR^-
=\fg_\bR^{-,\br}$.
Furthermore, if we write
\begin{align*}
&g_\lam
=\exp\big(\ts_{h\in\bN^j}g_{\lam,-h}\big),\quad
g_{\lam,-h}\in\Lie(G_{W^{(1)},\dots,W^{(j)},\bR}),
\tag5\\
&(\text{$(\Ad\circ\tau_k)_
{1\le k\le j}$-weight of $g_{\lam,-h}$})=-h,\quad
g_{\lam,-h}\to0,
\end{align*}
we have
\begin{align*}
\Int\Big(\tp_{k\le j}\tau_k\Big(
\tsqrt{\tfrac{y_{\lam,c_k}}{y_{\lam,c_{k+1}}}}
\Big)\Big)(g_\lam)
=\exp\Big(\ts_{h\in\bN^j}g_{\lam,-h}\tp_{k\le j}
\sqrt{\tfrac{y_{\lam,c_k}}{y_{\lam,c_{k+1}}}}
^{\,-h(k)}\Big).\tag6
\end{align*}
By Proposition \ref{KU6.2.2},
\begin{align*}
|\text{$(\Ad\circ\tau_k)_{1\le k\le j}$-weight
of $\Ad(g_\lam)^{-1}(a_{\lam,h})$}|<h.\tag7
\end{align*}
 From (4)--(7), we obtain
\begin{align*}
&\tau\left(\sqrt{\frac{y_{\lam,c_1}}{y_{\lam,c_2}}},\dots,
      \sqrt{\frac{y_{\lam,c_n}}{y_{\lam,c_{n+1}}}}\right)
\exp\big(\ts_{s\in S}iy_{\lam,s}N_s\big)
F_\lam\tag8\\
=\;&\exp\Big(\ts_{h\in\bN^j}b_{\lam,h}
\tp_{k\le j}
\sqrt{\tfrac{y_{\lam,c_k}}{y_{\lam,c_{k+1}}}}
^{\,-h(k)}\Big)k_\lam\cdot\br,
\quad\text{where}\\
&k_\lam:=\Int(g_\lam)^{-1}(k_{1,\lam})k_{2,\lam}
\in K_\br,\;\;k_\lam\to1,\\
&b_{\lam,h}\in\fg_\bR^-,\;\;
\text{$b_{\lam,h}$ converges for each $h$},\;\;
b_{\lam,h}(\pm h)\to0,\\
&|\text{$(\Ad\circ\tau_k)
_{1\le k\le j}$-weight of $b_{\lam,h}$}|\le h.
\end{align*}
Here $b_{\lam,h}(\pm h)$ denotes the parts of
$b_{\lam,h}$ of weight $\pm h$ with respect to
$(\Ad\circ\tau_k)_{1\le k\le j}$.

In the case $j=n$, (8) already completes the
proof of (C$_j$).

If $0\le j<n$, it remains to prove (3) of (C$_j$).
This is shown by downward induction on $j$, as
we have mentioned.
First we show the following. 

\medskip

\noindent {\bf Claim 2.} {\it 
To prove $(\text{\rm C}_j)$ $(3)$,
we may assume
$\exp\big(\ts_{t\in S_{\ge j+1}}
iy_{\lam,t}N_t\big)F_\lam\in\Dc_{j+1}$
for any $\lam$.
}

\medskip

The following is a simplification of the proof of the generalization of \cite{KU2} 6.4.4 Claim 2 by an observation that we can take $y^*_{\lam,t}=y_{\lam,t}$ there, so that the condition ibid.\ (11) holds trivially.

We prove Claim 2.
By the condition (5) of \ref{III3.3.3} for $j+1$, there exist
$e\ge0$, $F_\lam^*\in\Dc$
satisfying the following
(9)--(12).
\medskip

\noindent
(9) $\exp\big(\ts_{t\in S_{\ge j+1}}
iy_{\lam,t}N_t\big)F_\lam^*\in\Dc_{j+1}$.

\medskip

\noindent
(10) $y^{2e}_{\lam,c_{j+1}}d(F_\lam,F^*_\lam)
\to0$.
\medskip

\noindent
(11) Let $\a_\lam=\Big(\tp_{k\ge j+1}
\tau_k\Big(\sqrt{
\frac{y_{\lam,c_k}}{y_{\lam,c_{k+1}}}}\Big)\Big)
\exp\big(\ts_{t\in S_{\ge j+1}}
iy_{\lam,t}N_t\big)$.
(Note $U_\lam=\a_\lam F_\lam$.)
Then,

$$
y_{\lam,c_{j+1}}^{-e}\Ad(\a_\lam)\to0
$$
in the space of linear endomorphisms of $\fg_\bC$.
\medskip

\noindent
(12) Let $\b_\lam=\tp_{k\ge j+1}\tau_k\Big(
\sqrt{\frac{y_{\lam,c_k}}{y_{\lam,c_{k+1}}}}
\Big)$.
Then,
$$
y_{\lam,c_{j+1}}^{-e}\Ad(\b_\lam)^\nu\to
0\quad\text{for}\;\;\nu=0,\pm 1
$$
in the space of linear endomorphisms of $\fg_\bC$.
\medskip

Define $b_{\lam,h}^*$ ($h\in\bN^j)$ just as
$b_{\lam,h}$ by replacing 
$F_\lam$ by $F_\lam^*$.
To prove Claim 2, it is sufficient to prove
$$
\Ad(\b_\lam)^\nu(b_{\lam,h})
-\Ad({\b_\lam})^\nu(b_{\lam,h}^*)\to0\quad
\text{for}\;\;\nu=0,\pm1.
$$
The left-hand side of this is equal to
$$
y_{\lam,c_{j+1}}^{-e}
\Ad(\b_\lam)^\nu(y_{\lam,c_{j+1}}^e
(b_{\lam,h}-b^*_{\lam,h})).
$$
Hence, by (12), it is sufficient to prove
$$
y_{\lam,c_{j+1}}^e(b_{\lam,h}-b_{\lam,h}^*)\to0.
$$
Since $b_{\lam,h}$ is a real analytic function in
$((N_{\lam,k})_{1\le k\le j}, U_\lam)$ (\cite{KU2}, 6.1.6), it
is sufficient to prove that 
\begin{align*}
y_{\lam,c_{j+1}}^ed(U_\lam, U_\lam^*)\to0,\tag{13}
\end{align*}
where $U_\lam^*=\a_\lam F_\lam^*$.
By (10), we can write
$$
F_\lam=\exp(x_\lam)F_\lam^*\quad
\text{with}\quad x_\lam\in\fg_\bC,\quad
y^{2e}_{\lam,c_{j+1}}x_\lam\to0.
$$
We have
$$
U_\lam=\a_\lam F_\lam=\exp(\Ad(\a_\lam)
(x_\lam))U_\lam^*.
$$
By (11) and by
$$
y_{\lam,c_{j+1}}^e\Ad(\a_\lam)(x_\lam)
=(y_{\lam,c_{j+1}}^{-e}
\Ad(\a_\lam))(y_{\lam,c_{j+1}}^{2e}x_\lam)\to0,
$$
we obtain (13).
Thus Claim 2 is proved.
\medskip

By Claim 2, we assume $\exp\big(\ts_{t\in S_{\ge j+1}}
iy_{\lam,t}N_t\big)F_\lam\in\Dc_{j+1}$ for any
$\lam$.
Then, by \ref{KU6.2.1} at $\br$,
we have elements $b_{\lam,h'}$ for
$h'\in\bN^{j+1}$ and we have, for $h\in\bN^j$,
\begin{align*}
b_{\lam,h}=\ts_{k=0}^\infty b_{\lam,(h,k)}
\sqrt{\frac{y_{\lam,c_{j+1}}}
{y_{\lam,c_{j+2}}}}^{\,-k}.\tag{14}
\end{align*}
Using (8), applied to $b_{\lam,(h,k)}$ replacing $j$ by
$j+1$, we have
\medskip

\noindent
(15) $|(\Ad\circ\tau_{j+1})$-weight of
$b_{\lam,(h,k)}|\le k$.
If $h=0$, the parts of $b_{\lam,(h,k)}$ of
$(\Ad\circ\tau_{j+1})$-weight $\pm k$
converge to $0$.
\medskip

\noindent
By the hypothesis of induction,
\medskip

\noindent
(16) $\Ad\Big(\tp_{k\ge j+2}\tau_k\Big(
\tsqrt{\tfrac{y_{\lam,c_k}}{y_{\lam,c_{k+1}}}}
\Big)\Big)^\nu(b_{\lam,(h,k)})$ converges for
$\nu=0,\pm1$.
\medskip

\noindent
By (14)--(16), we have (3) of (C$_j$).
\qed
\bigskip

The proof of Proposition \ref{III3.3.4} (3) is completed.

\end{sbpara}

\begin{sbrem}
\label{c_to_IVApp} 
As mentioned at the top of this section, \cite{KU2} and \cite{KNU2} Part III
contain some mistakes which are corrected as in \cite{KNU2} Part IV Appendix.
But still the corrections contain some errors. 
  We review the outline of the mistakes of \cite{KU2} in (1) below, and explain the errors of \cite{KNU2} Part IV Appendix in (2) below. 

(1) The essence of the errors are in 6.4.12 and 7.1.2 (3) of \cite{KU2}, called the errors (1) and (2) in \cite{KNU2} Part IV A.1.1, respectively. 
Both errors base on a misuse of \cite{KU2} Proposition 3.1.6: 
In both situations, the authors thought that they could show the existence of real numbers $y^*_{\lam,t}$ satisfying some conditions including the convergence 
$$y^e_{\lam,s}|y_{\lam,t}-y^*_{\lam,t}|\to 0,$$
where $e$ is a fixed nonnegative integer.  
But, \cite{KU2} Proposition 3.1.6 (applied in the way explained there) gives only those satisfying 
$$ y^e_{\lam,s}(q_{\lam,t} - q_{\lam,t}^*) \to 0, $$
which is not enough.

We use \cite{KU2} Proposition 3.1.6 correctly in \ref{oneN} and Lemma \ref{lem712} in the present paper.
Then, we generalize (and correct, see (2) below) in this section (Section \ref{ss:CKS}) an alternative argument \cite{KNU2} Part IV A.1.3--A.1.9 %
explained in \cite{KNU2} Part IV A.1.2, %
which shows that, as is said in \cite{KNU2} Part IV A.1.11, in fact 
we can take just $y_{\lam,s}$, $y_{\lam,t}$ as 
$y^*_{\lam,s}$, $y^*_{\lam,t}$ satisfying all conditions, that is, the stronger statement including $y^e_{\lam,s}|y_{\lam,t}-y^*_{\lam,t}|=0$.

(2) The correction in \cite{KNU2} Part IV Appendix contains several errors including some critical typos.
  For example, 

1. In A.1.2, \lq\lq problem (1)'' and \lq\lq problem (2)'' should be interchanged. 

2. In A.1.5, the definition of $N_1$,...,$N_{m-1}$ is missing. %

3. In the last line of A.1.5, $X$ should be replaced by some log blowing-up of $\check E_{\sig}$ and $S$ should be the corresponding subspace coming from $\tilde E_{\sig}$ (not $E_{\sig}$; the tilde was missed).

4. In A.1.8, line 3 : $(A_n)$ does not follow from Lemma A.1.6 if $m-1=n$ because, then, Lemma A.1.6 is an empty statement. %

5. In A.1.8, line 4: $j \le n$ should be $j < n$. 

6. In A.1.10 and A.3.4, all \lq\lq $m=0$'' should be \lq\lq $m=1$.'' 

7. In A.3.4, all \lq\lq A.1.6'' should be \lq\lq A.1.7.''

\medskip

  The proof given in this section (Section \ref{ss:CKS}), in particular, in Proposition \ref{prop712} and Theorem \ref{Eopen}, corrects and generalizes, in the present context, that in \cite{KNU2} Part IV A.1.3--A.1.9.

\end{sbrem}


\subsection{Properties of the spaces of nilpotent orbits}
\label{ss:property}
  
    The aim of this Section \ref{ss:property} is to prove the following theorem.
  
\begin{sbthm}\label{t:property} 
  Let $\Sig$ be a weak fan in $\Lie(G')$ (${\rm \ref{fan}}$) and 
  let $\Gamma$ be a  semi-arithmetic subgroup ($\ref{Gamma}$) of $G'(\bQ)$.   

In  $(1)$ and $(2)$ below, we assume that $\Gamma$ is compatible with $\Sig$.  In $(3)$--$(5)$ below, we assume that $\Gamma$ is strongly compatible with $\Sig$.

\medskip

$(1)$ Let $X$ be one of $D^{\sharp}_{\Sig}$, 
$D^{\sharp}_{\Sig, [:]}$, $D^{\sharp}_{\Sig, \val}$. Then the action of $\Gamma$ on $X$ is proper, and the quotient space $\Gamma \bs X$ is Hausdorff.

$(2)$  Let $X$ be as in $(1)$. Assume that $\Gamma$ is torsion-free. 
Then  the action of $\Gamma$ on $X$ is free, and 
the projection $X\to \Gamma \bs X$ is a  local homeomorphism.
 
$(3)$ The quotient space $\Gamma \bs D_{\Sig}$ is Hausdorff.

$(4)$ Assume that  $\Gamma$ is neat.  Then $\Gamma \bs D_{\Sig}$ is a log manifold ({\rm \ref{logmfd}}). In particular, it belongs to $\cB(\log)$.  For each $\sig\in \Sig$, the map $\Gamma(\sig)^{\gp}\bs D_{\sig}\to \Gamma \bs D_{\Sig}$ is locally an isomorphism of log manifolds. 

$(5)$ Assume that $\Gamma$ is neat. Then there is a homeomorphism 
$$(\Gamma \bs D_{\Sig})^{\log} \simeq \Gamma \bs D^{\sharp}_{\Sig}$$
 over $\Gamma \bs D_{\Sig}$. 

\end{sbthm}

  See \ref{scomp} for the compatibility and the strong compatibility of $\Gamma$ and $\Sig$.

\begin{sbrem}\label{redrem3} 
  
 In Theorem \ref{t:property}, we can use  a semi-arithmetic subgroup of $G(\Q)$ (not of $G'(\Q)$)  in the following situation (1) and also in the following situation  (2).

(1)  If either $G$ is semisimple or  the condition (1) in Lemma \ref{pol2} is satisfied, Theorem \ref{t:property} holds for a semi-arithmetic subgroup $\Gamma$ of $G(\Q)$. In fact, $\Gamma \cap G'(\Q)$ is of finite index in $\Gamma$ (see Proposition \ref{A=A*} for the latter case). Hence by \ref{proper} (5), we can replace $\Gamma$ by the semi-arithmetic subgroup $\Gamma\cap G'(\Q)$ of $G'(\Q)$.

 (2)  Assume that $G$ is reductive. Then Theorem \ref{t:property} remains true for a semi-arithmetic subgroup of $G(\Q)$ if we make the following 
 modifications $(*)$ and $(**)$ below. Let $\overline \Gamma$ be the image of $\Gamma$ in $(G/Z)(\Q)$, where $Z$ denotes the center of $G$. 
 
$(*)$ In (2) (resp.\  (4) and (5)), the torsion free (resp.\ neat) property is assumed for $\overline \Gamma$, not for $\Gamma$. 
 
$(**)$ In (1) (resp.\ (2)), the proper action (resp.\ free action) is stated  for $\overline \Gamma$, not for $\Gamma$.

 For the proof, see  \ref{redpf}. 
\end{sbrem}

\begin{sbrem}
\label{r:KP}
  In \cite{KP} Theorem 6.1 and in its proof, the conclusions of Theorem \ref{t:property} except the parts on $X=D^{\sharp}_{\Sig,[:]}, 
D^{\sharp}_{\Sig,\val}$ in (1) and (2) for the Mumford--Tate domain $D$ 
(and its extensions) 
associated to a polarized pure Hodge structure $H$ are proved with $G$ being the (semisimple) Mumford--Tate group associated to $H$ (which is isomorphic to $M/Z$ in \ref{relMT}), $\Gamma$ being an arithmetic subgroup of $G(\bQ)$ which is in the connected component of $G(\bR)$ containing $1$, and $\Sig$ being a fan in $\Lie(G)$ which is strongly compatible with $\Gamma$. 
  The proof in \cite{KP} is a reduction to the case (\cite{KU2} 4.1.1 Theorem A)  
of the extended period domains of classical period domains for pure Hodge structures. 
  Our proof of Theorem \ref{t:property} bases on our studies of the space of Borel--Serre orbits, the space of SL(2)-orbits, and the CKS map, which were not considered in \cite{KP}, and  provides alternative proofs of their results.
  (Minor remarks: Precisely, their $D$ and its extensions are connected components of our period domain $D(G,\bar h_0)$, where $\bar h_0 \colon S_{\C/\R}\to G_\R$ is the homomorphism associated to $H$, and its extensions, respectively, and their results are deduced from Theorem \ref{t:property}. 
  Conversely, under the above assumption, their method can treat $\Gamma \bs D(M,h_0)$ ($h_0$ is as in \ref{relMT}) and its extensions (which are open and closed subspaces of $\Gamma \bs D(G,\bar h_0)$ and its extensions, as is seen by the latter half of \ref{ss:functoriality}). 
  See \ref{KP2} for some details. 
  Additionally, their $\Gamma$ is assumed to be neat, though it is easy to reduce the non-neat case to this case.)

\end{sbrem}

\begin{sbpara}

  The relation of this theorem (Theorem \ref{t:property}) with the results and their proofs in the former parts is as follows. 
      
The above (3), (4), and (5) of  Theorem \ref{t:property}  is the $G$-MHS versions of Theorem 2.5.5, Theorem 2.5.2 plus Theorem 2.5.4, and Theorem 2.5.6, respectively, of 
\cite{KNU2} Part III for the extended period domains of classical period domains (which are the mixed Hodge structure versions of (v), (ii) plus (iv), and  (vi), respectively, of \cite{KU2} 4.1.1 Theorem A for pure Hodge structures). 
  In the case of the extended period domains of classical period domains, these (3)--(5) of Theorem \ref{t:property} 
are proved in Section 4 of \cite{KNU2} Part III. 

  In there, 
also are proved the portion of (1) where 
$X=D^{\sharp}_{\Sig}, D^{\sharp}_{\Sig,\val}$ and 
the portion of (2) where 
$X=D^{\sharp}_{\Sig}, D^{\sharp}_{\Sig,\val}$ and $\Gamma$ is neat.
  The space $X=D^{\sharp}_{\Sig, [:]}$ was not considered in \cite{KNU2} Part III, 
and the portion of (1) where $X=D^{\sharp}_{\Sig, [:]}$ and the portion of (2) 
where $X=D^{\sharp}_{\Sig, [:]}$ and  $\Gamma$ is neat 
follow from \cite{KNU2} Part IV Theorem 6.1.1.  
  
  Most of this Section \ref{ss:property} is devoted to the proof of Theorem \ref{t:property} which is completed in \ref{end(i)}.  
  Almost descriptions in the proof in this section are some abridged versions of arguments in 
\cite{KNU2} Part III Section 4 etc. 
  We describe the main steps in the proof, but often omit the details if the arguments are the same as those in 
\cite{KNU2} Part III Section 4 etc.

\end{sbpara}

We start to explain the proof of Theorem \ref{t:property}.

Let $\Sig$ be a weak fan. 

\begin{sbprop}\label{claim3}
  Let $\sig, \sig' \in \Sig$. 
  Let 
$$w_{\lam}  =(y_{\lam},F_{\lam})\in \{(y, F)\in \sig_\R\times \Dc\;|\; \exp(iy)F\in D\}=(|\torus|_{\sig}\times \Dc)\cap E^{\sharp}_{\sig}$$
and 
$$w'_{\lam}  =(y'_{\lam},F'_{\lam})\in \{(y, F)\in \sig'_\R\times \Dc\;|\; \exp(iy)F\in D\}=(|\torus|_{\sig'}\times \Dc)\cap E^{\sharp}_{\sig'}$$ 
be directed families with the same index set. Let $*$ be $[:]$ or $\val$. Assume that $w_{\la} $ converges to $\alpha$ in $E^{\sharp}_{\sig,*}$, $w'_{\lam}$ converges to   $\alpha'$ in $E^{\sharp}_{\sig',*}$, and 
   $\exp(iy_{\lam})F_{\lam}= \exp(iy'_{\lam})F'_{\lam}$  in $D$  for all $\lam$. 
  Then{\rm :}
  
  $(1)$ The images of $\alpha$ and $\alpha'$  in $D^{\sharp}_{\Sig,*}$  coincide.
  
  $(2)$  $y_{\la}-y_{\la}'$ converges in $\Lie(G'_\R)$. 
  \end{sbprop}

\begin{pf}   In the case $*=\val$, this is a $G$-MHS version of \cite{KNU2} Part III Proposition 4.2.3. The proof of the latter works for the proof of Proposition \ref{claim3} for $*=\val$ and also for $*=[:]$. A key point of the proof is the continuity of the CKS map $D^{\sharp}_{\Sig, *}\to D^I_{\SL(2)}$.
\end{pf}

\begin{sbpara}\label{sigact1}

  Let $\sig \in \Sig$. 
  We consider the continuous actions of $i \sig_\R (\subset \sig_\C$) on $E^{\sharp}_{\sig}$, $E^{\sharp}_{\sig,[:]}$, $E^{\sharp}_{\sig, \val}$. 
 For $a\in \sig_\R$, $ia$ sends the class of $(\tau, b, F)$ in $E^{\sharp}_{\sig}$ (resp.\ $((\sig_j, N_j)_j, b, F)$ in $E^{\sharp}_{\sig,[:]}$, resp.\ $(\tau, V, b, F)$ in $E^{\sharp}_{\sig, \val}$) in the additive presentation to the class of $(\tau, b+a, \exp(-ia)F)$ (resp.\ $((\sig_j, N_j)_j, b+a, \exp(-ia)F)$, resp.\ $(\tau, V, b+a, \exp(-ia)F)$).

It may seem strange that we regard the above action as an action of $ia$ ($a\in \sig_\R$), not of $a$, but we do so because, in the situation \ref{sigact2} below, this action of $i\sig_\R$ on $E^{\sharp}_{\sig}$ will be compatible with the action of $\sig_\C(\supset i \sig_\R)$ on $E_{\sig}$.

  We have $i\sig_{\bR} \bs E^{\sharp}_{\sig} = D^{\sharp}_{\sig}$, $i\sig_{\bR} \bs E^{\sharp}_{\sig,[:]} = D^{\sharp}_{\sig,[:]}$, and $i\sig_{\bR} \bs E^{\sharp}_{\sig,\val} = D^{\sharp}_{\sig,\val}$.

\end{sbpara}

\begin{sbprop}\label{claim2}
  Let $\sig \in \Sig$. 
  Then $E^{\sharp}_{\sig} \to D^{\sharp}_{\sig}$, $E^{\sharp}_{\sig,[:]} \to  D^{\sharp}_{\sig,[:]}$, and $E^{\sharp}_{\sig,\val} \to D^{\sharp}_{\sig,\val}$ are $i\sig_{\bR}$-torsors in 
the category of topological spaces. 
\end{sbprop}

To prove Proposition \ref{claim2}, we use the following Lemma \ref{torsor6} and Lemma \ref{KU7.2.7}.

\begin{sblem}\label{torsor6}
  
Let $H$ be a topological group,  $X$ a topological space, and assume that we have a continuous action $H\times X\to X$. Assume the following {\rm (i)}--{\rm (iii)}.

{\rm (i)} This action is  free set-theoretically. 

{\rm (ii)} This action is proper topologically. 

{\rm (iii)} For each point $x \in X$, there exist a subset $S$ of $X$ which  contains $x$ and an open neighborhood $U$ of $1$ in $H$ such that $U \times  S\to  X$, $(h, s)  \mapsto  hs$, induces a homeomorphism from $U \times S$ onto an open set of $X$.

Then $ X\to H\bs X$ is an $H$-torsor in the category of topological spaces.
\end{sblem}

See \cite{KU2} Lemma 7.3.3 for the proof. 

\begin{sblem}\label{KU7.2.7}
  Assume that a Hausdorff topological group $H$ acts on a Hausdorff topological space $X$ continuously and freely. Let $X'$  be a dense subset of $X$. Then, the following two conditions {\rm (i)} and {\rm (ii)} are equivalent.

{\rm (i)} The action of $H$ on $X$ is proper.

{\rm (ii)} Let $(h_{\la}, x_{\la})_{\la}$  be a directed family of elements of $H\times X'$ such that $(x_\la)_\la$ and $(h_\la x_\la)_\la$  converge in $X$. Then $(h_\la)_\la$  converges in $H$.
\end{sblem}

See \cite{KU2} Lemma 7.2.7 for the proof. 

\begin{sbpara}\label{pfclaim2}  We prove Proposition \ref{claim2}. 
  We  check that the actions of $i \sig_\R$ on $E^{\sharp}_{\sig}$, $E^{\sharp}_{\sig, [:]}$, $E^{\sharp}_{\sig, \val}$ satisfy the conditions (i), (ii), (iii) of Lemma \ref{torsor6}, in Claim 1, Claim 2, Claim 3 below, respectively, which completes the proof.

 {\bf Claim 1.} The actions of $i \sig_\R$ on $E^{\sharp}_{\sig}$,  $E^{\sharp}_{\sig, [:]}$, $E^{\sharp}_{\sig, \val}$ are  free.
 
  Proof of Claim 1. This is a $G$-MHS version of \cite{KNU2} Part III Proposition 4.2.2 (ii). The proof of the latter works for the proof of Claim 1. 
  
  {\bf Claim 2.} The actions of $i\sig_{\bR}$ on  $E^{\sharp}_{\sig}$,  $E^{\sharp}_{\sig, [:]}$, $E^{\sharp}_{\sig, \val}$ are proper.  
  
 Proof of Claim 2.  For $E^{\sharp}_{\sig, [:]}$ and $E^{\sharp}_{\sig, \val}$,  this follows from (2) of Proposition \ref{claim3} by Lemma \ref{KU7.2.7}. 
  The result for $E^{\sharp}_{\sig}$ follows from this by \ref{proper} (3.2).

  Finally we have to check 
  
  {\bf Claim 3.}  For $X= E^{\sharp}_{\sig}$, $E^{\sharp}_{\sig,[:]}$, $E^{\sharp}_{\sig, \val}$ and for each  $x=(q,F') \in X$ , there are a subspace $S$ of 
$X$ passing through $x$ and an open neighborhood $U$ of $0$ in $i\sig_{\bR}$ such that the induced map $U \times S \to X$ is an open immersion. 

Proof of Claim 3. The construction of $S$ at the end of the proof of \cite{KNU2} Part III Proposition 4.4.3 works here. 
\end{sbpara}

\begin{sbprop}\label{claim1}
  Let $\sig \in \Sig$. 
  Then the inclusion maps $D^{\sharp}_{\sig} \to D^{\sharp}_{\Sig}$, $D^{\sharp}_{\sig,[:]} \to D^{\sharp}_{\Sig,[:]}$, and  $D^{\sharp}_{\sig,\val} \to D^{\sharp}_{\Sig,\val}$  are continuous open maps. 
\end{sbprop}

\begin{pf} This is a $G$-MHS version of \cite{KNU2} Part III Theorem 4.3.1. The proof of the latter works for the proof of Proposition \ref{claim1}. 
\end{pf}

\begin{sbprop}\label{shsh}  
  The canonical continuous surjections $D^{\sharp}_{\Sig, \val}\to D^{\sharp}_{\Sig, [:]}$ and $D^{\sharp}_{\Sig, [:]}\to D^{\sharp}_{\Sig}$ are proper.
\end{sbprop}

\begin{pf}
  This is a $G$-MHS version of \cite{KNU2} Part III Proposition 4.3.2, Part IV 4.4.3, and Part IV 4.4.6. 
  By Proposition \ref{claim2} and Proposition \ref{claim1},  this is  reduced to the fact that the maps $E^{\sharp}_{\sig, \val}\to E^{\sharp}_{\sig, [:]}$ and $E^{\sharp}_{\sig, [:]}\to E^{\sharp}_{\sig}$ are proper, and hence to the fact that $|\toric|_{\sig, \val}\to |\toric|_{\sig, [:]}$ and $|\toric|_{\sig, [:]}\to |\toric|_{\sig}$ are proper.
\end{pf}

\begin{sbprop}\label{sharpH} The spaces $D^{\sharp}_{\Sig}$, $D^{\sharp}_{\Sig,[:]}$,  and $D^{\sharp}_{\Sig, \val}$  are  Hausdorff.   (This is the case $\Gamma=\{1\}$ of $(1)$ of Theorem $\ref{t:property}$.)
\end{sbprop}
  
 \begin{pf}  Let $*$ be $[:]$ or $\val$. By Proposition \ref{shsh} and by the case $H=\{1\}$ of (3.2) of \ref{proper}, the Hausdorffness of $D^{\sharp}_{\Sig}$ follows from that of $D^{\sharp}_{\Sig, *}$.

  To see that $D^{\sharp}_{\Sig,*}$ is Hausdorff, by Proposition \ref{claim1}, it is enough to show the following.  

\smallskip

$(1)$ Let $\sig, \sig' \in \Sig$ and let $\beta \in D^{\sharp}_{\sig,*}$ and $\beta' \in D^{\sharp}_{\sig',*}$. 
  Assume that $x_{\lam} \in D$ converges to $\beta$ in $D^{\sharp}_{\sig,*}$ and to $\beta'$ in $D^{\sharp}_{\sig',*}$. 
  Then $\beta=\beta'$ in $D^{\sharp}_{\Sig,*}$.
  
\smallskip
  
  We prove $(1)$.
  By Proposition \ref{claim2}, there exist an open neighborhood $U$ of $\beta$ in $D^{\sharp}_{\sig,*}$ 
(resp.\ $U'$ of $\beta'$ in $D^{\sharp}_{\sig',*}$) and a continuous section $s_{\sig}:U \to E^{\sharp}_{\sig,*}$ 
(resp.\ $s_{\sig'}:U' \to E^{\sharp}_{\sig',*}$) of the projection 
$E^{\sharp}_{\sig,*} \to D^{\sharp}_{\sig,*}$
(resp.\ 
$E^{\sharp}_{\sig',*} \to D^{\sharp}_{\sig',*}$). 
  Let 
$w_{\lam}  =s_{\sig}(x_{\lam})$, $w'_{\lam}  =s_{\sig'}(x_{\lam})$, $\alpha =s_{\sig}(\beta)$, and $\alpha' =s_{\sig'}(\beta')$. 
  Then we can apply Proposition \ref{claim3} (1) and conclude $\beta=\beta'$.
  \end{pf}

\begin{sbpara}\label{compa1}   We prove (1) and (2) of Theorem \ref{t:property}.  Let $\Gamma$ be a semi-arithmetic subgroup of $G'(\bQ)$ which is compatible (\ref{scomp}) with $\Sig$. Under this assumption, the action of $\Gamma$ on $D^I_{\SL(2)}$ is proper by Theorem \ref{SL2gl}. 
  Together with Theorem \ref{CKS1} and the fact that $D^{\sharp}_{\Sig,[:]}$ is Hausdorff by Proposition \ref{sharpH}, we see that the action of $\Gamma$ on $D^{\sharp}_{\Sig,[:]}$ is proper (\ref{proper} (3.1)).
  Hence, again by Proposition \ref{sharpH} and \ref{proper} (3.1), 
the action of $\Gamma$ on $D^{\sharp}_{\Sig, \val}$ is proper.
  Since $D^{\sharp}_{\Sig, [:]}\to D^{\sharp}_{\Sig}$ is proper and surjective, the action of $\Gamma$ on $D^{\sharp}_{\Sig}$ is also proper (\ref{proper} (3.2)).
By \ref{proper} (1), the quotient spaces by these proper actions are Hausdorff. Thus we proved  (1) of Theorem \ref{t:property}.  
  
  We prove  (2) of Theorem \ref{t:property}. Assume that 
 $\Gamma$ is torsion-free.  Then the action of $\Gamma$ on $D^{\sharp}_{\Sig}$ is free. In fact, by using a neat subgroup of $\Gamma$ of finite index, this is reduced to the case where $\Gamma$ is neat. Then it becomes the $G$-MHS version of \cite{KNU2}  Part III Theorem 4.3.5 (i) whose proof also works in the present situation.
By \ref{proper} (2), we have the stated local homeomorphism.

\end{sbpara}
In \ref{sigact2}--\ref{end(i)}, we assume that $\Gamma$ is a semi-arithmetic subgroup of $G'(\Q)$ and is strongly compatible with $\Sig$.

\begin{sbpara}\label{sigact2}

  Let $\sig \in \Sig$. 
  We consider the action 
$$\sig_{\bC} \times E_{\sig} \to E_{\sig};\; (a, (q,F)) \mapsto (\bold e(a)q, \exp(-a)F)$$
of $\sig_\C$ on $E_{\sig}$,  where $a \in \sig_{\bC}$, $q \in \toric_{\sig}$, $F \in \Dc$ and 
$(q,F)\in E_{\sig}$. 
  This is an action in the category of log manifolds (we endow $\sig_\C$ with the trivial log structure). This action is compatible with the action of $i\sig_{\bR}$ on $E^{\sharp}_{\sig}$ in \ref{sigact1}. 
  We have $\sig_{\bC} \bs E_{\sig} = \Gamma(\sig)^{\gp}\bs D_{\sig}$. 
\end{sbpara}

The following Proposition \ref{Ctor} and Proposition \ref{Ctor0} are proved together. 

\begin{sbprop}\label{Ctor} In the category of locally ringed spaces over $\C$ with log structures, $E_{\sig}\to \Gamma(\sig)^{\gp}\bs D_{\sig}$ is a $\sig_\C$-torsor.

\end{sbprop}

\begin{sbprop}\label{Ctor0} 
The space $\Gamma(\sig)^{\gp}\bs D_{\sig}$ is a log manifold. 
\end{sbprop}

\begin{sblem}\label{torsor8} 
 Let $H$ be a complex analytic group, $X$ a log manifold, and assume that we have an action $H\times X\to X$ in the category of log manifolds (we regard $H$ as having the trivial log structure).  Assume the following {\rm (i)}--{\rm (iii)}.
 
{\rm (i)} This action is free set-theoretically.

{\rm (ii)} This action is proper topologically.   

{\rm (iii)} For each point $x \in X$, there exist a log manifold $S$ and a  morphism $\iota: S\to  X$ of log manifolds  whose image contains $x$ and an open neighborhood $U$ of $1$ in $H$ such that $U \times  S \to  X$, $(h, s) \to  h\iota(s)$, induces an isomorphism of log manifolds from $U \times S$ onto an open set of $X$. 

Then{\rm :}

$(1)$  The quotient topological space $H\bs X$ has a unique structure of a log manifold  such that, for an open set $V$ of $H\bs X$, $\cO_{H\bs X}(V)$ (resp.\ $M_{H\bs X}(V)$) is the set of all functions on $V$ whose pullbacks to the inverse image $V '$ of $V $ in $X$ belong to $\cO_X(V')$ (resp.\ $M_X (V')$).  (Here $M_*$ denotes the log structure of $*$.)

$(2)$ $X \to H\bs X$ is an $H$-torsor in the category of log manifolds.

\end{sblem}

See \cite{KU2} Lemma 7.3.3 for the proof. 

\begin{sbpara}\label{pfCtor}  We prove Proposition \ref{Ctor} and Proposition \ref{Ctor0}. 
  By Lemma \ref{torsor8}, it is sufficient to prove the following Claim 1--Claim 3. 
  We prove these claims one by one, which completes the proofs. 

{\bf Claim 1.} The action of $\sig_\C$ on $E_{\sig}$ is free. 

Proof of Claim 1. This is a $G$-MHS version of \cite{KNU2} Part III 4.2.2 (i) whose proof also works as a proof in the present situation.

{\bf Claim 2.} The action of $\sig_\C$ on $E_{\sig}$ is proper.

Proof of Claim 2. This is a $G$-MHS version of  \cite{KNU2} Part III 4.4.2 whose proof also works as a proof in the present situation. (The properness of the action of $i\sig_\R$ on $E^{\sharp}_{\sig}$ in Claim 2 in \ref{pfclaim2} is used in the proof through the projection $E_{\sig}\to E^{\sharp}_{\sig}$.)

{\bf Claim 3.} For each point $x \in E_{\sig}$, there exist a log manifold $S$ and a  morphism $\iota: S\to  E_{\sig}$ of log manifolds  whose image contains $x$ and an open neighborhood $U$ of $0$
in $\sig_\C$ such that $U \times  S \to  E_{\sig}$, $(h, s) \to  h\iota(s)$, induces an isomorphism of log manifolds from $U \times S$ onto an open set of $E_{\sig}$.

Proof of Claim 3. The construction of $S$ at the end of the proof of \cite{KNU2} Part III 4.4.3 works here. 
\end{sbpara}

\begin{sbpara}\label{Dlog1}  We have the special case $(\Gamma(\sig)^{\gp}\bs D_{\sig})^{\log}\simeq \Gamma(\sig)^{\gp}\bs D^{\sharp}_{\sig}$ of (5) of Theorem \ref{t:property} as follows. 

The map $E_{\sig}\to \Gamma(\sig)^{\gp}\bs D_{\sig}$ induces a continuous map $(E_{\sig})^{\log}\to (\Gamma(\sig)^{\gp}\bs D_{\sig})^{\log}$. On the other hand, the 
 map $|\toric|_{\sig}= \Hom(P(\sig), \R_{\geq 0}^{\mult}) \to (\toric_{\sig})^{\log}= \Hom(P(\sig), \R^{\mult}_{\geq 0}\times {\bold S}^1)$ ($P(\sig)$ is as in \ref{toric} and ${\bold S}^1=\{u\in \C^\times\;|\; |u|=1\}$) induces a continuous map $E^{\sharp}_{\sig} \to (E_{\sig})^{\log}$. The composition $E^{\sharp}_{\sig} \to  (E_{\sig})^{\log}\to (\Gamma(\sig)^{\gp}\bs D_{\sig})^{\log}$ induces $\Gamma(\sig)^{\gp}\bs D^{\sharp}_{\sig} \to (\Gamma(\sig)^{\gp}\bs D_{\sig})^{\log}$. The last map is a homeomorphism by Proposition \ref{Ctor} as in \cite{KNU2} Part III  4.4.4. 
\end{sbpara}

\begin{sbprop}\label{fixsig} Assume that $\Gamma$ is neat. Let $x=(\sig, Z)\in D_{\Sig}$,   $\gamma\in \Gamma$, and assume $\gamma x =x$. Then 
 $\gamma\in \Gamma(\sig)^{\gp}$.  

\end{sbprop}

\begin{pf} This is a $G$-MHS version of \cite{KNU2} Part III 4.3.5 (ii), and the proof of the latter works for the proof of Proposition \ref{fixsig}. \end{pf}

\begin{sblem}\label{strlem}

  Let $X$ be a topological space with a continuous action of a discrete group $\Gamma$,   let $Y$ be a set with an action of $\Gamma$, and let $f:X\to Y$ be a $\Gamma$-equivariant surjective map.  Let $\Gamma_0$  be a subgroup of $\Gamma$. We introduce the quotient topologies of $X$ on $\Gamma_0\bs Y$ and on $\Gamma\bs Y$. Let $V$ be an open set of $\Gamma_0\bs Y$ and let $U$ be the inverse image of $V$ in $\Gamma_0\bs X$. We assume moreover the three conditions {\rm (i)}--{\rm (iii)} below. Then,  the map $V \to  \Gamma\bs Y$ is a local homeomorphism.

{\rm (i)} $X \to  \Gamma\bs X$ is a local homeomorphism and $\Gamma\bs X$ is Hausdorff. 

{\rm (ii)} The map $U \to  V$ is proper.

{\rm (iii)} If $x\in X$ and $\gamma\in \Gamma$, and if the images of $\gamma x$ and $x$ in $\Gamma_0\bs Y$ are contained in $V$ and they coincide, then $\gamma \in \Gamma_0$. 
\end{sblem}

\begin{pf} This is \cite{KU2} Lemma 7.4.7. \end{pf}

\begin{sbpara}\label{CsS}
  Assume that $\Gamma$ is neat. 
  We prove that the map $\Gamma(\sig)^{\gp}\bs D_{\sig}\to \Gamma \bs D_{\Sig}$ is a local homeomorphism. 
  
We use Lemma \ref{strlem} for $X =D^{\sharp}_{\Sig}$, $Y =D_{\Sig}$, $\Gamma=\Gamma$, $\Gamma_0 =\Gamma(\sig)^{\gp}$,  $V =\Gamma(\sig)^{\gp}\bs D_{\sig}$ and $U = \Gamma(\sig)^{\gp}\bs D^{\sharp}_{\sig}$. 
  The (1) and  (2)  of Theorem \ref{t:property}  show that the condition (i) in Lemma \ref{strlem} is satisfied. 
  By \ref{Dlog1}, the condition (ii) in Lemma \ref{strlem} is satisfied. 
  Proposition \ref{fixsig}  shows that the condition (iii) in Lemma \ref{strlem} is satisfied.
\end{sbpara}
  
  \begin{sbpara}  We obtain (4) of Theorem \ref{t:property}  by Proposition \ref{Ctor0} and by \ref{CsS}. 
  
   We obtain (5) of Theorem \ref{t:property}  by \ref{Dlog1}  and by \ref{CsS}. 
  
  \end{sbpara} 
  
  \begin{sbprop}\label{fromsh} The map $\Gamma \bs D^{\sharp}_{\Sig} \to \Gamma \bs D_{\Sig}$ is proper and surjective. 
  \end{sbprop} 
  
  \begin{pf} Replacing $\Gamma$ by its 
 neat subgroup of finite index, we may assume that $\Gamma$ is 
 neat. Then this follows from 
   (5) of Theorem \ref{t:property}.    
       \end{pf} 
 
\begin{sbpara}\label{end(i)} By Proposition \ref{fromsh} and by (1) of Theorem \ref{t:property}, we obtain (3) of Theorem \ref{t:property}. 
\end{sbpara}

The proof of Theorem  \ref{t:property} is completed.

\begin{sbprop}\label{t:pro6} Let the assumption be as in $(4)$ of Theorem $\ref{t:property}$. 
  Then we have canonical homeomorphisms 
\begin{align*}
&(\Gamma\bs D_{\Sig})_{[:]}\simeq \Gamma\bs D_{\Sig, [:]}, \quad (\Gamma\bs D_{\Sig})_{\val}\simeq \Gamma\bs D_{\Sig, \val}, \\
&((\Gamma \bs D_{\Sig})^{\log})_{[:]} \simeq \Gamma \bs D^{\sharp}_{\Sig,[:]}, \quad \text{\rm{and} }
((\Gamma \bs D_{\Sig})^{\log})_{\val} \simeq \Gamma \bs D^{\sharp}_{\Sig,\val}.
\end{align*}
Here topologies of the spaces on the right hand side of these homeomorphisms are as in $\ref{rat8}$, $\ref{stval1}$, and $(1)$ of Theorem $\ref{t:property}$.
\end{sbprop}

\begin{pf} This follows from (4) and (5) of Theorem \ref{t:property}.  
\end{pf}

\begin{sbrem}\label{t:pro6rem}
The conclusion of Proposition \ref{t:pro6} holds if $G$ is reductive, $\Gamma$ is a semi-arithmetic subgroup of $G(\Q)$, and if the image of $\Gamma$ in $(G/Z)(\Q)$ is neat, where $Z$ is the center of $G$. See \ref{redpf} for the proof. 
\end{sbrem}

\subsection{Valuative spaces, III} 
\label{ss:val3}

The spaces $S_{[:]}$ of ratios are endowed with new log structures and the associated valuative spaces $S_{[\val]}$ are extensively studied in \cite{KNU2} Part IV Section 4. 
In this section, we review this subject and obtain a space $D^{\sharp}_{\Sig, [\val]}$ over $D^{\sharp}_{\Sig, [:]}$.

\begin{sbpara}\label{[:][val]} 
  Let $E$, $S$ and $M_S$ be as in \ref{rat1}. 
  Let $S_{[:]}$ be the topological space defined in \ref{ratSs} and  \ref{ratSt}.

We review the definition of the {\it new log structure on $S_{[:]}$} (\cite{KNU2} Part IV 4.3.3).
We endow $S_{[:]}$ with the sheaf $O_{S_{[:]}}$ of all $\R$-valued continuous functions.
Assume that we are given a chart $\cS\to M_S$ with $\cS$ being a sharp fs monoid such that $|f(s)| <1$ for any $f\in \cS\smallsetminus \{1\}$ and any $s\in S$.  
Let $\Phi=\{\cS^{(j)}\;|\; 0\leq j\leq n\}$ with 
$\cS=\cS^{(0)}\supsetneq \cS^{(1)}\supsetneq \dots \supsetneq \cS^{(n)}$ 
be as in \ref{rat4}.
Take $q_j\in \cS^{(j-1)}\smallsetminus \cS^{(j)}$ for $1\leq j\leq n$. 
We consider the log structure on $S_{[:]}(\Phi)$ (\ref{rat4}) associated to a chart
$$\N^n \to \cO_{S_{[:]}}\;;\; m\mapsto 
(\prod_{j=1}^{n-1} r(q_{j+1}, q_j)^{m(j)/2})\cdot (-1/\log(|q_n|))^{m(n)/2}.$$
These log structures on $S_{[:]}(\Phi)$ are glued to an fs log structure 
$M_{S_{[:]}}^{\text{new}}\to \cO_{S_{[:]}}$ on $S_{[:]}$ which is independent of any choices. 

We denote by $S_{[\val]}$ the valuative space $(S_{[:]})_{\val}$ associated to $S_{[:]}$ endowed with this new log structure. %
We have a proper and surjective map $S_{[\val]}\to S_{[:]}$  (\cite{KNU2} Part IV Corollary 3.1.10). %
\end{sbpara}

\begin{sbpara}\label{sig[val]}  Let $\Sigma$ be a weak fan in $\Lie(G')$.  
Let $\sig\in\Sig$.

There is a log structure on $E^{\sharp}_{\sig}$ which is the inverse image of that of $E_{\sig}$ but it depends on $\Gamma$ (precisely speaking, on $\Gamma(\sig)^{\gp}$).

We define the {\it new log structure $M^{\text{\rm new}}=M_{D^{\sharp}_{\Sig,[:]}}^{\text{\rm new}}$ of $D^{\sharp}_{\Sig,[:]}$} by using the new log structure $M_{S_{[:]}}^{\text{new}}\to \cO_{S_{[:]}}$ for $S=E^{\sharp}_\sig$ in \ref{[:][val]} as follows.

Let $U$ be an open set of $D^{\sharp}_{\Sig,[:]}$.
An $\R$-valued continuous  function on $U$ belongs to the new log structure of $D^{\sharp}_{\Sig,[:]}$ if its pullback on $E^{\sharp}_{\sig,[:]}$ belongs to the new log structure of $E^{\sharp}_{\sig,[:]}$ for all $\sig\in\Sig$. 

Locally on $D^{\sharp}_{\Sig,[:]}$, for $\sig\in \Sig$, for an open set $U$ of $D^{\sharp}_{\sig,[:]}$, and for a continuous section $s: U \to E^{\sharp}_{\sig,[:]}$ of $E^{\sharp}_{\sig,[:]}\to D^{\sharp}_{\sig, [:]}$ given on $U$, the restriction of this new log structure to $U$ coincides with the inverse image of the new log structure of $E^{\sharp}_{\sig, [:]}$ by $s$. 

Hence the new log structure $M^{\text{\rm new}}$ on $D^{\sharp}_{\Sig, [:]}$ is an fs log structure which is independent of the choice of $\Gamma$.
Denote by $D^{\sharp}_{\Sig, [\val]}$ the valuative space associated to the new log structure $M^{\text{\rm new}}$ of $D^{\sharp}_{\Sig, [:]}$.

\end{sbpara}

\begin{sbpara}\label{CKSval}

Consider the log structure of $D^I_{\SL(2)}$ which is defined in Proposition \ref{str3} by using distance to the boundary.
Similarly as in \cite{KNU2} Part IV 4.5.12, we can prove that the continuous map $D^{\sharp}_{\Sigma,[:]}\to D^I_{\SL(2)}$ (Theorem \ref{CKS1}) respects these log structures. 
  (In \cite{KNU2} Part IV 4.5.12, the reference 4.3.3 is a little ambiguous.
  The precise meaning is as above.)

Thus the map  $D^{\sharp}_{\Sigma,[:]} \to D^I_{\SL(2)}$ induces the continuous map  $D^{\sharp}_{\Sig,\lval}\to D^I_{\SL(2),\val}$ of associated valuative spaces  (cf.\ \cite{KNU2} Part IV 4.5.13).
\end{sbpara}

\begin{sbprop} 
  Let the notation and the assumption be as in $(1)$ and $(2)$ of Theorem $\ref{t:property}$. Then the conclusions of $(1)$ and $(2)$ of Theorem $\ref{t:property}$ and their variants in Remark $\ref{redrem3}$ are true also for  $X=D^{\sharp}_{\Sig, [\val]}$. 
\end{sbprop}

\begin{pf} 
This follows from the corresponding results for $X=D^{\sharp}_{\Sig,[:]}$ because the map $D^{\sharp}_{\Sig,[\val]}\to D^{\sharp}_{\Sig,[:]}$ is separated. 
\end{pf}

\subsection{Mild degeneration}
\label{ss:mild}

We consider the $G$-MHS version of \cite{KNU2} Part IV Section 5.1 in which we studied mild degenerations.
\begin{sbpara}

We define the {\it mild part $D^{\mild}_{\Sig}$ of $D_{\Sig}$} as the part of 
points $(\sig, Z)$ which satisfy the following condition.

(C) For each $N$ in the cone $\sig$, there is an $\R$-splitting of $W$ (which can depend on $N$) that is compatible with $N$.

For the other spaces of nilpotent orbits $D^{\sharp}_{\Sig}$, $D_{\Sig, [:]}$, $D^{\sharp}_{\Sig, [:]}$, $D_{\Sig, \val}$, and so on, we
define their mild parts $D^{\sharp,\mild}_{\Sig}$, $D^{\mild}_{\Sig, [:]}$, $D^{\sharp,\mild}_{\Sig, [:]}$, $D^{\mild}_{\Sig, \val}$, and so on as the inverse images
of $D^{\mild}_{\Sig}$.
\end{sbpara}

\begin{sbpara}
\label{diadef}
Let $D^{\diamond}_{\SL(2)}$ be the subset of $D^{\star,\mild}_{\SL(2)}\times \cL$   consisting 
of all elements $(p, Z, \delta)$ ($(p, Z) \in  D^{\star, \mild}_{\SL(2)}$ with $p \in  D_{\red,\SL(2)}$ and $Z \subset D$ (%
\ref{SL2AB}), $\delta \in \cL=W_{-2}\Lie(G_{u,\R}$)) satisfying the following conditions (i) and (ii).

(i) Let $n$ be the rank of $p$, and let ${\bold 0} := (0,...,0) \in \Z^n$. Then  the $\Ad(\tau_p^{\star})$-weights of $\delta$ are $\leq \bold 0$.

(ii) For every $x \in Z$, $\delta_W(x)$ coincides with the component of $\delta$ of $\Ad(\tau^{\star}_p)$-weight $\bold 0$.

We define the structure of $D^{\diamond}_{\SL(2)}$ as an object of $\cB_{\R}' (\log)$ by the embedding $D^{\diamond}_{\SL(2)}\overset{\subset}\to D^{\star,\mild}_{\SL(2)}\times \cL$.

We regard $D$ as a subspace of $D^{\diamond}_{\SL(2)}$ via the embedding $x\mapsto (x, \delta_W(x))$.

This is the $G$-MHS version of \cite{KNU2} Part IV 5.1.8, 5.1.9.

\end{sbpara}

\begin{sbprop}
The canonical map $D^{\star}_{\SL(2)}\to D_{\red, \SL(2)}\times \spl(W)$ induces a bijection from $D^{\diamond}_{\SL(2)}$ to the subset of $D_{\red,\SL(2)}\times \spl(W) \times \cL$ consisting of $(p, s, \delta)$ satisfying the following conditions {\rm (i)} and {\rm (ii)}. 

{\rm (i)} The $\Ad(\tau^{\star}_p)$-weights of $\delta$ are $\leq \bold 0$.

{\rm (ii)} Let $(\rho, \varphi)$ be an $\SL(2)$-orbit for $G_{\red}$ which represents $p$. Then the component of $\delta$ of $\Ad(\tau^{\star}_p)$-weight $\bold 0$  is of Hodge type $(\leq -1,\leq -1)$  with respect to $\varphi(\bi)$.
\end{sbprop}

This is a $G$-MHS version of \cite{KNU2} Part IV Proposition 5.1.11 and is proved in the same way. 

\begin{sbthm}\label{mildCKS} The identity map of $D$ extends uniquely to a continuous map $D^{\sharp,\mild}_{\Sig,[:]}\to D^{\diamond}_{\SL(2)}$. 
The last map is compatible with the new log structure of $D^{\sharp,\mild}_{\Sig, [:]}$ and the log structure of $D^{\diamond}_{\SL(2)}$, and induces a continuous map between the associated valuative spaces $D^{\sharp,\mild}_{\Sig,[\val]}\to D^{\diamond}_{\SL(2),\val}$. 
\end{sbthm}

This is a $G$-MHS version of \cite{KNU2} Part IV Theorem 5.1.10  and is proved in the same way.

\begin{sbpara}\label{afd2}
For mild degenerations,
we can replace the upper right part of the  fundamental diagram in Introduction by the following 
commutative diagram (maps respect structures of the spaces) which contain the space $D^{\diamond}_{\SL(2)}$ and its associated valuative space $D^{\diamond}_{\SL(2),\val}$.
$$\begin{matrix}
&&D^{\sharp,\mild}_{\Sigma,\lval}&
\overset{\psi}\to  & D^{\diamond}_{\SL(2),\val} & \to &D^{\star,\mild}_{\SL(2),\val}&\overset{\eta^{\star}}{\underset{\subset}\to} & D^{\mild}_{\BS, \val}\\
&&\downarrow &&\downarrow&&\downarrow &&\downarrow \\
&&D^{\sharp,\mild}_{\Sig,[:]} &\overset{\psi}\to & D^{\diamond}_{\SL(2)} &\to & D^{\star,\mild}_{\SL(2)}&& D^{\mild}_{\BS}\\
&&\downarrow&&&&\downarrow&&\\
\Gamma\bs D^{\mild}_{\Sigma}&\leftarrow&D^{\sharp,\mild}_{\Sigma}&&&&D_{\SL(2)}&&
\end{matrix}$$
\end{sbpara}

\begin{sbprop} 

The conclusions of Theorem $\ref{SL2gl}$ and their variants in $\ref{SL2SA1}$ and Proposition $\ref{SL2SA2}$ are true also for $X=D^{\diamond}_{\SL(2)}$.

\end{sbprop}

\begin{pf} This follows from the corresponding results for $X=D^{\star}_{\SL(2)}$ because the map $D^{\diamond}_{\SL(2)}\to D^{\star}_{\SL(2)}$ is continuous and $D^{\diamond}_{\SL(2)}$ is Hausdorff. 
\end{pf}

\subsection{The fundamental diagram in examples}\label{ss:fund}

We explain what our fundamental diagram tells in special examples. 

In particular, by using the fundamental diagram, we give a complement to the work of Goresky and Tai \cite{GT} on the relation of the toroidal compactification and the reductive Borel--Serre compactification (see \ref{Shimu}). 

In \ref{Ex1}, \ref{Ex2}, \ref{Shimu}, \ref{noSLBS}, we consider cases in which $G$ is reductive. For a reductive $G$, the main part of the fundamental diagram in Introduction  becomes 
$$\begin{matrix}     & & D^{\sharp}_{\Sig, [\val]} & \overset{\psi}\to & D_{\SL(2),\val} & \overset{\eta}\to & D_{\BS,\val}\\
&& \downarrow && \downarrow && \downarrow\\
&& D^{\sharp}_{\Sig,[:]} & \overset{\psi}\to & D_{\SL(2)} && D_{\BS}\\
&& \downarrow &&&&\\
\Gamma \bs D_{\Sig} &\leftarrow & D^{\sharp}_{\Sig} &&&&
\end{matrix},$$
and we consider this part. 

\begin{sbpara}\label{Ex1} {\it Example.} Let $G=\GL(2)$ and  $h: S_{\C/\R}\to G$ be the standard one (cf.\ \ref{clEx}).

Then $D=D(G,h)$ is the complex analytic manifold $ \frak H^{\pm}=\{\tau\;|\; \text{Im}(\tau)\neq 0\}=\C\smallsetminus \R$ on which $G(\R)=\GL(2,\R)$ acts naturally.

We describe the space  $D_{\BS}$ of Borel--Serre orbits. Let
$$P=\begin{pmatrix} *& *\\ 0&*\end{pmatrix} \subset \GL(2)_\Q=G.$$ 
All parabolic subgroups of $G$ other than $G$ are conjugate to $P$ under $G(\Q)$, and hence the whole $D_{\BS}$ can be understood as the union of the following picture of the open set $D_{\BS}(P)$. 
Let $P_1:=P\cap G'$. 
For $iy\in D$ ($y\in \R^\times$), the Borel--Serre lifting of $S_{P_1}= P_1/P_{1,u}$ is the group of matrices 
$\begin{pmatrix} a&0\\ 0& b \end{pmatrix}$ such that $ab=1$. The adjoint action of this matrix on $\begin{pmatrix} 0&1 \\0& 0\end{pmatrix} \in \Lie(P_{1,_u})$ is the multiplication by $ab^{-1}$, and hence the fundamental root sends this matrix to $a^{-1}b$. Hence the isomorphism $A_P\simeq \R_{>0}$ given by the fundamental root sends $r\in \R_{>0}$ to the matrix $\begin{pmatrix} 1/\sqrt{r} & 0\\ 0 & \sqrt{r}\end{pmatrix}$. Hence the Borel--Serre action of $r\in \R_{>0}\simeq A_P$ sends $x+iy$ to $x+ir^{-1}y$. From this, we have a commutative diagram of spaces 
$$\begin{matrix}   D  &\simeq & \R \times  \R_{>0}\times \{\pm 1\}\\ \cap&& \cap\\ D_{\BS}(P) &\simeq & \R\times \R_{\geq 0}\times \{\pm 1\},\end{matrix}$$
in which the upper horizontal arrow sends $x+iy$ ($x\in \R, y\in \R^\times$) to $(x, 1/|y|, \text{sgn}(y))$ and the lower isomorphism preserves the structure of real analytic manifolds with corners. We have $D_{\BS,\val}=D_{\BS}$. 

Next we describe the space $D_{\SL(2)}$ of $\SL(2)$-orbits. Let $W' \in \frak W(G)$ be the filtration associated to the homomorphism 
$$\alpha: {\bf G}_m \to G\; ;\; t \mapsto \begin{pmatrix} 1/t & 0\\0 & t\end{pmatrix},$$
and let $\Phi=\{W'\}$. Since $D_{\SL(2)}= \bigcup_{g\in G(\Q)}\; gD_{\SL(2)}(\Phi)$, the whole $D_{\SL(2)}$ can be understood as the union of the following picture of the open set $D_{\SL(2)}(\Phi)$. 
We have the distance $\beta: D\to \R_{>0}$ to $\Phi$-boundary defined as $x+iy\mapsto 1/\sqrt{|y|}$. Then the injective real analytic map  
$$\nu_{\alpha, \beta}: D_{\SL(2)}(\Phi) \overset{\subset}\to \R_{\geq 0} \times D \times \spl(W')$$
in Proposition \ref{emb1} sends $x+iy$ ($x\in \R, y\in \R^\times$) to $(1/\sqrt{|y|}, x|y|^{-1}+i\cdot \text{sgn}(y), x)$, where we identify $\spl(W')$ with $\R$ by sending $x\in \R$ to $\begin{pmatrix} 1 & x\\0&1\end{pmatrix}\cdot s$ with $s\in \spl(W')$  given by $\alpha$.  
The closure of $\nu_{\alpha, \beta}(D)$ in the target space is $C=\{(r, z, x)\;|\; z= xr^2\pm i\}$. Since $D$ is dense in $D_{\SL(2)}(\Phi)$, $\nu_{\alpha, \beta}$ induces an injective map $D_{\SL(2)}(\Phi)\to C$. Let $p\in D_{\SL(2)}(\Phi)$ be the class of the $\SL(2)$-orbit in \ref{clEx2}. Then $\nu_{\alpha, \beta}$ sends $\begin{pmatrix}\epsilon & x\\0 & 1 \end{pmatrix}\cdot p\in D_{\SL(2)}(\Phi)$ for $\epsilon \in \{\pm 1\}$  to $(0, i\epsilon, x)$. This proves the surjectivity of $D_{\SL(2)}(\Phi)\to C$ and hence  we have an isomorphism $D_{\SL(2)}(\Phi) \overset{\sim}\to C$  of objects $\cB'_\R(\log)$. 
From this, we have a commutative diagram of spaces 
$$\begin{matrix}   D  &\simeq & \R \times  \R_{>0}\times \{\pm 1\}\\ \cap&& \cap\\ D_{\SL(2)}(\Phi) &\simeq & \R\times \R_{\geq 0}\times \{\pm 1\},\end{matrix}$$
in which the upper horizontal arrow sends $x+iy$ ($x\in \R, y\in \R^\times$) to $(x, 1/\sqrt{|y|}, \text{sgn}(y))$ and the lower isomorphism preserves the real analytic structure and the log structure with sign.  The image of $p\in D_{\SL(2)}(\Phi)$ under the lower horizontal arrow is $(0, 0, 1)$. We have $D_{\SL(2),\val}=D_{\SL(2)}$. 

By Theorem \ref{Shim}, the identity map of $D$ extends uniquely to a morphism 
$D_{\SL(2)}\to D_{\BS}$ in $\cB'_\R(\log)$. This induces $D_{\SL(2)}(\Phi) \to D_{\BS}(P)$ for which the following diagram is commutative.
$$\begin{matrix}   D_{\SL(2)}(\Phi) & \overset{\sim}\to & \R \times \R_{\geq 0} \times \{\pm 1\}\\
\downarrow && \downarrow\\
D_{\BS}(P) &\overset{\sim}\to & \R \times\R_{\geq 0}\times \{\pm 1\}\end{matrix}$$ Here the right vertical arrow is $(x,r, \epsilon)\mapsto (x, r^2, \epsilon)$. Thus the map $D_{\SL(2)}\to D_{\BS}$ is a homeomorphism, but  their real analytic structures are slightly different.

Next we describe the spaces of nilpotent orbits for the fan $\Sig$ consisting of all cones $\R_{\geq 0}N$ with $N\in \Lie(G')=\frak{sl}(2, \Q)$ such that $N^2=0$ as a $(2,2)$-matrix. Then for a congruence subgroup $\Gamma$ of $\SL(2, \Z)$, $\Gamma \bs D_{\Sig}$ is a compactified  modular curve. We have $D^{\sharp}_{\Sig, [\val]}=D^{\sharp}_{\Sig, [:]}= D^{\sharp}_{\Sig}$, and 
 the CKS map induces a homeomorphism $D^{\sharp}_{\Sig}=D^{\sharp}_{\Sig, [:]}\overset{\sim}\to D_{\SL(2)}$. 

 Define $\sig,\tau\in \Sig$ by   $$\sig:=\R_{\geq 0}N, \quad \tau:= \R_{\geq 0}(-N)\quad \text{with}\;\; N=\begin{pmatrix} 0&1\\ 0& 0\end{pmatrix}.$$ The CKS map induces a homeomorphism $D^{\sharp}_\sig \cup D^{\sharp}_\tau \overset{\sim}\to D_{\SL(2)}(\Phi)$. For $a\in \R$, this map sends the nilpotent $i$-orbit $(\sig, a+i\R) \in D^{\sharp}_{\sig}$ (resp.\ $(\tau, a+i\R)\in D^{\sharp}_\tau$) to the element of $D_{\SL(2)}(\Phi)$ corresponding to the element $(a, 0, 1)$ (resp.\ $(a,0,-1)$) of $\R\times \R_{\geq 0}\times \{\pm 1\}$.
Here $a+i\R$ is regarded as a subset of $\Dc$ by identifying $\Dc$ with ${\bf P}^1(\C)\supset a+i\R$.

 The part about $D_{\BS}$ and $D_{\SL(2)}$ as topological spaces in this \ref{Ex1} is essentially described in \cite{KU1} 6.2.

\end{sbpara}

\begin{sbpara}\label{Ex3} {\it Example.}
  Let $G= \begin{pmatrix} \bG_m& \bG_a\\0& 1\end{pmatrix}$ and let $h : S_{\C/\R}\to G_{\red,\R}={\bf G}_{m,\R}$ be the homomorphism which induces $S_{\C/\R}(\R)=\C^\times\to G_{\red}(\R)\::\:  z\mapsto  z^{-1}\bar z^{-1}$ ($z\in \C^\times$). 
Let $x\in D=D(G, h)$ be the point defined by the homomorphism $S_{\C/\R} \to G_{\R}$ which induces 
$$S_{\C/\R}(\R)=\C^{\times}\to G(\R) = \begin{pmatrix}  \R^\times& \R\\ 0& 1 \end{pmatrix}\;;\; z\mapsto \begin{pmatrix} z^{-1} \bar z^{-1}&0\\0&1\end{pmatrix}
\quad (z\in \C^\times).
$$ Then we have an isomorphism of complex analytic manifolds
$$\C \overset{\sim}\to D\;;\; c \mapsto \begin{pmatrix} 1& c \\0 & 1\end{pmatrix}\cdot x.$$
 Let $p$ be the unique element of $D_{\red}$. The real analytic isomorphism $D \overset{\sim}\to \spl(W)\times \cL(p)$ in Proposition \ref{Dandgr} sends $c=a+ib\in \C=D$ ($a,b\in \R$) to $(s(a), \delta(b))$, where 
 $$s(a)= \begin{pmatrix} 1& a\\ 0 & 1\end{pmatrix} \spl_W(x), \quad \delta(b)=  \begin{pmatrix} 0& b \\ 0& 0\end{pmatrix} \in W_{-2}\Lie(G_\R)=\gr^W_{-2}\Lie(G_\R).$$

Let $P$ be the unique parabolic subgroup of $G_{\red}$, that is,  $G_{\red}$ itself. 
 Since $A_P=\{1\}$, $B_P=\R_{>0}$. By \ref{Bact}, the Borel--Serre action of $t\in B_P$ on $D= \C$ sends $c=a+ib$ ($a,b\in \R$) to $a+it^{-2}b$. 
As topological spaces with sheaves of real analytic spaces and with log  structures with sign, all  the spaces $D_{\BS}$, $D_{\BS, \val}$, $D^{\star}_{\SL(2)}$, $D^{\star}_{\SL(2),\val}$, $D^{II}_{\SL(2)}$, $D^{II}_{\SL(2), \val}$, $D^I_{\SL(2)}$, $D^I_{\SL(2),\val}$ coincide with the real analytic manifold with corners $\R \times [-\infty, \infty] \supset \R \times \R \simeq D$,  where $\R \times \R \simeq D$ sends $(a,b) \in \R\times \R$ to $a+ib\in \C=D$.   For $a\in \R$, $(a, \infty)\in \R \times [-\infty, \infty]$ corresponds to the element $(p, a+i\R_{>0})$ of $D_{\SL(2)}$ and to the element $(P, a+i\R_{>0})$ of $D_{\BS}$, and $(a, -\infty)\in \R \times [-\infty, \infty]$ corresponds to the element $(p, a+i\R_{<0})$ of $D_{\SL(2)}$ and to the element $(P, a+i\R_{<0})$ of $D_{\BS}$.

Define $$\sig:=\R_{\geq 0}N, \quad \tau:= \R_{\geq 0}(-N)\quad \text{with}\;\; N=\begin{pmatrix} 0&1\\ 0& 0\end{pmatrix}\in \Lie(G'),$$
 $$ \Sig:= \{\{0\}, \sig, \tau\},\quad \Gamma:=  \begin{pmatrix} 1 & \Z\\ 0 & 1\end{pmatrix}.$$
Then $\Gamma$ and $\Sig$ are strongly compatible. We have isomorphisms of complex analytic manifolds 
$$\Gamma \bs D\simeq \C^\times, \quad \Gamma \bs D_{\Sig} \simeq {\bf P}^1(\C),$$
where the class of  $c\in \C=D$ in $\Gamma \bs D$ is identified with $\exp(2\pi ic)\in \C^\times$. The point $0\in {\bf P}^1(\C)$ corresponds to the class of the nilpotent orbit $(\sig, \C)$ and the point $\infty\in {\bf P}^1(\C)$ corresponds to the class of the nilpotent orbit $(\tau, \C)$. 
We have $D^{\sharp}_{\Sig, [\val]}= D^{\sharp}_{\Sig,[:]}=D^{\sharp}_{\Sig}$, and the CKS map induces  
a homeomorphism $D^{\sharp}_{\Sig}= D^{\sharp}_{\Sig, [:]}\overset{\sim}\to D^I_{\SL(2)}$. For $a\in \R$, this map sends the  nilpotent $i$-orbit $(\sig, a+i\R)\in D^{\sharp}_{\sig}$ to the element of $D_{\SL(2)}$ corresponding to $(a, \infty)\in \R \times [-\infty, \infty]$, and the nilpotent $i$-orbit $(\tau, a+i\R)\in D^{\sharp}_{\tau}$  to the element of $D_{\SL(2)}$ corresponding to $(a, -\infty)\in \R \times [-\infty, \infty]$. Here $a+i\R$ is regarded as a subset of $\Dc=D=\C$.

\end{sbpara}

\begin{sbpara}\label{Ex2} {\it Example.} We omit the details in this \ref{Ex2}. 

Let $G=\GL(2) \times \GL(2)$ and let $h:S_{\C/\R}\to \GL(2)_\R \times \GL(2)_\R$ be the diagonal embedding of the homomorphism in \ref{Ex1}. This $h$  is $\bR$-polarizable.
We have
$$D=\frak H^{\pm} \times \frak H^{\pm}.$$
Let $\Sig$ be the fan consisting of all cones of the form $\R_{\geq 0}N_1\times \R_{\geq 0}N_2\subset \Lie(G'_\R)={\frak {sl}}(2, \R) \times {\frak {sl}}(2, \R)$, where $N_j$ is an element of 
the $j$-th ${\frak {sl}}(2,\Q)$ in $\Lie(G')$ such that $N_j^2=0$ as a $(2,2)$-matrix ($j=1,2$). 
Then in the category of topological spaces, the fundamental diagram presents 
$$\begin{matrix} D^{\sharp}_{\Sig, [\val]} & \overset{\sim}\to & D_{\SL(2),\val} &\overset{\sim}\to & D_{\BS,\val}\\
\downarrow && \downarrow && \\
D^{\sharp}_{\Sig, [:]} & \overset{\sim}\to & D_{\SL(2)}&&\downarrow\\
\downarrow && &\searrow  &\\
D^{\sharp}_{\Sig} && \overset{\sim}\longrightarrow && D_{\BS}.\end{matrix}$$
Furthermore, $D^{\sharp}_{\Sig}$ is canonically isomorphic to the product of two copies of $D^{\sharp}_{\Sig}$ of \ref{Ex1}, and $D_{\BS}$ is canonically isomorphic to the product of two copies of $D_{\BS}$ of \ref{Ex1} (see Proposition \ref{prodD}).

The proper surjective map $D^{\sharp}_{\Sig, [:]} \to D^{\sharp}_{\Sig}$ is not injective as the following property of the convergences show, and hence the map $D_{\SL(2)}\to D_{\BS}$ is not injective. 

Consider the point $p=(iy_1, iy_2)\in \frak H \times \frak H\subset D$ ($y_1, y_2\in \R_{>0}$). In the following (1), (2), (3), we give examples of the convergence of $p$ in $D^{\sharp}
_{\Sig}$, $D^{\sharp}_{\Sig, [:]}$, and $D^{\sharp}_{\Sig, [\val]}$, respectively, to show how the topologies of these three spaces are different.

(1) If  $y_j\to \infty$ for $j=1,2$, $p$ converges to the class of the nilpotent $i$-orbit $(\sig, Z)$, where $\sig=\R_{\geq 0}N_1\times \R_{\geq 0}N_2$ with $N_j$  the matrix $\begin{pmatrix} 0&1\\0&0\end{pmatrix}$ in the $j$-th ${\frak {sl}}(2, \Q)$ in $\Lie(G')$ ($j=1,2$) and $Z$ is the $\exp(i\sig_\R)$ orbit which passes $(i, i)\in D$. 

(2) If  $y_1, y_2 \to \infty$ and $y_1/y_2\to \infty$, then $p$ converges in $D^{\sharp}_{\Sig,[:]}$ to a point $a$, and if $y_1,y_2\to \infty$ and $y_1/y_2\to 1$, then $p$ converges in $D^{\sharp}_{\Sig,[:]}$ to a point $b$, and $a\neq b$. These $a,b\in D^{\sharp}_{\Sig, [:]}$ lie over the above class of $(\sig, Z)$ in $D^{\sharp}_{\Sig}$.

(3)  If $y_2 \to \infty$ and $y_1/y_2^2\to \infty$, then $p$ converges in $D^{\sharp}_{\Sig, [\val]}$ to a point  $c$, and if $y_2\to \infty$ and $y_1/y_2^3\to \infty$, $p$ converges in $D^{\sharp}_{\Sig, [:]}$ to a point $d$, and $c\neq d$. These $c,d\in D^{\sharp}_{\Sig, [\val]}$ lie over the point $a\in D^{\sharp}_{\Sig, [:]}$.

\end{sbpara}

\begin{sbpara}\label{Shimu}  {\it Shimura varieties.} 

Assume that $G$ is reductive and that $h_0:S_{\C/\R}\to G_\R$ satisfies the condition that the Hodge type of $\Lie(G_\R)$ via $h_0$ is in $\{(1,-1), (0, 0), (-1,1)\}$ (as in \ref{Shi0}). 

Then for an arithmetic subgroup (that is, a subgroup satisfying the condition $(A)$  in \ref{Gamma}) $\Gamma$ of $G'(\Q)$, there is a  fan $\Sig$ which is strongly compatible with $\Gamma$ such that   $\Gamma \bs D_{\Sig}$ is compact. This compact space is  called a Mumford (or toroidal) compactification of $\Gamma \bs D$. 

As in Theorem \ref{Shim}, we have a morphism $D_{\SL(2)}\to D_{\BS}$ which extends the identity morphism of $D$. We can prove $D_{\SL(2),\val}\overset{\sim}\to D_{\BS,\val}$, but we do not give the proof here.

As an application of the fundamental diagram, we have the following complement to the work \cite{GT} of Goresky and Tai on the relation of toroidal compactifications and reductive Borel--Serre compactification. 
 
 The reductive Borel--Serre space, which we denote by $D_{\BS}^{\flat}$ here, is defined to be the quotient of $D_{\BS}$ by the following equivalence relation. 
For $p_1=(P_1, Z_1), p_2=(P_2, Z_2)\in D_{\BS}$, $p_1\sim p_2$ if and only if $P_1=P_2$ and $P_{1,u}Z_1=P_{2,u}Z_2$, where $(\cdot)_u$ denotes the unipotent radical. %
The quotient  $\Gamma\bs D^{\flat}_{\BS}$ is compact, and is called the reductive Borel--Serre compactification of $\Gamma \bs D$.

Let $\Gamma$ be a neat arithmetic subgroup of $G'(\Q)$ and let $\Sig$ be strongly compactible with $\Gamma$ such that $\Gamma\bs D_{\Sig}$ is compact. 

Then concerning the relation of the compactifications $\Gamma \bs D_{\Sig}$ and $\Gamma \bs D^{\flat}_{\BS}$ of $\Gamma \bs D$, Goresky and Tai obtained the following result. The identity map of $\Gamma \bs D$ extends to a \lq\lq continuous map modulo homotopy'' from $\Gamma \bs D_{\Sig}\to \Gamma \bs D^{\flat}_{\BS}$ if we replace $\Sig$ by a sufficiently finer subdivision.  Precisely speaking, if we replace $\Sig$ by a sufficiently finer subdivision, there are a compact topological space $T$ which contains $\Gamma \bs D$ as a dense open subspace,  and  continuous surjective maps $f: T\to \Gamma\bs D_{\Sig}$ and  $g: T\to \Gamma \bs D^{\flat}_{\BS}$ such that $f$ is a homotopy equivalence and such that $f$ and $g$ induce the identity map of $\Gamma \bs D$. 

We have the following result. %
  The map $\Gamma \bs D_{\Sig, [:]}\to \Gamma \bs D_{\Sig}$ is proper surjective and a weak homotopy equivalence, and we have the continuous surjective  map $\Gamma \bs D_{\Sig, [:]}\to \Gamma \bs D^{\flat}_{\BS}$ induced by the continuous maps $D^{\sharp}_{\Sig, [:]}\to D_{\SL(2)}\to D_{\BS}$ by passing to the quotients. That is, compared to \cite{GT}, 
we do not need a subdivision of $\Sig$ here and we present a standard space $\Gamma \bs D_{\Sig, [:]}$ which connects $\Gamma \bs D_{\Sig}$ and $\Gamma \bs D^{\flat}_{\BS}$. 
  This gives an alternative proof for the existence of the canonical maps 
$H^m(\Gamma \bs D^{\flat}_{\BS},A) \to H^m(\Gamma \bs D_{\Sig}, A)$ 
with $A$ being an abelian group for any $m$, $A$ being a group for $m=1$, and $A$ being a set for $m=0$, obtained by Goresky--Tai. 
  We plan to discuss the details in a forthcoming paper.
\end{sbpara}
 
 \begin{sbpara}\label{noSLBS} An example of $D=D(G,h)$ with $G$ reductive  for which the identity map of $D$ does not extend to a continuous map $D_{\SL(2)}\to D_{\BS}$ is given in \ref{Ex6} basing on \cite{KU1} and \cite{KU2}.

 \end{sbpara}

\begin{sbrem} The spaces $D^{\sharp}_{\Sig, \val}$ and $D^{\sharp}_{\Sig, [:]}$ in the fundamental diagram play similar roles in our work. The former appears in \cite{KU2} and in \cite{KNU2} Part II, Part III and so on of our series of papers, and the latter appears in \cite{KNU2} Part IV and in this part. 
  Both have canonical proper continuous maps to $D^{\sharp}_{\Sig}$ and also  continuous maps (the CKS maps) to $D_{\SL(2)}$.  In this Part V , we are using $D^{\sharp}_{\Sig, [:]}$ (and the related space $\Gamma \bs D_{\Sig, [:]}$) more than $D^{\sharp}_{\Sig, \val}$ (and the related space $\Gamma\bs D_{\Sig,\val}$). The advantages of the space $D^{\sharp}_{\Sig, [:]}$ and $\Gamma \bs D_{\Sig, [:]}$ are: 

(1) The space   $D^{\sharp}_{\Sig, [:]}$ naturally produces the space $D^{\sharp}_{\Sig, [\val]}$ from which we can go to $D^I_{\SL(2),\val}$ (and in the case $G$ is reductive, to $D_{\BS,\val}$ and $D_{\BS}$).
  
(2) The definition of the CKS map by using $D^{\sharp}_{\Sig, [:]}$ is simpler and more natural than that by using $D^{\sharp}_{\Sig, \val}$.
The convergences of ratios such as $y_j/y_{j+1}\to a_j$ ($a_j\in \R_{>0}$), $y_j/y_{j+1}\to \infty$ appear in SL(2)-orbit theorem and these are the convergences in the space of ratios, and so the relation to the SL(2)-orbit theorem of $D^{\sharp}_{\Sig, [:]}$ seems stronger than that of $D^{\sharp}_{\Sig, \val}$. 

(3) The spaces $D^{\sharp}_{\Sig, [:]}$ and $\Gamma \bs D_{\Sig, [:]}$ have the application described in \ref{Shimu}.

On the other hand, the space $\Gamma \bs D_{\Sig, \val}$ has the sheaf of holomorphic functions which the space $\Gamma \bs D_{\Sig,[:]}$ does not have.

\end{sbrem}

\subsection{Functoriality in $G$}
\label{ss:functoriality}

\begin{sbpara}\label{G1G2} 
  Assume that we are given a homomorphism $f: G_1\to G_2$.
  We describe how we can relate the period domains and extended period domains for $G_1$ and those for $G_2$.

  We have to introduce some conditions for this functoriality. 
 
 We assume $f(G_{1,u})\subset G_{2,u}$ and hence $f$ induces $f_{\red}: G_{1,\red}\to G_{2,\red}$. 
We assume that we are given $h_1: S_{\C/\R}\to G_{1,\red, \R}$  as in  \ref{D}, and we further assume  that the induced homomorphism $h_2:=f_{\red}\circ h_1: S_{\C/\R}\to G_{2, \red, \R}$ satisfies  the condition 
 that the composition $\bG_{m,\R}\to S_{\C/\R} \overset{h_2}\longrightarrow   G_{2, \red, \R}$ is central. 

Then we have holomorphic maps 
$$D(G_1,h_1)\to D(G_2, h_2), \quad D(G_{1,\red},h_1)\to D(G_{2,\red},h_2).$$

 In \ref{ft1}--\ref{ft3}, we further assume that $h_1$ and $h_2$ are $\R$-polarizable.

\end{sbpara}

\begin{sbpara}\label{ft1} We first consider the spaces of Borel--Serre orbits.

Assume the the following condition (i) is satisfied.

\medskip

(i) The map $\Lie(G'_1)\to \Lie(G'_2)$ is surjective.

\medskip

Here as usual, $(-)'$ denote the commutator groups. Then we have 
a morphism of real analytic manifolds with corners
 $$D(G_1, h_1)^{\mild}_{\BS} \to D(G_2,h_2)^{\mild}_{\BS}\;;\; (P_1, Z_1)\mapsto (P_2, Z_2),$$
 where $P_1$ is a parabolic subgroup of $G_{1,\red}$, $Z_1$ is an $A_{P_1}$-orbit in $D(G_1,h_1)$, $P_2$ is the algebraic subgroup of $G_2$ generated by the image of $P_1$ and  the center of $G_2^{\circ}$, which is a parabolic subgroup of $G_2$, and $Z_2$ is the unique $A_{P_2}$-orbit in $D(G_2,h_2)$ which contains the image of $Z_1$. 
 
 Assume that  the above condition (i) and the following condition (ii) is satisfied.
 
 \medskip
 
 (ii) The map  $\Lie(G_{1,u})\to \Lie(G_{2,u})$ is injective.
 
 \medskip
 
 Then we have a morphism of real analytic manifolds with corners
 $$D(G_1,h_1)_{\BS}\to D(G_2, h_2)_{\BS}\;;\;(P_1, Z_1) \mapsto (P_2, Z_2),$$
 where $P_1$ and $P_2$ are as above, $Z_2$ is as above if $Z_1$ is an $A_{P_1}$-orbit, and in the case where $Z_1$ is a $B_{P_1}$-orbit, $Z_2$ is the unique $B_{P_2}$-orbit which contains the image of $Z_1$. 
 
 Here in the case where $Z_1$ is a $B_{P_1}$-orbit, the image of $Z_1$ in  $D(G_2,h_2)$ does not meet $D(G_2, h_2)_{\spl}$ by the injectivity of $\Lie(G_{1,u})\to \Lie(G_{2,u})$.

 \end{sbpara}

 \begin{sbpara}\label{ft2} Next we consider the spaces of $\SL(2)$-orbits.

 In the case $G_1$ and $G_2$ are reductive, we have a morphism $D(G_1, h_1)_{\SL(2)}\to D(G_2, h_2)_{\SL(2)}$ of locally ringed spaces with log structures with sign, which sends the class of an $\SL(2)$-orbit $(\rho_1, \varphi_1)$ ($\rho_1: \SL(2)^n_\R\to G_\R$, $\varphi_1: \frak H^n \to D(G_2, h_2)$) to the class of $(\rho_2, \varphi_2)$, where $\rho_2$ is the composition $\SL(2)_\R^n \overset{\rho_1}\to G_{1,\R}\to G_{2,\R}$ and $\varphi_2$ is the composition $\frak H^n \to D(G_1, h_1)\to D(G_2, h_2)$.
 
 In general, we have morphisms 
\begin{align*}
D(G_1, h_1)^A_{\SL(2)} & \to D(G_2, h_2)^A_{\SL(2)}\;\;  (\text{for the structures}\; I,II), \\
D(G_1, h_1)^{\star,\mild}_{\SL(2)} & \to D(G_2, h_2)^{\star,\mild}_{\SL(2)},
\end{align*}
where $(-)^A$ denotes the part consisting of $A$-orbits, which sends $(p_1, Z_1)$ to $(p_2, Z_2)$ ($p_1\in D(G_1, h_1)_{\red,\SL(2)}$, $Z_1$ is an $A$-orbit in $D(G_1, h_1)$, $p_2$ is the image of $p_1$ in $D(G_2, h_2)_{\red,\SL(2)}$ and $Z_2$ is the unique $A$-orbit in $D(G_2, h_2)$ containing the image of $Z_1$). 

If $\Lie(G_{1,u})\to \Lie(G_{2,u})$ is injective, we have morphisms
\begin{align*}
D(G_1, h_1)_{\SL(2)} & \to D(G_2, h_2)_{\SL(2)}\;\;  (\text{for the structures}\; I,II), \\
D(G_1, h_1)^{\star}_{\SL(2)} & \to D(G_2, h_2)^{\star}_{\SL(2)},
\end{align*}
which sends $(p_1, Z_1)$ to $(p_2, Z_2)$ ($p_1\in D(G_1, h_1)_{\red,\SL(2)}$, $p_2$ is the image of $p_1$ in $D(G_2, h_2)_{\red,\SL(2)}$, either $Z_1$ is an $A$-orbit in $D(G_1, h_1)$ 
 and $Z_2$ is the unique $A$-orbit in $D(G_2, h_2)$ containing the image of $Z_1$, or  
$Z_1$ is a $B$-orbit in $D(G_1, h_1)$ 
 and $Z_2$ is the unique $B$-orbit in $D(G_2, h_2)$ containing the image of $Z_1$). 

\end{sbpara}

\begin{sbpara}\label{ft3} Lastly, we consider the spaces of nilpotent orbits.

  Assume that we are given $\Sig_1$ for $G_1$. 
  Assume that the images of elements of $\Sig_1$ in $\Lie(G'_{2,\R})$ form a weak fan $\Sig_2$ in $\Lie(G_2')$. 
  (For example, this is satisfied if $\Lie(G_1')\to \Lie(G_2')$ is injective.)
  Then we have 
a map 
$$D(G_1, h_1)_{\Sig_1}\to D(G_2, h_2)_{\Sig_2},$$ 
a continuous map $$D(G_1,h_1)^{\sharp}_{\Sig_1, *}\to D(G_2, h_2)^{\sharp}_{\Sig_2, *}\;\;\;\text{for $*=[:], \val, [\val]$},$$ 
and a morphism of log manifolds  
$$\Gamma_1 \bs D(G_1, h_1)_{\Sig_1}\to \Gamma_2 \bs D(G_2, h_2)_{\Sig_2},$$
where $\Gamma_j\subset G'_j(\Q)$ $(j=1,2)$ are neat semi-arithmetic subgroups such that $\Gamma_j$ is strongly compatible with 
$\Sig_j$ and the image of $\Gamma_1$ in $G_2(\Q)$ is contained in $\Gamma_2$.

\end{sbpara}

\begin{sbprop}\label{prod0} Let $G_1$ and $G_2$ be linear algebraic groups over $\Q$ and let $h_j: S_{\C/\R}\to G_{j,\R}$ ($j=1,2$) be homomorphisms as in $\ref{D}$. 
  Let $h=(h_1, h_2): S_{\C/\R}\to (G_1\times G_2)_\R$. 
  Then $h$ satisfies the condition as in $\ref{D}$. The morphisms $D(G_1\times G_2, h) \to D(G_j, h_j)$ associated to the projections $G_1\times G_2\to G_j$ ($j=1,2$) induce an isomorphism $$D(G_1\times G_2, h)\overset{\sim}\to D(G_1, h_1) \times D(G_2, h_2).$$

\end{sbprop}

\begin{pf}
In the case where $G_1$ and $G_2$ are reductive, this is proved as 
$$D(G_1\times G_2, h)= \{(G_1\times G_2)(\R)\text{-conjugate of}\; h\}$$
$$= \prod_{j=1}^2 \{G_j(\R)\text{-conjugate of}\; h_j\}= D(G_1,h_1)\times D(G_2, h_2).$$
The general case is reduced to the reductive case by Proposition \ref{Dandgr} using 
$(G_1\times G_2)_u(\R)=G_{1,u}(\R)\times G_{2,u}(\R)$ (this proves that $\spl(W)$ for $G_1\times G_2$ is the product of $\spl(W)$ for $G_j$) and  $\cL(p)= \cL(p_1) \times \cL(p_2)$ for $p_j\in D(G_j, h_j)_{\red}$ and $p=(p_1,p_2)\in D(G_1\times G_2,h)_{\red}$. 
\end{pf}

 \begin{sbprop}\label{prodD} Let the notation be as in Proposition $\ref{prod0}$ and assume that $h_j$ ($j=1,2$) are $\R$-polarizable. 
 
 Then $h$ is $\R$-polarizable and we have$:$

  $(1)$  The canonical maps $D(G_1\times G_2, h)^{\mild}_{\BS} \to D(G_j, h_j)^{\mild}_{\BS}$ ($j=1,2$) induce an isomorphism of real analytic manifolds with corners 
 $$D(G_1\times G_2, h)^{\mild}_{\BS}\overset{\sim}\to D(G_1, h_1)^{\mild}_{\BS} \times D(G_2, h_2)^{\mild}_{\BS}.$$  
 In particular, if $G_1$ and $G_2$ are reductive, we have 
 $$D(G_1\times G_2, h)_{\BS}\overset{\sim}\to D(G_1, h_1)_{\BS} \times D(G_2, h_2)_{\BS}.$$

 $(2)$ Assume that $\Sig_j$ for $G_j$ are given ($j=1,2$). 
  Let $\Sig= \{\sig_1\times \sig_2\;|\; \sig_i\in \Sig_i\}$. 
Then $$D(G_1\times G_2,h)^{\sharp}_{\Sig_1\times\Sig_2} \overset{\sim}\to  D(G_1,h_1)^{\sharp}_{\Sig_1}\times D(G_2, h_2)^{\sharp}_{\Sig_2}$$ as topological spaces, and for neat semi-arithmetic subgroups $\Gamma_j$ of $G'_j(\Q)$ which are strongly compatible with $\Sig_j$  ($j=1,2$), we have an isomorphism of log manifolds
 $$(\Gamma_1\times \Gamma_2)\bs D(G_1\times G_2, h)_{\Sig}\overset{\sim}\to\Gamma_1\bs D(G_1,h_1)_{\Sig_1} \times \Gamma_2 \bs D(G_2, h_2)_{\Sig_2}.$$
 
 \end{sbprop} 

\begin{pf} (1) Parabolic subgroups of $G_1\times G_2$ are $P_1\times P_2$ for parabolic subgroups $P_j$ of $G_j$. (This fact is deduced from the surjectivity of the map $\Hom({\bf G}_m, \cG) \to \{\text{parabolic subgroups of}\;\cG\}$ in \ref{P_is_para} and from  $\Hom(\Gm, G_1 \times G_2) = \Hom(\Gm, G_1) \times \Hom(\Gm, G_2)$.)  We have $A_{P_1\times P_2}=A_{P_1}\times A_{P_2}$. 
 Hence the converse map is given by $((P_1, Z_1), (P_2, Z_2))\mapsto (P_1\times P_2, Z_1\times Z_2)$.

(2) The converse map is given by $((\sig_1, Z_1), (\sig_2, Z_2))\mapsto (\sig_1\times \sig_2, Z_1\times Z_2)$. 
\end{pf}

\begin{sbrem} On the other hand, even if $G$ is reductive, $D(G_1\times G_2, h)_{\SL(2)}$ need not be the product of $D(G_j, h_j)_{\SL(2)}$  $(j=1,2)$ as examples in \ref{Ex1} and \ref{Ex2} show. 

\end{sbrem}

In \ref{prod2}--\ref{Ex6}, we give examples of $G_1\to G_2$ whose associated morphisms of spaces of Borel--Serre orbits do not exist,  looking at examples of  convergence and divergence in the extended period domains. 

 \begin{sbpara}\label{prod2} {\it Example.}
  Let $(G_1,h_1)$  be the $(G, h)$ of \ref{Ex3} and let $G_1=G\times G$, $h_1=(h,h): S_{\C/\R}\to G\times G$, let $G_2=G$, let $G_1\to G_2$ be the second projection. Then the induced morphism $D(G_1,h_1)\to D(G_2,h_2)$ is understood as $\C^2\to \C\;;\; (c_1,c_2)\mapsto c_2$. This morphism does not extend to a continuous map $D(G_1,h_1)_{\BS} \to D(G_2, h_2)_{\BS}$. Note that in this case, $\Lie(G_{1,u})\to \Lie(G_{2,u})$ is not injective. Let $P_1=G_{1,\red}$,  the unique parabolic subgroup of $G_{1,\red}$. 
 For $\theta\in \R$, let $Z(\theta)$ be the $B_{P_1}$-orbit $\{(ir\cos(\theta), ir\sin(\theta)\;|\; r\in \R_{>0}\}$ in $\C^2=D(G_1,h_1)$.
  Then when $\theta\to 0$ and $r\to \infty$, 
 $(ir\cos(\theta), ir\sin(\theta))\in \C^2=D(G_1, h_1)$ converges to $(P_1, Z(0))\in D(G_1,h_1)_{\BS}$, but the image $ir\sin(\theta)\in \C=D(G_2, h_2)$ does not converge in  $D(G_2,h_2)_{\BS}$. 
 
 The point is that for finite dimensional graded $\R$-vector spaces $V_1$ and $V_2$ of weigh $\leq -1$ and for a linear map $f: V_1\to V_2$ which is compatible with the gradings, the following three conditions are equivalent. (i) Either $f$ is injective or $V_2=0$. (ii) $f$ extends to a continuous map $\overline{V}_1\to \overline{V}_2$. If these conditions are satisfied, $f$ extends to a morphism $\overline{V}_1\to \overline{V}_2$ of real analytic manifolds with corners uniquely.  
 
 \end{sbpara}

 \begin{sbpara}\label{Ex5} 
{\it Example.} Let $G_1=\GL(2) \times \GL(2)$, $G_2=\text{GSp}(4)$, and let $f:G_1\to G_2$ be the natural embedding. Let $h_1:S_{\C/\R}\to G_{1,\R}$ be as in \ref{Ex2} and let $h_2:S_{\C/\R}\to G_{2,\R}$ be $f\circ h_1$.   Then both $h_1$ and $h_2$ are $\R$-polarizable. We show that the canonical map $D(G_1, h_1)\to D(G_2, h_2)$, which we also denote by $f$, does not extend to a continuous map $D(G_1, h_1)_{\BS}\to D(G_2, h_2)_{\BS}$. 
 
 In fact, we have:
 
(1)  When $y_1, y_2 \to \infty$, $(iy_1, iy_2)\in \frak H^{\pm} \times \frak H^{\pm}=D(G_1, h_1)$ converges in $D(G_1, h_1)_{\BS}$.  
 
 On the other hand, the identity map of $D(G_2, h_2)$ extends to a homeomorphism $D(G_2, h_2)_{\SL(2)}\overset{\sim}\to D(G_2, h_2)_{\BS}$ (\cite{KU1}  Theorem 6.7). We have
 
 (2) When $y_2, y_1/y_2\to \infty$, $f(iy_1, iy_2)$ converges in $D(G_2, h_2)_{\SL(2)}$ to a point $a$, and  when $y_1, y_2/y_1\to \infty$, $f(iy_1, iy_2)$ converges in $D(G_2, h_2)$ to a point $b$, and $a\neq b$.  
 
 Hence in $D(G_2, h_2)_{\BS}$, $f(iy_1, iy_2)$ for $y_1, y_2\to \infty$ with $y_1/y_2\to \infty$ and that with $y_2/y_1\to \infty$ have different limits. 
  
  \end{sbpara}

\begin{sbpara}\label{Ex6} 
{\it Example.} 
 Let $G_1=\GL(2) \times \GL(2)$, $G_2=\text{GSp}(6)$, and let $f: G_1\to G_2$ be the homomorphism $(g_1, g_2)\mapsto g_1\otimes \text{Sym}^2(g_2)$. Let $h_1:S_{\C/\R}\to G_{1,\R}$ be as in \ref{Ex2} and let $h_2:S_{\C/\R}\to G_{2,\R}$ be $f\circ h_1$.  Both $h_1$ and $h_2$ are $\R$-polarizable. 

We explain the following (1) and (2).

(1) The canonical map $D(G_1, h_1)\to D(G_2, h_2)$, which we also denote by $f$,  does not extend to a continuous map $D(G_1, h_1)_{\BS} \to D(G_2, h_2)_{\BS}$.

(2) The identity map of $D(G_2,h_2)$ does not extend to a continuous map $D(G_2, h_2)_{\SL(2)} \to D(G_2, h_2)_{\BS}$. 

These (1) and (2) follow from (3)--(6) below. Note that $D(G_1, h_1)=\frak H^{\pm} \times \frak H^{\pm}$.

(3) When $y_1, y_2\to \infty$, $(iy_1, iy_2)$ converges in $D(G_1, h_1)_{\BS}$.

(4) When $y_2, y_2/y_1\to \infty$, $f(iy_1, iy_2)$ converges in $D(G_2, h_2)_{\SL(2)}$. 

Let $P$ and $Q$ be the parabolic subgroup of $G_2$ associated to the following homomorphisms $\mu: {\bf G}_m\to G_2$ and $\nu: {\bf G}_m\to G_2$, respectively.  $$\mu(t):=  
f\Bigl(
\begin{pmatrix} t^{-3}&0\\0&t^3\end{pmatrix}, \begin{pmatrix} t^{-1}&0\\0&t\end{pmatrix}\Bigr),\quad 
\nu(t):= 
f\Bigl(\begin{pmatrix} t^{-3}&0\\0&t^3\end{pmatrix}, \begin{pmatrix} t^{-2}&0\\0&t^2\end{pmatrix}\Bigr).$$ 
  Then $P\neq Q$ because the adjoint action of ${\bf G}_m$ on $\Lie(G_2)$ defined by $\mu$ (resp.\ $\nu$) multiplies 
$$N= f\Bigl(\begin{pmatrix} 0&1\\0&0\end{pmatrix}, \begin{pmatrix} 0&0\\1&0\end{pmatrix}\Bigr)$$
by $t^{-2}$ (resp.\ $t^2$).

 (5) When $y\to \infty$, $f(iy^3, iy)$ converges in $D(G_2, h_2)_{\BS}$ to a point whose associated parabolic subgroup of $G_2$ is $P$.

 (6) When $y\to \infty$, $f(iy^3, iy^2)$ converges in $D(G_2, h_2)_{\BS}$ to a point whose associated parabolic subgroup of $G_2$ is $Q$.

(These (5) and (6) are essentially contained in
\cite{KU1} Proposition 6.10 and also in \cite{KU2} Section 12.4.)

\end{sbpara}

In the rest of this %
Section \ref{ss:functoriality}, we assume that $G$ is reductive and let $Z$ be the center of $G$. We consider the case $G_1=G$ and $G_2=G/Z$. 
We assume that $h_0: S_{\C/\R}\to G_\R$ is $\R$-polarizable. 
We will see that the extended period domains for $G$ are understood from those of the semisimple group $G/Z$.  
 
 \begin{sblem}\label{redssP}  Let $\bar h_0$ be the homomorphism $S_{\C/\R}\to (G/Z)_\R$ induced by $h_0: S_{\C/\R}\to G_\R$.
  Then $\bar h_0$ is $\R$-polarizable. 
\end{sblem}

\begin{pf}
This is because $\Lie((G/Z)')=\Lie(G')$ and hence $\Ad(\bar h_0(i))$ on $\Lie((G/Z)'_\R)$ is a Cartan involution.
\end{pf}

\begin{sbprop}\label{redss3} 

The complex analytic manifold $D(G, h_0)$ is an open and closed submanifold of the complex analytic manifold $D(G/Z, \bar h_0)$. 

\end{sbprop}
\begin{pf} We have a canonical morphism $D(G, h_0) \to D(G/Z, \bar h_0)$ of complex analytic manifolds. Hence, to prove Proposition \ref{redss3}, it is sufficient to prove that via this map, $D(G, h_0)$ is an open and closed real analytic submanifold of $D(G/Z, \bar h_0)$.

We first prove

 \smallskip
 
 {\bf Claim.} Let $g\in G(\C)$. Then  $\Int (g)(h_0)=h_0$ in $\Hom(S_{\C/\R,\C}, G_\C)$ if and only if $\Int(g)(\bar h_0)= \bar h_0$ in $\Hom(S_{\C/\R, \C}, (G/Z)_\C)$.
 Here $\text{Int}(g)$ denotes the inner-automorphism given by $g$. 
  \smallskip
  
  Proof of Claim. Assume  $\text{Int}(g)(\bar h_0)= \bar h_0$. Then there is a homomorphism $z: S_{\C/\R,\C}\to Z_\C$ such that $\text{Int}(g)(h_0(s))=z(s)h_0(s)$ for all $s\in S_{{\C/\R},\C}$.  Since $z(s)= \text{Int}(g)(h_0(s))h_0(s)^{-1}$ belongs to the commutator subgroup $G'_\C$ of $G_\C$ and $G'_\C\cap Z_\C$ is finite, the image of $z: S_{\C/\R.\C}\to Z_\C$ is finite. Since $S_{\C/\R,\C}$ is connected, $z$ is the trivial homomorphism. 
Claim  is proved.

  Let $N= \{g\in G(\C)\;|\; \text{Int}(g)(h_0)=h_0\}= \{g\in G(\C)\;|\; \text{Int}(g)(\bar h_0)=\bar h_0\}$ and let $\bar N= N/Z(\C)$. Since  $G(\R)/Z(\R)$ is an open and closed real analytic submanifold of $(G/Z)(\R)$,
  $D(G, h_0)= G(\R)/(G(\R)\cap N)= (G(\R)/Z(\R))/((G(\R)/Z(\R))\cap \bar N)$ is an open and closed real analytic submanifold of  $(G/Z)(\R)/((G/Z)(\R)\cap \bar N)= D(G/Z, \bar h_0)$ as a real analytic manifold. 
  \end{pf}

\begin{sbpara}\label{redss4} {\it Example.} The map $D(G, h_0) \to D(G/Z, \bar h_0)$ for a reductive $G$ need not be bijective. Let $(\Z/4\Z)(1)$ be the algebraic group of $4$-th roots of $1$ over $\R$ and let $G$ be the semi-direct product of $S_{\C/\R}$ and $\Z/4\Z(1)$ in which  $S_{\C/\R}$ is the normal subgroup and the action of the generator of $\Z/4\Z(1)$ on  $S_{\C/\R}$ via the inner-automorphism is $z\mapsto z^{-1}$. Let $h_0: S_{\C/\R}\to G$ be the inclusion map. Then $G(\R)= \C^\times \times \{\pm 1\}$ (here $\{\pm 1\}\subset \Z/4\Z(1)$)
and $D(G, h_0)$ is a one-point set. On the other hand, $G/Z$ is the semi-direct product of $S_{\C/\R}/\{\pm 1\}$ and $\Z/4\Z(1)/\{\pm 1\}\simeq \{\pm 1\}$ in which $S_{\C/\R}/\{\pm 1\}$ is normal and the generator of $\Z/4\Z(1)/\{\pm 1\}$ acts on it via the inner-automorphism by $z\mapsto z^{-1}$. Hence $D(G/Z, \bar h_0)$ consists of two points. 

Both $\Dc(G, h_0)$ and $\Dc(G/Z, \bar h_0)$ consist of two points. 

 \end{sbpara}
  
\begin{sbprop}\label{redss5}  As a real analytic manifold with corners, $D(G, h_0)_{\BS}$  is canonically isomorphic to an open and closed subspace of 
$D(G/Z, \bar h_0)_{\BS}$. 
\end{sbprop}

\begin{pf}  By Proposition \ref{redss3}, we have that $Y:=D(G/Z, \bar h_0)$ is the disjoint union $Y_1\coprod Y_2$ of  open closed subspaces $Y_1$ and $Y_2$ of $Y$ such that the map $X:=D(G, h_0)\to D(G/Z, \bar h_0)$ induces an isomorphism $X \overset{\sim}\to Y_1$. As is easily seen, $Y_{\BS}$ is the disjoint union $Y_{\BS,1}\coprod Y_{\BS, 2}$, where $Y_{\BS, j}$ ($j=1,2$) is the open and closed subset of $Y_{\BS}$ consisting of all elements $(P, Z)$ such that $Z\subset Y_j$. 
The morphism $X \to Y$  induces a  morphism 
 $X_{\BS}\to Y_{\BS}$ and this induces a morphism $X_{\BS}\to Y_{\BS,1}$. We show that the last morphism is an isomorphism.

We have a bijection $P\mapsto P/Z$ from the set of all parabolic subgroups of $G$ to that of $G/Z$. It is sufficient to prove that we have an isomorphism $X_{\BS}(P) \overset{\sim}\to Y_{\BS,1}(P/Z):= Y_{\BS, 1}\cap Y_{\BS}(P/Z)$. Since $A_P \overset{\sim}\to A_{P/Z}$ and $\bar A_P \overset{\sim}\to {\bar A}_{P/Z}$, we have
$$X_{\BS}(P)= D(G, h_0)\times^{A_P} \bar A_P \overset{\sim}\to Y_1 \times^{A_{P/Z}} \bar A_{P/Z}=Y_{\BS,1}(P/Z).$$
\end{pf}

\begin{sbprop}\label{redss6}  As a locally ringed space with log structure with sign, $D(G, h_0)_{\SL(2)}$ (resp.\ $D(G, h_0)_{\BS, \val}$, resp.\ $D(G, h_0)_{\SL(2),\val}$)  is canonically isomorphic to an open and closed subspace of 
$D(G/Z, \bar h_0)_{\SL(2)}$ (resp.\  $D(G/Z, \bar h_0)_{\BS, \val}$, resp.\ $D(G/Z, \bar h_0)_{\SL(2),\val}$). 
\end{sbprop}

\begin{pf} We use the notation in the proof of Proposition \ref{redss5}.

The case of $D_{\SL(2)}$ is proved as follows. As is easily seen, $Y_{\SL(2)}$ is the disjoint union of $Y_{\SL(2),1}\coprod Y_{\SL(2), 2}$, where $Y_{\SL(2), j}$ ($j=1,2$) is the open and closed subspace of $Y_{\SL(2)}$ consisting of all elements whose torus orbits are contained in $Y_j$. 
The map $X_{\SL(2)}\to Y_{\SL(2)}$ induces a morphism $X_{\SL(2)}\to Y_{\SL(2),1}$. We show that the last morphism is an isomorphism. 
We use

\medskip

{\bf Claim.} For any field $E\supset \Q$, the map $\Hom(\SL(2)^n_E,G_E) \to \Hom(\SL(2)^n_E, (G/Z)_E)$ is a bijection. 

\medskip

We prove Claim.  We first prove the injectivity. Assume that $h_1, h_2\in \Hom(\SL(2)^n_E, G_E)$ have the same image in $\Hom(\SL(2)^n_E, (G/Z)_E)$. Then there is a homomorphism $a: \SL(2)^n_E \to Z_E$ such that $h_2=ah_1$. But $a$ is trivial because $\SL(2)^n= [\SL(2)^n, \SL(2)^n]$. Next we prove the surjectivity. Let $h\in \Hom(\SL(2)^n_E, (G/Z)_E)$. Then the image of $h$ is contained in the commutator subgroup $G'_E/(G'\cap Z)_E$ of $(G/Z)_E$, where $G'=[G, G]$. Since $G'\to G'/(G'\cap Z)$ is an isogeny and since $\SL(2)^n$ is simply connected,  this homomorphism $\SL(2)^n_E\to G'_E/(G'\cap Z)_E$ comes from a homomorphism $\SL(2)^n_E\to G'_E$. Claim is proved. 

By Claim and by using the description of the set of $\SL(2)$-orbits in (ii) or (iii) in Lemma \ref{redp1}, we see that the map $X_{\SL(2)}\to Y_{\SL(2),1}$ is bijective. It remains to compare the real analytic structures and the log structures with sign. By using the local descriptions Proposition \ref{emb1} of these structures, it is sufficient to prove that the canonical map from space $\spl(W')$ in Proposition \ref{emb1} for $(G,h_0)$ to the corresponding  space $\spl(\bar W')$ for $(G/Z, \bar h_0)$ is an isomorphism of real analytic manifolds. This map is identified with the map $G_{\R, W',u}\to G_{\R, \bar W', u}$ and is identified with the isomorphism $\Lie(G_{\R, W', u}) \overset{\sim}\to \Lie(G_{\R, \bar W', u})$. This completes the proof for $D_{\SL(2)}$. 

The proof for $D_{\BS,\val}$ (resp.\ $D_{\SL(2),\val}$) is similar to that for $D_{\BS}$ (Proposition \ref{redss5}) (resp.\ for $D_{\SL(2)}$). 
\end{pf}

\begin{sbprop}\label{redss7} Let $\Sig$ be a weak fan in $\Lie(G')$. Then 
$D(G, h_0)^{\sharp}_{\Sig}$ (resp.\  $D(G, h_0)^{\sharp}_{\Sig,[:]}$, $D(G, h_0)^{\sharp}_{\Sig,\val}$, $D(G, h_0)^{\sharp}_{\Sig,[\val]}$)     is canonically homeomorphic to an open and closed subset of $ D(G/Z, \bar h_0)^{\sharp}_{\Sig}$ (resp.\ $D(G/Z, h_0)^{\sharp}_{\Sig,[:]}$, $D(G/Z, h_0)^{\sharp}_{\Sig,\val}$, $D(G/Z, h_0)^{\sharp}_{\Sig,[\val]}$).

Here we denote the image of $\Sig$ under the isomorphism $\Lie(G'_\R)\overset{\sim}\to \Lie((G/Z)'_\R)$  by the same letter $\Sig$. 
\end{sbprop}

\begin{pf} We use the notation in the proof of Proposition \ref{redss6}.

For $j=1,2$, let $Y^{\sharp}_{\Sig, [:], j}$ be the inverse image  of $Y_{\SL(2), [:], j}$ under the CKS map $Y^{\sharp}_{\Sig, [:]}\to Y_{\SL(2)}$. Then  $Y^{\sharp}_{\Sig, [:]}$ is the disjoint union of an open and closed subsets $Y^{\sharp}_{\Sig,[:], 1}$ and $Y^{\sharp}_{\Sig, [:]. ,2}$. Since the map $Y^{\sharp}_{\Sig, [:]}\to Y^{\sharp}_{\Sig}$ is proper and all fibers of this map are connected, $Y^{\sharp}_{\Sig}$ is the disjoint union of the open and closed subsets $Y^{\sharp}_{\Sig,1}$ and $Y^{\sharp}_{\Sig, 2}$, where $Y^{\sharp}_{\Sig, j}$ ($j=1,2$) denotes the image of $Y^{\sharp}_{\Sig,[:],j}$ in $Y^{\sharp}_{\Sig}$.  The set $Y^{\sharp}_{\Sig, j}$ is the subset of $Y^{\sharp}_{\Sig}$ consisting of nilpotent $i$-orbits $(\sig, Z)$ such that there is an $F\in Z$ having the property that if $N_1, \dots, N_n$ generate $\sig$, then $\exp(\ts_{k=1}^n iy_kN_k)F\in Y_j$ if $y_k\geq 0$ for all $k$. From this description, we see that the map $X^{\sharp}_{\Sig}\to Y^{\sharp}_{\Sig,1}$ is bijective. The coincidence of the topologies can be seen by the fact that for each $\sig\in \Sig$, both $X^{\sharp}_{\sig}$ and $Y^{\sharp}_{\sig,1}$ have  the quotient topologies of the topology of $\{(q, F)\in |\toric|_{\sig} \times X\;|\; (q,F)\;\text{belongs to}\; E_{\sig}\;\text{of}\; (G,h_0)\}\simeq 
\{(q, F)\in |\toric|_{\sig} \times Y_1\;|\; (q,F)\;\text{belongs to}\; E_{\sig}\;\text{of}\; (G/Z,\bar h_0)\}$ (this fact follows from Proposition \ref{claim2}). This proves Proposition \ref{redss7} for $D^{\sharp}_{\Sig}$. The proofs for $D^{\sharp}_{\Sig, [:]}$, $D^{\sharp}_{\Sig, \val}$ and $D^{\sharp}_{\Sig, [\val]}$ are similar. 
\end{pf}

\begin{sbprop}\label{redss8} Let $\Sig$ be a weak fan in $\Lie(G')$, and let $\Gamma$ be a semi-arithmetic subgroup of $G(\Q)$ which is strongly compatible with $\Sig$ and such that the image $\overline \Gamma$ of $\Gamma$ in $(G/Z)(\Q)$ is neat. Then as a locally ringed space with log structure,
$\Gamma\bs D(G, h_0)_{\Sig}$ is canonically isomorphic to an open and closed subspace of $\overline \Gamma\bs D(G/Z, \bar h_0)_{\Sig}$. 
\end{sbprop}

\begin{pf} We use the notation in the proof of Proposition \ref{redss7}. For $j=1,2$, let $Y_{\Sig,j}$ be the subset of $Y_{\Sig}$ consisting of nilpotent orbits $(\sig, Z)$ such that for some $F\in Z$, if $N_k$ ($1\leq k\leq n$) generate $\sig$, then $\exp(\ts_{k=1}^n z_kN_k)F\in Y_j$ if $\text{Im}(z_k)\geq 0$ for all $k$. By the similar descriptions of $Y^{\sharp}_{\Sig,j}$ for $j=1,2$ in the proof of Proposition \ref{redss7}, we have that $Y_{\Sig}$ is the disjoint union of $Y_{\Sig,j}$ for $j=1,2$.  From the bijectivity of $X^{\sharp}_{\Sig}\to Y^{\sharp}_{\Sig,1}$, we obtain the bijectivity of  $X_{\Sig}\to Y_{\Sig, 1}$. Hence the map $\Gamma\bs X_{\Sig}\to {\overline \Gamma}\bs Y_{\Sig,1}$ is bijective. The coincidence of the sheaf of rings of holomorphic functions and the coincidence of the log structure can be seen by the following facts (i) and (ii) concerning  both $\Gamma(\sig)^{\gp}\bs X_{\sig}$ and ${\overline \Gamma}(\sig)^{\gp}\bs Y_{\sig, 1}$ for $\sig\in \Sig$. 

\medskip

(i) A  subset $U$ is open if and only if for any fs log analytic space $S$ and for any morphism $S\to \toric_{\sig} \times X\simeq \toric_{\sig}\times Y_1$ whose image is contained in $\{(q, F)\in \toric_{\sig}\times X\;|\; (q,F)\in E_{\sig}\;\text{of}\; (G, h_0)\}\simeq \{(q, F)\in \toric_{\sig}\times Y_1\;|\; (q,F)\in E_{\sig}\;\text{of}\; (G/Z, \bar h_0)\}$, the inverse image of $U$ in $S$ is  open. 

(ii) For an open set $U$ and for a function $f:U\to \C$, $f$ belongs to $\cO(U)$ (resp.\ $M(U)$)  if and only if for any $S$ and $S\to \toric_{\sig}\times X\simeq \toric_{\sig} \times Y_1$ as in (i) such that $\cO_S$ is a sheaf of reduced rings (i.e., rings without non-zero nilpotent elements), the pullback of $f$ on the inverse image of $U$ in $S$ belongs to $\cO_S(U)$ (resp.\ $M_S(U)$). 
 
 \medskip
 
 These (i) and (ii) follow from Proposition \ref{Ctor}.
\end{pf}

\begin{sbpara}\label{redpf} 
  We prove Proposition \ref{BSSA2}, Proposition \ref{SL2SA2}, 
(2) of Remark \ref{redrem3}, and Remark \ref{t:pro6rem}.
If $\Gamma$ is a semi-arithmetic subgroup of $G(\Q)$,  the image $\overline \Gamma$ of $\Gamma$ in $(G/Z)(\Q)$ is a semi-arithmetic subgroup, and $\overline \Gamma \cap (G/Z)'(\Q)$ is of finite index in $\overline \Gamma$.  
  By \ref{proper} (5) and by Propositions \ref{redss5}, \ref{redss6}, \ref{redss7}, \ref{redss8}, we can replace $G$ by $G/Z$ and replace $\Gamma$ by the semi-arithmetic subgroup $\overline \Gamma \cap (G/Z)'(\Q)$ of $(G/Z)'(\Q)$. 
  Thus Proposition \ref{BSSA2}, Proposition \ref{SL2SA2}, (2) of Remark \ref{redrem3}, and Remark \ref{t:pro6rem} are reduced to 
Theorem \ref{BSgl}, Theorem \ref{SL2gl}, Theorem \ref{t:property}, and Proposition \ref{t:pro6}, respectively.
\end{sbpara}

\begin{sbpara}
\label{KP2}

  We describe some details of the relation of this paper with the work \cite{KP} explained in Remark \ref{r:KP}. 

Let $H$, $M$, $h_0:S_{\C/\R}\to M_\R$, $\bar h_0: S_{\C/\R}\to (M/Z)_{\bR}$ be as in \ref{relMT}. Let $\Gamma$ be a neat arithmetic subgroup of $(M/Z)(\Q)$ which is contained in the connected component of $(M/Z)(\R)$ containing $1$ and let $\Sig$ be a fan in $\Lie((M/Z)')=\Lie(M')$ which is strongly compatible with $\Gamma$. Then $\Gamma$ acts on $D(M, h_0)$ and on $D(M, h_0)_{\Sig}$. The work  \cite{KP} shows that $\Gamma \bs D(M, h_0)_{\Sig}$ is a logarithmic manifold. The method in \cite{KP} is to use the inclusion map $D(M, h_0)\to D(\La)$ ($\La$ is as in \ref{relMT}) and to use the work \cite{KU2} on the toroidal partial compactification for $D(\La)$. 

This result can be deduced  also from Theorem \ref{t:property} for $G=M/Z$ as follows. Take a neat semi-arithmetic subgroup $\Gamma_1$ of $M'(\Q)$ whose image $\Gamma_2$ in $(M/Z)(\Q)$ is a normal subgroup of $\Gamma$ of finite index. Then %
by Proposition \ref{redss8}, $\Gamma_1 \bs D(M, h_0)_{\Sig}$ is %
an open and closed subspace of %
$\Gamma_2 \bs D(M/Z, \bar h_0)_{\Sig}$. 
  By taking the quotients by $\Gamma/\Gamma_2$, we have that $\Gamma \bs D(M, h_0)_{\Sig}$ is an open and closed subspace of $\Gamma \bs D(M/Z, \bar h_0)_{\Sig}$, which is a logarithmic manifold by Theorem \ref{t:property} and by the part of Remark \ref{redrem3} (1) for semisimple algebraic groups, and hence is a logarithmic manifold.

\end{sbpara}

\subsection{$G$-log mixed Hodge structures}\label{ss:GLMH}

We consider the $G$-MHS version of the notion log mixed Hodge structure.

In this Section \ref{ss:GLMH}, $\Gamma$ denotes a semi-arithmetic  subgroup (\ref{Gamma}) of $G(\Q)$. 

\begin{sbpara}\label{GLMH}
  Let $S$ be an object of the category $\cB(\log)$ (\ref{logmfd}). Recall that the topological space $S^{\log}$ is endowed with a proper surjective continuous map $\tau: S^{\log}\to S$ and a sheaf of rings $\cO_S^{\log}$ over $\tau^{-1}(\cO_S)$. A log $\Q$-mixed Hodge structure on $S$ is a triple $(H_\Q, W, H_\cO)$, where $H_\Q$ is a locally constant sheaf on $S^{\log}$ of finite-dimensional $\Q$-vector spaces, $W$ is an increasing filtration on $H_\Q$, $H_\cO$ is a vector bundle on $S$ endowed with an isomorphism $\cO_S^{\log}\otimes_\Q H_\Q\simeq \cO_S^{\log} \otimes_{\tau^{-1}(\cO_S)} \tau^{-1}(H_\cO)$ and with a decreasing filtration $F$,   satisfying certain conditions (see \cite{KNU2} Part III 1.3).  
   We denote by $\LMH (S)$ the category of log $\Q$-mixed Hodge structures over $S$.
  
  A {\it $G$-log mixed Hodge structure}  ({\it $G$-LMH}, for short) over $S$ is an exact $\otimes$-functor from $\Rep(G)$ 
   to $\LMH (S)$.
    
  A $G$-LMH on $S$ with a {\it $\Gamma$-level structure} is a $G$-LMH $H$ over $S$ endowed with a global section of the quotient sheaf $\Gamma\bs {\cal I}$, where $\cal I$ is the following sheaf on $S^{\log}$. For an open set $U$ of $S^{\log}$, ${\cal I}(U)$ is the set of all isomorphisms $H_\Q|_U\overset{\sim}\to \text{id}$ of $\otimes$-functors from $\Rep(G)$  to the category of local systems of $\Q$-modules over $U$. 
\end{sbpara}

\begin{sbpara}
\label{type}
     Let $\Sig$ be a weak fan in $\Lie(G')$ which is strongly compatible with $\Gamma$ (\ref{scomp}). 
  A $G$-LMH $H$ over $S$ with a $\Gamma$-level structure $\lambda$  is said to be {\it of type $(h_0,\Sig)$} if for any $s\in S$, any $t\in s^{\log}$, and 
any $\otimes$-isomorphism $\tilde \lambda_t: H_{\Q,t}\simeq  \;\text{id}$ which belongs to $\lambda_t$, 
there is a $\sig\in \Sig$ satisfying the following (i) and (ii). 

  (i) The logarithm of the action of $\Hom((M_S/\cO^\times_S)_s, \N)\subset \pi_1(s^{\log})$ on $H_{\Q,t}$ is contained, via $\tilde \lambda_t$, in $\sigma \subset \Lie(G_\R)$. 
  
  (ii) Let  $a: \cO_{S,t}^{\log}\to \C$ be a ring homomorphism which induces the evaluation $\cO_{S,s}\to \C$ at $s$ and consider the element $F: V\mapsto {\tilde \lambda}_t a(H(V))$ of $Y$ (\ref{Y}). Then this element belongs to $\Dc$ and $(\sig, F)$ generates a nilpotent orbit (\ref{nilp2}).
\end{sbpara}

\begin{sbrem}
  The definition of the type $(h_0,\Sig)$ in \cite{KNU3} 4.2.2 should be modified as above because \lq\lq the smallest cone satisfying (i)'' in the condition (ii) there may not be well-defined when $\Sig$ is not a fan.
\end{sbrem}

\begin{sbpara}
  If $(H, \lambda)$ is a $G$-LMH with a $\Gamma$-level structure of type $(h_0,\Sig)$, we have a map $S \to \Gamma \bs D_{\Sig}$, called the {\it period map} associated to $(H, \lambda)$, which sends $s\in S$ to the class of the nilpotent orbit $(\sig, Z)\in D_{\Sig}$. 
  Here $\sig$ is the smallest cone of $\Sig$ satisfying (i) and (ii) in \ref{type}, which exists by a variant of \cite{KNU2} Part III Lemma 2.2.4 (see also Appendix of this paper), and $Z$ is the associated $\exp(\sig_{\bC})$-orbit obtained in (ii) in \ref{type}.
\end{sbpara}

\begin{sbpara}\label{GMHSS}
   Let $S$ be an object of $\cB(\log)$. Let $S^{\circ}$ be the underlying locally ringed space over $\C$ of $S$ with the trivial log structure. 
  By a $G$-MHS on $S$ with a $\Gamma$-level structure, we mean a 
  $G$-LMH on $S^{\circ}$ with a $\Gamma$-level structure. By a {\it $G$-MHS on $S$ with a $\Gamma$-level structure of type $h_0$}, we mean a $G$-LMH on $S^{\circ}$ with a $\Gamma$-level structure of type $(h_0,\Sig)$, where $\Sig$ is the fan consisting of the one cone $\{0\}$. 
  \end{sbpara}

\subsection{Moduli of $G$-log mixed Hodge structures and period maps}\label{ss:mGLMH}
We show that $\Gamma \bs D_{\Sig}$ is a moduli space of $G$-LMH.

In this Section \ref{ss:mGLMH}, let $\Gamma$ be a subgroup of $G(\Q)$, and assume that either one of the following two conditions is satisfied. 

(i) $\Gamma$  is a neat semi-arithmetic subgroup (\ref{neat}, \ref{Gamma}) of $G'(\Q)$.

(ii) $G$ is reductive, $\Gamma$ is a semi-arithmetic subgroup of $G(\Q)$, and the image of $\Gamma$ in $(G/Z)(\Q)$ is neat, where $Z$ is the center of $G$. 

Note first the following. 

\begin{sbprop} The complex analytic manifold $\Gamma \bs D$ represents the functor 
$$S\mapsto \{\text{isomorphism class of }
G\text{-MHS  on $S$ with a $\Gamma$-level structure of type $h_0$}\}$$
from $\cB(\log)$ to the category of sets. 
\end{sbprop}

  The main result here is the following. 

\begin{sbthm}\label{t:main} Let $\Sig$ be as in $\ref{type}$. 
  Then 
$\Gamma \bs D_{\Sig}$ represents the functor 
$$S\mapsto \{\text{isomorphism class of }
G\text{-LMH on }S\text{ with a }\Gamma\text{-level structure of type $(h_0,\Sig)$}\}\text{.}$$
from $\cB(\log)$ to the category of sets. 
\end{sbthm} 

\begin{sbpara}
The proof of Theorem \ref{t:main} is similar to the proof of \cite{KNU2} Part III Theorem 2.6.6.

The first part of the proof is to understand the functor which $E_{\sigma}$ represents.
Then take the quotient $\Gamma(\sigma)\bs D_{\sigma}$ of $E_{\sigma}$ by $\sigma_{\C}$. 
\end{sbpara}

\begin{sbpara}\label{adelic} There is a  variant of Theorem \ref{t:main} for $G$-LMH with adelic level structure.

Let $G_1$ be a closed algebraic subgroup of $G$ and let $K$ be an open compact subgroup of $G_1({\bf A}_\Q^{\infty})$, where ${\bf A}^{\infty}_\Q$ is the adele ring of $\Q$ without the $\infty$-component. We show that under certain assumptions, the space
$$G_1(\Q) \bs (D \times G_1({\bf A}^{\infty}_\Q)/K)$$
is a moduli space of $G$-MHS with $K$-level structure and its toroidal partial compactification is a moduli space of $G$-LMH with $K$-level structure.

For each $g\in G_1({\bf A}_\Q^{\infty})/K$, let $\Gamma(g)= G_1(\Q)\cap {\tilde g}K{\tilde g}^{-1}$, where $\tilde g$ denotes a lifting of $g$ to $G_1({\bf A}^{\infty}_\Q)$. Then $\Gamma(g)$ is an arithmetic subgroup of $G_1(\Q)$. We have $\Gamma(\gamma g)=\gamma \Gamma(g)\gamma^{-1}$ for $\gamma\in G_1(\Q)$ and $g\in G_1({\bf A}^{\infty}_\Q)/K$.  

Let $R$ be a representative of $G_1(\Q)\bs G_1({\bf A}^{\infty}_\Q)/K$ in $G_1({\bf A}^{\infty}_\Q)/K$. 

We assume that for every $g\in G_1({\bf A}^{\infty}_\Q)/K$ (equivalently, for each $g\in R$),  the subgroup of $\Gamma(g)$ of $G(\Q)$ is  neat and satisfies either one of the conditions (i) and (ii) at the beginning of Section \ref{ss:mGLMH}.  

We also assume that for each $g\in R$, we are given a weak fan $\Sig(g)$ in $\Lie(G')$ which is strongly compatible with $\Gamma(g)$. 
For each $g\in G_1({\bf A}^{\infty}_\Q)/K$,  define $\Sig(g):= \Ad(\gamma)\Sig(g_0)$, where $g=\gamma g_0$ with $\gamma\in G_1(\Q)$ and $g_0\in R$ (then $\Sig(g)$ is independent of the choices of such $\gamma$ and $g_0$).  We have $\Sig(\gamma g)= \Ad(\gamma)\Sig(g)$ for all $\gamma\in G_1(\Q)$ and $g\in G_1({\bf A}^{\infty}_\Q)/K$. 

Then we have a log manifold 
$$G_1(\Q)\bs \coprod_{g\in G_1({\bf A}^{\infty})/K} D_{\Sig(g)} =G_1(\Q)\bs (\bigcup_{g\in G_1({\bf A}^{\infty}_\Q)/K}\; D_{\Sig(g)} \times \{g\})/K =  \coprod_{g\in R}  \;\Gamma(g)\bs D_{\Sig(g)}$$
(here the action of $\gamma \in G_1(\Q)$ sends an element $(x,g)$ of $\coprod_{g\in G_1({\bf A}^{\infty})/K} D_{\Sig(g)}$ with $x \in D_{\Sig(g)}$ and $g\in G_1({\bf A}^{\infty}_\Q)/K$  to the element $(\gamma x, \gamma g)$). This log manifold  contains  $$G_1(\Q) \bs (D \times G_1({\bf A}^{\infty}_\Q)/K)=\coprod_{g\in R}\; \Gamma(g)\bs D$$ as an open set. 

By Theorem \ref{t:main}, this log manifold represents the functor 
$$S\mapsto \{\text{isomorphism class of \;}
G\text{-LMH on $S$ with a}\;K\text{-level structure of type $(h_0,\Sig)$}\}$$
and the above open set represents the functor 
$$S\mapsto \{\text{isomorphism class of \;}
G\text{-MHS on $S$ with a}\;K\text{-level structure of type $h_0$}\}.$$
Here $\Sig$ denotes the family $(\Sig(g))_g$. A $K$-level structure on a $G$-LMH (or its special case $G$-MHS) $H$ on $S$ means a global section of the quotient sheaf $K\bs {\cal J}_1$, where ${\cal J}_1$ is the following sheaf on $S^{\log}$. Let $\cal J$ be the following sheaf on $S^{\log}$.
 For an open set $U$ of $S^{\log}$, ${\cal J}(U)$ is the set of all isomorphisms $H_\Q|_U\otimes_{\Q} {\bf A}^{\infty}_\Q\overset{\sim}\to \text{id}\otimes_{\Q} {\bf A}^{\infty}_\Q$ of $\otimes$-functors from $\Rep(G)$  to the category of sheaves  of ${\bf A}^{\infty}_\Q$-modules over $U$. Then $G_1({\bf A}^{\infty}_{\Q})$ acts on $\cal J$ and we have a canonical injective morphism ${\cal I} \to {\cal J}$ for $\cal I$ as in \ref{GLMH}. Let ${\cal J}_1=G_1({\bf A}^{\infty}_\Q)\cal I \subset {\cal J}$. For  a $G$-LMH (resp.\ $G$-MHS) $H$ on $S$ with a $K$-level structure, we say that $H$ is of type $(h_0, \Sig)$ if for each $s\in S$, if the $K$-level structure at the point of $S^{\log}$ lying over $s\in S$ is $g^{-1}\la$ with $g\in G_1({\bf A}^{\infty}_{\Q})$ and $\la\in {\cal I}$, then $H$ with the $\Gamma(g)$-level structure $\la$ is of type $(h_0, \Sig(g))$ (resp.\ $h_0$). The period morphism from $S$ to this log manifold (resp.\ the above open set) associated to $H$ with this $K$-level structure of type $(h_0, \Sig)$ (resp.\ $h_0$)  is as follows. On an open neighborhood of $s$,  taking the above $g$ in $R$, it is the period map from $S$ to  $\Gamma(g)\bs D_{\Sig(g)}$ (resp.\ $\Gamma(g)\bs D$) associated to  $(H, \la)$. 
 
 \end{sbpara}

Next we consider extensions of the associated period maps in Theorem \ref{t:epm1} and Theorem \ref{t:epm2}.

\begin{sbthm}\label{t:epm1}
  
  Let $S$ be a connected, log smooth, fs log analytic space, and let $U$ be the open subspace of $S$ consisting of all points of $S$ at which the log structure of $S$ is trivial.
Let $(H,\lam)$ be a $G$-MHS on $U$ with a $\Gamma$-level structure 
of  type $h_0$ ($\ref{GMHSS}$).
Let $\varphi: U\to \G\bs D$ be the associated period map. 
Assume that $(H,\lam)$ extends to a $G$-LMH on $S$ with a $\Gamma$-level structure ($\ref{GLMH}$). 
Then{\rm:}

\medskip

{\rm(1)} For any point $s\in S$, there exist an open neighborhood
$V$ of $s$, a log modification $V'$ of $V$ ({\rm{\cite{KU2}}}\ $3.6.12$), a subgroup $\G_1$ of $\G$, and a fan (we do not need a weak fan here) $\Sig$ in $\Lie(G')$
which is strongly compatible with $\G_1$ such that the period map
$\varphi|_{U\cap V}$ lifts to a morphism $U\cap V \to \G_1\bs D$ which
extends uniquely to a morphism $V'\to\G_1\bs D_\Sig$ of log manifolds.
Furthermore, we can take a 
commutative group $\Gamma_1$.
$$
\CD
U&\supset&U\cap V&\quad\subset\quad& V'\\
@V{\varphi}VV@VVV@VVV\\
\G \bs D@<<< \G_1 \bs D&\subset& \G_1 \bs D_\Sig.
\endCD
$$

{\rm(2)} Assume that $S\smallsetminus U$ is a smooth divisor. 
Then we can take $V=V'=S$ and $\G_1=\G$ in $(1)$. 
That is, we have a commutative diagram 
$$
\CD
U&\subset&S\\
@V{\varphi}VV@VVV\\
\G \bs D&\quad\subset\quad& \G \bs D_\Sig.
\endCD
$$

{\rm(3)} Assume that $\G$ is commutative.
Then we can take $\G_1=\G$ in $(1)$.

\medskip

{\rm (4)} 
Assume that $\G$ is commutative and that the following condition {\rm(i)} is satisfied.
\smallskip

{\rm(i)} There is a finite family $(S_j)_{1\leq j\leq n}$ of connected
locally closed analytic subspaces of $S$ such that $S=\bigcup_{j=1}^n S_j$
as a set and such that, for each $j$, the inverse image of the sheaf
$M_S/\cO^\times_S$ on $S_j$ is locally constant.
\smallskip

Then we can take $\G_1=\G$ and $V=S$ in $(1)$. 
\end{sbthm}

\medskip

This is the $G$-MHS version of \cite{KNU2} Part III Theorem 7.5.1 in mixed Hodge case and of \cite{KU2} Theorem 4.3.1 in pure Hodge case. 
The proof goes exactly in the same way as in the pure case treated in \cite{KU2}. 

\begin{sbthm}\label{t:epm2}
Let the notation $S$, $U$, $(H,\lam)$ and the assumptions be as in Theorem $\ref{t:epm1}$.
Let $\varphi:U\to\Gamma\bs D$ be the associated period map.
Let $S^{\log}_{[:]}= S^{\log} \times_S S_{[:]}$ and let $S^{\log}_{[\val]}=S^{\log}\times_S S_{[\val]}$, and regard $U$ as open sets of these spaces. 
Then{\rm :}

$(1)$ 
The map $\varphi:U\to\G\bs D$ extends uniquely to continuous maps
$$
S_{[:]}^{\log}\to\G\bs D_{\SL(2)}^I, \qquad S_{[\val]}^{\log} \to \G \bs D^I_{\SL(2),\val}.
$$

$(2)$ Assume that the complement $S\smallsetminus U$ of $U$ is a smooth divisor on $S$. Then 
the map $\varphi:U\to\G\bs D$ extends uniquely to a continuous map
$$
S^{\log} \to \G \bs D^{\star}_{\SL(2),\val}
$$
and hence extends uniquely to a continuous map $S^{\log} \to \G \bs D_{\BS,\val}$ and to $S^{\log} \to \G \bs D_{\BS}$. 

\end{sbthm}

This is the $G$-MHS version of \cite{KNU2} Part IV Theorem 6.3.1. 
The proof goes exactly in the same way as there.
In the proof of %
loc.\ cit., there are typos, i.e., the two $S_{[:]}$ in the middle of the proof should be changed to $S^{\log}_{[:]}$.

\subsection{Infinitesimal study}\label{ss:inf}

In this section, we consider the logarithmic tangent bundle of the moduli space of $G$-LMH. 
  We will prove

\begin{sbprop}\label{p:tangent}
Let $Z=\G\bs D_\Sig$ be as in Theorem $\ref{t:main}$ and let $(H_{\text{\rm univ}},\lambda)$ be the universal object on $Z$.
Let $\theta_Z$ be the logarithmic tangent bundle of $Z$ (i.e., the $\cO_Z$-dual of 
the sheaf 
$\omega^1_Z$ of differential forms with log poles). Then, we have a canonical isomorphism of $\cO_Z$-modules
$$
\theta_Z\simeq H_{\text{\rm univ}}(\Lie(G))_\cO/F^0H_{\text{\rm univ}}(\Lie(G))_\cO.
$$
\end{sbprop}

This is an analogue of \cite{KU2} Proposition 4.4.3 and proved in a similar way as below.

\begin{sbpara}\label{tan1}  In \ref{tan1}--\ref{tan4}, let $S$ be an object of  
$\Cal B(\log)$ and assume that $S$ is log smooth (this means that $S$  is locally a strong subspace (\ref{logmfd})  of a log smooth fs log analytic space). 

Let $H$ be a $G$-LMH on $S$. 
We will construct a commutative diagram of $\cO_S$-modules
$$\begin{matrix}&& \theta_S &\overset{\sim}\to & %
\tilde\theta_S\\
&&\downarrow && \downarrow\\
H(\Lie(G))_\cO/F^0&\overset{\sim}\to & \cE(H) &\overset{\sim}\to & %
\tilde\cE(H)
\end{matrix}$$
whose horizontal arrows are isomorphisms. Here $\theta_S$ is the logarithmic tangent bundle of $S$ %
and the sheaves %
$\tilde\theta_S$, $\cE(H)$ and %
$\tilde\cE(H)$ are defined below.

Furthermore, in the case $S=\Gamma \bs D_{\Sig}$, we will show that the right vertical arrow is an isomorphism. This will give the isomorphism in Proposition \ref{p:tangent}. 

\end{sbpara}

\begin{sbpara}\label{theta} Define the sheaf
 $\tilde\theta_S$ as follows.

Let $U$ be an open set of $S$, and let $\tilde U=U[T]/(T^2)= (U, \cO_U[T]/(T^2))$. 
Then  %
$\tilde\theta_S(U)$ is 
 the set of all morphisms $\tilde U\to S$ which extend
the inclusion morphism $U \to S$.

We have a canonical isomorphism  $\theta_S \overset{\sim}\to \tilde\theta_S$. 

\end{sbpara}

\begin{sbpara}\label{manyE1}  We define  the sheaf $\cE(H)$.

A section  of $\cE(H)$ is a collection of an $\cO_S$-homomorphism $\delta_{p,V}: F^pH(V)_{\cO}\to H(V)_\cO/F^p$  for $p\in \Z$ and $V\in \Rep(G)$ 
satisfying the following conditions (i)--(iii).

(i) Functoriality in $V$. For a morphism $V_1\to V_2$ in $\Rep(G)$, the square 
$$\begin{matrix}  F^pH(V_1)_\cO &\to & H(V_1)_\cO/F^p\\
\downarrow&&\downarrow\\
F^pH(V_2)_\cO &\to & H(V_2)_\cO/F^p
\end{matrix}$$
is commutative. 

(ii) The diagram
$$\begin{matrix}  F^{p+1}H(V)_\cO &\overset{\delta_{p+1,V}}\longrightarrow & H(V)_{\cO}/F^{p+1}\\
\downarrow &&\downarrow\\
F^pH(V)_\cO & \overset{\delta_{p,V}}\longrightarrow& H(V)_\cO/F^{p}\end{matrix}$$
is commutative for every $p$.

(iii) For $V_1,V_2\in \Rep(G)$ and for $p,q\in \Z$, the following diagram is commutative.
$$\begin{matrix}  F^pH(V_1)_\cO \otimes F^qH(V_2)_\cO &\to & (H(V_1)_\cO/F^p\otimes F^qH(V_2)_\cO)\oplus (F^pH(V_1)_\cO \otimes H(V_2)_\cO/F^q)\\
\downarrow &&\downarrow\\
F^{p+q}H(V_1\otimes V_2)_\cO &\to& H(V_1\otimes V_2)_\cO/F^{p+q}\end{matrix}$$
Here the vertical arrows are the evident ones.
The upper horizontal row is $x\otimes y\mapsto (\delta_{p,V_1}(x) \otimes y,\, x\otimes \delta_{q,V_2}(y))$.
The lower horizontal arrow is $\delta_{p+q, V_1\otimes V_2}$. 

A section of $\cE(H)$ on an open set $U$ of $S$ is defined in the same way by replacing $S$ by $U$.

We have the evident homomorphism $H(\Lie(G))_\cO/F^0 \to \cE(H)$.

\end{sbpara}

\begin{sbpara}\label{E'(H)} We define the sheaf %
$\tilde\cE(H)$.

Let $U$ be an open set of $S$, and let $\tilde U$ be as in \ref{theta}.
Then %
$\tilde\cE(H)(U)$ is 
the set of all isomorphism classes of $G$-LMH $\tilde H$ on $\tilde U$ whose
pullbacks to $U$ coincide with the restriction of $H$ to $U$.

\end{sbpara}

\begin{sbpara}\label{tangent} We define the map $\theta_S\to \cE(H)$.

Note that for $V\in \Rep(G)$, $H(V)_\cO 
=\tau_*(\cO_S^{\log}\otimes_\bQ H(V)_\bQ)$. 
Then, $d\otimes1_{H(V)_\bC}:\Cal O_S^{\log}\otimes_\bC H(V)_\bC
\to\omega_S^{1,\log}\otimes_\bC H(V)_\bC$ induces
a connection
$$
\nabla:\Cal H(V)_\cO\to\omega_S^1\otimes_{\Cal O_S}\cH(V)_\cO.
$$

We define a map
$$
\theta_S\to\bigoplus_p\Cal Hom(F^pH(V)_\cO, H(V)_\cO/F^p)
$$
by assigning $\delta\in\theta_S$ the element of
$\bigoplus_p\Cal Hom(F^pH(V)_\cO,H(V)_\cO/F^p)$ induced by the composite
map
$$
H(V)_\cO \overset\nabla\longrightarrow\omega_S^1\otimes_{\Cal O_S}H(V)_\cO
\overset{\delta}\longrightarrow H(V)_\cO.
$$
Note that this map is not $\Cal O_S$-linear but that
it induces an $\Cal O_S$-linear map
$F^pH(V)_\cO\to H(V)_\cO/F^p$.
Thus, we have a homomorphism of $\Cal O_S$-modules
$$
\theta_S\to\cE(H).
$$
\end{sbpara}

\begin{sbpara} We have a map %
$\tilde\theta_S\to \tilde\cE(H)$
 which sends a morphism $f: \tilde U\to S$ to the class of $f^*H$.

\end{sbpara}

\begin{sbpara}\label{tan2}
We have a canonical isomorphism 
$\cE(H)\to %
\tilde\cE(H)$ defied as follows. 

Let $(\delta_{p,V})_{p,V}$ be a section of $\cE(H)$ on $U$. Define the corresponding $G$-LMH $\tilde H$ on $\tilde U$ as follows. We identify the topological spaces $U^{\log}$ and $(\tilde U)^{\log}$. We define $\tilde H(V)_\Q$ with the weight filtration and the level structure as the same as the restriction of $H(V)_\Q$ to $U^{\log}$. Since $\cO_{\tilde U}^{\log}= \cO_{\tilde U} \otimes_{\cO_U} \cO_U^{\log}$,  $\tilde H(V)_\cO= \tau_*(\cO_{\tilde U}^{\log}\otimes \tilde H(V)_\Q)$ is identified with $\cO_{\tilde U} \otimes_{\cO_U}  H(V)_{\cO}|_U$. We define the $p$-th Hodge filter on $\tilde H(V)_\cO$ as the $\cO_{\tilde U}$-submodule of $\tilde H_{\cO}= \cO_{\tilde U} \otimes_{\cO_U} H(V)_\cO|_U$ generated by $x+Ty$, where $x\in F^pH(V)_{\cO}|_U$ and $y\in H(V)_\cO|_U$ such that $y \bmod F^pH(V)_\cO|_U$ coincides with $\delta_{p,V}(x)$. It is easy to see that this $\cE(H)\to %
\tilde\cE(H)$ is an isomorphism.

\end{sbpara}
\begin{sbpara}\label{tan3}
The diagram 
$$\begin{matrix} \theta_S &\overset{\sim}\to & %
\tilde\theta_S\\
\downarrow && \downarrow\\
\cE(H) &\overset{\sim}\to & %
\tilde\cE(H)
\end{matrix}$$
is commutative.

\end{sbpara}

\begin{sbpara}\label{manyE3}  

By Tannaka duality (see below),  $H(\Lie(G))_\cO$ is identified with the collection of $\delta_V\in \cE nd_{\cO_S}(H(V)_\cO)$ for $V\in \Rep(G)$ satisfying the 
following conditions (i) and (ii).

\smallskip

(i) It is 
 functorial in $V$. That is, for a morphism $h:V_1\to V_2$ in $\Rep(G)$, we have  $h\circ \delta_{V_1}=\delta_{V_2}\circ h$. 
 
(ii) For all $V_1,V_2\in \Rep(G)$, $\delta_{V_1\otimes V_2}$ coincides with $\delta_{V_1}\otimes 1 +1\otimes\delta_{V_2}$. 

\smallskip

Here Tannaka duality is used as follows. By \cite{SR}, locally on $S$, the functors $\Rep(G)\ni V \mapsto \cO_S \otimes_{\Q} V$ and $\Rep(G)\ni V \mapsto H(V)_{\cO}$ are isomorphic as $\otimes$-functors. By this and by the Tannaka duality (\cite{SR}) applied to the functor $\Rep(G)\ni V \mapsto \cO_S[T]/(T^2) \otimes_\Q V$, we have the above understanding of $H(\Lie(G))_\cO$.

\end{sbpara}

\begin{sblem}\label{manyE4} The map $H(\Lie(G))_{\cO}/F^0H(\Lie(G))_\cO \to \cE(H)$ is an isomorphism.

\end{sblem}

\begin{pf} By \cite{SR} Ch.\ IV 2.4, the $\otimes$-functors with filtrations $V\mapsto H(V)_\cO$ and $V\mapsto \gr_F(H(V)_\cO)$ are isomorphic locally on  $S$. Hence we may assume that they are isomorphic. Fix an isomorphism. Then via it, $H(\Lie(G))_\cO/F^0$ is identified with $\bigoplus_{p<0} \;\gr_F^p(H(\Lie(G))_\cO)$ and $\cE(H)$ is identified with the sheaf of collections $\delta_V : \gr_F(H(V)) \to \gr_F(H(V))$ ($V\in \Rep(G)$) satisfying the following conditions (i)--(iii). 

(i) $V\mapsto \delta_V$ is functorial in $V$. 

(ii)  $\delta_{V_1\otimes V_2}=\delta_{V_1}\otimes 1 + 1\otimes  \delta_{V_2}$ for all $V_1, V_2\in \Rep(G)$. 

(iii)   $\text{Image}(\delta_V)\subset \bigoplus_{p<0}\; \gr^p_F(H(V)_\cO)$. 

\noindent 
By \ref{manyE3}, the sheaf of those $(\delta_V)_V$ satisfying (i) and (ii) (not necessarily satisfying (iii)) is isomorphic to $\gr_F(H(\Lie(G))_\cO)$. 
  This proves Lemma \ref{manyE4}. 
\end{pf}

\begin{sbpara}\label{tan4}  If $S=\Gamma \bs D_{\Sig}$ and  $H = H_{\text{\rm univ}}$, %
$\tilde\theta_S\to \tilde\cE(S)$ is an isomorphism because $S$ is the moduli space. 

Consequently, we have $\theta_S\overset{\sim}\to \cE(H)\simeq H(\Lie(G))_\cO/F^0$.  This completes the proof of Proposition \ref{p:tangent}. 
\end{sbpara}

\begin{sbpara}\label{htang}
Let $Z=\Gamma\bs D_\Sig$ be as in Proposition \ref{p:tangent}.
We define %
$\theta_Z^1=\gr_F^{-1}H_{\text{\rm univ}}(\Lie(G))_\cO$.

\end{sbpara}

\begin{sbpara}\label{diffprd}
Let $Z=\Gamma\bs D_\Sig$ be as above.
Let $S$ be a logarithmically smooth object of $\Cal B(\log)$, let
$(H,\lambda)$ be a $G$-$\LMH$ on $S$ of type $(h_0,\Sigma)$, and let $\varphi:S\to Z$ be the corresponding period map.
Then $H$ is the pullback of $H_{\text{univ}}$ on $Z$ by
$\varphi$, and 
 the map $\theta_S\to \cE(H)\simeq H(\Lie(G))_\cO/F^0$ is identified with the
canonical map $\theta_S\to\varphi^*(\theta_Z)$ and
$\gr_F^{-1}(H(\Lie(G))_\cO)$ is identified with
$\varphi^*(%
\theta_Z^1)$.

The connection of $H(V)_\cO$ satisfies the Griffiths transversality for every $V\in \Rep(G)$  if and only if the map $\theta_S\to \varphi^*\theta_Z$ factors through $\varphi^*%
\theta_Z^1$.
\end{sbpara}

\subsection{Generalizations, I}
\label{ss:gen1}
  In this section and the next, we give two generalizations of the theory in this paper.

\begin{sbpara} Let $E$ be a subfield of $\R$.  
  We have the following generalization whose case $E=\bQ$ is the theory explained so far, as far as $\Gamma$ is not involved. 
  
  Let $\cG$ be a linear algebraic group over $E$. Assume that we are given a homomorphism 
  $h_0: S_{\C/\R}\to \cG_{\red}\otimes_E \R$ of algebraic groups over $\R$ such that the composition  $\bG_{m,\R}\to S_{\C/\R} \to  \cG_{\red}\otimes_E \R$ comes from a homomorphism $k_0: \bG_{m,E}\to \cG_{\red}$  whose image is contained in the center of $\cG_{\red}$ such that for some (and hence for every) lifting $\tilde k_0: \bG_{m,E} \to \cG$ of $k_0$, the adjoint action of $\bG_{m,E}$ on $\Lie(\cG_u)$ via $\tilde k_0$ is of weights $\leq -1$. Let $\Rep_E(\cG)$ be the category of finite-dimensional representations of $\cG$ over $E$. Note that every $V\in \Rep_E(\cG)$ has a $\cG$-stable weight filtration $W_{\bullet}V$.

  Let $D=D(\cG, h_0)$ be the set of isomorphism classes of exact $\otimes$-functors $H: \Rep_E(\cG) \to E\MHS$ over $E$ which keeps the underlying $E$-vector spaces and their weight filtrations satisfying the following condition. The homomorphism $S_{\C/\R}\to \cG_{\red}\otimes_E\R$ associated to the restriction of $H$ to $\Rep_E(\cG_{\red})$ is $\cG_{\red}(\R)$-conjugate to $h_0$. 
  
  By the method of Section \ref{ss:cxanastrD}, $D$ is regarded as  a complex analytic manifold.

\end{sbpara}

\begin{sbpara} Define $D_{\BS}$ in the same way as in Section \ref{s:DBS} except that for a parabolic subgroup $P$ of $\cG_{\red}$ (defined over $E$), $S_P$ is defined this time to be the maximal $E$-split torus in the center of $P_{\red}$ and $A_P:= \Hom(X(S_P)^+, \R^{\mult}_{>0})$. ($X(S_P)^+$ is defined in the same way as in  \ref{Drt}.) We have the Borel--Serre action of $A_P$ on $D$ and the Borel--Serre action of $\R_{>0}\times A_P$ on $D_{\nspl}$.

Thus $D_{\BS}$ is defined to be the set of all pairs $(P, Z)$, where $P$ is a parabolic subgroup of $\cG_{\red}$ and $Z$ is either an $A_P$-orbit in $D$ or a $B_P$-orbit in $D_{\nspl}$ for the Borel--Serre action. 

By the method of Section \ref{ss:BSan}, $D_{\BS}$ is regarded as a real analytic manifold with corners.
\end{sbpara}

\begin{sbpara} Define the objects $D^{\star}_{\SL(2)}$, $D^I_{\SL(2)}$, and $D^{II}_{\SL(2)}$ of $\cB'_\R(\log)$ in the same way as in Section \ref{s:DSL} except that we change the condition of the rationality of the weight filtrations $W_{\bullet}^{(j)}$ to the $E$-rationality. 

\end{sbpara}

\begin{sbpara} We define the topological spaces 
 $D^{\sharp}_{\Sig}$, $D^{\sharp}_{\Sig,[:]}$ and $D^{\sharp}_{\Sig,[\val]}$ as explained in Section \ref{s:DSig} except that the rationality of cones in $\Sig$ are replaced by $E$-rationality. 
 
\end{sbpara}

\begin{sbpara}

We have the following fundamental diagram (without $\Gamma$):
$$\begin{matrix}
&& D^{\sharp}_{\Sig,[\val]} & \to & D_{\SL(2),\val} & \leftarrow & D^{\star}_{\SL(2),\val}& \to & D_{\BS,\val}\\
&& \downarrow && \downarrow && \downarrow && \downarrow\\
D^{\sharp}_{\Sig} &\leftarrow& D^{\sharp}_{\Sig,[:]}& \to & D_{\SL(2)} &&D^{\star}_{\SL(2)}&& D_{\BS},
\end{matrix}$$
where the arrows respect the structures (the structure of $D^{\star}_{\SL(2)}$ as an object of $\cB'_\R(\log)$, etc.) of these spaces. 

These extended period domains are Hausdorff spaces. 

\end{sbpara}

\begin{sbpara}

Assume that $\cG= G \otimes_\Q E$ for a linear algebraic group $G$ over $\Q$ and that $k_0: \bG_{m,E}\to \cG_{\red}$ comes from $\bG_m \to G_{\red}$. Then $D$ is the same as the period domain  $D(G,h_0)$ for $G$, and the above extended period domains contain the ones for $G$.

\end{sbpara}

\begin{sbpara} In this generalization, however, we can not have a nice theory of the quotients by $\Gamma$.

We show an example in which  $\cG=G\otimes_\Q E$, where $G$ is a reductive  algebraic group over $\Q$, $\Gamma$ is a semi-arithmetic subgroup of $G'(\Q)$, and $\Gamma \bs D_{\BS}$, $\Gamma \bs D_{\SL(2)}$, and $\Gamma \bs D_{\Sig}$ for some $E$-rational fan $\Sig$ which is strongly compatible with $\Gamma$, are not Hausdorff. 

Let $L$ be a totally real field of degree $>2$ and let $E$ be a subfield of $\R$ which contains all conjugates of $L$. Let $G:=$ Res$_{L/\Q} (\GL(2)_L)$, $\cG:=G \otimes_\Q E=\prod_{j=1}^n  \nu_j^*(\GL(2)_L)=\GL(2)_E^n$, where $n=[L:\Q]$ and  $\nu_1, \dots, \nu_n$ are  all the different field homomorphisms  $L\to E$. Let $h_0: S_{\C/\R}\to \cG\otimes_k \R=\GL(2)_\R^n$ be the homomorphism $z\mapsto (\langle z\rangle, \dots, \langle z \rangle)$ (\ref{clEx}).  Then $D$ is canonically isomorphic to $\frak H^n$, the $n$-fold product of the upper half plane, and $G(\Q)=\GL(2, L)$ acts on it via $(\nu_1, \dots, \nu_n)$. 

Let  $\Gamma$ be a subgroup of $\SL(2, O_L)$ of finite index.  
We show that the quotient spaces  $\Gamma \bs D_{\BS}$, $\Gamma \bs D_{\SL(2)}$, and $\Gamma \bs D_{\Sig}$ are not Hausdorff, where $\Sig$ is the set of all nilpotent cones in $\Lie(\cG_\R)= \frak{gl}(2, \R)^n$ of the form $\bigoplus_{j=1}^n \R N_j$ with $N_j$  a nilpotent $(2,2)$-matrix over $E$. 

First we show that $\Gamma \bs D_{\BS}$ is not Hausdorff. For $z\in \frak H$, let $p_{\BS}(z)\in D_{\BS}$ be the limit of $(iy_1, \dots, iy_{n-1}, z)\in D=\frak H^n$, where $y_j\in \R_{>0}$ and $y_j\to \infty$. That is, $p_{\BS}(z)$ is the $A_P$-orbit containing $(i, \dots, i, z)$, where 
$P$ is the $E$-parabolic subgroup $\prod_{j=1}^n P_j$ of $\cG=\GL(2)_E^n$ with $P_j= \begin{pmatrix} *&*\\0&* \end{pmatrix}\subset \GL(2)_E$ for $1\leq j\leq n-1$ and $P_n=\GL(2)_E$. Let $S= \{p_{\BS}(ia)\;|\; a\in \R_{>0}\}$. Then  the diagonal matrix $(u, u^{-1})$ with $u\in O_L^\times$ in $\SL(2, L)$ acts on $S$ as $(i\infty, \dots, i\infty, ia) \mapsto (i\infty, \dots, i\infty, ia \nu_n(u)^2)$. But $\{\nu_n(\gamma)^2\;|\; \gamma\in \Gamma_1\}\bs \R_{>0}$ for a subgroup $\Gamma_1$ of $O_L^\times$ of finite index is not Hausdorff. 

Next we show that $\Gamma \bs D_{\SL(2)}$ is not Hausdorff. For $z\in \frak H$, let $p_{\SL(2)}(z)$ be the limit point of $(iy_1, \dots, iy_{n-1}, z)$, where $y_j\in \R_{>0}$, $y_j/y_{j+1}\to \infty$ ($1\leq j\leq n-1$) with $y_n:=1$. 
That is, $p_{\SL(2)}(z)$ is the class of $(\rho, \br)$, where $\rho: \SL(2)_\R^{n-1}\to \cG_\R=\GL(2)_\R^n$ is the homomorphism $(g_q, \dots, g_{n-1})\mapsto (g_1, \dots, g_{n-1}, 1)$ and $\br=(i, \dots, i, z)$. 
Let $S$ be the set $\{p_{\SL(2)}(ia)\; |\; a\in \R_{>0}\}$. The rest of the proof is similar to the above.

Lastly, we show that $\Gamma \bs D_{\Sig}$ is not Hausdorff. 
For $z\in \frak H$, let $p_{\Sig}(z)\in \Gamma \bs D_{\Sig}$ be the limit point of $(z_1, \dots, z_{n-1}, z)\bmod \Gamma$, where $z_j\in \frak H$ and $\text{Im}(z_j)\to \infty$. That is, $p_{\SL(2)}(z)$ is the class mod $\Gamma$ of the 
nilpotent orbit $(\sig, Z)$, where $\sig=\sig_1\times \dots \times \sig_n \subset \frak{gl}(2, \R)^n$ 
with $\sig_j= \begin{pmatrix} 0&*\\0&0\end{pmatrix}$ for $1\leq j\leq n-1$ and $\sig_n=\{0\}$ and $Z$ is 
the $\exp(\sig_\C)$-orbit in $\Dc$ which contains $\frak H^{n-1} \times \{z\}$. Let $S= \{p_{\Sig}(z)\;|\; z\in \frak H\}$. The rest of the proof is similar to the above (we use the fact that for a subgroup $\Gamma_1$ of $O_L^\times$ of finite index, $\{\nu_n(\gamma)^2\;|\; \gamma\in \Gamma_1\}\bs \frak H$ is not Hausdorff).

\end{sbpara}

\subsection{Generalizations, II}
\label{ss:gen2}

Here we give a generalization which contains partial toroidal compactifications of higher Albanese manifolds treated in \cite{KNU3}. This generalization is used in \cite{Ka} for applications to number theory. 

\begin{sbpara}
 Let $\cG$ be a normal algebraic subgroup of $G$ and let $\cQ=G/\cG$. We fix an element $b$ of $D(G, h_0)$. Let $D(G, h_0, \cG, b)\subset D(G, h_0)$ be the inverse image in $D(G,h_0)$ 
 of the image of $b$ in $D(\cQ, h_{0,\cQ})$, where $h_{0,\cQ}$ denotes the composite homomorphism $S_{\C/\R}\overset{h_0}\to G_{\red,\R}\to \cQ_{\red,\R}$. 

\end{sbpara}

\begin{sbprop} The morphism $D(G, h_0)\to D(\cQ, h_{0,\cQ})$ is smooth. 

\end{sbprop}
This follows from the surjectivity of the map of tangent spaces. 
\begin{sbcor}

$D(G, h_0, \cG, b)$ is smooth. 

\end{sbcor}

\begin{sbpara}
We consider the quotient space $\Gamma\bs D(G, h_0, \cG, b)$ and its toroidal partial compactification for a neat semi-arithmetic subgroup $\Gamma$ of $\cG'(\Q)$.  
\end{sbpara}

\begin{sbpara}\label{gen2e1}  {\it Example} 1. A higher Albanese manifold is regarded as an example of $\Gamma \bs D(G, h_0, \cG, b)$ as is explained in \cite{KNU3}.   

\end{sbpara}

\begin{sbpara}\label{gen2e2} {\it Example} 2. Let $H_0$ be a $\Z$-MHS, that is, a $\Q$-MHS endowed with a $\Z$-lattice $H_{0,\Z}$ in $H_{0,\Q}$. Assume that we have a polarization $p_w: \gr^W_w H_0 \otimes \gr^W_w H_0\to \Q(-w)$ in the sense of Deligne for each $w\in \Z$. Then  for $G$, $\cG$, $h_0$ and $b$ as below, $\Gamma \bs D(G,h_0,  \cG, b)$ is identified with 

$(*)$ the set of all isomorphism classes of $\Z$-MHS $H$ endowed with isomorphisms $\gr^W_wH \simeq \gr^W_wH_0$ of $\Z$-HS for all $w$.

 Let 
$$G= \{(g, t)\in \Aut(H_{0, \Q}, W) \times \bG_m\;|\; p_w(gx\otimes gy) = t^wp_w(x, y) \text{ for all }w\in \Z,  x,y\in \gr^W_wH_{0, \Q}\}, $$ 
$$\cG=G_u,$$
and hence 
$$\cQ= G/\cG= G_{\red}=$$ $$ \{(g, t) \in \bigl(\prod_w \Aut(\gr^W_wH_{0, \Q})\bigr) \times \bG_m \;|\; p_w(g_wx\otimes g_wy)= t^wp_w(x\otimes y)\text{ for all }
w\in \Z, x,y\in \gr^W_wH_{0,\Q}\}.$$ 
Then we have $h_0: S_{\C/\R}\to G_{\red,\R}=\cQ_\R$ and $b\in D(G, h_0)$ corresponding to $H_0$ as follows. We have the homomorphism $h_0$ which is associated to the Hodge structure $\gr^WH_0$. Then the condition (1) in Lemma \ref{pol2}  is satisfied. Let $\cC_0$ (resp.\ $\cC_{0, \red}$) be the smallest full subcategory of $\Q\MHS$ (resp.\ $\Q\HS$) which contains $H_0$ and $\Q(1)$ (resp.\ $\gr^WH_0$ and $\Q(1)$) and which is stable under taking $\oplus$, $\otimes$, duals and subquotients. Let $\Cal T_0$ be the Tannakian group of $\cC_0$ associated to the fiber functor $H\mapsto H_\Q$. Then ${\Cal T}_{0, \red}$ is identified with the Tannakian group of $\cC_{0, \red}$ associated to the fiber functor  $H \mapsto H_\Q$. We have a homomorphism $S_{\C/\R}\to {\Cal T}_{0, \red, \R}$ associated to $\gr^W H_0$, which we denote also by $h_0$. We have a canonical homomorphism ${\Cal T}_0\to G$ and this induces a morphism 
$D({\Cal T}_0, h_0)\to D(G, h_0)$. We have the canonical element of $D({\Cal T_0}, h_0)$ which sends $V\in \Rep({\Cal T}_0)$ to the corresponding $\Q\MHS$. 
Let $b$ be the image of this canonical element under  $D({\Cal T}_0, h_0)\to D(G, h_0)$.

 Let $$\Gamma= \{g\in \Aut(H_{0, \Z}, W)\; |\; \gr^W_w(g)=1 \; \text{for all}\; w\in \Z\}.$$

We prove that $\Gamma \bs D(G,h_0,  \cG, b)$ is identified with the above set $(*)$. For an element $H\in D(G, h_0, \cG, b)$, $H(V)$ for $V= H_{0,\Q}\in \Rep(G)$ with the $\Z$-lattice $H_{0, \Z}$ gives an element of the set $(*)$.  
  This gives a map  $\Gamma \bs D(G, h_0, \cG, b)\to (*)$. Conversely let  $H$ be a $\Z\MHS$ with $\gr^W_wH\simeq \gr^W_w H_0$ for all $w$. Let $\cC$ be the smallest full subcategory of $\Q\MHS$ which contains $H$ and $\Q(1)$ and which is stable under taking $\oplus$, $\otimes$, duals and subquotients. Let $\Cal T$ be the Tannakian group of $\cC$ associated to the fiber functor $H'\mapsto H'_\Q$. We have ${\Cal T}_{\red}= {\Cal T}_{0,\red}$ and 
we have a canonical element $D({\Cal T}, h_0)$. We have a canonical homomorphism $\Cal T\to G$, and the image of this canonical element in $D(G, h_0)$ belongs to  $D(G, h_0, \cG, b)$.  This gives the converse map $(*)\to \Gamma \bs D(G, h_0, \cG, b)$.

\end{sbpara}

\begin{sbpara}
Toroidal partial compactifications of $\Gamma\bs D(G, h_0,  \cG, b)$ are obtained by using a weak fan $\Sig$ in $\Lie(G)$ which is strongly compatible with $\Gamma$ such that $\sig\subset \Lie(\cG_\R)$ for all $\sig\in \Sig$. In fact, for such a $\Sig$, the map $D(G, h_0) \to D(\cQ, h_{0, \cQ})$ induces a morphism $\Gamma \bs D(G, h_0)_{\Sig}\to D(\cQ, h_{0, \cQ})$ of log manifolds. Let $D(G, h_0, \cG, b)_{\Sig}$ be the inverse image of $b_{\cQ}\in D(\cQ, h_{0,\cQ})$ in $D(G, h_0)_{\Sig}$. The fiber $\Gamma\bs D(G, h_0, \cG, b)_{\Sig}$ of $b_{\cQ}$ in $\Gamma \bs D(G, h_0)_{\Sig}$  is our toroidal partial compactification  of $\Gamma \bs D(G, h_0, \cG, b)$. 

\end{sbpara}

\begin{sbprop}\label{gen2ls} The space $\Gamma \bs D(G, b_0, \cG, b)_{\Sig}$ is a 
 log manifold which represents the  following functor on $\cB(\log)${\rm :} It  sends $S\in \cB(\log)$ to the set of all morphisms $S\to D(G, h_0)_{\Sig}$ such that the composition $S\to D(G, h_0)_{\Sig} \to D(\cQ, (h_0)_{\cQ})$ is the constant function $b_{\cQ}$. 

\end{sbprop}

\begin{pf} Let $X=\Gamma \bs D(G, h_0)_{\Sig}$, $Y=D(\cQ, b_{\cQ})$,  and let $x\in X$. Then the proof of Claim 3 in \ref{pfCtor} shows that there are an open neighborhood $U$ of $x$ in $X$, 
a log smooth fs log analytic space $Z$ over $\C$, and log differential forms $\omega_1, \dots, \omega_n$ on $Z$  such that the morphism $U\to Y$ factors as $U\to Z \to Y$ satisfying the following conditions 
(i) and (ii). Let $Z'=\{z\in Z\;|\; \omega_j(z)=0\;\text{for}\;1\leq j\leq n\}$. Here $\omega_j(z)$ denotes the log differential form on the log point $z$ obtained from $\omega_j$  (\ref{logmfd}).

(i) $U$ is isomorphic over $Z$ to an open subspace of $Z'$  for the strong topology of $Z'$ in $Z$.

(ii) The morphism $Z\to Y$ is log smooth.

Hence the fiber $Z_b$ of
$b_{\cQ}$ in $Z$ is log smooth,  and the fiber of  $b_{\cQ}$ in $U$ is an open subspace (for the strong topology)
of $\{z\in Z_b\;|\; \omega_j(z)=0\;\text{for}\;1\leq j\leq n\}$ and hence is  a log manifold. 
\end{pf}

\begin{sbpara}
In Example 1 in \ref{gen2e1}, $\Gamma \bs D(G, h_0, \cG, b)_{\Sig}$ is the toroidal partial compactification of the higher Albanese manifold discussed in \cite{KNU3} Section 5. 

\end{sbpara}

\begin{sbpara}
In  Example 2 in \ref{gen2e2}, $\Gamma \bs D(G, h_0, \cG, b)_{\Sig}$ represents the functor on $\cB(\log)$ which sends $S\in \cB(\log)$ to the set of all isomorphism classes of $\Z$-LMH $H$ of type $(h_0, \Sig)$  on $S$ endowed with isomorphisms $\gr^W_w H \simeq  \gr^W_w H_0$ of $\Z$-LMH for all $w\in \Z$ (this tells that $\gr^W_wH$ are constant $\Z$-HS for all $w$). See \cite{KNU2} Part III Section 5 for a more general moduli space of $\Z$-LMH on $S$ with given  graded quotients for the weight filtration. 

\end{sbpara}

\begin{sbpara}
On this toroidal partial compactification $\Gamma\bs D(G, h_0, \cG, b)_{\Sig}$, $\cH_{\text{univ}}(\Lie(\cG))_{\cO}/F^0$ is the log tangent bundle (this is a generalization of Proposition \ref{p:tangent}). 
\end{sbpara}


\setcounter{section}{0}
\renewcommand{\thesection}{\Alph{section}}
\section{Appendix}
  Here we give a complement to \cite{KNU2} Part III. 

  In the proof of (2) $\Rightarrow$ (2)${}''$ in \cite{KNU2} Part III Lemma 2.2.4, it is stated that \lq\lq $(\sig_0,F)$ generates a nilpotent orbit because both $(\tau,F)$ and $(\sig, F)$ generate nilpotent orbits.''
  It turns to be valid, though not trivial, and follows from the following proposition.

\begin{prop}
  Let $\Lambda=(H_0,W,(\langle\cdot,\cdot\rangle_w)_w, (h^{p,q})_{p,q})$ be as in {\rm \cite{KNU2} Part III 2.1.1}. 
  Let $\sig$ be a nilpotent cone, 
  let $\sig_0$ be its face, let $N$ be an interior of $\sig_0$, and let $F \in \check D$. 
  Assume that both $(\sig, F)$ and $(N,F)$ generate nilpotent orbits. 
  Then $(\sig_0,F)$ also generates a nilpotent orbit. 
\end{prop}

  First, the admissibility and the Griffiths transversality for $(\sig_0,F)$ follow from those for $(\sig,F)$. 
  Thus the problem is reduced to the following proposition by forgetting $\sig$. 
  (Notation is changed.) 

\begin{prop}
\label{p:defno}
  Let $\sig$ be an admissible nilpotent cone, let $N_0$ be an interior of $\sig$, and let $F \in \check D$ 
such that $F$ satisfies the Griffiths transversality with respect to $\sig$. 
  Assume that $(N_0,F)$ generates a nilpotent orbit. 
  Then $(\sig, F)$ also generates a nilpotent orbit. 
\end{prop}

\begin{sbrem}
  The proof below shows that the assumption of admissibility is weakened to that the 
weight filtration with respect to an interior is constant. 
  (Cf.\ the first and the second lines of \cite{CKS} p.505.) 
\end{sbrem}

  We prove Proposition \ref{p:defno}. 
  Let $N_1,\ldots,N_n$ generate $\sig$. 
  Then the conclusion is equivalent to that $\exp(\sum i y_j N_j)F \in D$ for any $y_j \gg 0$. 
  Since this condition can be checked on each $\gr^W_w$, we can reduce to the pure case. 
  In the rest of this proof, we assume that we are in the situation of pure weight $k$. 

  Let $M$ be the weight filtration of one (and hence for any) interior of $\sig$. 
  By Remark in Part III 1.3.2 (which is seen by \cite{CKS} (4.66)), 
the conclusion of Proposition \ref{p:defno} is equivalent to that $(H_0, \langle\cdot,\cdot\rangle_0,M[-k],F)$ is a 
mixed Hodge structure polarized in the sense of \cite{CKS} (2.26) by any interior $N$ of $\sig$. 
  By the assumption that $(N_0,F)$ generates a nilpotent orbit, $(H_0,M[-k],F)$ is at least a mixed Hodge structure.
  Hence for an interior $N$, the condition for $N$ is equivalent to that for any $l \ge0$, 
the primitive part $P_{k+l}:=\Ker(N^{l+1}\colon \gr^M_{k+l} \to \gr^M_{k-l-2})$ for $N$ is a Hodge structure polarized by 
$S_N:=\langle \cdot,N^l\cdot\rangle_0$. 
  Note here that by the assumption, each $\gr^M_w$ is already a Hodge structure, 
and by 
Griffiths transversality, any interior $N$ is a homomorphism of Hodge structures so that 
the primitive part $P_{k+l}$ for any $N$ is a Hodge structure (of weight $k+l$) whose Hodge numbers are independent of $N$ 
(determined by those for $\gr^M$).

  Consider the Hodge decomposition of the Hodge structure on $P_{k+l}$ for $N$ and denote its $(p,q)$-component by $P_N^{p,q}$. 
  Since the dimension of $P_N^{p,q}$ is constant with respect to $N$ (as noted in the above) 
and the pair of the subspace $P_N^{p,q}$ of $(\gr^M_{k+l})_{\bC}$ and 
the induced %
Hermitian form $S_N^{p,q}$ by $S_N(\cdot,\overline \cdot)$ on $P_N^{p,q}$ varies continuously with respect to $N$, 
if we prove that $S_N^{p,q}$ is always nondegenerate, the positivity for one $N$ inherits to all $N$. 

  We prove the nondegeneracy. 
  Assume that $S_N^{p,q}$ degenerates for some $N, p, q$. 
  Then there is a nonzero vector $v \in P_N^{p,q}$ such that $S_N(\cdot,\overline v)$ is zero on $P_N^{p,q}$. 
  Since the Hodge components are orthogonal with respect to $S_N(\cdot,\overline \cdot)$, 
$S_N(\cdot,\overline v)$ is zero on the whole $P_{k+l}$ for this $N$. 
  Further, consider the decomposition of $\gr^M_{k+l}$ into the images of primitive parts of weights $\ge k+l$.
  Then this decomposition is orthogonal with respect to 
$S_N=\langle\cdot,N^l\cdot\rangle_0$ (\cite{Sc} Lemma (6.4)).  
  Hence 
$\langle\cdot,N^l\overline v\rangle_0$ 
is zero on $\gr^M_{k+l,\bC}$.
  Since $M_{k+l}$ is the orthogonal complement to $M_{k-l-1}$ (the same lemma of ibid.), 
$N^l\overline v = 0$ in $\gr^M_{k-l,\bC}$ and $v=0$ in $\gr^M_{k+l,\bC}$, a contradiction.

\section*{Correction to \cite{KNU2} Part IV}

There are some typos in the proof of Theorem 6.1.1 in \cite{KNU2} Part IV.
Here we give corrections to them.

In the line 3 of the proof of loc.\ cit., $D_{\Sig,[\val]}\to D_{\Sig,[:]}\to\Gamma\bs D_\Sig$ should be changed to $D_{\Sig,[\val]}^\sharp\to D_{\Sig,[:]}^\sharp\to D_\Sig^\sharp$.

In the lines 4, 5 of the proof of loc.\ cit., 
$D_{\Sig,[\val]}$, $D_{\Sig,[:]}$, and $\Gamma\bs D_\Sig$ should be changed to 
$D_{\Sig,[\val]}^{\sharp}$, 
$D_{\Sig,[:]}^\sharp$, and $D_\Sig^\sharp$, respectively.

There are some typos in the proof of Theorem 6.1.3 in \cite{KNU2} Part IV.

In the line 2, $D^I_{\SL(2)}$ should be $\Gamma \bs D^I_{\SL(2)}$. 

In the lines 4, 5, $S_{[:]}$ should be $S_{[:]}^{\log}$.

\section*{Correction to \cite{KNU3}}
  Here we give a correction to \cite{KNU3} 6.1.2.
  The construction there does not necessarily give an MHS on 
$\Lie(\cG_{\Gamma})$.

  Change the part 

\bigskip

We define 
the weight filtration on $\Lie(\cG_{\Gamma})$ (resp.\ the Hodge 
filtration on $\Lie(\cG_{\Gamma})_\C$) as the image of that of 
$\Lie(\cG_{\Gamma_r})$ (resp.\ $\Lie(\cG_{\Gamma_r})_\C$) (2.2.4, 2.3). 
This gives a 
structure of an MHS on $\Lie(\cG_{\Gamma})$ which is independent of the 
choice of $r$. 

\bigskip

into 

\bigskip

Hereafter we assume that $\Lie(\cG_{\Gamma})$ has an MHS 
which is a quotient of MHS on $\Lie(\cG_{\Gamma_r})$ 
induced from the canonical variation of MHS in 2.2.4. 
Note that this MHS on $\Lie(\cG_{\Gamma})$ is independent of the choice of $r$. 

\section*{List of notation}

\noindent
 {\bf Period domains}
 
 $D=D(G,h_0)$  \quad  \ref{perd}, \ref{D}

 $D_{\red}$ \quad  \ref{Dred}
 
 $\Dc$ \quad   \ref{Y}

 \noindent
 {\bf Extended period domains}

 $D_{\BS}$ \quad \ref{DBSs}
 
 $D_{\SL(2)}$ \quad \ref{redSL} (case of reductive groups),  \ref{DSL}  (general case)
 
Structures of $D_{\SL(2)}$:  $D^I_{\SL(2)}$  \quad \ref{str3}, \quad  $D^{II}_{\SL(2)}$  \quad \ref{str2}
 
 $D^{\star}_{\SL(2)}$ \quad  \ref{DSL}
 
 $D^{\star,\pa}_{\SL(2)}$ \quad  \ref{lgmdW}

 $D^{\star,+}_{\SL(2)}$ \quad  \ref{lgmd+}

 $D_{\BS, \val}$,  $D^{\star}_{\SL(2), \val}$, $D^I_{\SL(2),\val}$, $D^{II}_{\SL(2), \val}$ \quad  \ref{val2}

 $D_{\Sig}$, $D^{\sharp}_{\Sig}$ \quad \ref{DSigdef} 
 
$D_{\Sig,[:]}$, $D^{\sharp}_{\Sig,[:]}$ \quad  \ref{rat3}
 
 $D_{\Sig,\val}$,  $D^{\sharp}_{\Sig, \val}$ \quad  \ref{Dval}

 $D^{\sharp}_{\Sig, [\val]}$ \quad   \ref{sig[val]}
 
 $D^{\diamond}_{\SL(2)}$  \quad \ref{diadef}

 \noindent
 {\bf Categories}

  $\cB'_\R(\log)$ \quad\ref{BR1}

$\cB(\log)$ \quad \ref{logmfd}

\bigskip

\noindent {\rm Kazuya KATO
\\
Department of mathematics
\\
University of Chicago
\\
Chicago, Illinois, 60637, USA}
\\
{\tt kkato@math.uchicago.edu}

\bigskip

\noindent {\rm Chikara NAKAYAMA
\\
Department of Economics 
\\
Hitotsubashi University 
\\
2-1 Naka, Kunitachi, Tokyo 186-8601, Japan}
\\
{\tt c.nakayama@r.hit-u.ac.jp}

\bigskip

\noindent
{\rm Sampei USUI
\\
Graduate School of Science
\\
Osaka University
\\
Toyonaka, Osaka, 560-0043, Japan}
\\
{\tt usui@math.sci.osaka-u.ac.jp}
\end{document}